\documentclass[12pt, a4paper]{scrartcl}

\usepackage[T1]{fontenc}
\usepackage[utf8]{inputenc}
\usepackage[english]{babel}

\usepackage{amsmath}
\usepackage{amssymb}
\usepackage{amsthm}

   \usepackage[all]{xy}
\usepackage{cancel}
\usepackage{mathtools}
\usepackage{bbm}
\usepackage[inline]{enumitem}

\usepackage{tikz}
\usetikzlibrary{arrows}
\usetikzlibrary{arrows.meta, calc}
\usepackage{tikz-cd}
\tikzset{symbol/.style={draw=none,
        every to/.append style={edge node={node [sloped, allow upside down, auto=false]{$#1$}}}
    }
}
\usepackage{hyperref}
\hypersetup{ colorlinks=true,
linkcolor=blue,
pdfpagemode=UseOutlines,
    pdfencoding=unicode          }

\usepackage{aliascnt}
\usepackage{microtype}
\usepackage{cleveref}

\theoremstyle{plain}

\newaliascnt{Thm}{mycounter}
\newtheorem{Thm}[Thm]{Theorem}
\aliascntresetthe{Thm}
\newaliascnt{Lem}{mycounter}
\newtheorem{Lem}[Lem]{Lemma}
\aliascntresetthe{Lem}
\newaliascnt{Prp}{mycounter}
\newtheorem{Prp}[Prp]{Proposition}
\aliascntresetthe{Prp}
\newaliascnt{Cor}{mycounter}
\newtheorem{Cor}[Cor]{Corollary}
\aliascntresetthe{Cor}
\newaliascnt{Con}{mycounter}

\aliascntresetthe{Con}

\newaliascnt{Def}{mycounter}
\newtheorem{Def}[Def]{Definition}
\aliascntresetthe{Def}
\newaliascnt{DefRem}{mycounter}

\aliascntresetthe{DefRem}
\newaliascnt{DefThm}{mycounter}

\aliascntresetthe{DefThm}
\newaliascnt{DefLem}{mycounter}
\newtheorem{DefLem}[DefLem]{Definition/Lemma}
\aliascntresetthe{DefLem}
\newaliascnt{Axm}{mycounter}

\aliascntresetthe{Axm}
\newaliascnt{Not}{mycounter}
\newtheorem{Not}[Not]{Notation}
\aliascntresetthe{Not}
\newaliascnt{Asm}{mycounter}
\newtheorem{Asm}[Asm]{Assumption}
\aliascntresetthe{Asm}
\newaliascnt{Asms}{mycounter}

\aliascntresetthe{Asms}

\newaliascnt{Rem}{mycounter}
\newtheorem{Rem}[Rem]{Remark}
\aliascntresetthe{Rem}
\newaliascnt{Eg}{mycounter}
\newtheorem{Eg}[Eg]{Example}
\aliascntresetthe{Eg}

\crefname{Thm}{Theorem}{Theorems}
\Crefname{Thm}{Theorem}{Theorems}
\crefname{Lem}{Lemma}{Lemmas}
\Crefname{Lem}{Lemma}{Lemmas}
\crefname{Prp}{Proposition}{Propositions}
\Crefname{Prp}{Proposition}{Propositions}
\crefname{Cor}{Corollary}{Corollaries}
\Crefname{Cor}{Corollary}{Corollaries}
\crefname{Con}{Conjecture}{Conjectures}
\Crefname{Con}{Conjecture}{Conjectures}
\crefname{Def}{Definition}{Definitions}
\Crefname{Def}{Definition}{Definitions}
\crefname{DefRem}{Definition/Remark}{Definitions/Remarks}
\Crefname{DefRem}{Definition/Remark}{Definitions/Remarks}
\crefname{DefThm}{Definition/Theorem}{Definitions/Theorems}
\Crefname{DefThm}{Definition/Theorem}{Definitions/Theorems}
\crefname{DefLem}{Definition/Lemma}{Definitions/Lemmas}
\Crefname{DefLem}{Definition/Lemma}{Definitions/Lemmas}
\crefname{Axm}{Axiom}{Axioms}
\Crefname{Axm}{Axiom}{Axioms}
\crefname{Not}{Notation}{Notations}
\Crefname{Not}{Notation}{Notations}
\crefname{Asm}{Assumption}{Assumptions}
\Crefname{Asm}{Assumption}{Assumptions}
\crefname{Asms}{Assumptions}{Assumptions}
\Crefname{Asms}{Assumptions}{Assumptions}
\crefname{Rem}{Remark}{Remarks}
\Crefname{Rem}{Remark}{Remarks}
\crefname{Eg}{Example}{Examples}
\Crefname{Eg}{Example}{Examples}

\crefname{proof}{proof}{proofs}
\newenvironment{Prf}[2][]{\par\bigskip\noindent \refstepcounter{mycounter}\label[proof]{#2:prf}\textbf{Proof \themycounter\ for \Cref{#2}}\ifx&#1&\else\space(#1)\fi \textbf{.}}{\hfill\qed\par\bigskip}

 \newcommand{\N}{\mathbb{N}}
\newcommand{\Z}{\mathbb{Z}}
\newcommand{\R}{\mathbb{R}}
\newcommand{\Q}{\mathbb{Q}}

\newcommand{\id}{\mathrm{id}}
\newcommand{\Acal}{\mathcal{A}}
\newcommand{\Bcal}{\mathcal{B}}
\newcommand{\Ccal}{\mathcal{C}}
\newcommand{\Dcal}{\mathcal{D}}
\newcommand{\Ecal}{\mathcal{E}}
\newcommand{\Fcal}{\mathcal{F}}
\newcommand{\Gcal}{\mathcal{G}}
\newcommand{\Hcal}{\mathcal{H}}

\newcommand{\Lcal}{\mathcal{L}}
\newcommand{\Mcal}{\mathcal{M}}

\newcommand{\Pcal}{\mathcal{P}}
\newcommand{\Qcal}{\mathcal{Q}}

\newcommand{\Scal}{\mathcal{S}}
\newcommand{\Tcal}{\mathcal{T}}
\newcommand{\Ucal}{\mathcal{U}}
\newcommand{\Vcal}{\mathcal{V}}
\newcommand{\Wcal}{\mathcal{W}}
\newcommand{\Xcal}{\mathcal{X}}
\newcommand{\Ycal}{\mathcal{Y}}
\newcommand{\Zcal}{\mathcal{Z}}

\newcommand{\kerstyle}[1]{\mathbf{#1}}
\DeclareRobustCommand{\ids}{\textup{id}}

\newcommand{\Gk}{\kerstyle{G}}

\newcommand{\Kk}{\kerstyle{K}}

\newcommand{\Mk}{\kerstyle{M}}

\newcommand{\Pk}{\kerstyle{P}}
\newcommand{\Qk}{\kerstyle{Q}}

\newcommand{\Tk}{\kerstyle{T}}
\newcommand{\Uk}{\kerstyle{U}}
\newcommand{\Vk}{\kerstyle{V}}

\newcommand{\Xk}{\kerstyle{X}}
\newcommand{\Yk}{\kerstyle{Y}}
\newcommand{\Zk}{\kerstyle{Z}}

\newcommand{\prstyle}[1]{\mathsf{#1}}
\newcommand{\PrPf}{\prstyle{P}}
\newcommand{\Prpf}{\prstyle{S}}
\newcommand{\PrpF}{\prstyle{R}}
\newcommand{\PrPF}{\prstyle{K}}
\newcommand{\PrQm}{\prstyle{Q}}
\newcommand{\Giry}{\prstyle{G}}
\newcommand{\PrT}{\prstyle{T}}

\newcommand{\Mbb}{\mathbb{M}}
\newcommand{\Pch}{\mathrm{Pch}}

\newcommand{\ihom}{{\mathcal{QMS}}}

\newcommand{\deltabf}{{\boldsymbol \delta}}

\newcommand{\ubar}[1]{\text{\b{$#1$}}}

\newcommand{\two}{\mathbf{2}}
\newcommand{\twob}{{\ubar{\mathbf{2}}}}
\newcommand{\twoc}{\breve{\mathbf{2}}}
\newcommand{\twop}{{\bar{\mathbf{2}}}}
\newcommand{\one}{\mathbf{1}}
\newcommand{\zero}{\mathbf{0}}

\newcommand{\ismapof}{\precsim}

\newcommand{\sm}{\setminus}
\newcommand{\ins}{\subseteq}
\newcommand{\sni}{\supseteq}
\newcommand{\cmpl}{\mathsf{c}}

\DeclareMathOperator*{\bigand}{\bigwedge}
\DeclareMathOperator*{\bigor}{\bigvee}
\newcommand{\himp}{\mathbin{\Rightarrow}}

\newcommand{\pr}{\mathrm{pr}}
\newcommand{\ev}{\mathrm{ev}}
\newcommand{\Eq}{\mathrm{Eq}}
\newcommand{\CoEq}{\mathrm{CoEq}}
\newcommand{\pf}{\mathrm{pf}}
\newcommand{\incl}{\mathrm{incl}}
\newcommand{\res}{\mathrm{res}}
\newcommand{\Kleisli}{\mathrm{Kl}}

\newcommand{\srj}{\twoheadrightarrow}
\newcommand{\inj}{\hookrightarrow}
\newcommand{\bij}{\stackrel{\sim}{\longrightarrow}}

\newcommand{\patchable}{patchable}
\newcommand{\Patchable}{Patchable}
\newcommand{\patchability}{patchability}
\newcommand{\Patchability}{Patchability}
\newcommand{\sturdy}{sturdy}

\newcommand{\Sturdy}{Sturdy}

\DeclareMathOperator{\Sub}{\mathrm{Sub}}

\newcommand{\I}{\mathbbm{1}}
\newcommand{\Pa}{\mathrm{Pa}}

\newcommand{\Ch}{\mathrm{Ch}}

\newcommand{\Anc}{\mathrm{Anc}}

\newcommand{\Desc}{\mathrm{Desc}}

\DeclareMathOperator*{\Indep}{\perp\!\!\!\perp}
\DeclareMathOperator*{\Perp}{\perp}

\DeclareMathOperator*{\nPerp}{\cancel\perp}
\DeclareMathOperator*{\given}{|}

\newcommand{\dcup}{\,\dot{\cup}\,}
\newcommand{\bigdcup}{\mathop{\dot{\bigcup}}}

\DeclareMathOperator*{\ftimes}{\times}

\newcommand{\cop}{\operatorname{copy}}
\newcommand{\del}{\operatorname{del}}
\newcommand{\swap}{\operatorname{swap}}

\newcommand{\mult}{\operatorname{mult}}
\newcommand{\Inp}{I}

\newcommand{\QBS}{\mathbf{QBS}}
\newcommand{\QUS}{\mathbf{QUS}}
\newcommand{\StdQUS}{\mathbf{StdQUS}}
\newcommand{\SQUS}{\mathbf{SQUS}}
\newcommand{\PQUS}{\mathbf{PQUS}}
\newcommand{\QMS}{\mathbf{QMS}}
\newcommand{\SQMS}{\mathbf{SQMS}}
\newcommand{\PQMS}{\mathbf{PQMS}}
\newcommand{\Meas}{\mathbf{Meas}}
\newcommand{\SMeas}{\mathbf{SMeas}}
\newcommand{\UMeas}{\mathbf{UMeas}}
\newcommand{\Top}{\mathbf{Top}}
\newcommand{\Sets}{\mathbf{Sets}}

\newcommand{\lp}{\left ( }
\newcommand{\rp}{\right ) }

\newcommand{\lC}{\left \{ }
\newcommand{\rC}{\right \} }

\newcommand{\st}{\,\middle|\,}

\newcommand{\arrhead}{{latex}}
\newcommand{\arrtail}{{stealth reversed}}
\newcommand{\arrstar}{{Rays[n=6]}}

\newcommand*{\hut}[1][]{\mathrel{\tikz [baseline=-0.25ex,\arrhead-\arrtail, #1] \draw [#1] (0pt,0.5ex) -- (1.3em,0.5ex);}}

\newcommand*{\tuh}[1][]{\mathrel{\tikz [baseline=-0.25ex,\arrtail-\arrhead, #1] \draw [#1] (0pt,0.5ex) -- (1.3em,0.5ex);}}

\newcommand*{\sus}[1][]{\mathrel{\tikz [baseline=-0.25ex,\arrstar-\arrstar, #1] \draw [#1] (0pt,0.5ex) -- (1.3em,0.5ex);}}

\let\oldFootnote\footnote
\newcommand\nextToken\relax

\renewcommand\footnote[1]{\oldFootnote{#1}\futurelet\nextToken\isFootnote}

\newcommand\isFootnote{\ifx\footnote\nextToken\textsuperscript{,}\fi}

\providecommand{\doi}[1]{\href{https://doi.org/#1}{\textsf{doi:#1}}}
\providecommand{\eprintlink}[1]{\href{https://arxiv.org/abs/#1}{\textsf{arXiv:#1}}}
\providecommand{\urllink}[1]{\href{#1}{\textsf{#1}}}

\begin{document}

\begin{titlepage}

\title{\bfseries Quasi-Measurable Spaces\\[2mm]
\large A Convenient Foundation of Probability Theory}
\author{Patrick Forr\'e\\[+10pt]
    \small{AI4Science Lab}\\[-5pt]
 \small{Korteweg-de Vries Institute for Mathematics}\\[-5pt]
    \small{University of Amsterdam}}%
\date{}\maketitle
\footnotetext{\texttt{p.d.forre@uva.nl}}
\begin{abstract}
We introduce the categories of quasi-measurable spaces, mild generalizations of quasi-Borel spaces that allow for general sample spaces and less restrictive random variables. Each such category is complete, cocomplete, cartesian closed and, as we show, even locally cartesian closed, a \emph{quasitopos} and a \emph{Heyting category}, and so models an extensional dependent type theory with a typed intuitionistic first-order logic. We put one probability monad in the foreground, the \emph{push-forward probability monad} $\PrPf$, and isolate the two properties of the sample space its theory needs: a \emph{Fubini property}, making $\PrPf$ commutative, and a \emph{product isomorphism} $\Omega \times \Omega \cong \Omega$, making it strong. It is affine always.

    Our main interest is the category of \emph{quasi-universal spaces}, where the sample space is the \emph{universal Hilbert cube} $[0,1]^\N$ with the $\sigma$-algebra of all universally measurable subsets. That completion is what buys explicit descriptions of the induced $\sigma$-algebras --- as intersections of Lebesgue-complete $\sigma$-algebras, and for countably separated spaces as the completion of any countable separating family --- which quasi-Borel spaces do not admit. All push-forward probability monads collapse into one here, and on \emph{standard} quasi-universal spaces the Giry monad joins them. We prove a Fubini theorem, a disintegration theorem for Markov kernels, Kolmogorov extension theorems and a conditional de Finetti theorem, and the properties of the associated Markov category.
    
    Finally we formalize causal Bayesian networks over quasi-universal spaces and prove a global Markov property, translating transitional conditional independence into this setting and establishing its asymmetric separoid rules. Exponential objects then put variables and causal mechanisms on the same graphical footing, in \emph{partially generic} causal Bayesian networks.
\end{abstract}


\DeclareTOCStyleEntries[numwidth=2.9em]{tocline}{subsection}
\DeclareTOCStyleEntries[numwidth=3.6em]{tocline}{subsubsection}
{\small \tableofcontents}
\thispagestyle{empty}
\end{titlepage}

\pagestyle{headings}

\pagenumbering{arabic}
\pagestyle{headings}

\section{Introduction}
\label{sec:intro}

\subsection{Motivation}

Even though mathematicians have been on the quest for finding the right framework for pursuing topology for more than 50 years and been rewarded with the discovery of the \emph{convenient category} of compactly generated weakly Hausdorff spaces, quasi-topologies and the category of condensed sets, see \cite{Spa63,Ste67,McC69,EH99,EH02,Str09,SC19,SC20}, the same endeavour for pursuing \emph{higher probability theory} just started rather recently. However, by mirroring the found principles for those categories of topological spaces, a major milestone was quickly achieved by the construction of the category of \emph{quasi-Borel spaces} $\QBS$, see \cite{Heu17,SKV17,SSSW21}.

The main motivation for the latter was to enable probabilistic programming to incorporate higher-order functions, see
\cite{Sta16,Heu17}. Even though the probability monad of Giry, see \cite{Law62,Gir82}, is strong and provides a semantic foundation for first-order probabilistic programming language, see \cite{Sta16}, it does not  support higher-order functions, because the category of measurable spaces $\Meas$ is not cartesian closed, see \cite{Aum61}. The category of quasi-Borel spaces was then constructed to allow for such higher-order functions and was shown to be cartesian closed, see \cite{Heu17}, which thus solved those problems.

Since the lack of exponential objects inside the category of measurable spaces $\Meas$ also provides obstacles and complications in other application domains like probabilistic graphical models and causal reasoning, see below, the author of the paper at hand was motivated to further investigate those topics in slightly more generality under the umbrella term of the categories of \emph{quasi-measurable spaces}.
\begin{figure}[ht]
    \centering
    \begin{tikzpicture}[scale=1, transform shape]
        \tikzstyle{every node} = [draw,shape=circle,color=blue]
        \node (X) at (0,0) {$X$};
        \node (Y) at (4,0) {$Y$};
        \node (Z) at (8,0) {$Z$};
   \foreach \from/\to in {X/Y, Y/Z}
    \draw[-{Latex[length=3mm,width=2mm]}, color=blue] (\from) -- (\to);
     \end{tikzpicture}
     \caption{A Markov chain with output variables $X$, $Y$, $Z$.
         The graph allows us to read off the conditional independence:
         $Z \Indep X \given Y$.}
    \label{fig:markov-chain}
\end{figure}

To motivate further, consider a Markov chain with three variables $X$, $Y$, $Z$, see Figure \ref{fig:markov-chain}. Their joint distribution is then given by the factorization:
\begin{align*} P(X,Y,Z) = P(Z|Y) \otimes P(Y|X) \otimes P(X). \tag{I} \label{eq:markov-chain}
\end{align*}
This factorization then shows the conditional independence:
\[ Z \Indep X\given Y.\]
One could have as well read this off the graph in Figure \ref{fig:markov-chain} via the \ids-separation statement $Z \Perp X \given Y$ and the
global Markov property, see e.g.\ \cite{Lau90,Lau96,Ric03,KF09,Pearl09,FM17,FM18,FM20,For21}.
$Z \Indep X\given Y$ means that when the value of $Y$ is provided then $X$ has no additional information about the state of $Z$.
From the factorization in Equation \ref{eq:markov-chain} we can even read off
that the value of $Z$ is only determined through the Markov kernel $P(Z|Y)$ when $Y$ is provided.
So the distribution $P(X)$ of $X$ and the Markov kernel $P(Y|X)$ of $Y$ both have no influence on $Z$ when conditioned on $Y$. So in some sense $Z$ is independent not only of $X$, but also of $P(X)$ and $P(Y|X)$
when conditioned on $Y$.
If we use $Q(X)$, $Q(Y|X)$ as variables for the values $P(X)$, $P(Y|X)$, resp., then
we would like to be able, first, to formalize a conditional independence statement of form:
    \[ Z \Indep \lp X,Q(X),Q(Y|X) \rp \given Y,\]
and, second, to read this off a graph like in Figure \ref{fig:partially-generic-markov-chain} via some global Markov property. This comes with several challenges.

First, the variables $Q(X)$ and $Q(Y|X)$ are not random variables in the usual sense as there is no distribution defined over them. Conditional independence statements mixing stochastic and non-stochastic variables have a history in statistics: the decision-theoretic treatment of causal inference in \cite{Daw02b} puts non-random regime and intervention nodes into a graph, and \emph{extended conditional independence} of \cite{CD17b} is a calculus of conditional independence for exactly such a mixture, with separoid-style axioms. The problem was solved in the generality we need in \cite{For21}, in the category of measurable spaces $\Meas$, by introducing the notion of \emph{transitional conditional independence}, which is the notion we transport to $\QUS$ below; \cite{For21} Section 1 positions it against \cite{CD17b}.

\begin{figure}[ht]
    \centering
    \begin{tikzpicture}[scale=1, transform shape]
      \tikzstyle{every node} = [draw,shape=rectangle,color=teal, minimum size = 0.75cm]
      \node (PX) at (-1,2) {$Q(X)$};
      \node (PY) at (3,2) {$Q(Y|X)$};
        \tikzstyle{every node} = [draw,shape=circle,color=blue]
        \node (X) at (0,0) {$X$};
        \node (Y) at (4,0) {$Y$};
        \node (Z) at (8,0) {$Z$};
   \foreach \from/\to in {PX/X, PY/Y}
    \draw[-{Latex[length=3mm,width=2mm]}, color=teal] (\from) -- (\to);
   \foreach \from/\to in {X/Y, Y/Z}
    \draw[-{Latex[length=3mm,width=2mm]}, color=blue] (\from) -- (\to);
     \end{tikzpicture}
     \caption{The partially generic causal Bayesian network with output variables $X$ in $\Xcal$ and $Y$ in $\Ycal$ and $Z$ in $\Zcal$
         and input variables $Q(X)$ in $\PrPf(\Xcal)$ and $Q(Y|X)$ in $\PrPf(\Ycal)^\Xcal$.
         The graph allows us to read off the transitional conditional independence:
         $Z \Indep \lp X,Q(X),Q(Y|X)\rp \given Y$.}
    \label{fig:partially-generic-markov-chain}
\end{figure}

Second, for a global Markov property to hold we need to incorporate variables $Q(Y|X)$ into a (causal) Bayesian network.
For this, note that $X$ and $Y$ take values in a measurable spaces $(\Xcal,\Bcal_\Xcal)$ and $(\Ycal,\Bcal_\Ycal)$, resp., while $Q(Y|X)$ takes values in the space $\Lcal$ of Markov kernels from $\Xcal$ to $\Ycal$:
\[\Lcal:=\Meas\lp(\Xcal,\Bcal_\Xcal),(\Giry(\Ycal,\Bcal_\Ycal),\Bcal_{\Giry(\Ycal,\Bcal_\Ycal)})\rp,\]
where $\Giry(\Ycal,\Bcal_\Ycal)$ is the space of all probability measures on $(\Ycal,\Bcal_\Ycal)$ endowed with the smallest $\sigma$-algebra $\Bcal_{\Giry(\Ycal,\Bcal_\Ycal)}$ that makes all evaluation maps $P \mapsto P(B)$, $B \in \Bcal_\Ycal$, measurable.
So, to incorporate $Q(Y|X)$ into a Bayesian network, we would need to construct a (measurable) Markov kernel from the product space $\Xcal \times \Lcal$ to $\Ycal$.
More precisely, we would need to construct a $\sigma$-algebra $\Bcal_\Lcal$ on $\Lcal$ such that the following map is measurable:
\[ \ev:\;  \lp \Xcal\times\Lcal, \Bcal_\Xcal\otimes\Bcal_\Lcal\rp \to  \lp\Giry(\Ycal,\Bcal_\Ycal),\Bcal_{\Giry(\Ycal,\Bcal_\Ycal)}\rp,
\qquad (x,Q(Y|X)) \mapsto Q(Y|X=x). \]
Unfortunately and related to the fact that the category of measurable spaces $\Meas$ is not cartesian closed, in general, there does not exist any $\sigma$-algebra $\Bcal_\Lcal$ on $\Lcal$ that makes the above map $\ev$
measurable, where $\Bcal_\Xcal \otimes \Bcal_\Lcal$ is the product $\sigma$-algebra, which is generated by the cylinder sets $A \times D$, with $A \in \Bcal_\Xcal$ and $D \in \Bcal_\Lcal$; not even in the case where the measurable spaces $\Xcal$ and $\Ycal$ are (uncountable) standard Borel measurable spaces, see \cite{Aum61}.

To remedy this shortcoming we engage in studying quasi-Borel spaces and, slightly more general, \emph{quasi-measurable spaces}.
The main idea behind quasi-measurable spaces can best be described when contrasted with
the classical setup of measure-theoretic probability theory:

In the \emph{classical setup of measure-theoretic probability theory} one usually starts with the
assumption of the presence of a (typically not further specified)
measurable space $(\Omega,\Bcal_\Omega)$, the \emph{sample space}, and a probability measure $P$ on it.
Then one considers the \emph{data spaces} $\Xcal$, typically $\R$ or $\R^D$, where the real world measurements happen. One then specifies a $\sigma$-algebra $\Bcal_\Xcal$ of subsets of $\Xcal$ as the collection of all admissible \emph{events} on $\Xcal$. When $\Bcal_\Xcal$ is given then the admissible \emph{random variables}
on $\Xcal$ that are used to describe the data points or random measurements on $\Xcal$
are all measurable maps $X$ from $(\Omega,\Bcal_\Omega)$ to $(\Xcal,\Bcal_\Xcal)$.
The set of all admissible random variables is then:
\[\Xcal^\Omega:=\Fcal(\Bcal_\Xcal):= \Meas\lp(\Omega,\Bcal_\Omega),(\Xcal,\Bcal_\Xcal)\rp. \]
Note that the set of random variables is dependent on the choice of the $\sigma$-algebra.
So if we change the $\sigma$-algebra $\Bcal_\Xcal$ we implicitly change the set of admissible random variables $\Xcal^\Omega$ with it.

In the \emph{theory of quasi-measurable spaces} we swap the roles of the set of admissible events $\Bcal_\Xcal$ and set of admissible random variables $\Xcal^\Omega$.
We again start with a sample space $(\Omega,\Bcal_\Omega)$, but then directly specify the admissible set of random variables $\Xcal^\Omega$
on the  data space $\Xcal$.
This then \emph{induces} a $\sigma$-algebra of admissible events on $\Xcal$ via:
\[ \Bcal_\Xcal:=\Bcal(\Xcal^\Omega):= \lC A \ins \Xcal \st \forall X \in \Xcal^\Omega.\, X^{-1}(A) \in \Bcal_\Omega \rC.\]
Note that if we now change the set of admissible random variables $\Xcal^\Omega$ we would implicitly change the $\sigma$-algebra $\Bcal_\Xcal$. Also note that the $\sigma$-algebra $\Bcal_\Xcal$ is now
dependent on the choice of $(\Omega,\Bcal_\Omega)$ and the set of random variables $\Xcal^\Omega$.

The tuple $(\Xcal,\Xcal^\Omega)$ will then be called a \emph{quasi-measurable space w.r.t.\ $(\Omega,\Bcal_\Omega)$}. The properties of the corresponding category of quasi-measurable spaces $\QMS$ will then be dependent on the choice of the sample space $(\Omega,\Bcal_\Omega)$,
which thus needs to be specified and studied more closely, in contrast to the classical setup.

If we are now coming back to our example of Markov chains we will need to
construct/formalize probabilistic graphical models in this setting, here \emph{causal Bayesian networks}.
 Furthermore, we need to study
\emph{probability monads} and \emph{Markov kernels} between quasi-measurable spaces.
Also, to translate the mentioned notion of \emph{transitional conditional independence} from the category of measurable spaces $\Meas$ to the category of quasi-measurable spaces $\QMS$ we need to follow the steps in \cite{For21} and prove a \emph{disintegration theorem} for Markov kernels in $\QMS$.

For such disintegration results one usually needs close control over the $\sigma$-algebras of the involved spaces, see \cite{For21} Cor.\ C.8. This is the point where the quasi-Borel spaces from \cite{Heu17} become rather inconvenient to work with. The reason is that the induced $\sigma$-algebras come, as described above, as push-forward $\sigma$-algebras of the sample space along the quasi-measurable functions. This roughly means that the induced $\sigma$-algebras will be pushed onto the images of such maps and then extended by null-sets outside, without necessarily including null-sets on that image. So the induced $\sigma$-algebras will, in general and in some sense, be partially complete and partially incomplete. To overcome this inconsistency and harmonize the situation it seems natural to use sample spaces that are complete with respect to all relevant probability measures. Since, furthermore, for the disintegration theorem one requires
\emph{perfect} probability measures and \emph{universally countably generated} $\sigma$-algebras, see \cite{Fad85,Pac78,Fremlin,For21}, one then almost necessarily\footnote{See Lemma \ref{deflem:univ-meas}.} ends up requiring a sample space that is an uncountable Polish space endowed with its $\sigma$-algebra of universally measurable subsets, which then satisfies all those constraints and requirements.
The corresponding category of quasi-measurable spaces with such a sample space will be called the category of \emph{quasi-universal spaces} $\QUS$. These are the reasons we strongly engage in developing the theory for quasi-universal spaces in this paper.

\subsection{Contributions}

Closely following the construction of quasi-Borel spaces $\QBS$, see \cite{Heu17}, we introduce the categories of \emph{quasi-measurable spaces} $\QMS$. For this we mainly make two slight changes:

First, we set up the theory of quasi-measurable spaces to allow for different sample spaces $\Omega$, see \Cref{qms:def:sample-space-act-1}, other than the real numbers $\R$ with its Borel $\sigma$-algebra.
Since the properties of the resulting category depend on that choice, the sample space is not introduced in one go, but \emph{act by act}, each act adding exactly the structure that the next batch of results needs: the set of admissible random variables $\Omega^\Omega$ in \Cref{qms:def:sample-space-act-1}, a compatible $\sigma$-algebra $\Bcal_\Omega$ in \Cref{qms:def:sample-space-act-2}, a set of product-compatible probability measures $\PrPf$ in \Cref{qms:def:sample-space-act-3}, and finally one concrete choice of all three --- the \emph{universal Hilbert cube} --- in \Cref{qms:def:sample-space-act-4}. This staging makes visible which theorem needs which assumption, and it is what allows us in \Cref{sec:sample-spaces} to \emph{construct} sample spaces with prescribed good behaviour rather than to postulate them.

Second, for \Cref{qms:def:quasi-measurable-spaces} of quasi-measurable spaces we drop the 3rd property in the definition of quasi-Borel spaces, see \cite{Heu17} Def.\ 7. This simplifies and generalizes the theory of quasi-measurable spaces $\QMS$ at the same time. For our purposes, we (almost) lose none of the good properties of $\QBS$, see below:

The categories of quasi-measurable spaces $\QMS$ are shown to have all (small) limits and colimits, see \Cref{qms:thm:co-complete}. Arbitrary (small) products and coproducts distribute, see \Cref{qms:thm:coprod-prod-distr}.
They are \emph{cartesian closed}, see \Cref{qms:thm:function-space} and \Cref{qms:cor:qms-cartesian-closed}, which is what allows higher-order functions and a simply typed $\lambda$-calculus, and which is where the categorical development in \cite{Heu17} stops. We take it three steps further, because each step buys something that a probabilistic language can use.

\medskip
\noindent\textbf{More than cartesian closed.}
First, $\QMS$ is \emph{locally} cartesian closed, see \Cref{qms:thm:lcc}: every slice category $\QMS_\Tcal$ is cartesian closed, the fibre product $-\ftimes_\Tcal-$ having an internal hom as its right adjoint. Read type-theoretically, a quasi-measurable map $\Xcal \to \Tcal$ is a \emph{family} of spaces indexed by $\Tcal$, dependent sums and dependent products of such families exist and are stable under substitution, and $\QMS$ is therefore a model of an \emph{extensional dependent type theory}, see \cite{See84} and \cite{Joh02} A1.5.3. This is what lets the \emph{type} of a term depend on a value that has itself been sampled --- a distribution over sequences whose length is random, a model whose parameter space is drawn first --- rather than only its value.

Second, the indiscrete two-point space $\twop$ classifies embeddings, i.e.\ strong monomorphisms, and together with completeness, cocompleteness and local cartesian closedness this makes $\QMS$ a \emph{quasitopos}, see \Cref{qms:thm:quasitopos}. It is not a topos: $\twop$ does not classify arbitrary monomorphisms, and the one-point space $\one$ is a generator, so $\QMS$ behaves like a category of concrete structured sets rather than like a category of sheaves.

Third, the subspaces $\Sub(\Xcal)$ of a quasi-measurable space form a complete \emph{Heyting algebra}, the assignment $\Xcal \mapsto \Sub(\Xcal)$ is functorial, and every pullback $f^*:\,\Sub(\Ycal) \to \Sub(\Xcal)$ is a homomorphism of complete Heyting algebras with a left adjoint $\exists_f$ and a right adjoint $\forall_f$. So $\QMS$ is a \emph{Heyting category}, see \Cref{qms:thm:heyting}, and carries the internal logic of a typed \emph{intuitionistic first-order logic}: predicates over a space, implication, and quantification along any quasi-measurable map, all interpreted by subspaces.

It is worth being precise about what is new here, and the honest answer is: not the theorems, but the treatment. A quasi-measurable space is exactly a \emph{concrete presheaf} on a two-object site built from $\Omega$, and concrete sheaves on a concrete site form a Grothendieck quasitopos --- complete, cocomplete, locally cartesian closed --- by \cite{BH11} Thm.\ 5.25. That one citation subsumes all six results just quoted, for every sample space; we spell the identification out in \Cref{qms:rem:concrete-presheaf}. The same observation has been made for quasi-Borel spaces, which are the concrete sheaves on standard Borel spaces, see \cite{VO18} Section 11 and \cite{MMS22}, so $\QBS$ carries this structure too; and \cite{AKM26} calls the dependent-type reading of $\QBS$ ``well-known''. What we contribute is a direct and elementary treatment, aimed at readers who want to compute with the objects rather than cite their existence: it holds for an \emph{arbitrary} sample space and for the weaker \Cref{qms:def:quasi-measurable-spaces}; it produces the function space, the internal hom, the subspace lattice with its Heyting implication and both quantifier adjoints, and a coproduct that is more rigid than the one in $\QBS$, none of which the general theorem exhibits; and it lives in the same category as the probability monad and the disintegration theory below, so that the logic and the measure theory are available at once, on the same objects. The genuinely new mathematics of this paper is elsewhere --- in the universal completion and what it buys, see below.

\medskip
We also maintain the adjunction between the category of measurable spaces $\Meas$ and each of the categories of quasi-measurable spaces $\QMS$, see \Cref{thm:qms-meas-adjunction}. The fixed points of that adjunction, which we call \emph{\sturdy} spaces, form the common part of $\Meas$ and $\QMS$, see \Cref{qms:cor:sturdy-meas}.

By not insisting on the mentioned 3rd property that quasi-Borel spaces have, \cite{Heu17} Def.\ 7, the right-adjoint $\Fcal$ of that adjunction does not preserve countable coproducts of measurable spaces anymore. Quasi-measurable spaces with that property are studied on their own in \Cref{sec:patchable-qms} under the name \emph{\patchable}. We show in \Cref{thm:patchable-reflexive} that the category of \patchable\ quasi-measurable spaces $\PQMS$ is a reflective subcategory of $\QMS$, and in fact an exponential ideal, so that one may reflect back to $\PQMS$ whenever countable coproducts are needed, see \Cref{lem:patch-coprod}; we also study which constructions preserve \patchability. Since none of the main development depends on it, this material is placed in the appendix.

\medskip
\noindent\textbf{One monad in the foreground.}
A quasi-measurable space carries several candidate spaces of probability measures, and correspondingly several probability monads: the \emph{quasi-measurable probability monad} $\PrQm$, whose members are the measures that are quasi-measurable as maps on the space of events, and four monads of push-forward probability measures $\PrPF$, $\PrpF$, $\PrPf$, $\Prpf$, which are variations of the monad on $\QBS$ from \cite{Heu17} --- one more general ($\PrPF$), one more restrictive ($\Prpf$), and one complementary ($\PrPf$) to $\PrpF$. Rather than carrying all five through the paper, we put \emph{one} of them in the foreground: the \emph{push-forward probability monad} $\PrPf$ of \Cref{sec:pf-monad}. The other four are developed in full generality, in the same notation, in \Cref{qms:sec:app-other-monads}, where \Cref{qms:tab:monads} compares them side by side: they are distinguished purely by which ingredient of a push-forward $P(X)=X_*P$ --- the random variable $X$, the measure $P$, both, or neither --- may depend on the sample point. They are nested, see \Cref{lem:S-P-R-K-incl}, they collapse pairwise under mild hypotheses on the sample space, see \Cref{lem:K-P-S-R-eq}, and on quasi-universal spaces all four coincide, see \Cref{thm:qus-S=P-in-R}, so nothing is lost by the choice.

For $\PrPf$ we isolate the two assumptions that actually do the work, and state them explicitly rather than leaving them implicit in the choice of $\Omega$: the \emph{Fubini property}, \Cref{qms:asm:fubini}, and the \emph{product isomorphism property} $\Omega \times \Omega \cong \Omega$, \Cref{qms:asm:prod-iso}.
The first yields a Fubini theorem valid for \emph{all} quasi-measurable spaces and all non-negative functions, see \Cref{qms:thm:fubini}, together with a workable sufficient criterion, \Cref{qms:prp:fubini-crit}; the second yields the product of Markov kernels, \Cref{qms:thm:product-PrPf}, and with it the strong monad structure, \Cref{qms:thm:PrPf-monad}. This is what makes $\QMS$ support semantics for a probabilistic programming language in the monadic style with higher-order functions, similar to how the category of quasi-Borel spaces does, see \cite{Mog91,Heu17,SKV17,SSSW21}.

\Cref{sec:sample-spaces} then answers the question of where such sample spaces come from. We give two general construction theorems --- one starting from a measurable space with $(\Omega,\Bcal_\Omega) \cong (\Omega\times\Omega,\Bcal_\Omega\otimes\Bcal_\Omega)$, see \Cref{qms:thm:spl-spc-gen}, and one that additionally completes along a chosen class of probability measures, see \Cref{qms:thm:spl-spc-cmpl} --- and show that countably separated product spaces with perfect or countably compact probability measures are covered by them, see \Cref{qms:cor:sample-space-count-sep-perf}. Control over the induced $\sigma$-algebras comes from \emph{countable separatedness} together with \emph{perfectness}, see \Cref{qms:thm:count-sep-perf} and \Cref{qms:cor:count-sep-perf}.

\medskip
\noindent\textbf{The universal Hilbert cube.}
In \Cref{sec:qus} we commit to one concrete sample space and see what the theory looks like once we do. The fourth act, \Cref{qms:def:sample-space-act-4}, takes $\Omega$ to be the \emph{universal Hilbert cube} $[0,1]^\N$ endowed with the $\sigma$-algebra of all universally measurable subsets; \Cref{qms:thm:act4-well-behaved} verifies that this is an instance of the construction of \Cref{qms:thm:spl-spc-cmpl}, with the ``zip-locker'' map supplying $\Omega \cong \Omega \times \Omega$ elementarily. Up to Borel isomorphism this is the same as $\R^\N$ or $\R$ with their universally measurable subsets; the resulting category is the \emph{category of quasi-universal spaces} $\QUS$, and it differs from $\QBS$ of \cite{Heu17} exactly in the use of the universal completion.

That difference pays off. The induced $\sigma$-algebras $\Bcal_\Xcal$ now have a description as an intersection of Lebesgue-complete $\sigma$-algebras, see \Cref{lem:qus-universally-complete}, which becomes most pronounced for \emph{countably separated} quasi-universal spaces, see \Cref{thm:countably-separated}; giving such a description was not possible for quasi-Borel spaces. On $\QUS$ the four push-forward monads agree, see \Cref{thm:qus-S=P-in-R}, and on \emph{standard} quasi-universal spaces even $\PrQm$ and the Giry monad join them, see \Cref{thm:univ-qus-GPQKR}. We then prove a \emph{Fubini Theorem} \Cref{thm:fubini}, a \emph{Disintegration Theorem} \Cref{thm:disint-X} for Markov kernels, \emph{Kolmogorov Extension Theorems} \Cref{thm:kol-ext-univ} and \Cref{thm:kol-ext-cs}, and a \emph{Conditional De Finetti Theorem} \Cref{cor:de-finetti}.

Quasi-universal spaces with properties similar to those of standard Borel spaces are called \emph{standard quasi-universal spaces}, with category $\StdQUS$; they are exactly the retracts of the sample space $\Omega$ in $\QUS$, see \Cref{def:univ-qus}. We express the resulting structure in the language of \emph{categorical probability theory}, see \Cref{thm:Markov-category}:

    The category of standard quasi-universal spaces $(\StdQUS,\times,\one)$ is a
    symmetric cartesian (but not closed) monoidal category  with countable products and countable coproducts, the latter computed in $\PQUS$ rather than in $\QUS$. The triple $(\PrPf,\delta,\Mbb)$ is a strong affine symmetric monoidal/commutative monad on $(\StdQUS,\times,\one)$. Its Kleisli category $\Kleisli(\PrPf)$ on $\StdQUS$ is an a.s.-compatibly representable Markov category with conditionals and Kolmogorov powers, see \cite{Kle65,Cho17,Fri19,Fri20,FriG20,Fri21}; these notions are recalled in \Cref{qms:sec:app-markov-cat}.

In \Cref{sec:causal-models} we study conditional independence relations and causal models in $\QUS$.
There we translate the notion of \emph{transitional conditional independence} in \cite{For21}
from $\Meas$ to $\QUS$, see \Cref{def:tr-cond-ind}, and show all the corresponding
\emph{(asymmetric) separoid rules} in \Cref{thm:separoid_axioms-tci}.
These are then used to prove the \emph{global Markov property} for \emph{causal Bayesian networks},
see \Cref{thm:gmp-mI-CBN}, which relates the underlying graphical structure to the transitional conditional independencies of its variables via a \ids-separation criterion.

We then introduce the notion of \emph{partially generic causal Bayesian networks}, see \Cref{def:gen-cbn}.  Finally, with this and together with the strong probability monad $\PrPf$,  the exponential objects
$\PrPf(\Xcal)^\Zcal$ in $\QUS$, and the global Markov property we now arrive at our goal and can formulate conditional independence relations between variables and causal mechanisms on equal footing and reason with them graphically as wanted and explained in the beginning with the example of Markov chains, see \Cref{fig:partially-generic-markov-chain}.

\medskip
\noindent\textbf{Structure of this paper.}
\Cref{sec:qms} introduces quasi-measurable spaces over an abstract sample space and establishes the categorical structure just described. \Cref{sec:meas-vs-qms} adds the $\sigma$-algebra to the sample space and sets up the adjunction with $\Meas$, together with the induced $\sigma$-algebras that everything measure-theoretic afterwards is phrased in. \Cref{sec:giry-completions} recalls the Giry monad and completions of $\sigma$-algebras and assembles from them the quasi-measurable probability monad $\PrQm$. \Cref{sec:pf-monad} adds the probability measures to the sample space and develops the push-forward probability monad $\PrPf$ together with the Fubini theorem and the product of Markov kernels. \Cref{sec:sample-spaces} answers where such sample spaces come from. \Cref{sec:qus} fixes the universal Hilbert cube once and for all and proves the disintegration, Kolmogorov extension and conditional de Finetti theorems for quasi-universal spaces. \Cref{sec:causal-models} builds causal Bayesian networks in that setting and proves the global Markov property. \Cref{sec:discussion} summarizes and discusses what changes if one chooses a different sample space.

A reader who wants the categorical structure alone can stop after \Cref{sec:qms}; one who wants the probability theory but not the causal models can stop after \Cref{sec:qus}; one who is here for the graphical models will want \Cref{sec:pf-monad}, \Cref{sec:qus} and \Cref{sec:causal-models}.

Proofs of the results stated in the main text, together with the more technical lemmas they rest on, are collected in the appendix; each such statement carries a pointer to its proof, and each proof carries a pointer back.

\section{The Category of Quasi-Measurable Spaces}
\label{sec:qms}

In this section we will construct the \emph{category of quasi-measurable spaces} $\QMS$ w.r.t.\ a fixed sample space $\Omega$. These construction closely follow the constructions of quasi-Borel spaces from \cite{Heu17},
 also see \cite{SKV17,SSSW21}, but where we now allow for different sample spaces $\Omega$ other than the set of real numbers $\R$. We will have products, coproducts, limits, colimits and function spaces/exponential objects, etc., in $\QMS$.

In the next subsection we will go through the basic definitions of the sample space $\Omega$, quasi-measurable spaces and quasi-measurable maps. This will then constitute a category, the category of
\emph{quasi-measurable spaces} $\QMS$ w.r.t.\ $\Omega$.

\subsection{The Sample Space - Act 1 - The Random Variables}

\begin{Def}[The sample space - act 1]
    \label{qms:def:sample-space-act-1}
    In the following we fix a \emph{sample space} $(\Omega,\Omega^\Omega)$, which consists of:
    \begin{enumerate*}
        \item a non-empty set $\Omega$,
        \item a set $\Omega^\Omega$ of maps $\Phi:\,\Omega \to \Omega$ that contains:
            \begin{enumerate*}
                \item the identity map $\id_\Omega:\,\Omega \to \Omega$,
                \item all constant maps,
                \item all compositions $\Phi_1\circ \Phi_2$ for all $\Phi_1,\Phi_2 \in \Omega^\Omega$.
            \end{enumerate*}
    \end{enumerate*}
\end{Def}

The archetypical example of a sample space that follows \Cref{qms:def:sample-space-act-1} is the following:

\begin{Eg}
    \label{qms:eg:meas-sample-space}
    If $\Omega$ is a set endowed with any set $\Ecal_\Omega$ of subsets of $\Omega$ such that $\Omega, \emptyset \in \Ecal_\Omega$,
    e.g.\ a measurable space with $\sigma$-algebra $\Ecal_\Omega$ or a topological space with topology $\Ecal_\Omega$, then the following choice:
    \begin{align}
        \Omega^\Omega &:= \lC \Phi:\, \Omega \to \Omega \st \forall D \in \Ecal_\Omega.\, \Phi^{-1}(D) \in \Ecal_\Omega \rC,
    \end{align}
    turns $(\Omega,\Omega^\Omega)$ into a valid sample space according to \Cref{qms:def:sample-space-act-1}.
\end{Eg}

\subsection{Quasi-Measurable Spaces and Maps}

\begin{Def}[Quasi-measurable spaces]
    \label{qms:def:quasi-measurable-spaces}
    A \emph{quasi-measurable space} $(\Xcal,\Xcal^\Omega)$ (w.r.t.\ $(\Omega,\Omega^\Omega)$)  consists of:
   \begin{enumerate*}
        \item a set $\Xcal$,
        \item a set $\Xcal^\Omega$ of maps $X:\,\Omega \to \Xcal$ that contains:
            \begin{enumerate*}
                \item all constant maps,
                \item all compositions $X\circ \Phi$ for all $X \in \Xcal^\Omega$, $\Phi \in \Omega^\Omega$.
            \end{enumerate*}
   \end{enumerate*}
\end{Def}

For brevity and by abuse of notation we will often just refer to $\Xcal$ as a quasi-measurable space and implicitly assume that it is endowed with the set of maps $\Xcal^\Omega$, which we might refer to as the set of \emph{(admissible) random variables} $X:\, \Omega \to \Xcal$ on $\Xcal$.

\begin{Def}[Quasi-measurable maps]
    A \emph{quasi-measurable map} $g:\,(\Xcal,\Xcal^\Omega) \to (\Ycal,\Ycal^\Omega)$ between quasi-measurable spaces $(\Xcal,\Xcal^\Omega)$ and $(\Ycal,\Ycal^\Omega)$ is a map $g:\, \Xcal \to \Ycal$ such that for all $X \in \Xcal^\Omega$ we have that: $g(X):=g \circ X \in \Ycal^\Omega$.
\end{Def}

We denote the set of all quasi-measurable maps w.r.t.\ the sample space $(\Omega,\Omega^\Omega)$
as $\QMS\lp (\Xcal,\Xcal^\Omega),(\Ycal,\Ycal^\Omega)\rp$, or just $\QMS\lp\Xcal,\Ycal\rp$ for brevity. We then have the identification: $\Xcal^\Omega=\QMS(\Omega,\Xcal)$. One also immediately checks that the composition of quasi-measurable maps is again a quasi-measurable map. The following is then immediate:

\begin{DefLem}[The category of quasi-measurable spaces]
    The class of all quasi-measurable spaces and quasi-measurable maps
    w.r.t.\ a fixed sample space $(\Omega,\Omega^\Omega)$ forms a category,
    the \emph{category of quasi-measurable spaces} $\QMS$ with sample space $(\Omega,\Omega^\Omega)$.
\end{DefLem}

\begin{Eg}[Measurable spaces as quasi-measurable spaces]
       Let $(\Omega,\Bcal_\Omega)$, $(\Xcal,\Bcal_\Xcal)$ and $(\Ycal,\Bcal_\Ycal)$ be measurable spaces.
            Put:
            \begin{align}
                \Omega^\Omega &:= \Meas\lp(\Omega,\Bcal_\Omega),(\Omega,\Bcal_\Omega)\rp,\\
                \Xcal^\Omega &:=\Meas\lp(\Omega,\Bcal_\Omega),(\Xcal,\Bcal_\Xcal)\rp, \\
                \Ycal^\Omega &:=\Meas\lp(\Omega,\Bcal_\Omega),(\Ycal,\Bcal_\Ycal)\rp.
            \end{align}
        Then $(\Xcal,\Xcal^\Omega)$ and $(\Ycal,\Ycal^\Omega)$ are quasi-measurable spaces and every measurable map:
        $g:\, (\Xcal,\Bcal_\Xcal) \to (\Ycal,\Bcal_\Ycal)$ is also a quasi-measurable map $g:\, (\Xcal,\Xcal^\Omega) \to (\Ycal,\Ycal^\Omega)$. Thus:
        \begin{align}
           \Meas\lp(\Xcal,\Bcal_\Xcal), (\Ycal,\Bcal_\Ycal) \rp \ins \QMS\lp(\Xcal,\Xcal^\Omega), (\Ycal,\Ycal^\Omega)\rp.
        \end{align}
        The inclusion in the other direction might not hold, which already highlights the first difference between the theory of measurable spaces and quasi-measurable spaces.
\end{Eg}

\begin{Def}[Indiscrete quasi-measurable spaces]
    \label{qms:def:indiscrete-qms}
    A quasi-measurable space $(\Xcal,\Xcal^\Omega)$ is called \emph{indiscrete} if $\Xcal^\Omega = \Sets(\Omega,\Xcal)$, i.e.\
    $\Xcal^\Omega$ consists of all set-theoretic maps from $\Omega$ to $\Xcal$.
\end{Def}

Note that any set $\Xcal$ can be be turned into quasi-measurable space only allowing for constant random variables, or into an \emph{indiscrete} quasi-measurable space allowing all maps to be random variables via putting $\Xcal^\Omega := \Sets(\Omega,\Xcal)$.
Later in \Cref{qms:eg:discrete-space} we also define discrete quasi-measurable spaces.

\begin{Eg}[Important corner cases]
    \begin{enumerate}
        \item The \emph{empty space} $(\zero,\zero^\Omega):=(\emptyset,\emptyset)$ is a quasi-measurable space.
            For every quasi-measurable space $\Xcal$ we have:
            \begin{align}
                \QMS(\zero,\Xcal)=\lC \incl:\,\zero \inj \Xcal \rC,
            \end{align}
            where $\incl$ is the inclusion map.
            So $(\zero,\zero^\Omega)$ is an \emph{initial object} of $\QMS$.
        \item The \emph{one-point space} together with the constant map: $(\one,\one^\Omega):=(\lC 0 \rC, \lC !:\,\Omega \to \one\rC)$
            is a quasi-measurable space. For every quasi-measurable space $\Xcal$ we have:
            \begin{align}
                \QMS(\Xcal,\one)=\lC !:\,\Xcal \to \one \rC,
            \end{align}
            where $!$ is the constant map. So $(\one,\one^\Omega)$ is a \emph{terminal object} of $\QMS$.
        \item The \emph{indiscrete two-point space} is denoted by $(\twop,\twop^\Omega)$ with $\twop:=\lC 0,1 \rC$ and $\twop^\Omega:=\Sets\lp \Omega, \twop\rp$. For every quasi-measurable space $\Xcal$ we have:
            \begin{align}
              \twop^\Xcal := \QMS(\Xcal,\twop)=\Sets(\Xcal,\twop),
            \end{align}
            which can thus be identified with the power set of $\Xcal$.
    \end{enumerate}
\end{Eg}

Note that in \Cref{qms:eg:discrete-space} we also define the discrete two-point space $(\twob,\twob^\Omega)$, which is different from
$(\twop,\twop^\Omega)$.

In the following, we will for any set $\Xcal$ identify $\Xcal^\one$ with the set of all constant maps $g:\,\Omega \to \one \to \Xcal$.

\subsection{Quotients and Subspaces}

\begin{Def}[Quotients and embeddings]
    \label{qms:def:quot-emb}
    A quasi-measurable map $f:\,(\Xcal,\Xcal^\Omega) \to (\Ycal,\Ycal^\Omega)$ is called:
    \begin{enumerate}
\item \emph{embedding} if $f$ is injective and $\Xcal^\Omega = \lC X:\, \Omega \to \Xcal \st f \circ X \in \Ycal^\Omega\rC$.
        \item \emph{quotient} if $f$ is surjective and $\Ycal^\Omega = \lC f \circ X \st X \in \Xcal^\Omega \rC$.
    \end{enumerate}
\end{Def}

\begin{Rem}[Subspace]
    We can consider an injective quasi-measurable map  $f:\,(\Xcal,\Xcal^\Omega) \inj (\Ycal,\Ycal^\Omega)$ as a \emph{subspace} of $(\Ycal,\Ycal^\Omega)$ by identifying $\Xcal$ with the subset $f(\Xcal) \ins \Ycal$, and $\Xcal^\Omega$ with $f \circ \Xcal^\Omega \ins \Ycal^\Omega$ and $f$ with the inclusion map $f(\Xcal) \ins \Ycal$. \emph{Embedded subspaces} are then special subspaces with a maximal subset of admissible random variables relative to $\Ycal$.
\end{Rem}

\begin{Def}[Retracts]
    A quasi-measurable space $(\Xcal,\Xcal^\Omega)$ is called a \emph{retract} of another quasi-measurable space $(\Ycal,\Ycal^\Omega)$
    if there exist quasi-measurable maps
    $i:\, (\Xcal,\Xcal^\Omega) \to (\Ycal,\Ycal^\Omega)$ and $r:\, (\Ycal,\Ycal^\Omega) \to (\Xcal,\Xcal^\Omega)$ such that:
    $ r \circ i = \id_\Xcal$.
\end{Def}

For a retract the map $i$ is then automatically an embedding and $r$ a quotient map.

\begin{Def}[Standard quasi-measurable spaces]
    \label{qms:def:standard-qms}
    A quasi-measurable space $(\Xcal,\Xcal^\Omega)$ is called a \emph{standard quasi-measurable space} if it is a retract of
    the sample space $(\Omega, \Omega^\Omega)$.
    For the sample space of \Cref{qms:def:sample-space-act-4} these are the \emph{standard quasi-universal spaces} of \Cref{def:univ-qus}, which play the role that standard Borel spaces play in classical measure theory.
\end{Def}

More details can be found in \Cref{qms:sec:quot-sub}.

\subsection{The Product Space}

In this subsection we define the categorical \emph{product} of quasi-measurable spaces.
The construction again closely follows \cite{Heu17}.

\begin{Def}[The product space]
    Let $(\Xcal_i,\Xcal_i^\Omega)$, $i\in I$, be a family of quasi-measurable spaces indexed by any set $I$.
    We then define their \emph{product space} as:
    \begin{align}
        \prod_{i \in I} \lp \Xcal_i,\Xcal_i^\Omega \rp &:= \lp \prod_{i\in I} \Xcal_i,\prod_{i\in I} \Xcal_i^\Omega\rp,
    \end{align}
    with the set-theoretic products on the right.
\end{Def}

So a map $X=(X_i)_{i \in I}:\,\Omega \to \prod_{i\in I} \Xcal_i$ is an admissible random variable on the product space
if and only if every component $X_i$ is an admissible random variable on $\Xcal_i$, $i \in I$.
It is easily checked that $\prod_{i \in I} \lp \Xcal_i,\Xcal_i^\Omega \rp$ is a quasi-measurable space.
In fact we have the following result:

\begin{Thm}
   The product $\prod_{i \in I} \lp \Xcal_i,\Xcal_i^\Omega \rp$ defines a categorical product of $\QMS$,
i.e.\ we have a natural isomorphism:
\begin{align}
    \QMS\lp\Zcal, \prod_{i\in I} \Xcal_i\rp &\bij \prod_{i \in I} \QMS\lp\Zcal,\Xcal_i\rp,
    & g &\mapsto (\pr_i \circ g)_{i \in I},
\end{align}
where $\pr_j:\,\prod_{i\in I} \Xcal_i \to \Xcal_j$, $j \in I$, are the (quasi-measurable) projection maps.
\end{Thm}

The proofs can be found in \Cref{qms:sec:product}.

\subsection{The Function Space}

In this subsection we will introduce the \emph{function space} of quasi-measurable spaces. This will be a construction that is not possible in the category of measurable spaces. Functions spaces will serve as \emph{exponential objects/internal homs} in the category of quasi-measurable spaces $\QMS$, making $\QMS$ cartesian closed.
The construction again closely follows \cite{Heu17} and proofs can be found in \Cref{qms:sec:fct-space}.

\begin{Def}[The function space]
    \label{qms:def:function-space}
    Let $(\Xcal,\Xcal^\Omega)$ and $(\Ycal,\Ycal^\Omega)$ be quasi-measurable spaces.
    We then define their \emph{function space} as the set:
    \begin{align} \Ycal^\Xcal := \QMS(\Xcal,\Ycal),
    \end{align}
    together with the following set of admissible random variables:
    \begin{align} \lp \Ycal^\Xcal \rp^\Omega
        &:= \lC Y:\, \Omega \to  \Ycal^\Xcal\st \lp (\omega,x)  \mapsto Y(\omega)(x) \rp \in \QMS\lp\Omega \times \Xcal,\Ycal\rp\rC\\
        &= \lC Y:\, \Omega \to  \Ycal^\Xcal\st \forall \Phi \in \Omega^\Omega.\,\forall X \in \Xcal^\Omega.\,
        \lp \omega \mapsto Y(\Phi(\omega))(X(\omega)) \rp \in \Ycal^\Omega \rC.
    \end{align}
\end{Def}

It is easily checked that $\lp\Ycal^\Xcal, \lp\Ycal^\Xcal\rp^\Omega\rp$ is again a quasi-measurable space.
 \Cref{qms:def:function-space} is made in a way such that we can canonically identify: $(\Ycal^\Xcal)^\Omega = \Ycal^{\Omega \times \Xcal}$.
Even more, we now have the following result:

\begin{Thm}
    \label{qms:thm:function-space}
    The function space $\lp\Ycal^\Xcal,\lp \Ycal^\Xcal \rp^\Omega\rp$ defines a categorical function space in $\QMS$, i.e.\
    we have the natural \emph{curry} isomorphism:
    \begin{align} \QMS(\Zcal \times \Xcal,\Ycal) \bij \QMS(\Zcal,\Ycal^\Xcal),\qquad g \mapsto (z \mapsto (x\mapsto g(z,x))).  \end{align}
    In particular, using $\Zcal=\Ycal^\Xcal$, the \emph{evaluation map} is quasi-measurable:
    \begin{align}\ev:\, \Ycal^\Xcal \times \Xcal \to \Ycal,\qquad (Y,x) \mapsto Y(x).\end{align}
\end{Thm}

Note that our suggestive notation for the set of admissible random variables $\Xcal^\Omega$ of a quasi-measurable space $(\Xcal,\Xcal^\Omega)$
does not lead to much ambiguity:
$\Xcal^\Omega$ from \Cref{qms:def:quasi-measurable-spaces} and $\Xcal^\Omega:=\QMS\lp(\Omega,\Omega^\Omega),(\Xcal,\Xcal^\Omega)\rp$
from \Cref{qms:def:function-space} agree as sets.

In more categorical terms, see \cite{Joh02} A1.5, we can rephrase \Cref{qms:thm:function-space} as:

\begin{Cor}
    \label{qms:cor:qms-cartesian-closed}
    The category $\QMS$ of quasi-measurable spaces and quasi-measurable maps together with (finite) products $\times$ and the terminal object $\one$, the one-point space,
    is a \emph{cartesian closed} category.
\end{Cor}

\begin{Eg}[Power set]
    Let $(\Xcal,\Xcal^\Omega)$ be a quasi-measurable space and $(\twop,\twop^\Omega)$ the indiscrete two-point space.
    Then $\twop^\Xcal =  \Sets(\Xcal,\twop)$
    can be identified with the \emph{power set} of $\Xcal$ (with additional admissible random variables).
\end{Eg}

Note that the results from \Cref{qms:thm:function-space} and \Cref{qms:cor:qms-cartesian-closed} for the category of quasi-measurable spaces $\QMS$ cannot be overstated. \Cref{qms:thm:function-space} implies that for quasi-measurable spaces and maps we have a theory of computation that allows for \emph{simply typed $\lambda$-calculus}, see \cite{Chu40,Rey98,Joh02} D4.2.
The analogue for the category of measurable spaces $\Meas$ does not hold.
So in these results we already see that $\QMS$ categorically is better behaved than $\Meas$.

\subsection{Further Categorical Constructions}

\subsubsection{All Other Limits and Colimits}
Similar to the \emph{product} $\prod_{i \in I} \lp \Xcal_i,\Xcal_i^\Omega\rp$ of a family of quasi-measurable spaces
one can construct their \emph{coproduct} $\coprod_{i \in I} \lp \Xcal_i,\Xcal_i^\Omega\rp$ in $\QMS$ dually. It turns out that products and coproducts follow a well-behaved \emph{distributive law}. It is important to note that the coproduct in $\QMS$ behaves more ``rigidly'' than the coproduct constructed for $\QBS$ in \cite{Heu17}. Details can be found in \Cref{qms:sec:coproduct}.

For quasi-measurable maps $f_1,f_2:\, (\Xcal,\Xcal^\Omega) \to (\Ycal,\Ycal^\Omega)$ we can easily construct
their \emph{equalizer} $\Eq(f_1,f_2)$ and \emph{coequalizer} $\CoEq(f_1,f_2)$ in $\QMS$, see \Cref{qms:sec:equalizers}.
Since we then constructed all  products and coproducts indexed over any sets, and also all equalizers and coequalizers we
can, in fact, construct every (small\footnote{\emph{Small} means indexed over sets.\label{qms:fn:small}}) \emph{limit} and every (small\footref{qms:fn:small}) \emph{colimit} by general category theoretic arguments, see \cite{Mac98} Thm.\ V.2.1. In other words, we get:

\begin{Thm}
    \label{qms:thm:co-complete}
    The category $\QMS$  of quasi-measurable spaces and quasi-measurable maps is \emph{complete} and \emph{cocomplete},
    i.e.\ it has all (small) limits and all (small) colimits.
\end{Thm}

\begin{Eg} We have the following examples of \emph{limits} in $\QMS$:
    \begin{enumerate*}
        \item the terminal object $\one$,
        \item equalizers (with possibly more than two maps),
        \item products,
        \item pull-backs aka fibre products, e.g.\ $\Xcal \ftimes_\Tcal \Ycal$,
        \item inverse/projective limits (over directed sets).
    \end{enumerate*}

    We have the following examples of \emph{colimits} in $\QMS$:
\begin{enumerate*}[resume]
        \item coproducts,
        \item coequalizers (over possibly more than two maps),
        \item push-outs,
        \item direct limits (over directed sets).
    \end{enumerate*}
\end{Eg}

\subsubsection{Fibre Products and Internal Homs}

For a fixed quasi-measurable space $(\Tcal,\Tcal^\Omega)$ we can consider quasi-measurable spaces $\Xcal$ \emph{over $\Tcal$}, i.e.\ quasi-measurable maps $T_\Xcal:\,\Xcal \to \Tcal$. A quasi-measurable map \emph{over $\Tcal$} is then a quasi-measurable map $g:\,\Xcal \to \Ycal$ such that $T_\Ycal \circ g = T_\Xcal$. The corresponding category $\QMS_\Tcal$ is then called the \emph{slice category of $\QMS$ over $\Tcal$}.
The \emph{fibre product} $\Xcal \ftimes_\Tcal \Ycal$ in $\QMS$ is then the categorical product in $\QMS_\Tcal$.
Similar to the relation between the product $\Xcal \times \Ycal$ and function space $\Zcal^\Ycal$ in $\QMS$ the fibre product $\Xcal \ftimes_\Tcal \Ycal$ has a right-adjoint in $\QMS_\Tcal$, the \emph{internal hom}: $\ihom_\Tcal(\Ycal,\Zcal)$, which is naturally a quasi-measurable space over $\Tcal$ and satisfies the natural bijection:
    \begin{align}
        \QMS_\Tcal\lp \Xcal \ftimes_\Tcal \Ycal, \Zcal \rp & \cong \QMS_\Tcal\lp \Xcal, \ihom_\Tcal\lp\Ycal,\Zcal\rp\rp.
    \end{align}
Details can be found in \Cref{qms:sec:fibre-prod} and \Cref{qms:sec:int-homs}.
Since these results hold for every $\Tcal$, we can summarize them in categorical terms as:

\begin{Thm}
    \label{qms:thm:lcc}
    The category $\QMS$ of quasi-measurable spaces and quasi-measurable maps is \emph{locally cartesian closed} (with terminal object $\one$).
\end{Thm}

Type-theoretically this means that $\QMS$ allows for an \emph{extensional\footnote{I.e.\ identity types are propositions.} dependent type theory}, see \cite{See84,Joh02} A1.5.3.

\subsubsection{Structure from Subobjects}

If  $\Wcal \to \Xcal$ is an embedding of quasi-measurable spaces, $\twop$ the indiscrete two-point space, $\I_\Wcal:\, \Xcal \to \twop$  the indicator function, $1:\,\one \to \twop$ the constant $1$-map, then we can canonically identify: $\Wcal \cong \Xcal \ftimes_\twop \one$.
This means that $\twop$ can classify strong monomorphisms (=embeddings) in $\QMS$.

The facts that $\QMS$ is locally cartesian closed, has all (finite) limits, all (finite) colimits and a subobject classifier $\twop$ for strong monomorphisms, using \cite{Joh02} A1.6, A2.6, then culminate into the following:

\begin{Thm}
    \label{qms:thm:quasitopos}
    The category $\QMS$ of quasi-measurable spaces and quasi-measurable maps forms a \emph{quasitopos}.
\end{Thm}

Also note that $\one$ separates quasi-measurable maps as they are already determined on points,
i.e.\ $\QMS$ has the small generating set $\lC \one \rC$.
We cannot use $\twop$ to classify non-strong monomorphisms, and indeed $\QMS$ is in general
neither a \emph{topos} nor a \emph{pretopos}:

\begin{Prp}
    \label{qms:prp:not-a-topos}
    If $|\Omega| \ge 2$ then $\QMS$ is not balanced, and hence neither a topos nor a pretopos.
    If $|\Omega|=1$ then $\QMS \simeq \Sets$, which is a topos.
    \begin{proof}
        Let $|\Omega| \ge 2$ and pick $\omega_1 \neq \omega_2$ in $\Omega$.
        Let $\twoc$ be the two-point space carrying only its two constant random variables and
        $\twob$ the discrete two-point space, with $\twob^\Omega=\Meas(\Omega,\twob)$, see
        \Cref{qms:eg:discrete-space} and the Example after \Cref{qms:thm:PrPf-monad}.
        Both are quasi-measurable spaces and the identity $\id:\, \twoc \to \twob$ is
        quasi-measurable. It is injective, hence a monomorphism, and surjective; since $\one$ is a
        generator of $\QMS$ every surjective quasi-measurable map is an epimorphism, so it is also
        an epimorphism. It is not an isomorphism: the map $X:\,\Omega \to \lC 0,1\rC$ with
        $X(\omega_1)=0$ and $X(\omega)=1$ otherwise lies in $\twob^\Omega$ but is not constant, so
        it is not in $\twoc^\Omega$, and the inverse of $\id$ would have to send it there.
        So $\QMS$ has a morphism that is monic and epic but not invertible, i.e.\ $\QMS$ is not
        \emph{balanced}. Every topos and every pretopos is balanced, see \cite{Joh02} A1.4,
        so $\QMS$ is neither.\\
        If $|\Omega|=1$ then every map $\Omega \to \Xcal$ is constant, so
        $\Xcal^\Omega = \Sets(\Omega,\Xcal)$ is forced for every object, $\QMS \simeq \Sets$, and
        $\twop$ classifies all monomorphisms.
    \end{proof}
\end{Prp}

So $\QMS$ sits strictly between the quasitoposes and the toposes, and the two degenerate sample
spaces $|\Omega| \le 1$ are the only ones for which it collapses to $\Sets$.

However, we show in \Cref{qms:sec:heyting} that the set of subspaces $\Sub(\Xcal)$ for every quasi-measurable space forms a
\emph{Heyting algebra}, an algebraic structure to formalize \emph{intuitionistic logic}, see \cite{Hey30}. We then show that the construction is functorial in $\Xcal$. If $f:\,\Xcal \to \Ycal$ is a quasi-measurable map then we get a pullback map $f^*:\,\Sub(\Ycal) \to \Sub(\Xcal)$, which turns out to be a Heyting algebra homomorphism.
We also construct a left-adjoint $f_!:\, \Sub(\Xcal) \to \Sub(\Ycal)$ and a right-adjoint $f_*:\,\Sub(\Xcal) \to \Sub(\Ycal)$ to $f^*$.
$f_!$ allows for an interpretion as an \emph{existential quantifier} $\exists_f$, while $f_*$ as a \emph{universal quantifier} $\forall_f$.
With these constructions $\QMS$ inherits the structure of a \emph{Heyting category}, also called \emph{logos} in \cite{Fre90}:

\begin{Thm}
   \label{qms:thm:heyting}
   The category $\QMS$ of quasi-measurable spaces and quasi-measurable maps forms a \emph{Heyting category}.
\end{Thm}

So $\QMS$ has an internal logic of a (typed) \emph{intuitionistic first-order logic}.

\begin{Rem}[Quasi-measurable spaces are concrete presheaves]
    \label{qms:rem:concrete-presheaf}
    The results of this section are, taken together, an instance of a known theorem, and it is
    worth saying which.
    Let $\Ccal_\Omega$ be the category with the two objects $\one$ and $\Omega$ and
    \begin{align}
        \Ccal_\Omega(\one,\one) &= \lC \id \rC, &
        \Ccal_\Omega(\Omega,\one) &= \lC \, ! \, \rC, &
        \Ccal_\Omega(\one,\Omega) &= \Omega, &
        \Ccal_\Omega(\Omega,\Omega) &= \Omega^\Omega.
    \end{align}
    \Cref{qms:def:sample-space-act-1} says precisely that this \emph{is} a category --- closure
    under composition gives associativity, the identity map gives the unit, and the constant maps
    are the composites $\Omega \stackrel{!}{\to} \one \stackrel{\omega}{\to} \Omega$ --- with
    $\one$ terminal and $\Ccal_\Omega(\one,-)$ faithful. Endowed with the trivial coverage, in
    which only identities cover, it is a \emph{concrete site} in the sense of \cite{BH11}
    Def.\ 4.9: the coverage is subcanonical because every presheaf is a sheaf for it, and the
    joint-surjectivity condition is vacuous.
    A \emph{concrete presheaf} on $\Ccal_\Omega$ is then a set $\Xcal := F(\one)$ together with a
    subset $F(\Omega) \ins \Sets(\Ccal_\Omega(\one,\Omega),\Xcal)=\Sets(\Omega,\Xcal)$ closed under
    restriction along $\Omega \to \one$, which is exactly containing the constant maps, and along
    each $\Phi \in \Omega^\Omega$, which is exactly closure under $X \circ \Phi$. That is
    \Cref{qms:def:quasi-measurable-spaces} verbatim, and the morphisms agree. So
    \begin{align}
        \QMS &\simeq \mathrm{Conc}(\Ccal_\Omega,\text{trivial coverage}).
    \end{align}
    By \cite{BH11} Thm.\ 5.25 the category of concrete sheaves on a concrete site is a
    quasitopos with all small limits and colimits, and locally cartesian closed; see also
    \cite{DE06} and \cite{Joh02} A2.6, C2.2. This subsumes \Cref{qms:thm:co-complete},
    \Cref{qms:thm:function-space}, \Cref{qms:cor:qms-cartesian-closed}, \Cref{qms:thm:lcc},
    \Cref{qms:thm:quasitopos} and \Cref{qms:thm:heyting}, for every sample space $\Omega$.
    The same observation has been made for quasi-Borel spaces, which are the concrete sheaves on
    the standard Borel spaces with countable measurable covers, see \cite{VO18} Section 11 and
    \cite{MMS22}.

    We nevertheless give the constructions explicitly and elementarily, because the objects they
    produce are the ones the rest of the paper computes with: the function space $\Ycal^\Xcal$,
    the fibre product, the internal hom, the subspace lattice $\Sub(\Xcal)$ with its Heyting
    implication and its two quantifier adjoints $\exists_f \dashv f^* \dashv \forall_f$, and the
    coproduct --- whose \emph{rigidity} compared to $\QBS$, see \Cref{qms:sec:coproduct}, is a
    difference the general theorem does not see. A reader who wants only the abstract statements
    may take them from \cite{BH11} and skip \Cref{sec:cat-qms}.
\end{Rem}

\begin{Eg}[A two-point sample space: simple graphs]
    \label{qms:eg:simpgph}
    Take $\Omega = \lC 0,1 \rC$ and $\Omega^\Omega$ all four self-maps, which is a legal choice by
    \Cref{qms:def:sample-space-act-1}. Then $\Xcal^\Omega \ins \Sets(\lC0,1\rC,\Xcal)=\Xcal \times
    \Xcal$ must contain the diagonal, since it contains the constants, and be closed under
    precomposition with the swap; that is, $\Xcal^\Omega$ is exactly a \emph{reflexive symmetric
    relation} on $\Xcal$, and quasi-measurable maps are exactly graph homomorphisms. So $\QMS$ is
    then equivalent to the category of simple graphs, which is a known Grothendieck quasitopos
    and known not to be a topos, see \cite{Vig03}. This is a small independent check on
    \Cref{qms:thm:quasitopos}, \Cref{qms:thm:heyting} and \Cref{qms:prp:not-a-topos}.
    Its probabilistic half is degenerate, see \Cref{qms:rem:act3-infinite}.
\end{Eg}

\section{Measurable Spaces vs.\ Quasi-Measurable Spaces}
\label{sec:meas-vs-qms}
\subsection{The Sample Space - Act 2 - The $\sigma$-Algebra}

We now endow the sample space $(\Omega,\Omega^\Omega)$ from \Cref{qms:def:sample-space-act-1} with an additional compatible $\sigma$-algebra $\Bcal_\Omega$.
This will allow us to turn measurable spaces into quasi-measurable spaces and vice versa. It will later also allow us to have probability distributions on $\Omega$.

\begin{Def}[The sample space - act 2]
    \label{qms:def:sample-space-act-2}
    In the following we fix a \emph{sample space} $(\Omega,\Omega^\Omega,\Bcal_\Omega)$, where the tuple $(\Omega,\Omega^\Omega)$ satisfies the points of
    \Cref{qms:def:sample-space-act-1} and where we now also fix a $\sigma$-algebra  $\Bcal_\Omega$  of subsets of $\Omega$ such that every $\Phi \in \Omega^\Omega$ is $\Bcal_\Omega$-$\Bcal_\Omega$-measurable, i.e.\ such that:
    \begin{align}
        \Omega^\Omega & \ins \Meas\lp (\Omega,\Bcal_\Omega),(\Omega,\Bcal_\Omega)\rp. \label{qms:eq:sample-space-act-2}
    \end{align}
\end{Def}

Besides \Cref{qms:eg:meas-sample-space} we have:

\begin{Eg}
    The following choices for $(\Omega,\Omega^\Omega,\Bcal_\Omega)$ satisfy \Cref{qms:def:sample-space-act-1} and \Cref{qms:def:sample-space-act-2}.
    \begin{enumerate}
        \item If $(\Omega,\Bcal_\Omega)$ is a measurable space we can take the set of all measurable maps:
            \begin{align} \Omega^\Omega:=\Meas\lp (\Omega,\Bcal_\Omega),(\Omega,\Bcal_\Omega)\rp.
            \end{align}
\item If $(\Omega,\Omega^\Omega)$ satisfies \Cref{qms:def:sample-space-act-1} and $\Ecal_\Omega$ is any set of subsets of $\Omega$,
            e.g.\ $\Ecal_\Omega := \lC  \Phi(\Omega) \st  \Phi \in \Omega^\Omega \rC$, then one can take:
            \begin{align}
                \Bcal_\Omega&:= \sigma\lp\lC \Phi^{-1}\lp E\rp\st E \in \Ecal_\Omega, \Phi \in \Omega^\Omega \rC\rp,
        \end{align}
        or its universal completion (or any other suitable form of completion).
        \item $\Bcal_\Omega:= \lC\emptyset,\Omega\rC$ or the power set $\Bcal_\Omega:=\lC D \st D \ins \Omega\rC$ are always valid choices.
    \end{enumerate}
\end{Eg}

\subsection{Measurable and Topological Spaces}

If $(\Omega,\Omega^\Omega,\Bcal_\Omega)$ is a sample space satisfying \Cref{qms:def:sample-space-act-2} then we
can use the fixed $\sigma$-algebra $\Bcal_\Omega$ to turn every topological space or measurable space into quasi-measurable space by the following construction for the set of admissible random variables,
which basically allows for all measurable maps like in the classical measure-theoretic setting.

\begin{Def}[The admissible random variables for measurable and topological spaces]
    \label{qms:eg:meas-qms-space}
    Let $\Xcal$ be a set endowed with any set $\Ecal_\Xcal$ of subsets of $\Xcal$,
    e.g.\ a measurable space with $\sigma$-algebra $\Ecal_\Xcal$ or a topological space with topology $\Ecal_\Xcal$, then the following choice of admissible random variables:
    \begin{align}
        \Xcal^\Omega &:= \Fcal\lp \Ecal_\Xcal\rp := \lC X:\, \Omega \to \Xcal \st \forall A \in \Ecal_\Xcal.\, X^{-1}(A) \in \Bcal_\Omega \rC,
    \end{align}
    turns $(\Xcal,\Xcal^\Omega)$ into a quasi-measurable space w.r.t.\ $(\Omega,\Omega^\Omega)$.
\end{Def}

By default we endow every topological space or measurable space with this set of admissible random variables. A special case of importance is the following:

\begin{Def}[Discrete quasi-measurable spaces]
    \label{qms:eg:discrete-space}
    Let $\Xcal$ be any set. The \emph{discrete quasi-measurable space structure} on $\Xcal$ w.r.t.\  $(\Omega,\Omega^\Omega,\Bcal_\Omega)$ is given by:
    \begin{align}
        \Xcal^\Omega & := \lC X:\, \Omega \to \Xcal \st \forall x \in \Xcal.\, X^{-1}(x) \in \Bcal_\Omega \rC.
    \end{align}
\end{Def}

Unless stated otherwise, we, by default, endow every \emph{countable} space with this discrete quasi-measurable space structure, e.g.\
for $\twob:=\lC 0,1 \rC$, $\N$, $\Z$ or $\Q$. Note that $\twop$ from \Cref{qms:def:indiscrete-qms} does not carry the discrete quasi-measurable space structure.

\subsection{The Induced Sigma-Algebra}

In this subsection we define the \emph{induced $\sigma$-algebra} $\Bcal_\Xcal$ of a quasi-measurable space $(\Xcal,\Xcal^\Omega)$.
This $\sigma$-algebra makes all admissible random variables $X \in \Xcal^\Omega$ $\Bcal_\Omega$-$\Bcal_\Xcal$-measurable.
We show that the induced $\sigma$-algebra can be turned into a quasi-measurable space $(\Bcal_\Xcal, (\Bcal_\Xcal)^\Omega)$
in its own right. This allows us to compare $\sigma$-algebras and quasi-measurable spaces on the same footing inside $\QMS$.
Note that this was not possible inside the category of measurable spaces $\Meas$ (without difficulties or topological considerations).

\begin{Def}[The induced $\sigma$-algebra of a quasi-measurable space]
 Let $(\Xcal,\Xcal^\Omega)$ be a quasi-measurable space.
 Then $\Xcal$ can be endowed with the following $\sigma$-algebra:
 \begin{align} \Bcal_\Xcal:=\Bcal(\Xcal^\Omega)&:=\lC A\ins \Xcal\st \forall X \in \Xcal^\Omega.\,X^{-1}(A) \in \Bcal_\Omega \rC. \end{align}
 Similarly, we get the analogous construction on the product space $\Omega \times \Xcal$:
 \begin{align} \Bcal_{\Omega \times \Xcal}&:=\Bcal(\Omega^\Omega \times \Xcal^\Omega)= \lC D \ins \Omega \times \Xcal\st \forall \Phi \in \Omega^\Omega.\, \forall X \in \Xcal^\Omega.\, (\Phi\times X)^{-1}(D) \in \Bcal_\Omega \rC.\end{align}
 Using these notations we define the set of admissible random variables on $\Bcal_\Xcal$ as follows:
 \begin{align} (\Bcal_\Xcal)^\Omega := \lC \Psi:\, \Omega \to \Bcal_\Xcal\st \exists D \in \Bcal_{\Omega \times \Xcal}.\,
           \forall \omega \in \Omega. \, \Psi(\omega)=D_\omega \rC, \end{align}
  where we define the section $D_\omega$ as:
  \begin{align} D_\omega := \lC x \in \Xcal\st (\omega, x) \in D \rC \in \Bcal_\Xcal.\end{align}
\end{Def}

It can easily be checked, see \Cref{qms:sec:sigma-algebra}. that the induced $\sigma$-algebra $\lp \Bcal_\Xcal,\lp\Bcal_\Xcal \rp^\Omega \rp$ is again a quasi-measurable space and
we have a canonical identification $\lp \Bcal_\Xcal \rp^\Omega = \Bcal_{\Omega \times \Xcal}$.
By definition it is clear that every admissible random variable $X \in \Xcal^\Omega$ is $\Bcal_\Omega$-$\Bcal_\Xcal$-measurable.
Furthermore, the following result holds:

\begin{Lem}[The $\sigma$-algebra identifies with the space of indicator functions]
    \label{qms:lem:inclusion-check}
    Let $(\Xcal,\Xcal^\Omega)$, $(\Ycal,\Ycal^\Omega)$ be quasi-measurable spaces and $\twob:=\lC0,1 \rC$ the discrete two-point space as discussed in \Cref{qms:eg:discrete-space}.
    Then the indicator function (inclusion check) induces an isomorphism of quasi-measurable spaces:
    \begin{align}
        \I:\; (\Bcal_\Ycal,(\Bcal_\Ycal)^\Omega) &\bij (\twob^\Ycal,(\twob^\Ycal)^\Omega),& B &\mapsto \I_B,
\end{align}
where $\I_B(y):=1$ if $ y \in B$ and $\I_B(y):=0$ if $y \notin B$. Furthermore, the following maps:
 \begin{align}   \Bcal_\Ycal \times \Ycal^\Xcal  &\to \Bcal_\Xcal, &(B,f) &\mapsto f^{-1}(B),  \\
     \twob^\Ycal \times \Ycal^\Xcal & \to \twob^\Xcal, & (\I_B,f) & \mapsto \I_B \circ f = \I_{f^{-1}(B)},
 \end{align}
 are quasi-measurable and can be identified under the isomorphism $\I$.
\end{Lem}

The proof can be found in \Cref{qms:sec:sigma-algebra}.

\subsection{The Adjunction}

If $(\Omega,\Omega^\Omega,\Bcal_\Omega)$ is a sample space according to
\Cref{qms:def:sample-space-act-2}
then we can study how the category of measurable spaces $\Meas$ and the category of quasi-measurable spaces $\QMS$ are related to each other.
The main point can be summarized in the following \emph{adjunction} between measurable spaces and quasi-measurable spaces:

\begin{Thm}[Measurability equals quasi-measurability]
    \label{qms:thm:adjunction}
    A map $g:\, \Xcal \to \Ycal$ from a quasi-measurable space $(\Xcal,\Xcal^\Omega)$ to a measurable space $(\Ycal,\Bcal_\Ycal)$ is
    measurable if and only if it is quasi-measurable, provided one uses for $\Bcal_\Xcal$ and $\Ycal^\Omega$ the corresponding choices:
    \begin{align}
        \Bcal_\Xcal := \Bcal(\Xcal^\Omega)  &:= \lC A \ins \Xcal\st \forall X \in \Xcal^\Omega.\, X^{-1}(A) \in \Bcal_\Omega \rC, \\
        \Ycal^\Omega := \Fcal(\Bcal_\Ycal) &:= \Meas\lp (\Omega,\Bcal_\Omega),(\Ycal,\Bcal_\Ycal)\rp.
 \end{align}
 In other words, we have the natural identification of sets:
     \begin{align} \Meas\lp(\Xcal,\Bcal(\Xcal^\Omega)),(\Ycal,\Bcal_\Ycal)\rp = \QMS\lp(\Xcal,\Xcal^\Omega), (\Ycal,\Fcal(\Bcal_\Ycal)) \rp.   \end{align}
     Furthermore, iterating $\Fcal$ and $\Bcal$ always gives the following inclusions:
     \begin{align}
         \Fcal\Bcal(\Xcal^\Omega) &\sni \Xcal^\Omega, &\Bcal \Fcal\Bcal(\Xcal^\Omega) &= \Bcal(\Xcal^\Omega),\\
         \Bcal\Fcal(\Bcal_\Ycal) &\sni \Bcal_\Ycal,  &\Fcal\Bcal\Fcal(\Bcal_\Ycal) &= \Fcal(\Bcal_\Ycal).
    \end{align}
\end{Thm}

Note that both functors $\Bcal$ and $\Fcal$ are dependent on the choice of the $\sigma$-algebra $\Bcal_\Omega$.
Also note that we always have: $\Bcal(\Omega^\Omega) = \Bcal_\Omega$, but we might have:  $\Fcal(\Bcal_\Omega) \supsetneq \Omega^\Omega$.

\begin{Eg}
    Let $(\Ycal,\Ecal_\Ycal)$ be a topological space, e.g.\ $\Ycal=\R$.
    Then a map $g:\, (\Xcal,\Xcal^\Omega) \to (\Ycal,\Fcal(\Ecal_\Ycal))$ is quasi-measurable if and only if
    for all $B \in \Ecal_\Ycal$ we have that $g^{-1}(B) \in \Bcal(\Xcal^\Omega)$,
    which amounts to checking if $X^{-1}(g^{-1}(B)) \in \Bcal_\Omega$ for all $X \in \Xcal^\Omega$.
\end{Eg}

The adjunction \Cref{qms:thm:adjunction} also tells us that if we wanted to restrict our attention to measurable spaces, but still wanted to work
with the categorical constructions from $\QMS$, then we could just apply the composite functor $\Bcal\Fcal$. This turns every measurable map $h:\, (\Zcal,\Bcal_\Zcal) \to (\Ycal,\Bcal_\Ycal)$
into the measurable map $h:\, (\Zcal,\Bcal\Fcal(\Bcal_\Zcal)) \to (\Ycal,\Bcal\Fcal(\Bcal_\Ycal))$. Here the $\sigma$-algebra $\Bcal\Fcal(\Bcal_\Ycal)$ can be interpreted as a form of a completion of $\Bcal_\Ycal$ \emph{relative to $\Bcal_\Omega$}. The adjunction tells us that measurable maps between such complete measurable spaces are the same as quasi-measurable maps. So the corresponding category $\SQMS$ of such spaces can fully be embedded into both, $\Meas$
(with coreflector $\Bcal\Fcal$) and $\QMS$ (with reflector $\Fcal\Bcal$). $\SQMS$ could thus be regarded as the ``intersection'' of $\Meas$ and $\QMS$. More details can be found in \Cref{qms:sec:sturdy-qms}. A reason not to solely work in $\SQMS$ is that we would lose function spaces $\Ycal^\Xcal$, which exist in $\QMS$, but usually not in $\SQMS$.

It is worth mentioning, and a basic fact from category theory, that, as a right-adjoint, $\Fcal$ preserves arbitrary limits, like products, and, as a left-adjoint, $\Bcal$ preserves arbitrary colimits, like coproducts. For instance, if $(\Xcal_i,\Bcal_i)$ is a family of measurable spaces, $i \in I$,
then we have for the product $\sigma$-algebra $\bigotimes_{i \in I} \Bcal_i $, which is generated by finite cylinder sets, the following identification:
\begin{align} \Fcal\lp \bigotimes_{i \in I} \Bcal_i \rp = \prod_{i\in I} \Fcal(\Bcal_i). \label{qms:eq:adj-prod}
    \end{align}
More details and proofs can be found in \Cref{qms:sec:adjunction}.

\section{The Quasi-Measurable Probability Monad}
\label{sec:giry-completions}

For literature about monads see \cite{Koc70,Koc72,Str72,Mac98,Mog91}; the definitions we use are recalled in \Cref{qms:sec:app-monads}. As a first example we recall the strong probability monad of Giry, see \cite{Law62,Gir82}. It is not the monad this paper works with --- the category of measurable spaces $\Meas$ is not cartesian closed --- but it is the yardstick against which the monads on $\QMS$ are measured, and it is what the induced $\sigma$-algebras are built from.

This paper puts \emph{one} probability monad on $\QMS$ in the foreground: the push-forward probability monad $\PrPf$ of \Cref{sec:pf-monad}, which is the one that carries the interpretation of collecting the distributions of admissible random variables. Two others appear here in support: the Giry monad $\Giry$ just below, as the measure-theoretic yardstick, and the \emph{quasi-measurable probability monad} $\PrQm$ of \Cref{sec:general-monad}, which is the well-behaved companion of $\Giry$ on $\QMS$. The name is meant literally, and it is the reason this section carries it: a probability measure $P$ on the induced $\sigma$-algebra $\Bcal_\Xcal$ belongs to $\PrQm(\Xcal)$ precisely when the map $A \mapsto P(A)$ is itself \emph{quasi-measurable} on the space of events $\Bcal_\Xcal$, exactly as Giry\'s measures are the ones that are \emph{measurable} on it, see \Cref{qms:def:PrQm-maps}. The Giry monad and the completions of $\sigma$-algebras are treated first because they are the two ingredients $\PrQm$ is assembled from. Three further monads $\PrPF$, $\PrpF$, $\Prpf$, siblings of $\PrPf$, are developed in full generality in \Cref{qms:sec:app-other-monads}, where \Cref{qms:tab:monads} compares all of them and explains why $\PrPf$ is the one to work with.

\subsection{Giry's Probability Monad}

\begin{Def}[The space of probability measures for measurable spaces]
    \label{qms:def:prob-meas}
    Let $(\Xcal,\Bcal_\Xcal)$ be a measurable space. Then we abbreviate with $\Giry(\Xcal,\Bcal_\Xcal)$ the set of all probability measures $P:\, \Bcal_\Xcal \to [0,1]$. Unless explicitly stated otherwise, we consider $\Giry(\Xcal,\Bcal_\Xcal)$
    to be the measurable space that is endowed with the smallest $\sigma$-algebra $\Bcal_{\Giry(\Xcal,\Bcal_\Xcal)}$ such that the evaluation maps:
    \begin{align}
        \ev_A:\, \Giry(\Xcal,\Bcal_\Xcal) &\to [0,1],&  P &\mapsto P(A),
    \end{align}
    are measurable for all $A \in \Bcal_\Xcal$, i.e.: \begin{align}\Bcal_{\Giry(\Xcal,\Bcal_\Xcal)} &:= \sigma\lp\lC \ev_A^{-1}((t,1])\st A \in \Bcal_\Xcal, t \in \R \rC\rp. \end{align}
   When turned into a quasi-measurable space, we, as usual, endow the measurable space $\Giry(\Xcal,\Bcal_\Xcal)$ with the following
   set of admissible random variables:
   \begin{align}
        \Giry(\Xcal,\Bcal_\Xcal)^\Omega &:= \Fcal\lp \Bcal_{\Giry(\Xcal,\Bcal_\Xcal)}\rp\\
                                         &=\lC K(X|\Inp):\,\Omega \to \Giry(\Xcal,\Bcal_\Xcal)\st \forall A \in \Bcal_\Xcal.\,
                                         (\omega \mapsto K(X \in A|\Inp=\omega)) \in [0,1]^\Omega \rC.
   \end{align}
\end{Def}

\begin{Thm}[The strong probability monad of Giry, \cite{Law62,Gir82,Sat18}]
    Let $(\Meas,\times,\one)$ be the category of all measurable spaces and measurable maps, where the cartesian monoidal structure
    $\times$ is given
    by the usual product-$\sigma$-algebra:
    \begin{align} (\Xcal,\Bcal_\Xcal)\times(\Ycal,\Bcal_\Ycal) &= (\Xcal \times \Ycal,\Bcal_\Xcal\otimes\Bcal_\Ycal),\\
     \Bcal_\Xcal\otimes\Bcal_\Ycal &= \sigma\lp\lC A \times B\st A \in \Bcal_\Xcal, B \in \Bcal_\Ycal \rC\rp.  \end{align}
    The \emph{Giry monad} $(\Giry,\delta,\Mbb)$ is given, as described in \Cref{qms:def:prob-meas}, by:
    \begin{align} \Giry(\Xcal,\Bcal_\Xcal) &= \lC P:\; \Bcal_\Xcal \to [0,1] \text{ probability measure}\rC, \\
    \Bcal_{\Giry(\Xcal,\Bcal_\Xcal)} &= \sigma\lp\lC \ev_A^{-1}((t,1])\st A \in \Bcal_\Xcal, t \in \R \rC\rp.
    \end{align}
The \emph{monad unit} $\delta$ is given by mapping points to their Dirac delta point measures:
    \begin{align} \delta:\; (\Xcal,\Bcal_\Xcal) \to \lp\Giry(\Xcal,\Bcal_\Xcal),\Bcal_{\Giry(\Xcal,\Bcal_\Xcal)} \rp, \qquad
    x \mapsto \delta_x=(A \mapsto \I_A(x)).\end{align}
    The \emph{monad action} $\Mbb$ is given (while suppressing the $\sigma$-algebras for readability) via:
    \begin{align} \Mbb:\; \Giry(\Giry(\Xcal)) \to \Giry(\Xcal),\qquad \Mbb(\Pi)(A) := \int_{\Giry(\Xcal)} \ev_A(P)\,d\Pi(P).\end{align}
    The Giry monad is \emph{affine}, \emph{strong} and \emph{commutative} --- the last being the classical Fubini--Tonelli theorem on a product of $\sigma$-algebras --- see \cite{Gir82,Pra09,Bar16}; affinity is the case $T=\Giry$ of \Cref{qms:lem:affine}. Its \emph{costrength} is given by:
    \begin{align} \rho_{\Xcal,\Ycal}:\; \Giry(\Xcal) \times \Ycal \to \Giry(\Xcal \times \Ycal),\qquad (P,y) \mapsto P \otimes \delta_y,\end{align}
\end{Thm}

Also note that it is known that $\Meas$ is complete and cocomplete,
i.e.\ it has all small limits and all small colimits, but it is not cartesian closed, see \cite{Aum61}.
The Giry monad supports first order semantics of continuous probabilistic programming languages,
    but not higher-order ones, as $\Meas$ is not cartesian closed.

\begin{Rem}[Strength or costrength: which orientation we use]
    \label{qms:rem:costrength}
    A strong monad on a symmetric monoidal category carries two natural transformations, related
    by the braiding:
    \begin{align}
        \tau_{\Xcal,\Ycal}:\; \Xcal \times T(\Ycal) &\to T(\Xcal \times \Ycal), &
        \rho_{\Xcal,\Ycal}:\; T(\Xcal) \times \Ycal &\to T(\Xcal \times \Ycal),
    \end{align}
    the \emph{strength} and the \emph{costrength}, see \Cref{qms:sec:app-monads}. Either one
    determines the other, so the choice is one of convenience --- but it is not a free one here,
    because this paper writes conditional distributions as $P(X|Y)$, with the \emph{measured}
    variable first and the \emph{conditioning} variable second.

    That convention fixes the orientation of everything downstream. The product of Markov kernels,
    \Cref{qms:thm:product-PrPf}, is
    \begin{align}
        \otimes:\; \PrPf(\Xcal)^{\Ycal \times \Zcal} \times \PrPf(\Ycal)^\Zcal
        &\to \PrPf(\Xcal \times \Ycal)^\Zcal, &
        K(X|Y,Z) \otimes Q(Y|Z) &\mapsto P(X,Y|Z),
    \end{align}
    with the measure-valued argument on the left; and the corresponding degenerate case, obtained
    by taking the second factor to be a point mass, is
    \begin{align}
        \rho_{\Xcal,\Ycal}:\; \PrPf(\Xcal) \times \Ycal \to \PrPf(\Xcal \times \Ycal), \qquad
        (P,y) \mapsto P \otimes \delta_y,
    \end{align}
    which is the \emph{costrength}. Using $\tau$ instead would put the plain object on the left and
    force a braiding $\swap_{\Xcal,\Ycal}$ into the statement of essentially every product, Fubini
    and disintegration result in \Cref{sec:qus}. So we consistently work with $\rho$, and quote
    $\tau$ only where the general theory of \Cref{qms:sec:app-monads} needs both --- for instance in
    the definition of a commutative monad, which is a statement about the two of them agreeing.
\end{Rem}

\begin{Rem}[Which properties of a probability monad we track]
    \label{qms:rem:monad-properties}
    Several probability monads appear in this paper, and for each of them we record the same
    four properties, in the same order. Each one licenses a specific construction downstream, so
    it is worth saying once what they buy. Precise definitions are collected in
    \Cref{qms:sec:app-monads}.
    \begin{enumerate}
        \item \emph{Affine}: $T(\one) \cong \one$, i.e.\ there is exactly one ``probability
            measure'' on the one-point space. This is what makes $T$ a monad of \emph{probability}
            measures rather than of sub-probability measures or of measures, and it is what makes
            the marginalisation maps $T(\Xcal \times \Ycal) \to T(\Xcal)$ behave.
        \item \emph{Strong}: $T$ comes with a natural \emph{costrength}
            $\rho_{\Xcal,\Ycal}:\, T(\Xcal) \times \Ycal \to T(\Xcal \times \Ycal)$, compatible
            with the unit $\delta$ and the multiplication $\Mbb$. The strength is what lets a
            distribution be carried along with a \emph{parameter}, so that $T$ can be used in a
            context: without it there is no Kleisli composition of maps that also depend on side
            information, and no monadic semantics for a probabilistic programming language with
            variables. On the choice of $\rho$ rather than the strength
            $\tau_{\Xcal,\Ycal}:\, \Xcal \times T(\Ycal) \to T(\Xcal \times \Ycal)$ see
            \Cref{qms:rem:costrength}.
        \item \emph{Commutative} (equivalently, for a strong monad on a symmetric monoidal
            category, \emph{symmetric monoidal}, see \cite{Koc70,Koc72}): the two maps
            $T(\Xcal) \times T(\Ycal) \to T(\Xcal \times \Ycal)$ obtained by applying the left and
            the right strength in the two possible orders agree. Written out in integrals, this
            says exactly that the order of integration may be exchanged --- it \emph{is} the
            Fubini property. This is the property that makes the product of two independent
            distributions, and more generally the product of Markov kernels, well defined
            independently of the order in which the factors are taken.
        \item \emph{Product of Markov kernels}: the derived operation
            \[ \otimes:\; T(\Xcal)^{\Ycal \times \Zcal} \times T(\Ycal)^\Zcal \to
               T(\Xcal \times \Ycal)^\Zcal, \]
            which composes a conditional distribution with a marginal one. It is the concrete form
            commutativity takes in this paper, and it is what the causal models of
            \Cref{sec:causal-models} are built from.
    \end{enumerate}
    These are not independent bookkeeping: for a monad on a cartesian monoidal category, being
    \emph{affine and commutative} is precisely what makes its Kleisli category a
    \emph{Markov category} in the sense of \cite{Fri20}, which is the statement we reach in
    \Cref{thm:Markov-category}. So the three properties above are exactly the hypotheses of the
    structural result at the end of \Cref{sec:qus}, and each monad is measured against them.

    For the Giry monad all four hold classically. For the monads on $\QMS$ they hold under
    explicit hypotheses on the sample space: $\PrQm$ is affine and strong with no further
    assumptions, whereas for $\PrPf$ strength needs $\Omega \cong \Omega \times \Omega$,
    \Cref{qms:asm:prod-iso}, and commutativity needs the Fubini property,
    \Cref{qms:asm:fubini}. Isolating these two assumptions, rather than burying them in a
    particular choice of $\Omega$, is one of the points of \Cref{sec:pf-monad}.
\end{Rem}

\subsection{Completions of Sigma-Algebras}

In this subsection we gather some results of how measurable maps and Markov kernels behave under completion of $\sigma$-algebras.
This will be used later, when deciding what $\sigma$-algebra a sample space should carry.

\begin{Def}[Completions]
    Let $\Xcal$ be a set and $\Ecal$ any set of subsets of $\Xcal$. We abbreviate: $\Giry(\Xcal,\Ecal):=\Giry(\Xcal,\sigma(\Ecal))$.
    \begin{enumerate}
        \item For $\emptyset \neq \Pcal \ins \Giry(\Xcal,\Bcal_\Xcal)$ with $\Ecal \ins \Bcal_\Xcal$
    we define the \emph{$\Pcal$-completion} of $\Ecal$ as:
    \begin{align} (\Ecal)_\Pcal &:= \bigcap_{P \in \Pcal} \lp\sigma(\Ecal)\rp_{P|_{\sigma(\Ecal)}},
        \label{qms:eq:def-compl}
    \end{align}
     i.e.\ the intersection of all Lebesgue completions of the $\sigma$-algebra $\sigma(\Ecal)$ generated by $\Ecal$ w.r.t.\
     the restrictions of every probability measure $P$ in $\Pcal$ to $\sigma(\Ecal)$.
\item We also put $(\Ecal)_\emptyset := \lC A \st A \ins \Xcal \rC$, the power set of $\Xcal$.
\item The \emph{universal completion} of $\Ecal$ is $(\Ecal)_\Giry$  for $\Giry:=\Giry(\Xcal,\Ecal)$.
\item We say that $\Ecal$ or $(\Xcal,\Ecal)$, resp., is \emph{$\Pcal$-complete}  if $\Ecal = (\Ecal)_\Pcal$.
\item We say that $\Ecal$ or $(\Xcal,\Ecal)$, resp., is \emph{universally complete} if $\Ecal=(\Ecal)_\Giry$.
    \end{enumerate}
\end{Def}

\begin{Lem}[Measurable maps under completions]
    \label{qms:lem:map-completetions}
    Let $f:\, (\Xcal,\Bcal_\Xcal) \to (\Ycal,\Bcal_\Ycal)$ be a measurable map with induced map:
    \begin{align}
    f_*:\, \Giry(\Xcal,\Bcal_\Xcal) \to \Giry(\Ycal,\Bcal_\Ycal), \qquad P \mapsto f_*P.
    \end{align}
    Consider subsets: $\Pcal \ins \Giry(\Xcal,\Bcal_\Xcal)$ and $\Qcal \ins \Giry(\Ycal,\Bcal_\Ycal)$,
    such that: $f_*(\Pcal) \ins \Qcal$.

    Then $f$ is also $(\Bcal_\Xcal)_\Pcal$-$(\Bcal_\Ycal)_\Qcal$-measurable:
    \begin{align}
        f:\, (\Xcal,(\Bcal_\Xcal)_\Pcal) \to (\Ycal,(\Bcal_\Ycal)_\Qcal).
    \end{align}
    In particular, this holds for the universal completions\footnote{Similarly, we could take $\Pcal$ and $\Qcal$, resp., as the sets of all \emph{perfect} probability measures, or, the sets of all \emph{countably compact} probability measures, because push-forwards of such probability measures along measurable maps stay in the same class, see \cite{Fremlin} 451E, 452R.}:
    \begin{align}
        f:\, (\Xcal,(\Bcal_\Xcal)_\Giry) &\to (\Ycal,(\Bcal_\Ycal)_\Giry).
    \end{align}
\end{Lem}

\begin{Thm}[Markov kernels under completions]
    \label{qms:thm:markov-kernel-completions}
    Let $K(Y|X):\,(\Xcal,\Bcal_\Xcal) \to \Giry(\Ycal,\Bcal_\Ycal)$ be a measurable map, i.e.\ a Markov kernel, with the induced map given
    for $B \in \Bcal_\Ycal$ and $P \in \Giry(\Xcal,\Bcal_\Xcal)$ via:
    \begin{align}K(Y|X_\circ):\; \Giry(\Xcal,\Bcal_\Xcal) &\to \Giry(\Ycal,\Bcal_\Ycal), &
    K(Y \in B|X_\circ=P) &:=\int K(Y \in B|X=x) \, P(dx).\end{align}
    Let $\Pcal \ins\Giry(\Xcal,\Bcal_\Xcal)$  and $\Qcal \ins \Giry(\Ycal,\Bcal_\Ycal)$ be subsets such that:
    \begin{align}
        K \circ \bar \Pcal :=  \lC K(Y|X_\circ = P) \st P \in \Pcal \rC \cup \lC K(Y|X=x) \st x \in \Xcal \rC \stackrel{!}{\ins} \Qcal.
    \end{align}
    Then the following induced Markov kernel is a well-defined and measurable map:
    \begin{align} K(Y|X):\, (\Xcal,(\Bcal_\Xcal)_\Pcal) \to \Giry(\Ycal,(\Bcal_\Ycal)_\Qcal).\end{align}
    In particular, this holds for $\Qcal = K \circ \bar \Pcal$ and for the universal completions:
    \begin{align}
        K(Y|X):\, (\Xcal,(\Bcal_\Xcal)_\Pcal) &\to \Giry(\Ycal,(\Bcal_\Ycal)_{K \circ \bar \Pcal}),\\
        K(Y|X):\, (\Xcal,(\Bcal_\Xcal)_\Giry) &\to \Giry(\Ycal,(\Bcal_\Ycal)_\Giry).
    \end{align}
\end{Thm}

An important example that occurs in proofs later on is the following.

\begin{Eg}
    \label{qms:eg:markov-kernel-completions}
    Let $(\Zcal,\Bcal_\Zcal)$ and $(\Xcal,\Bcal_\Xcal)$ be measurable spaces and
    \begin{align}
        Q(Z|X):\, (\Xcal,\Bcal_\Xcal) &\to \Giry(\Zcal,\Bcal_\Zcal),
    \end{align}
    be a Markov kernel.
    With help of Dynkin's Lemma we then see that the following map defines a Markov kernel as well:
    \begin{align}
        K(Y|X):\, (\Xcal,\Bcal_\Xcal) &\to \Giry(\Zcal \times \Xcal, \Bcal_\Zcal\otimes\Bcal_\Xcal), \\
        x & \mapsto Q(Z|X)\otimes \delta_x(X)  \\
          &\qquad =(D \mapsto Q(Z \in D^x|X=x)),
    \end{align}
    where $D^x := \lC z \in \Zcal \st (z,x) \in D \rC$.
    The induced map:
    \begin{align}
    K(Y|X_\circ):\, \Giry(\Xcal,\Bcal_\Xcal) \to \Giry(\Zcal \times \Xcal, \Bcal_\Zcal\otimes\Bcal_\Xcal),
    \end{align}
    is on $D \in \Bcal_\Zcal\otimes\Bcal_\Xcal$ and at $P \in \Giry(\Xcal,\Bcal_\Xcal)$ then given by:
    \begin{align}
        K(Y \in D|X_\circ=P)  &= \int Q(Z \in D^x|X=x)\, P(X \in dx)  \\
        &= (Q(Z|X) \otimes P(X))(D).
    \end{align}
    Now let $\Pcal \ins\Giry(\Xcal,\Bcal_\Xcal)$ be any subset. Then we get:
    \begin{align}
        \lC K(Y|X=x) \st x \in \Xcal \rC &= \lC Q(Z|X)\otimes \delta_x(X) \st x \in \Xcal \rC,\\
        \lC K(Y|X_\circ = P) \st P \in \Pcal \rC &= \lC Q(Z|X) \otimes P(X) \st P \in \Pcal \rC.
    \end{align}
    So if we put:
    \begin{align}
        \bar \Pcal &:= \Pcal \cup \lC \delta_x \st x \in \Xcal \rC, & \Qcal &:= Q \otimes \bar \Pcal \quad\text{(or bigger)},
    \end{align}
    then we get a well-defined Markov kernel on the completions:
    \begin{align}
        K(Y|X):\, (\Xcal,(\Bcal_\Xcal)_{\Pcal}) &\to \Giry(\Zcal \times \Xcal, (\Bcal_\Zcal\otimes\Bcal_\Xcal)_\Qcal), \\
        x & \mapsto Q(Z|X)\otimes \delta_x(X).
    \end{align}
\end{Eg}

\subsection{The Monad $\PrQm$ of Quasi-Measurable Probability Measures}
\label{sec:general-monad}

In this subsection we introduce the \emph{quasi-measurable probability monad} $\PrQm$, which is well-behaved w.r.t.\ the categorical constructions of the category of quasi-measurable spaces $\QMS$ like products, function spaces, induced $\sigma$-algebras as quasi-measurable spaces, etc. Its members are the probability measures that are quasi-measurable as maps on the space of events, which is where the name comes from: $\PrQm(\Xcal)$ is cut out of $\Giry(\Xcal)$ by asking that $A \mapsto P(A)$ lie in $[0,1]^{\Bcal_\Xcal}$, equivalently that $\omega \mapsto P(D^\omega)$ be an admissible random variable for every sample-point-indexed family of events $D \in \Bcal_{\Xcal \times \Omega}$.
$\PrQm$ plays a similar role for $\QMS$ as the Giry monad $\Giry$ for $\Meas$. It turns out that $\PrQm$ also constitutes an affine and strong probability monad.
However, $\PrQm$ does not yet have a proper interpretation as collecting the probability distributions of (admissible) random variables, which we will develop later. Probability measures $P \in \PrQm(\Xcal)$ are rather interpreted as those probability measures that are well-behaved under parameterizations via quasi-measurable maps.

In the following let $(\Omega,\Omega^\Omega,\Bcal_\Omega)$ be a sample space according to \Cref{qms:def:sample-space-act-2}.

\begin{Def}[The space of quasi-measurable probability measures]
    \label{qms:def:PrQm-maps}
    Let $(\Xcal,\Xcal^\Omega)$ be a quasi-measurable space and $\Bcal_\Xcal:=\Bcal(\Xcal^\Omega)$ its induced $\sigma$-algebra, which we will also consider
    as the quasi-measurable space $(\Bcal_\Xcal,(\Bcal_\Xcal)^\Omega)$.
    Then we define:
    \begin{align}
        \PrQm(\Xcal,\Xcal^\Omega) &:= \Giry(\Xcal,\Bcal_\Xcal) \cap [0,1]^{\Bcal_\Xcal} \\
                     &= \lC P \in \Giry(\Xcal,\Bcal_\Xcal) \st\forall D \in \Bcal_{\Xcal\times\Omega}.\,
                        (\omega \mapsto P(D^\omega)) \in [0,1]^\Omega \rC,\\
        \PrQm(\Xcal,\Xcal^\Omega)^\Omega &:= \iota^*([0,1]^{\Bcal_\Xcal})^\Omega\\
                                         &= \lC K(X|\Inp):\,\Omega \to \PrQm(\Xcal,\Xcal^\Omega)\st K(X|\Inp) \in ([0,1]^{\Bcal_\Xcal})^\Omega \rC\\
                                         &= \bigg\{  K(X|\Inp):\,\Omega \to \PrQm(\Xcal,\Xcal^\Omega)\,\bigg|\, \forall \Phi \in \Omega^\Omega.\,\forall D \in \Bcal_{\Xcal \times \Omega}.\\
                                         & \qquad\quad \lp\omega \mapsto K(X \in D^\omega|\Inp=\Phi(\omega))\rp \in [0,1]^\Omega\bigg\}.
   \end{align}
\end{Def}

Since $(\PrQm(\Xcal,\Xcal^\Omega),\PrQm(\Xcal,\Xcal^\Omega)^\Omega)$ inherited the subspace structure of $[0,1]^{\Bcal_\Xcal}$ it is clear that it is a quasi-measurable space and also that the evaluation map is quasi-measurable:
\begin{align}
   \ev:\, \PrQm(\Xcal) \times \Bcal_\Xcal &\to [0,1], & (P,A) \mapsto P(A).
\end{align}
From this we then get that the inclusion map: $\PrQm(\Xcal) \inj \Giry(\Xcal)$ is quasi-measurable.
Even more, we have the following:
\begin{Thm}
    \label{qms:thm:PrQm-maps}
    For quasi-measurable spaces $(\Xcal,\Xcal^\Omega)$, $(\Ycal,\Ycal^\Omega)$, $(\Zcal,\Zcal^\Omega)$   the following maps are well-defined and quasi-measurable:
    \begin{align}
        \int:\,&& [0,\infty]^\Xcal \times \PrQm(\Xcal) &\to [0,\infty], & (f,P) &\mapsto \int f \,dP,\\
        \pf:\;&& \Ycal^\Xcal \times \PrQm(\Xcal) &\to \PrQm(\Ycal), & (g,P) &\mapsto g_*P, \\
        \Mbb:\;&& \PrQm(\PrQm(\Xcal)) &\to \PrQm(\Xcal), & \Mbb(\Pi)(A) &:= \int \ev(P,A)\, d\Pi(P), \\
        \rho:\;&& \PrQm(\Xcal) \times \Ycal &\to \PrQm(\Xcal \times \Ycal), &
        (P,y) &\mapsto P \otimes \delta_y,\\
        \otimes:\;&& \PrQm(\Xcal)^{\Ycal \times \Zcal} \times \PrQm(\Ycal)^\Zcal &\to \PrQm(\Xcal \times \Ycal)^\Zcal,&
    \end{align}
    where the last map is given on $z \in \Zcal$ and $D \in \Bcal_{\Xcal \times \Ycal}$ by:
    \begin{align}
        \lp K(X|Y,Z) \otimes Q(Y|Z) \rp(D|z)  &:= \int K(X \in D^y|Y=y,Z=z)\, Q(Y \in dy|Z=z).
    \end{align}
\end{Thm}

Furthermore, one can show that these constructions are functorial and satisfy the monad axioms. We arrive at:

\begin{Thm}[The quasi-measurable probability monad, see \Cref{thm:Q-strong-monad}]
    \label{qms:thm:PrQm-monad}
    The triple $(\PrQm,\delta,\Mbb)$ is an affine and strong monad on the
    quasitopos $(\QMS,\times,\one)$, see \Cref{qms:thm:quasitopos}, and it admits a product of
    Markov kernels, see \Cref{qms:thm:PrQm-maps}. It is moreover commutative as
    soon as the sample space has the Fubini property \emph{and} $\PrQm(\Xcal)=\PrPf(\Xcal)$ for the
    spaces involved, see \Cref{qms:thm:commutative} and \Cref{qms:rem:comm-PrQm}; both hold on
    standard quasi-universal spaces.
\end{Thm}

\section{The Push-forward Probability Monad}
\label{sec:pf-monad}

\subsection{The Sample Space - Act 3 - The Probability Measures}

\begin{Def}[The sample space - act 3]
    \label{qms:def:sample-space-act-3}
    In the following we fix a \emph{sample space} $(\Omega,\Omega^\Omega,\Bcal_\Omega,\PrPf)$, where the triple $(\Omega,\Omega^\Omega,\Bcal_\Omega)$ satisfies the points of \Cref{qms:def:sample-space-act-1} and
    \Cref{qms:def:sample-space-act-2} and where we now also fix a non-empty set $\PrPf$  of \emph{product-compatible probability measures} $P:\, \Bcal_P \to [0,1]$ on $\Omega$, each possibly living on its own $\sigma$-algebra $\Bcal_P$ as its domain.\footnote{Both non-emptiness requirements are needed. If $\Omega=\emptyset$ then $\Omega^\Omega=\lC \id_\emptyset\rC$, every $\Xcal^\Omega$ is a singleton, $\QMS \simeq \Sets$, and acts 2 and 3 as well as \Cref{qms:asm:fubini} and \Cref{qms:asm:prod-iso} hold vacuously; but $\PrPf(\Xcal)=\emptyset$ for every $\Xcal$, so $\PrPf(\one) \not\cong \one$ and $\PrPf$ is not affine. If $\PrPf=\emptyset$ the same happens.}
    More formally, we require the following points to hold:
    \begin{enumerate}
        \item For every $\Phi \in \Omega^\Omega$ and every $P \in \PrPf$ we have the inclusion\footnote{Since $\Omega^\Omega \ins \Meas\lp(\Omega,\Bcal_\Omega),(\Omega,\Bcal_\Omega)\rp$ we actually only require: $\Bcal_\Omega \ins \Bcal_P$ for all $P \in \PrPf$; as the rest would then follow automatically. However, we sometimes construct $\PrPf$ and $\Omega^\Omega$ before $\Bcal_\Omega$. In such case, the above requirement is then a condition for $\Bcal_\Omega$, which then needs to be checked in full, as stated, see \Cref{qms:eg:sample-space-3} \Cref{qms:item:sample-space-3-2}.}:
     \begin{align}
        \Bcal_\Omega \ins \Phi_*\Bcal_P := \lC D \ins \Omega \st \Phi^{-1}(D) \in \Bcal_P \rC. \label{qms:eq:sample-space-act-3-a}
    \end{align}
    \item For every $D \in \Bcal_{\Omega \times \Omega}$ and $P \in \PrPf$ the following map:
            \begin{align}
                \Omega &\to [0,1], & \omega \mapsto P(D^\omega), \label{qms:eq:sample-space-act-3-b}
            \end{align}
  is $\Bcal_\Omega$-$\Bcal_{[0,1]}$-measurable\footnote{Note that by \Cref{qms:thm:adjunction} there is no ambiguity between measurability and quasi-measurability here. So we can check measurability just by looking on preimages of intervals of $[0,1]$.}, where $D^\omega := \lC \tilde \omega \in \Omega \st (\tilde \omega,\omega) \in D \rC$.
     \item  For every $P_1,P_2 \in \PrPf$ there exist $P \in \PrPf$ and $\Phi=(\Phi_1,\Phi_2) \in (\Omega\times\Omega)^\Omega$ such that for every $D \in \Bcal_{\Omega \times \Omega}$ we have the equality:
            \begin{align}
                (P_1 \otimes P_2)(D) :=       \int P_1(D^{\omega_2}) \, P_2(d\omega_2)
                       &\stackrel{!}{=} P\lp\lC \omega \in \Omega \st (\Phi_1(\omega),\Phi_2(\omega)) \in D \rC\rp  \label{qms:eq:sample-space-act-3-c} \\
&= P(\Phi \in D).  \notag
            \end{align}
  \end{enumerate}
\end{Def}

\begin{Rem}
    \label{qms:rem:sample-space-act-3}
    \begin{enumerate}
        \item Note that in \Cref{qms:def:sample-space-act-3} we (for now) allow every $P \in \PrPf$ to be defined on its own $\sigma$-algebra $\Bcal_P$ of subsets of $\Omega$.
        \item However, the requirements in \Cref{qms:def:sample-space-act-3} are already satisfied if:
        \begin{enumerate}
            \item $\PrPf \ins \PrQm(\Omega,\Omega^\Omega)$,  and:
            \item the restriction $\otimes:\, \PrPf \times \PrPf \to \PrQm(\Omega \times \Omega)$ already maps to $\PrPf(\Omega \times \Omega)$, the set of all push-forward measures of $P \in \PrPf$ along all $\Phi \in (\Omega \times \Omega)^\Omega$ inside $\PrQm(\Omega \times \Omega)$.
        \end{enumerate}
        Also see \Cref{qms:def:PrQm-maps} and \Cref{qms:thm:PrQm-maps}.
\item   We will later make several more technical assumptions about the sample space $(\Omega,\Omega^\Omega,\Bcal_\Omega,\PrPf)$ to ensure the validity of certain results.
    \end{enumerate}
\end{Rem}

\begin{Eg}
    Let $(\Omega,\Omega^\Omega,\Bcal_\Omega)$ be a sample space according to \Cref{qms:def:sample-space-act-1} and \Cref{qms:def:sample-space-act-2} then $\PrPf := \PrQm(\Omega,\Omega^\Omega)$ satisfies all points of \Cref{qms:def:sample-space-act-3} by \Cref{qms:rem:sample-space-act-3}.
\end{Eg}

\begin{Eg}
    \label{qms:eg:sample-space-3}
    The following examples satisfy the points of \Cref{qms:def:sample-space-act-1}, \Cref{qms:def:sample-space-act-2} and (only) the \Cref{qms:eq:sample-space-act-3-a} of \Cref{qms:def:sample-space-act-3}.
    \begin{enumerate}
        \item If $\Omega$ is a set and $\PrPf$ is a set of probability measures on $\Omega$, we can put:
            \begin{align} \Bcal_\Omega &:= \bigcap_{P \in \PrPf} \Bcal_P, & \Omega^\Omega &:=\Meas\lp(\Omega,\Bcal_\Omega),(\Omega,\Bcal_\Omega)\rp.
            \end{align}
        \item If $(\Omega,\Omega^\Omega)$ satisfies \Cref{qms:def:sample-space-act-1} and $\PrPf$ is a set of probability measures on $\Omega$, we can put: \label{qms:item:sample-space-3-2}
            \begin{align} \Bcal_\Omega &:= \bigcap_{\substack{\Phi \in \Omega^\Omega\\P \in \PrPf}} \Phi_*\Bcal_P
            = \lC D \ins \Omega \st \forall \Phi \in \Omega^\Omega.\,\forall P \in \PrPf.\, \Phi^{-1}(D) \in \Bcal_P \rC.\end{align}
        \item If $(\Omega,\Omega^\Omega,\Bcal_\Omega)$ satisfies \Cref{qms:def:sample-space-act-1}  and \Cref{qms:def:sample-space-act-2}
            we can take $\PrPf:=\Giry(\Omega,\Bcal_\Omega)$, the set of all probability measures $P$ defined on $\Bcal_P:=\Bcal_\Omega$.
\end{enumerate}
The \Cref{qms:eq:sample-space-act-3-b} and \Cref{qms:eq:sample-space-act-3-c} of \Cref{qms:def:sample-space-act-3} in those cases still would need checking.
\end{Eg}

\begin{Rem}[Act 3 forces the sample space to be infinite]
    \label{qms:rem:act3-infinite}
    \Cref{qms:def:sample-space-act-3} \Cref{qms:eq:sample-space-act-3-c} is a self-productivity
    requirement at the level of measures, and it already rules out every finite sample space with
    more than one point. Take $\Omega=\lC 0,1\rC$ with $\Bcal_\Omega$ the power set and
    $\Omega^\Omega$ all four self-maps, so that $\QMS$ is the category of simple graphs of
    \Cref{qms:eg:simpgph}, and let $\PrPf$ be all probability measures on $\Omega$. Then
    $\Bcal_\Xcal$ is the power set for every $\Xcal$ and $\PrPf(\Xcal)$ consists of the measures
    supported on at most two points that are joined by an edge. Let $\Xcal$ be the complete graph
    on $\lC a,b,c,d\rC$ and put
    \begin{align}
        P_0 &:= \tfrac12(\delta_a + \delta_b), &
        P_1 &:= \tfrac12(\delta_c + \delta_d), &
        \Pi &:= \tfrac12 \delta_{P_0} + \tfrac12 \delta_{P_1} \; \in \; \PrPf(\PrPf(\Xcal)).
    \end{align}
    Then $\Mbb(\Pi)=\tfrac14(\delta_a+\delta_b+\delta_c+\delta_d)$ is supported on four points and
    is therefore not in $\PrPf(\Xcal)$: the monad multiplication does not land in $\PrPf$. In
    accordance with this, act 3 \Cref{qms:eq:sample-space-act-3-c} fails here --- the product
    $P_1 \otimes P_2$ is uniform on four points, while every $P(\Phi \in -)$ with
    $\Phi \in (\Omega \times \Omega)^\Omega$ is supported on at most two --- and so does
    \Cref{qms:asm:prod-iso}, since $|\Omega \times \Omega| = 4 \neq 2 = |\Omega|$. This is the
    cheapest witness that the monad structure genuinely needs a hypothesis on $\Omega$.
\end{Rem}

\subsection{The Space of Push-forward Probability Measures}

\begin{Def}
    Let $(\Omega,\Omega^\Omega,\Bcal_\Omega,\PrPf)$ be a sample space according to \Cref{qms:def:sample-space-act-3}
    and $(\Xcal,\Xcal^\Omega)$ a quasi-measurable space with induced $\sigma$-algebra $\Bcal_\Xcal:=\Bcal(\Xcal^\Omega)$.
    We then define the \emph{space of push-forward probability measures} as:
    \begin{align}
        \PrPf(\Xcal) :=  \PrPf(\Xcal,\Xcal^\Omega) &:= \lC P(X):\, \Bcal_\Xcal \to [0,1] \st  X \in \Xcal^\Omega, P \in \PrPf \rC,
    \end{align}
    where the push-forward probability measure $P(X) := X_*P$ is defined for $X \in \Xcal^\Omega$, $P \in \PrPf$, $A \in \Bcal_\Xcal$ via:
    \begin{align}
         P(X \in A) := P(X)(A)  &:= P( X^{-1}(A)).
    \end{align}
    We define the set of admissible random variables on $\PrPf(\Xcal)$ as follows:
    \begin{align}
        \PrPf(\Xcal)^\Omega &:= \lC P(X|\Inp):\, \Omega \to \PrPf(\Xcal) \st  X \in \lp\Xcal^\Omega\rp^\Omega, P \in \PrPf \rC,
    \end{align}
    where we define $P(X|\Inp)$ for $X \in \lp\Xcal^\Omega\rp^\Omega$, $P \in \PrPf$,  $\omega \in \Omega$, $A \in \Bcal_\Xcal$ via:
    \begin{align}
        P(X \in A|\Inp=\omega) : =P(X|\Inp)(A|\omega)  &:= P( X(\omega)^{-1}(A)) = P\lp\lC \tilde \omega \in \Omega \st X(\omega)(\tilde \omega) \in A \rC\rp.
    \end{align}
\end{Def}

Here, $\Inp$ in $P(X|\Inp)$ is just a symbol indicating a hypothetical input variable.
Note that the probability measures $P(X)$ are restricted to $\Bcal_\Xcal$,
even though, as push-forward measures, they probably could measure more subsets of $\Xcal$.

\begin{Rem}
    Let $P \in \PrPf$ and $\Phi \in((\Omega \times \Omega)^\Omega)^\Omega$ defined by:
    $\Phi(\omega)(\tilde \omega):=(\tilde \omega, \omega)$. Then for $D \in \Bcal_{\Omega \times \Omega}$ and $\omega \in \Omega$ we get: $P(\Phi(\omega) \in D) = P(D^\omega)=(P\otimes\delta_\omega)(D)$.
    So we get:
    \begin{align}
     \lp \omega \mapsto  P(\Phi|\Inp=\omega) = P \otimes \delta_\omega \rp \in \PrPf(\Omega \times \Omega)^\Omega.
    \end{align}
\end{Rem}

Furthermore, it is easily checked that $\lp \PrPf(\Xcal),\PrPf(\Xcal)^\Omega \rp$ is again a quasi-measurable space, and we have the following:

\begin{Lem}[See \Cref{qms:thm:PrQm-maps}] Let $(\Omega,\Omega^\Omega,\Bcal_\Omega,\PrPf)$ be a sample space according to \Cref{qms:def:sample-space-act-3}.
    Then for every quasi-measurable spaces $(\Xcal,\Xcal^\Omega)$ and $(\Ycal,\Ycal^\Omega)$
    the inclusion map:
    \begin{align}
        \lp \PrPf(\Xcal), \PrPf(\Xcal)^\Omega \rp \inj \lp \PrQm(\Xcal), \PrQm(\Xcal)^\Omega \rp,
    \end{align}
    the push-forward map, the evaluation map and integration maps:
    \begin{align}
    \pf:&& \Ycal^\Xcal \times \PrPf(\Xcal) &\to \PrPf(\Ycal), & (f,P(X)) &\mapsto P(f(X))=f_*X_*P, \\
    \ev:&& \PrPf(\Xcal) \times \Bcal_\Xcal &\to [0,1], & (P(X),A) &\mapsto P(X \in A) = P(X^{-1}(A)),\\
    \int:&& [0,\infty]^\Xcal \times \PrPf(\Xcal) &\to [0,\infty],& (f,P(X)) &\mapsto \int f(X)\,dP = \int f(x)\,P(X \in dx),
\end{align}
are always well-defined and quasi-measurable.
\end{Lem}

\subsection{Fubini Theorem}

We first show how \Cref{qms:def:sample-space-act-3} plus some extra assumption on the sample space implies Fubini's theorem for all quasi-measurable spaces and non-negative functions.

\begin{Asm}[Fubini property]
    \label{qms:asm:fubini}
    Let $(\Omega,\Omega^\Omega,\Bcal_\Omega,\PrPf)$ be a sample space satisfying \Cref{qms:def:sample-space-act-3}. We say that it satisfies the \emph{Fubini property} if
    for every $P_1,P_2 \in \PrPf$ we have the equality on $\Bcal_{\Omega \times \Omega}$:
    \begin{align}
        s_*(P_2 \otimes P_1) &= P_1 \otimes P_2,
    \end{align}
    where $s$ is the (quasi-measurable) swap map:
     \begin{align*}
         s:\, \Omega \times \Omega &\to \Omega \times \Omega, & (\omega_1,\omega_2) &\mapsto (\omega_2,\omega_1).
     \end{align*}
    More explicitly, this means that for every $D \in \Bcal_{\Omega \times \Omega}$ we have:
            \begin{align}
                \int P_2(D_{\omega_1}) \, P_1(d\omega_1) &=  \int P_1(D^{\omega_2}) \, P_2(d\omega_2), \label{qms:eq:fubini}
            \end{align}
            where, following the convention used throughout --- a \emph{superscript} fixes the trailing coordinates of a set in a product and a \emph{subscript} fixes the leading ones --- we write $D_{\omega_1} := \lC \omega_2 \in \Omega \st (\omega_1,\omega_2) \in D \rC$ and
             $D^{\omega_2} := \lC \omega_1 \in \Omega \st (\omega_1, \omega_2) \in D \rC$.
\end{Asm}

With the Fubini property on the sample space we can now show that Fubini's theorem holds for all quasi-measurable spaces.

\begin{Rem}[On the independence of the two assumptions]
    \label{qms:rem:asm-independence}
    We do not know whether \Cref{qms:asm:fubini} and \Cref{qms:asm:prod-iso} are independent of
    \Cref{qms:def:sample-space-act-3}: we have no sample space that satisfies act 3 and fails
    either of them. Act 3 \Cref{qms:eq:sample-space-act-3-c} is itself a self-productivity
    condition at the level of measures, see \Cref{qms:rem:act3-infinite}, so the three are
    certainly related. We state the two properties separately because they are what the proofs
    below actually use, and because they are what one checks in practice --- not because we can
    separate them by an example.
\end{Rem}

\begin{Rem}[What the Fubini property excludes]
    \label{qms:rem:sierpinski}
    It is worth naming the classical obstruction that \Cref{qms:asm:fubini} rules out. Under the
    continuum hypothesis, well-order $[0,1]$ in order type $\omega_1$ and put
    $E:=\lC (x,y) \st x \preceq y \rC$. Every vertical section of $E$ is countable and every
    horizontal section is co-countable, so the two iterated integrals of $\I_E$ are $0$ and $1$.
    Sierpi\'nski's set $E$ is not measurable for the completion of the product of Lebesgue measure
    with itself, hence $\I_E \notin [0,\infty]^{\Omega \times \Omega}$ for the sample space of
    \Cref{qms:def:sample-space-act-4}, and the hypothesis $h \in [0,\infty]^{\Xcal \times \Ycal}$
    of the theorem below excludes it. The Fubini property is exactly the requirement that no such
    set becomes an admissible event of $\Omega \times \Omega$.
\end{Rem}

\begin{Thm}[Fubini Theorem, see \Cref{qms:thm:fubini:prf}]
    \label{qms:thm:fubini}
    Let $(\Omega,\Omega^\Omega,\Bcal_\Omega,\PrPf)$ be a sample space satisfying \Cref{qms:def:sample-space-act-3} and the Fubini property, \Cref{qms:asm:fubini}.
    Then for every quasi-measurable spaces $(\Xcal,\Xcal^\Omega)$ and $(\Ycal,\Ycal^\Omega)$, every
     $P(X) \in \PrPf(\Xcal)$, $Q(Y) \in \PrPf(\Ycal)$ and $h \in [0,\infty]^{\Xcal \times \Ycal}$,
     the following formula holds:
    \begin{align}
         \iint h(x,y)\,P(X \in dx) \, Q(Y \in dy) &= \iint h(x,y)\, Q(Y \in dy) \,P(X \in dx).
    \end{align}
\end{Thm}

We now investigate when a sample space satisfies the Fubini property.

\begin{Prp}[Fubini property criterion, see \Cref{qms:prp:fubini-crit:prf}]
    \label{qms:prp:fubini-crit}
    Let $(\Omega,\Omega^\Omega,\Bcal_\Omega,\PrPf)$ be a sample space satisfying \Cref{qms:def:sample-space-act-3}. Further assume that we have the inclusion:
    \begin{align}
        \Bcal_{\Omega\times\Omega} \ins (\Bcal_\Omega \otimes \Bcal_\Omega)_{\PrPf(\Omega\times\Omega)}.
    \end{align}
    Then $(\Omega,\Omega^\Omega,\Bcal_\Omega,\PrPf)$ satisfies the Fubini property, \Cref{qms:asm:fubini}.
\end{Prp}

\begin{Rem}
    In \Cref{qms:cor:count-sep-perf} it is shown that a sample space $(\Omega,\Omega^\Omega,\Bcal_\Omega,\PrPf)$ with countably separated $(\Omega,\Bcal_\Omega)$ and perfect $\PrPf$ satisfies the criterion in \Cref{qms:prp:fubini-crit}, which implies the Fubini property, \Cref{qms:asm:fubini}, and thus Fubini's \Cref{qms:thm:fubini}.
\end{Rem}

\begin{Def}[Composition of Markov kernels]
    \label{def:kleisli-comp}
    Let $(\Xcal,\Xcal^\Omega)$, $(\Ycal,\Ycal^\Omega)$ and $(\Zcal,\Zcal^\Omega)$ be
    quasi-measurable spaces and let $K(X|Y) \in \PrPf(\Xcal)^\Ycal$ and
    $Q(Y|Z) \in \PrPf(\Ycal)^\Zcal$ be Markov kernels. Their \emph{composition} is
    \begin{align}
        \circ:\;\PrPf(\Xcal)^\Ycal \times \PrPf(\Ycal)^\Zcal &\to \PrPf(\Xcal)^\Zcal, \\
        \lp K(X|Y) \circ Q(Y|Z)\rp (A|z) &:= \int K(X \in A|Y=y)\, Q(Y \in dy|Z=z),
    \end{align}
    for $A \in \Bcal_\Xcal$ and $z \in \Zcal$. It is the marginal of the product of Markov kernels
    of \Cref{qms:thm:product-PrPf} in the $\Ycal$-component, so it is quasi-measurable and lands in
    $\PrPf(\Xcal)$ whenever the product does. In categorical terms $\circ$ is the composition of
    the Kleisli category $\Kleisli(\PrPf)$, and $\deltabf$ is its identity; it is also written as
    the \emph{bind} of the monad, see \cite{Mog91}.
\end{Def}

\subsection{Countably Separated Spaces}

The $\sigma$-algebra $\Bcal_\Xcal$ of a quasi-measurable space $(\Xcal,\Bcal_\Xcal)$ can, in general,
become unmanageably big. Already the $\sigma$-algebra $\Bcal_{\Omega \times \Omega}$ of the sample space product might be much bigger than the product $\sigma$-algebra $\Bcal_\Omega \otimes \Bcal_\Omega$ of measurable spaces:
\begin{align*}
        \Bcal_\Omega \otimes \Bcal_\Omega \ins \Bcal_{\Omega \times \Omega}.
\end{align*}
However, we can regain control over such $\sigma$-algebras $\Bcal_\Xcal$ if they are \emph{countably separated} and the set of probability measures $\PrPf$ on the sample space $\Omega$ only contains \emph{perfect probability measures}. Furthermore, things become even more elegant if the $\sigma$-algebra $\Bcal_\Omega$ is \emph{$\PrPf(\Omega)$-complete} and/or itself countably separated.
We will state this in the following more precisely.

\begin{Def}
    A measurable space $(\Xcal,\Bcal_\Xcal)$ is called \emph{countably separated} if there exists a countable subset $\Ecal \ins \Bcal_\Xcal$
    that \emph{separates the points of $\Xcal$}, i.e.\ for all $x_1,x_2 \in \Xcal$ with $x_1\neq x_2$ there exists an $A \in \Ecal$ such that
    either $x_1 \in A$ and $x_2 \notin A$ or $x_1 \notin A$ and $x_2 \in A$.
    A quasi-measurable space $(\Xcal,\Xcal^\Omega)$ is called \emph{countably separated} if $(\Xcal,\Bcal(\Xcal^\Omega))$ is.
\end{Def}

\begin{Thm}
    \label{qms:thm:count-sep-perf}
    Let $(\Omega,\Omega^\Omega,\Bcal_\Omega,\PrPf)$ be a sample space according to \Cref{qms:def:sample-space-act-3}.
    Assume that every probability measure $P \in \PrPf$ is \emph{perfect}.
    Let $(\Xcal,\Xcal^\Omega)$ be a \emph{countably separated} quasi-measurable space and $\Ecal \ins \Bcal(\Xcal^\Omega)$ be any countable subset that separates the points of $\Xcal$.
    Then we have the inclusion of $\sigma$-algebras:
    \begin{align} \Bcal(\Xcal^\Omega)&\ins \lp \Ecal \rp_{\PrPf(\Xcal)}.   \end{align}
    Furthermore, equality holds if $\Bcal_\Omega$ is $\PrPf(\Omega)$-complete, i.e.\ if $\Bcal_\Omega=\lp\Bcal_\Omega\rp_{\PrPf(\Omega)}$.\\
\end{Thm}

\begin{Rem}
    Note that for countably separated $(\Xcal,\Xcal^\Omega)$ and perfect $\PrPf$ \Cref{qms:thm:count-sep-perf} gives us control over the induced $\sigma$-algebra on the product:
\begin{align}\Bcal_{\Xcal \times \Xcal} &\ins \lp\Bcal_\Xcal \otimes \Bcal_\Xcal\rp_{\PrPf(\Xcal \times \Xcal)}, \label{qms:eq:count-sep-perf}
\end{align}
with equality if $\Bcal_\Omega$ is $\PrPf(\Omega)$-complete.
\end{Rem}

A special case of interest emerges when $(\Omega,\Bcal_\Omega)$ is itself countably separated:

\begin{Cor}
    \label{qms:cor:count-sep-perf}
     Let $(\Omega,\Omega^\Omega,\Bcal_\Omega,\PrPf)$ be a sample space according to \Cref{qms:def:sample-space-act-3} such that $(\Omega,\Bcal_\Omega)$ is countably separated and  $\PrPf$ only consists of perfect probability measures. Then we have the inclusion:
\begin{align}\Bcal_{\Omega \times \Omega} &\ins \lp\Bcal_\Omega \otimes \Bcal_\Omega\rp_{\PrPf(\Omega \times \Omega)}, \label{qms:eq:count-sep-perf-om}
\end{align}
with equality if $\Bcal_\Omega$ is $\PrPf(\Omega)$-complete.
In particular, $(\Omega,\Omega^\Omega,\Bcal_\Omega,\PrPf)$ satisfies the Fubini property, \Cref{qms:asm:fubini}, by \Cref{qms:prp:fubini-crit}.
\end{Cor}

\begin{Rem}
    \label{qms:rem:count-sep-perf-cond}
    \Cref{qms:cor:count-sep-perf} indicates that, for the theory of quasi-measurable spaces, we might want to restrict to sample spaces $(\Omega,\Omega^\Omega,\Bcal_\Omega,\PrPf)$ such that:
    \begin{enumerate}
        \item $(\Omega,\Bcal_\Omega)$ is countably separated,
        \item $\PrPf$ is perfect, i.e.\ only contains perfect probability measures,
        \item $\Bcal_\Omega$ is $\PrPf(\Omega)$-complete, i.e.\ $\Bcal_\Omega = (\Bcal_\Omega)_{\PrPf(\Omega)}$.
    \end{enumerate}
\end{Rem}

\subsection{Products of Markov Kernels}

In the following we will state the most important implication that result from \Cref{qms:def:sample-space-act-3}. The corresponding proofs can be found in \Cref{qms:sec:push-for-prob-mon}.
After that we will investigate how to construct sample spaces that satisfy those assumptions.

\begin{Asm}[Product isomorphism property]
    \label{qms:asm:prod-iso}
     Let $(\Omega,\Omega^\Omega)$ be a sample space satisfying \Cref{qms:def:sample-space-act-1}.
     We say that $(\Omega,\Omega^\Omega)$ satisfies the \emph{product isomorphism property}
     if there exists an isomorphism of quasi-measurable spaces:
    \begin{align}
        \Theta: \; \Omega & \cong \Omega \times \Omega. \label{qms:eq:prod-iso}
     \end{align}
\end{Asm}

\begin{Thm}[Product of Markov kernels, see \Cref{qms:thm:product-PrPf:prf}]
    \label{qms:thm:product-PrPf}
    Let $(\Omega,\Omega^\Omega,\Bcal_\Omega,\PrPf)$ be a sample space satisfying \Cref{qms:def:sample-space-act-3} and also \Cref{qms:asm:prod-iso}: $\Omega  \cong \Omega \times \Omega$.

    Let  $(\Xcal,\Xcal^\Omega)$,  $(\Ycal,\Ycal^\Omega)$, $(\Zcal,\Zcal^\Omega)$ be any quasi-measurable spaces.
       Then the product of Markov kernels:
    \begin{align}
        \otimes:\; \PrPf(\Xcal)^{\Ycal \times \Zcal} \times \PrPf(\Ycal)^\Zcal & \to \PrPf(\Xcal \times \Ycal)^\Zcal,
    \end{align}
    given for $D \in \Bcal_{\Xcal \times \Ycal}$ and $z \in \Zcal$ via the formula:
    \begin{align}
        \lp K(X|Y,Z) \otimes Q(Y|Z)\rp(D|z) & := \int K(X \in D^y |Y=y,Z=z)\, Q(Y \in dy|Z=z),
    \end{align}
    is a well-defined and quasi-measurable map.
\end{Thm}

\begin{Thm}[The push-forward probability monad, see \Cref{qms:thm:PrPf-monad:prf}]
    \label{qms:thm:PrPf-monad}
    Let $(\Omega,\Omega^\Omega,\Bcal_\Omega,\PrPf)$ be a sample space satisfying \Cref{qms:def:sample-space-act-3} and also \Cref{qms:asm:prod-iso}: $\Omega  \cong \Omega \times \Omega$.

    Then the triple $(\PrPf,\delta,\Mbb)$ is an affine and  strong monad on the quasitopos $(\QMS,\times,\one)$, see \Cref{qms:thm:quasitopos}.
    If, furthermore, the sample space has the Fubini property \Cref{qms:asm:fubini}, then
    $(\PrPf,\delta,\Mbb)$ is also \emph{commutative}, equivalently \emph{symmetric monoidal}, see
    \Cref{qms:thm:commutative}.
\end{Thm}

The proofs can be found in \Cref{qms:sec:push-for-prob-mon}.

\begin{Eg} Let $\Omega=[0,1]$ with Borel $\sigma$-algebra $\Bcal_\Omega$ and
    $\Omega^\Omega:=\Meas(\Omega,\Omega)$
    and $\PrPf=\lC P \rC$, where $P$ is the uniform distribution on $[0,1]$.
    \begin{enumerate}
        \item Consider $\twop=\lC 0,1 \rC$ with $\twop^\Omega=\Sets(\Omega,\twop)$, which gives: $\Bcal_{\twop}=\lC \twop,\emptyset \rC$.
        \item[] Let $V \notin \Bcal_\Omega$ be the non-Lebesgue-measurable Vitali set and put $X:=\I_V \in \twop^\Omega$.
        \item[] Then $\PrPf(\twop) = \lC P(X) \rC$ is trivial, with $P(X \in \twop)=1$  and $P(X \in \emptyset)=0$.
        \item Consider $\twoc=\lC 0,1 \rC$ with $\twoc^\Omega=\lC X_0, X_1 \rC$, the two constant maps.
        \item[] With this we get: $\Bcal_{\twoc}=\lC \twoc,\lC0\rC,\lC1\rC, \emptyset \rC$ and
            $\PrPf(\twoc) = \lC \delta_0=P(X_0),\delta_1=P(X_1) \rC$.
        \item Consider $\twob=\lC 0,1 \rC$ with $\twob^\Omega=\Meas(\Omega,\twob)$, which implies: $\Bcal_{\twob}=\lC \twob,\lC0\rC,\lC1\rC, \emptyset \rC$.
        \item[] For $p \in [0,1]$ let $X_p:=\I_{[0,p]} \in \twob^\Omega$.
         Then $\PrPf(\twob) = \lC \mathrm{Ber}(p) = P(X_p) \st p \in [0,1] \rC$, the set of all Bernoulli distributions.
    \end{enumerate}

\end{Eg}

\section{Construction of Well-Behaved Sample Spaces}
\label{sec:sample-spaces}

\begin{Thm}[Sample spaces from measurable spaces, see \Cref{qms:thm:spl-spc-gen:prf}]
    \label{qms:thm:spl-spc-gen}
    Let $(\Omega,\Bcal_\Omega)$ be any measurable space such that we have an isomorphism of measurable spaces:
    \begin{align}
        \Theta:\, (\Omega,\Bcal_\Omega) \cong (\Omega \times \Omega, \Bcal_\Omega\otimes\Bcal_\Omega).
    \end{align}
    We then further put:
    \begin{align}
        \Omega^\Omega&:= \Meas\lp (\Omega,\Bcal_\Omega),(\Omega,\Bcal_\Omega) \rp, &
\PrPf &:= \Giry(\Omega,\Bcal_\Omega).
    \end{align}
    Then $(\Omega,\Omega^\Omega,\Bcal_\Omega,\PrPf)$ is a sample space satisfying all the requirements of
    \Cref{qms:def:sample-space-act-1}, \Cref{qms:def:sample-space-act-2}, \Cref{qms:def:sample-space-act-3} and \Cref{qms:asm:fubini} and \Cref{qms:asm:prod-iso}.
    Furthermore, we have that:
    \begin{align}
        \Bcal_{\Omega \times \Omega} &= \Bcal_\Omega \otimes \Bcal_\Omega.
    \end{align}
\end{Thm}

\begin{Thm}[Sample spaces from completing measurable spaces, see \Cref{qms:thm:spl-spc-cmpl:prf}]
    \label{qms:thm:spl-spc-cmpl}
    Let $(\Omega,\Bcal_\otimes)$ be any measurable space such that we have an isomorphism of measurable spaces:
    \begin{align}
        \Theta:\, (\Omega,\Bcal_\otimes) \cong (\Omega \times \Omega, \Bcal_\otimes\otimes\Bcal_\otimes).
    \end{align}
    Now choose any subset $\PrPf \ins \Giry(\Omega,\Bcal_\otimes)$ such that the following implications hold:
    \begin{align}
        P \in \PrPf \quad \land \quad  \Phi \in \Meas\lp (\Omega,(\Bcal_\otimes)_\PrPf),(\Omega,\Bcal_\otimes) \rp &\quad \implies \quad \bar P(\Phi)|_{\Bcal_\otimes} \in \PrPf, \label{qms:eq:PrPf-compl-cond} \\
        P_1, P_2 \in \PrPf &\quad \implies \quad \Theta_*^{-1}(P_1 \otimes P_2)|_{\Bcal_\otimes} \in \PrPf, \label{qms:eq:PrPf-prod-cond}
    \end{align}
    where we identify $P \in \PrPf \ins \Giry(\Omega,\Bcal_\otimes)$ with its unique completion $\bar P$ living on $(\Bcal_\otimes)_P \sni (\Bcal_\otimes)_\PrPf$.
    We then further choose:
    \begin{align}
        \Bcal_\Omega  &:= (\Bcal_\otimes)_\PrPf, &
        \Omega^\Omega &:= \Meas\lp (\Omega,\Bcal_\Omega),(\Omega,\Bcal_\Omega) \rp.
    \end{align}
    Then $(\Omega,\Omega^\Omega,\Bcal_\Omega,\PrPf)$ is a sample space satisfying all the requirements of
    \Cref{qms:def:sample-space-act-1}, \Cref{qms:def:sample-space-act-2}, \Cref{qms:def:sample-space-act-3} and \Cref{qms:asm:fubini} and \Cref{qms:asm:prod-iso}.
    Furthermore, we have that $\Bcal_\Omega$ is $\PrPf(\Omega)$-complete, and:
    \begin{align}
        \Fcal(\Bcal_\otimes) &= \Fcal(\Bcal_\Omega) = \Omega^\Omega, &
        \Bcal_{\Omega \times \Omega} &= (\Bcal_\Omega \otimes \Bcal_\Omega)_{\PrPf(\Omega \times \Omega)}.
    \end{align}
\end{Thm}

\begin{Eg}
    In the setting of \Cref{qms:thm:spl-spc-cmpl} the following choices for $\PrPf$ satisfy \Cref{qms:eq:PrPf-compl-cond} and \Cref{qms:eq:PrPf-prod-cond}:
    \begin{enumerate}
        \item $\PrPf := \Giry(\Omega,\Bcal_\otimes)$,
        \item $\PrPf := \lC P \in \Giry(\Omega,\Bcal_\otimes) \st P \text{ is perfect} \rC$,
        \item $\PrPf := \lC P \in \Giry(\Omega,\Bcal_\otimes) \st P \text{ is countably compact} \rC$.
    \end{enumerate}
    Indeed, the completions, push-forwards and products of perfect/countably compact probability measures
    are again perfect, countably compact, resp., by \cite{Fremlin} 451E, 451G, 451I, 452R.
    Note that countably compact probability measures are always perfect by \cite{Fremlin} 451C.
    So the latter two choices for $\PrPf$ also satisfy the conditions of \Cref{qms:thm:count-sep-perf}.
\end{Eg}

\begin{Cor}[Countably separated product sample space for perfect probability distributions, see \Cref{qms:cor:sample-space-count-sep-perf:prf}]
    \label{qms:cor:sample-space-count-sep-perf}
    Let $\Omega_0$ be a non-empty set and $\Ecal_0$ a countable set of subsets of $\Omega_0$ that separates the points of $\Omega_0$.
    We then put:
    \begin{align}
        \Omega &:= \prod_{n \in \N} \Omega_0,
        & \Ecal &:= \lC \pr_n^{-1}(A) \ins \Omega \st n \in \N, A \in \Ecal_0 \rC,
    \end{align}
 where the intersection is w.r.t.\ the set of all complete, perfect probability measures on $\Omega$ that can measure $\Ecal$:
    \begin{align}
        \PrPf &:=\lC P:\, \Bcal_P \to [0,1] \st P \text{ complete, perfect probability measure on } \Omega \text{ with } \Ecal \ins \Bcal_P \rC. \end{align}
    We further put:
    \begin{align}
        \Bcal_\Omega &:= \bigcap_{P \in \PrPf} \Bcal_P, &
        \Omega^\Omega&:= \Meas\lp (\Omega,\Bcal_\Omega),(\Omega,\Bcal_\Omega) \rp.
    \end{align}
    Then $\Ecal \ins \Bcal_\Omega$ is a countable set of subsets of $\Omega$ that separates the points of $\Omega$
    and $(\Omega,\Omega^\Omega,\Bcal_\Omega,\PrPf)$ is a sample space satisfying all the requirements of
    \Cref{qms:def:sample-space-act-1}, \Cref{qms:def:sample-space-act-2}, \Cref{qms:def:sample-space-act-3} and \Cref{qms:asm:fubini} and \Cref{qms:asm:prod-iso}. The ``zip-locker'' bijection provides the needed isomorphism of quasi-measurable spaces:
    \begin{align}
        \Theta=\Theta_0 \times \Theta_1:\, \Omega \cong \Omega \times \Omega, \qquad (\omega_n)_{n \in \N} \mapsto \lp (\omega_{2n})_{n \in \N}, (\omega_{2n+1})_{n \in \N} \rp.
    \end{align}
    Furthermore, all conditions mentioned in \Cref{qms:rem:count-sep-perf-cond} are satisfied as well.
\end{Cor}

We close with the example that the rest of this paper runs on. It builds two sample spaces out of the construction above, one of which is the choice made in \Cref{sec:qus}.

\begin{Eg}[The Hilbert cube]
    \label{qms:def:univ-hilbert-cube}\label{qms:eg:hilbert-cube}
    Let $\Omega := \prod_{ n \in \N} [0,1]$ be the \emph{Hilbert cube}, the countably infinite product  of unit intervals $[0,1]$, and $\Ecal_\Omega$ its Borel-$\sigma$-algebra,
    which is generated by (finite) cubes.
    Let $\PrPf:=\lC U \rC$, where $U := \bigotimes_{n \in \N} U[0,1]$ is the uniform distribution on $\Omega$
    defined on the completed $\sigma$-algebra $U:\, \Bcal_U:=\lp \Ecal_\Omega\rp_U \to [0,1]$.
We consider the following two choices for $\Bcal_\Omega$, namely:
 \begin{enumerate*}
     \item $\Bcal_\Omega:=\Ecal_\Omega$, \label{qms:item:hilbert-cube}
     \item $\Bcal_\Omega:=\lp\Ecal_\Omega\rp_\Giry$, the universal completion of $\Ecal_\Omega$. \label{qms:item:univ-hilbert-cube}
 \end{enumerate*}
 In both cases, we put $\Omega^\Omega := \Meas\lp (\Omega,\Bcal_\Omega),(\Omega,\Bcal_\Omega)\rp$.
 We will refer to the sample space given by the quadrupels $\lp \Omega,\Omega^\Omega,\Bcal_\Omega,\PrPf \rp$ in \Cref{qms:item:hilbert-cube} as the \emph{Borel Hilbert cube} and
 the one in \Cref{qms:item:univ-hilbert-cube} as the \emph{universal Hilbert cube}, resp.

 \emph{Hilbert's hotel}, the ``zip locker map'':
 \begin{align*}
     \Theta= \Theta_0 \times \Theta_1:\, \Omega &\to \Omega \times \Omega, &\Theta\lp(e_n)_{n \in \N}\rp := \lp (e_{2k})_{k \in \N},(e_{2m+1})_{m \in \N} \rp,
 \end{align*}
 is a Borel isomorphism, which by \Cref{qms:eq:adj-prod} induces an isomorphism of quasi-measurable spaces (for both choices of $\Bcal_\Omega$)
 with components: $\Theta_0,\Theta_1\in \Omega^\Omega$.
 For (finite) cylinder sets $A_0,A_1 \in \Ecal_\Omega$ and $\omega \in \Omega$ we get the equations:
 \begin{align}
     U\lp\lp A_0 \times A_1 \rp^\omega \rp &= U(A_0) \cdot \I_{A_1}(\omega),\\
     U(\Theta_0 \in A_0,\Theta_1 \in A_1) & = U(A_0) \cdot U(A_1) \\&= U(A_1) \cdot U(A_0).
 \end{align}
Together with Dynkin's Lemma, see \cite{Kle20} Thm.\ 1.19, this verifies for $\Bcal_\Omega=\Ecal_\Omega$ \Cref{qms:def:sample-space-act-3} \Cref{qms:eq:sample-space-act-3-b}, \Cref{qms:eq:sample-space-act-3-c}: $U(\Theta_0,\Theta_1)  = U \otimes U$, and, \Cref{qms:asm:fubini}.
Using \Cref{qms:thm:markov-kernel-completions} the same then holds also for its universal completion $\Bcal_\Omega=\lp \Ecal_\Omega\rp_\Giry$.

In both cases we get the equalities of sets:
\begin{align}
    \PrPf(\Omega) &= \Giry(\Omega,\Bcal_\Omega), & \PrPf(\Omega)^\Omega &= \Meas\lp (\Omega,\Bcal_\Omega),\Giry(\Omega,\Bcal_\Omega)\rp.
\end{align}
Those can be proved by standard arguments using a Borel isomorphism $\Omega \cong \R$ and conditional quantile functions for reparameterizing Markov kernels, see \cite{Cen82,Kal21,For21} Appendix G.
\end{Eg}

\section{Quasi-Universal Spaces: Probability Theory over the Universal Hilbert Cube}
\label{sec:qus}

Up to this point the sample space was kept general: \Cref{qms:def:sample-space-act-1} fixed its
random variables, \Cref{qms:def:sample-space-act-2} its $\sigma$-algebra, \Cref{qms:def:sample-space-act-3}
its probability measures, and \Cref{sec:sample-spaces} showed how to build sample spaces satisfying all
three acts together with the Fubini property, \Cref{qms:asm:fubini}, and the product isomorphism property,
\Cref{qms:asm:prod-iso}.

In this section we commit to one concrete choice --- the \emph{universal Hilbert cube} --- and show what
the theory looks like once we do. The resulting category is the category of \emph{quasi-universal spaces}
$\QUS$. It is the setting in which the induced $\sigma$-algebras become explicitly describable, in which
the four push-forward monads collapse into one, and in which disintegration, Kolmogorov extension and
de Finetti theorems become available. It is also the setting for the causal models of
\Cref{sec:causal-models}.

\subsection{The Sample Space - Act 4 - The Universal Hilbert Cube}
\label{qms:sec:act4}

\begin{Def}[The sample space - act 4]
    \label{qms:def:sample-space-act-4}
    The \emph{universal Hilbert cube} is the sample space $(\Omega,\Omega^\Omega,\Bcal_\Omega,\PrPf)$ given by:
    \begin{align}
        \Omega &:= [0,1]^\N, &
        \Bcal_\otimes &:= \Bcal_{[0,1]^\N} = \bigotimes_{n \in \N} \Bcal_{[0,1]}, \\
        \Bcal_\Omega &:= \lp \Bcal_\otimes \rp_\Giry, &
        \Omega^\Omega &:= \Meas\lp (\Omega,\Bcal_\Omega),(\Omega,\Bcal_\Omega)\rp,
    \end{align}
    where $\Bcal_\otimes$ is the Borel $\sigma$-algebra of the Hilbert cube and
    $\Bcal_\Omega=\lp\Bcal_\otimes\rp_\Giry$ its universal completion, i.e.\ the $\sigma$-algebra of all
    universally measurable subsets of $\Omega$, and where:
    \begin{align}
        \PrPf &:= \Giry(\Omega,\Bcal_\Omega),
    \end{align}
    the set of \emph{all} probability measures on $(\Omega,\Bcal_\Omega)$.
    Restriction along $\Bcal_\otimes \ins \Bcal_\Omega$ identifies these with the Borel probability
    measures on the Hilbert cube, each Borel measure having a unique extension to its universal
    completion; the latter is the form in which \Cref{qms:thm:spl-spc-cmpl} produces them, see
    \Cref{qms:eg:hilbert-cube}.
    The corresponding category of quasi-measurable spaces is called the category of
    \emph{quasi-universal spaces} $\QUS$; its objects are \emph{quasi-universal spaces}.
\end{Def}

\begin{Thm}[The universal Hilbert cube is a well-behaved sample space, see \Cref{qms:thm:act4-well-behaved:prf}]
    \label{qms:thm:act4-well-behaved}
    The universal Hilbert cube of \Cref{qms:def:sample-space-act-4} satisfies
    \Cref{qms:def:sample-space-act-1}, \Cref{qms:def:sample-space-act-2} and
    \Cref{qms:def:sample-space-act-3}, as well as the Fubini property \Cref{qms:asm:fubini}
    and the product isomorphism property \Cref{qms:asm:prod-iso}.
    Moreover $(\Omega,\Bcal_\Omega)$ is countably separated, every $P \in \PrPf$ is perfect,
    $\Bcal_\Omega$ is $\PrPf(\Omega)$-complete, and:
    \begin{align}
        \Fcal(\Bcal_\otimes) &= \Fcal(\Bcal_\Omega) = \Omega^\Omega, &
        \Bcal_{\Omega \times \Omega} &= \lp \Bcal_\Omega \otimes \Bcal_\Omega\rp_{\PrPf(\Omega\times\Omega)}.
    \end{align}
    So all conditions collected in \Cref{qms:rem:count-sep-perf-cond} hold.
\end{Thm}

This is the universal Hilbert cube of \Cref{qms:eg:hilbert-cube}, whose two variants were built there directly from \Cref{qms:thm:spl-spc-cmpl}. That example also works the choice out for a single probability measure, and records that nothing is lost by it.

\begin{Rem}[One measure is enough]
    \label{qms:rem:act4-one-measure}
    By \Cref{thm:qus-S=P-in-R} the whole of $\PrPf(\Xcal)$ is already reached from the single uniform
    distribution $U=\bigotimes_{n \in \N} U[0,1]$ on the Hilbert cube:
    \begin{align}
        \PrPf(\Xcal) &= \lC X_*U \st X \in \Xcal^\Omega \rC,
    \end{align}
    for every quasi-universal space $(\Xcal,\Xcal^\Omega)$.
    So one may equally well take $\PrPf:=\lC U \rC$ in \Cref{qms:def:sample-space-act-4}: the induced spaces
    of probability measures, and hence the whole theory, are the same. Which of the two presentations is
    more convenient depends on the argument at hand.
\end{Rem}

\begin{Rem}[Other choices of the sample space]
    \label{qms:rem:act4-other-choices}
    Since there are Borel isomorphisms $[0,1]^\N \cong \R^\N \cong \R$ between uncountable Polish spaces, we
    could equally well have taken:
    \begin{align}
        (\Omega,\Bcal_\Omega) &= \lp \R^\N, \lp \bigotimes_{n \in \N} \Bcal_\R\rp_\Giry \rp,
        &\text{or}&&
        (\Omega,\Bcal_\Omega) &= \lp \R, \lp\Bcal_\R\rp_\Giry\rp,
    \end{align}
    or any other uncountable Polish space endowed with its $\sigma$-algebra of all universally measurable
    subsets: all of these give \emph{isomorphic} categories $\QUS$, the isomorphism being
    transport along any Borel isomorphism between the two sample spaces. It is not canonical ---
    there is no preferred Borel isomorphism $[0,1]^\N \cong \R$ --- which is why we fix one sample
    space in \Cref{qms:def:sample-space-act-4} and use the others only as descriptions of it.
    The reason for presenting the Hilbert cube, or $\R^\N$, rather than $\R$ is that the isomorphism
    $\Omega \times \Omega \cong \Omega$ is then the elementary zip-locker map rather than a Borel
    isomorphism. A countable product also makes it easy to model countably many random experiments with this
    one sample space, which could capture all stochasticity in reality.

    The only difference to the sample space $(\R,\Bcal_\R,\Fcal(\Bcal_\R))$ of the category of
    \emph{quasi-Borel spaces} $\QBS$ of \cite{Heu17} is the use of the universal completion in place of the
    Borel $\sigma$-algebra. That difference is what buys us the explicit description of the induced
    $\sigma$-algebras in \Cref{lem:qus-universally-complete} and \Cref{thm:countably-separated}: the
    push-forward $X_*\Bcal_\R$ pushes $\Bcal_\R$ onto the image $X(\Omega)$ and then extends by null sets
    \emph{outside} the image, without necessarily containing the null sets \emph{inside} it. Working with
    complete $\sigma$-algebras $(\Bcal_\R)_P$ and their push-forwards $X_*(\Bcal_\R)_P$, which are then
    complete as well, harmonizes the situation.

    A second difference is that we do not require quasi-universal spaces to be \emph{\patchable}, see
    \Cref{def:patchable}, whereas quasi-Borel spaces are \patchable\ by definition, \cite{Heu17} Def.\ 7.
    The only reason to impose \patchability\ here would be to preserve countable coproducts in the passage
    from measurable spaces to quasi-universal spaces, see \Cref{lem:patch-coprod}; for everything else it is
    simpler and more general to proceed without it. \Patchability\ is treated on its own in
    \Cref{sec:patchable-qms}.
\end{Rem}

\begin{Rem}[Which monads survive]
    \label{qms:rem:act4-monads}
    The isomorphism $\Omega \times \Omega \cong \Omega$ of \Cref{qms:thm:act4-well-behaved} makes
    both $(\PrPf,\delta,\Mbb)$ and $(\PrQm,\delta,\Mbb)$ strong on the quasitopos
    $(\QUS,\times,\one)$, see \Cref{qms:thm:PrPf-monad} and \Cref{qms:thm:PrQm-monad}, and the
    Fubini property makes them commutative, see \Cref{qms:thm:commutative} --- for $\PrQm$ this
    uses $\PrQm(\Xcal)=\PrPf(\Xcal)$, which holds on standard quasi-universal spaces by
    \Cref{thm:univ-qus-GPQKR}. These two are the only
    probability monads that occur in what follows.\footnote{The three further monads of
    \Cref{qms:sec:app-other-monads} all collapse onto $\PrPf$ here, and on standard
    quasi-universal spaces even $\PrQm$ and $\Giry$ join them; see \Cref{thm:qus-S=P-in-R} and
    \Cref{thm:univ-qus-GPQKR}. This is why they play no role in the main text.}
\end{Rem}

\subsection{Countably Separated Quasi-Universal Spaces}

In this subsection we study \emph{countably separated} quasi-universal spaces.
The reason is that those spaces have a very convenient descriptions for their $\sigma$-algebra.
It turns out that such countably separated quasi-universal spaces $(\Xcal,\Xcal^\Omega)$ will satisfy:
\[ \Bcal(\Xcal^\Omega) = \lp \Ecal \rp_{\PrPf(\Xcal,\Xcal^\Omega)},  \]
for every countable subset $\Ecal \ins \Bcal(\Xcal^\Omega)$ that separates the points of $\Xcal$.
Surprisingly, that means that those $\sigma$-algebras are then automatically countably \emph{generated} up to some form of completion. Many things will simplify for those spaces.

\subsubsection{Countably Separated Spaces}

Here we will go through the definition of \emph{countably separated} spaces and variants thereof.

Countable separatedness was defined in \Cref{sec:pf-monad}; for quasi-universal spaces
we will need two variations of it.

\begin{Def}[Separated and universally countably separated]
    \label{qms:def:sep-variants}
    Let $\Xcal$ be a set and $\Ecal$ a set of subsets of $\Xcal$.
    \begin{enumerate}
        \item $(\Xcal,\Ecal)$, or just $\Ecal$, is called \emph{separated} if $\Ecal$
            separates the points of $\Xcal$, i.e.\ for all $x_1 \neq x_2$ in $\Xcal$ there
            is an $A \in \Ecal$ containing exactly one of them.
        \item $(\Xcal,\Ecal)$ is called \emph{countably separated} if some countable
            $\Ecal' \ins \Ecal$ separates the points of $\Xcal$. For $\Ecal=\Bcal_\Xcal$
            this is the notion of \Cref{sec:pf-monad}.
        \item $(\Xcal,\Ecal)$, or just $\Ecal$, is called \emph{universally countably
            separated} if $(\Xcal,(\Ecal)_\Giry)$ is countably separated, i.e.\ if the
            universal completion of $\Ecal$ contains a countable family that separates
            the points of $\Xcal$.
        \item A quasi-measurable space $(\Xcal,\Xcal^\Omega)$ is called
            \emph{(countably) separated}, \emph{universally countably separated}, resp.,
            if $(\Xcal,\Bcal(\Xcal^\Omega))$ is so in the sense above.
        \item A quasi-measurable space $(\Xcal,\Xcal^\Omega)$ is called
            \emph{quasi-separated} if its diagonal is an event of the product:
            \[\Delta_\Xcal:=\lC (x,x) \in \Xcal \times \Xcal\st x \in \Xcal \rC
              \in \Bcal_{\Xcal \times \Xcal}. \]
    \end{enumerate}
\end{Def}

\begin{Rem}[See \Cref{qms:qus:rem:auto-4:prf}]
    \label{qms:qus:rem:auto-4}
    \begin{enumerate}
        \item Every countably separated (quasi-)measurable space is separated and universally countably separated.
        \item Every countably separated quasi-measurable space is also quasi-separated.
        \item The countable product of countably separated (quasi-)measurable spaces is countably separated.
    \end{enumerate}

\end{Rem}

\begin{Def}[Separated quotient of a (quasi-)measurable space]
    \label{def:separated-quotient}
    Let $(\Xcal,\Bcal_\Xcal)$, $(\Xcal,\Xcal^\Omega)$, resp., be a (quasi-)measurable space.
    Then the quotient:
            \[ t:\; \Xcal \srj \tilde \Xcal := \Xcal/\mathord{\sim},\]
            together with the quotient $\sigma$-algebra $\Bcal_{\tilde \Xcal}:=t_*\Bcal_\Xcal$,
            the quotient functions $\tilde \Xcal^\Omega:=t \circ \Xcal^\Omega$, resp.,
    will be called the \emph{separated quotient} of $(\Xcal,\Bcal_\Xcal)$, $(\Xcal,\Xcal^\Omega)$, resp.
    Here the equivalence relation is given by:
    \[x_1 \sim x_2 \quad :\iff\quad \forall A \in \Bcal_\Xcal.\quad \lC x_1,x_2 \rC \ins A \; \lor\; \lC x_1,x_2 \rC \ins A^\cmpl. \]
    Also note for the quasi-measurable space $(\Xcal,\Xcal^\Omega)$ that we take $\Bcal_\Xcal:=\Bcal(\Xcal^\Omega)$ and that:
    \[ t_*\Bcal_\Xcal= t_*\Bcal(\Xcal^\Omega) = \Bcal(t \circ \Xcal^\Omega) = \Bcal(\tilde \Xcal^\Omega) = \Bcal_{\tilde \Xcal}. \]
\end{Def}

\begin{Rem}
    \label{rem:separated-quotient}
    Let $(\Xcal,\Bcal_\Xcal)$, $(\Xcal,\Xcal^\Omega)$, resp., be a (quasi-)measurable space.
    Then the quotient map of the separated quotient:
    \[t:\; \Xcal \srj \tilde \Xcal := \Xcal/\mathord{\sim},\]
    induces a bijection:
    \[ t^*:\; \Bcal_{\tilde \Xcal} \to \Bcal_\Xcal, \qquad \tilde A \mapsto t^{-1}(\tilde A).  \]
    Furthermore, every section $s$ of $t$, which exists by the axiom of choice:
    \[s:\; \tilde \Xcal \to \Xcal,\qquad t(s(\tilde x)) = \tilde x,\]
    is measurable.
\end{Rem}

The following Lemmas might be of interest to understand the role of
countably separated quasi-universal spaces better:

\subsubsection{The Sigma-Algebra of Countably Separated Quasi-Universal Spaces}

In this subsection we will highlight the structure of the induced $\sigma$-algebra for countably separated quasi-universal spaces.

The next theorem is \Cref{qms:thm:count-sep-perf} for the sample space of
\Cref{qms:def:sample-space-act-4}: by \Cref{qms:thm:act4-well-behaved} that sample space
is countably separated, its probability measures are perfect and $\Bcal_\Omega$ is
complete, so the inclusion there becomes an equality here.

\begin{Thm}[See \Cref{thm:countably-separated:prf}]
    \label{thm:countably-separated}
Let $(\Xcal,\Xcal^\Omega)$ be a countably separated quasi-universal space.
    Then we have:
    \[ \Bcal(\Xcal^\Omega) = \lp \Ecal \rp_{\PrPf(\Xcal,\Xcal^\Omega)},  \]
    for every countable subset $\Ecal \ins \Bcal(\Xcal^\Omega)$ that separates the points of $\Xcal$,
    where:   \[
       \PrPf(\Xcal,\Xcal^\Omega):=\lC P(X)=X_*P:\; \Bcal(\Xcal^\Omega) \to [0,1]\st  X \in \Xcal^\Omega, P \in
   \Giry(\Omega,\Bcal_\Omega)\footref{fn:G-Q-1}\rC.\]

\end{Thm}

\begin{Cor}[See \Cref{qms:qus:cor:auto-13:prf}]
    \label{qms:qus:cor:auto-13}
    Let $(\Xcal,\Xcal^\Omega)$ be a countably separated quasi-universal space.
    Then the quasi-universal space of probability measures:
    \[ \PrPf(\Xcal,\Xcal^\Omega)  \]
    is countably separated as well.

\end{Cor}

\subsection{Universal Measurable and Standard Quasi-Universal Spaces}

In this section we will gather properties of \emph{universal measurable spaces} and \emph{standard quasi-universal spaces}. The former were  studied in \cite{For21}.
They will be the most well-behaved quasi-universal spaces and play a similar role in the category of quasi-universal spaces $\QUS$ as standard Borel measurable spaces play inside the category of measurable spaces $\Meas$ or quasi-Borel spaces $\QBS$. They also allow us to translate (some) results from $\Meas$ to $\QUS$.
We will see that this mainly entails checking which properties of standard Borel measurable spaces are preserved under universal completion.

\subsubsection{Universal Measurable Spaces}

Here we quickly review \emph{universal measurable spaces}, see \cite{For21} Appendix B, which
are measurable spaces that are isomorphic to a universally measurable subset of $\R$, thus generalizing
standard Borel measurable spaces and analytic measurable spaces.
One design choice one has to make is if one explicitly wants their $\sigma$-algebra, per definition, to be universally complete or if one implicitly builds this into the definition.
Since we do not want to study more (co)reflective subcategories of $\Meas$, here given by the universal completion, and we want to be able to say: ``Standard and analytic measurable spaces are universal.'' and ``Countable products of universal measurable spaces are universal.'' (without changing the product-$\sigma$-algebra), etc., we allow for non-universally complete universal measurable spaces, which then get their convenient properties only after universal completion.

The class of spaces in the next Definition/Lemma is classical, and predates \cite{For21}: it is
Shortt's class of \emph{universally measurable} measurable spaces, see \cite{Sho84}, also called
\emph{metrically standard Borel spaces} after Mackey, see \cite{Kec95} 21.D for the underlying
descriptive set theory. Shortt's Lemma 4 characterises a countably generated separated measurable
space as universally measurable exactly when every finite measure on it gives full measure to a
subset that is itself standard Borel, which is the content of our conditions (2) and (3); the
perfectness that appears there is the one studied in \cite{Dar71,Saz62,Fremlin}. What the
following adds to \cite{Sho84} is the equivalence with the retract condition (4), which is the
form the rest of this paper uses, and the bookkeeping of the universal completions.

\begin{DefLem}[Universal measurable spaces, see \cite{Sho84} and \cite{For21} Appendix B, see \Cref{deflem:univ-meas:prf}]
    \label{deflem:univ-meas}
 Let $\Bcal_\R$ be the Borel $\sigma$-algebra on $\R$ and $(\Omega,\Bcal_\Omega)=(\R,(\Bcal_\R)_\Giry)$ its universal completion.
A measurable space $(\Xcal,\Bcal_\Xcal)$ is called \emph{universal measurable space} if it satisfies any of the following equivalent statements:
\begin{enumerate}
    \item There exist $\Zcal \in \Bcal_\Omega$ and an isomorphism of measurable spaces:
            \[ (\Xcal,(\Bcal_\Xcal)_{\Giry}) \cong (\Zcal,\Bcal_{\Omega|\Zcal}). \]
    \item $(\Xcal,\Bcal_\Xcal)$ is universally countably generated, (universally countably) separated and perfect
        (i.e.\ every probability measure $P \in \Giry(\Xcal,\Bcal_\Xcal)$ is perfect).
    \item There exists a countable subset $\Ecal \ins (\Bcal_\Xcal)_\Giry$ that separates the points of $\Xcal$ such that
        $(\Ecal)_\Giry = (\Bcal_\Xcal)_\Giry$ and every probability measure $P \in \Giry(\Xcal,\Ecal)$ is perfect.
    \item $\Xcal = \emptyset$, or $(\Xcal,(\Bcal_\Xcal)_\Giry)$ is a retract of $(\Omega,\Bcal_\Omega)$ in $\Meas$, i.e.\ there are measurable maps
        $i:\; (\Xcal,(\Bcal_\Xcal)_\Giry) \to (\Omega,\Bcal_\Omega)$
        and $r:\;  (\Omega,\Bcal_\Omega) \to (\Xcal,(\Bcal_\Xcal)_\Giry)$ such that $r \circ i = \id_\Xcal$.\footnote{The empty space has to be listed separately: it satisfies the first three conditions, but there is no map $\Omega \to \emptyset$, so it is not a retract of $\Omega$. Keeping $\emptyset$ a universal measurable space is what makes $\UMeas$ and $\StdQUS$ closed under countable coproducts, including the empty one, see \Cref{thm:Markov-category}.\label{fn:empty-retract}}
\end{enumerate}

\end{DefLem}

\begin{Rem}[See \cite{For21} Appendix B]
    \label{rem:umeas}
    \begin{enumerate}
        \item The countable product of universal measurable spaces is a universal measurable space (with the product $\sigma$-algebra, i.e.\ in $\Meas$).
        \item The countable coproduct (in $\Meas$) of universal measurable spaces is a universal measurable space.
        \item Every universally measurable subset of a universal measurable space
            is a universal measurable space.
        \item $(\Omega,\Bcal_\Omega)=(\R,(\Bcal_\R)_\Giry)$ is a universal measurable space.
        \item Every standard (Borel) measurable space and every analytic measurable space is a universal measurable space.
        \item Every countable discrete measurable space is a universal measurable space.
        \item Every Radon Hausdorff space that has a countable network of universally measurable subsets is a universal measurable space together with (the universal completion of) its Borel $\sigma$-algebra.
        \item Every Polish space is a universal measurable space
            together with (the universal completion of) its Borel $\sigma$-algebra.
    \end{enumerate}
\end{Rem}

\begin{Not}[The category of universally complete universal measurable spaces]
    We abbreviate the full subcategory of all universally complete universal measurable spaces inside $\Meas$ as $\UMeas$.
\end{Not}

\subsubsection{Standard Quasi-Universal Spaces}

Here we now give a definition for \emph{standard quasi-universal spaces}, which are now objects in the category of quasi-universal spaces $\QUS$. We then show how they correspond to universal measurable spaces.
This then allows us to translate properties of standard Borel measurable spaces to universal measurable spaces to standard quasi-universal spaces.

\begin{DefLem}[See \Cref{def:univ-qus:prf}]
    \label{def:univ-qus}
    A quasi-universal space $(\Xcal,\Xcal^\Omega)$ will be called \emph{standard quasi-universal space}
    if it satisfies any of the following equivalent conditions:
    \begin{enumerate}
        \item  $(\Xcal,\Bcal(\Xcal^\Omega))$ is a (universally complete) universal measurable space,
            see Definition/Lemma \ref{deflem:univ-meas}, and $\Fcal\Bcal(\Xcal^\Omega)=\Xcal^\Omega$.
            \item There exists a $\sigma$-algebra $\Bcal_\Xcal$ on $\Xcal$ such that $(\Xcal,\Bcal_\Xcal)$ is a universal measurable space, see Definition/Lemma \ref{deflem:univ-meas}, and $\Xcal^\Omega=\Fcal(\Bcal_\Xcal)$.
            \item $\Xcal=\emptyset$, or $(\Xcal,\Xcal^\Omega)$ is a retract of $(\Omega,\Omega^\Omega)$ in $\QUS$, i.e.\ there exist quasi-measurable maps $i:\; (\Xcal,\Xcal^\Omega) \to  (\Omega,\Omega^\Omega)$ and $r:\;(\Omega,\Omega^\Omega)\to (\Xcal,\Xcal^\Omega)$ such that $r \circ i = \id_\Xcal$, see \footref{fn:empty-retract}.
    \end{enumerate}

\end{DefLem}

\begin{Rem}[Countable coproducts of standard quasi-universal spaces]
    The countable coproduct of universal measurable spaces is a universal measurable space.
    The functor $\Fcal$ only preserves countable coproducts if one reflects back to the category
    of \patchable\ quasi-universal spaces $\PQUS$, see Lemma \ref{lem:patch-coprod}.
    So we get that if $(\Xcal_i,\Xcal_i^\Omega)$ is a countable family of standard quasi-universal spaces then
    \[ (\Xcal,\Xcal^\Omega) := \coprod_{i \in I}^\PQUS(\Xcal_i,\Xcal_i^\Omega) = \Pch \lp \coprod_{i \in I}^\QUS(\Xcal_i,\Xcal_i^\Omega)\rp
    = \lp \coprod_{i \in I} \Xcal_i, \Pch\lp \coprod_{i \in I} \Xcal_i^\Omega \rp \rp  \]
    is also a standard quasi-universal space.
    In other words, the countable $\PQUS$-coproduct of standard quasi-universal
    spaces is a standard quasi-universal space.
\end{Rem}

\begin{Cor}[See \Cref{cor:umeas-uqus-eq:prf}]
    \label{cor:umeas-uqus-eq}
    The functors $\Fcal$ and $\Bcal$ from Theorem \ref{thm:qms-meas-adjunction}
    establish an equivalence between the full subcategory of
    universally complete universal measurable spaces $\UMeas$ inside $\Meas$ and the full subcategory of standard quasi-universal spaces $\StdQUS$ inside $\QUS$.

\end{Cor}

\begin{Rem}[No \patchability\ is needed here]
    \label{rem:umeas-uqus-eq-patch}
    The equivalence needs no extra hypothesis, and in particular no \patchability. The reason is
    that it is built into \Cref{def:univ-qus}: a standard quasi-universal space satisfies
    $\Xcal^\Omega = \Fcal\lp\Bcal(\Xcal^\Omega)\rp$, so it is already a fixed point of
    $\Fcal\Bcal$, while a universally complete universal measurable space is a fixed point of
    $\Bcal\Fcal$. The two functors therefore restrict to mutually inverse equivalences on these
    full subcategories, and the morphisms agree by \Cref{qms:thm:adjunction}.
    Note also that every standard quasi-universal space is automatically \patchable: its
    admissible random variables are \emph{all} measurable maps $\Omega \to \Xcal$, and gluing
    countably many measurable maps along a countable measurable partition of $\Omega$ gives a
    measurable map again, see \Cref{def:patchable}.
    Where \patchability\ does matter is not the equivalence but the \emph{preservation of
    countable coproducts} by $\Fcal$, see \Cref{lem:patch-coprod}, which is why the coproducts in
    \Cref{thm:Markov-category} are those of $\PQUS$ and not those of $\QUS$.
\end{Rem}

\subsection{Spaces of Probability Measures for Quasi-Universal Spaces}

Here $\PrPf$ absorbs the other push-forward monads of \Cref{qms:sec:app-other-monads}, and on standard quasi-universal spaces even $\PrQm$ and $\Giry$ agree with it, see \Cref{thm:qus-S=P-in-R} and \Cref{thm:univ-qus-GPQKR}. So from here on $\PrPf$ is the only space of probability measures we need.

The push-forward probability monads of \Cref{qms:sec:app-other-monads} all agree on the category of quasi-universal spaces $\QUS$, see \Cref{thm:qus-S=P-in-R}, so $\PrPf$ alone carries the theory here.
We then go on to show in Theorem \ref{thm:univ-qus-GPQKR} that for \emph{standard} quasi-universal spaces all defined probability spaces even agree with $\Giry$ and $\PrQm$,
as sets and as quasi-measurable spaces, which allows us to translate properties of standard Borel measurable spaces to standard quasi-universal spaces.

\begin{Cor}[See \Cref{cor:random-fct:prf}]
    \label{cor:random-fct}
    Let $(\Xcal,\Xcal^\Omega)$ and $(\Zcal,\Zcal^\Omega)$  quasi-universal spaces.
    Then we have a commutative diagram  of quasi-measurable maps:
    \[\xymatrix{
            \lp\Xcal^\Zcal\rp^\Omega \ar@{=}_-{\swap}[d]\ar@{->>}^-{(\_)_*P}[rr] && \PrPf(\Xcal^\Zcal)\ar@{->}[d]  \\
            \lp\Xcal^\Omega\rp^\Zcal \ar@{->}_-{[(\_)_*P]^\Zcal}[rr] && \PrPf(\Xcal)^\Zcal,
    }\]
    where the horizontal maps are induced by $\psi \mapsto \psi_*P$, where $P$ is the uniform distribution on $\Omega\cong[0,1]$, and where the right vertical map is:
\[\psi_*P \mapsto ( z \mapsto \psi_*(z)P).\]
    Furthermore, if either $(\Xcal,\Xcal^\Omega)$ or $(\Zcal,\Zcal^\Omega)$  is a standard quasi-universal space then all maps are (surjective) quotient maps of quasi-universal spaces.

\end{Cor}

\begin{Rem}
    \label{rem:random-fct}
    To understand the difference between the objects
$\PrPf(\Xcal^\Zcal)$ and $\PrPf(\Xcal)^\Zcal$ from Corollary \ref{cor:random-fct} let $\Zcal=\lC 0,1,2\rC$. Then an element $P(X) \in \PrPf(\Xcal^\Zcal)$ corresponds to encoding the joint distribution:
    \[ P(X_0,X_1,X_2).\]
 This is in contrast to an element $P(X|Z) \in \PrPf(\Xcal)^\Zcal$, which would correspond to encoding the tuple of marginal distributions:
 \[ \lp P(X_0), P(X_1), P(X_2) \rp.  \]
 The quasi-measurable map: $\PrPf(\Xcal^\Zcal) \srj \PrPf(\Xcal)^\Zcal$ then corresponds to mapping the joint distribution to the tuple of marginal distributions:
 \[ P(X_0,X_1,X_2)\mapsto \lp P(X_0), P(X_1), P(X_2) \rp.\]
\end{Rem}

\subsection{Fubini Theorem for Probability Monad $\PrPf$ on Quasi-Universal Spaces}

We show that integration w.r.t.\ product probability measures on products of quasi-universal spaces commute.
Such results are known as \emph{Fubini theorems}.

\begin{Thm}[Fubini Theorem for Markov kernels, see \Cref{thm:fubini:prf}]
    \label{thm:fubini}
    Let $(\Xcal,\Xcal^\Omega)$, $(\Ycal,\Ycal^\Omega)$, $(\Ucal,\Ucal^\Omega)$, $(\Zcal,\Zcal^\Omega)$
    be quasi-universal spaces then the following
    diagram is commutative:
    \[\xymatrix{
            \PrPf(\Xcal)^\Ucal \times \PrPf(\Ycal)^\Zcal \ar_-{\swap_{\PrPf(\Xcal)^\Ucal,\PrPf(\Ycal)^\Zcal}}[d] \ar^-{\otimes}[rrrr] &&&& \PrPf(\Xcal \times \Ycal)^{\Ucal \times \Zcal} \ar^-{\swap_{(\Xcal,\Ycal),*}}[d]\\
            \PrPf(\Ycal)^\Zcal \times \PrPf(\Xcal)^\Ucal \ar_-{\otimes}[rr] &&
            \PrPf(\Ycal \times \Xcal)^{\Zcal \times \Ucal} \ar_-{\swap_{\Ucal,\Zcal}^*}[rr]
                && \PrPf(\Ycal \times \Xcal)^{\Ucal \times \Zcal}.
                                                        }\]
    More concretely, for $K(X|U) \in  \PrPf(\Xcal)^\Ucal$, $Q(Y|Z) \in \PrPf(\Ycal)^\Zcal$, $u \in \Ucal$, $z \in \Zcal$ and $f \in [0,\infty]^{\Xcal \times \Ycal}$ we have:
    \begin{align*}
        &\int \int f(x,y)\,K(X \in dx|U=u)\,Q(Y \in dy|Z=z)\\
        &\qquad\qquad = \int \int f(x,y)\footnotemark\,Q(Y \in dy|Z=z)\,K(X \in dx|U=u).
    \end{align*}

\end{Thm}

\subsection{Disintegration of Markov Kernels on Quasi-Universal Spaces}

In this subsection we will show under which conditions a Markov kernel on a product space of quasi-universal spaces can be \emph{disintegrated} into a marginal part and a conditional part such that the chain rule holds.
For this we use the a bit more suggestive notation $K(X,Y|Z)$ for Markov kernels to make it easier
to indicate the marginals $K(Y|Z)$, which usually would need to be written as a push-forward,
and the chain rule for conditioning like $K(X,Y|Z)=K(X|Y,Z) \otimes K(Y|Z)$,
which quickly became unreadable otherwise. We then discuss the (essential) uniqueness of such factorizations.
We first recall the strongest disintegration theorem for measurable spaces from \cite{For21}.

Before the theorem we record the elementary fact that makes the last of its cases work without any
hypothesis on $\Xcal$: for a discrete $\Ycal$ the induced $\sigma$-algebra of the product is
already the product $\sigma$-algebra, so the gap described in \Cref{rem:disint-why-Y} does not open
at all.

\begin{Lem}[Products with a discrete factor, see \Cref{lem:prod-discrete:prf}]
    \label{lem:prod-discrete}
    Let $(\Xcal,\Xcal^\Omega)$ be a quasi-measurable space and let $(\Ycal,\Ycal^\Omega)$ be
    \emph{discrete}, i.e.\ $\Ycal$ is countable and $\Ycal^\Omega = \Fcal(\two^\Ycal)$. Then
    \begin{align}
        \Bcal_{\Xcal \times \Ycal} &= \Bcal_\Xcal \otimes \Bcal_\Ycal
        = \lC \bigdcup_{y \in \Ycal} A_y \times \lC y \rC \st A_y \in \Bcal_\Xcal \rC.
    \end{align}
    In particular $\Delta_\Ycal \in \Bcal_{\Ycal \times \Ycal}$, i.e.\ every discrete
    quasi-measurable space is quasi-separated.
\end{Lem}
\begin{Def}[Disintegration triple]
    \label{def:disint-triple-qus}
    A triple $(\Xcal,\Xcal^\Omega)$, $(\Ycal,\Ycal^\Omega)$, $(\Zcal,\Zcal^\Omega)$ of
    quasi-universal spaces is called a \emph{disintegration triple} if \emph{every} Markov kernel
    \[ K(X,Y|Z) \in \PrPf(\Xcal \times \Ycal)^\Zcal \]
    admits a \emph{conditional Markov kernel}
    \[ K(X|Y,Z) \in \PrPf(\Xcal)^{\Ycal \times \Zcal}, \]
    i.e.\ one with
    \[ K(X,Y|Z) = K(X|Y,Z) \otimes K(Y|Z), \]
    where $K(Y|Z)$ denotes the marginal Markov kernel. Equivalently, the product of Markov kernels
    \[ \otimes:\; \PrPf(\Xcal)^{\Ycal\times\Zcal} \times \PrPf(\Ycal)^\Zcal \srj
       \PrPf(\Xcal \times \Ycal)^\Zcal \]
    is surjective. If, moreover, $(\Xcal,\Ycal,\Zcal \times \Omega)$ is a disintegration triple,
    then that map is even a \emph{quotient map}, i.e.\ surjective on the level of admissible
    random variables.
    As in \Cref{def:disintegration-triple}, this is a property of the triple of spaces alone.
\end{Def}

\begin{Thm}[Disintegration of Markov kernels between quasi-universal spaces, see \Cref{thm:disint-X:prf}]
    \label{thm:disint-X}
    Let $(\Xcal,\Xcal^\Omega)$, $(\Ycal,\Ycal^\Omega)$ and $(\Zcal,\Zcal^\Omega)$  be quasi-universal spaces.
    Assume one of the following points:
    \begin{enumerate}
        \item $(\Ycal,\Ycal^\Omega)$ is \emph{countably separated} and $(\Xcal,\Xcal^\Omega)$ is a
            \emph{standard} quasi-universal space, or:
        \item $(\Ycal,\Ycal^\Omega)$ is \emph{countably separated} and $(\Zcal,\Zcal^\Omega)$ is a
            \emph{standard} quasi-universal space
            (e.g.\ $\Zcal \cong \one$, $\N$, $\Omega$, etc.), or:
        \item $(\Ycal,\Ycal^\Omega)$ is \emph{discrete}, i.e.\ $\Ycal$ is countable and
            $\Ycal^\Omega=\Fcal(\two^\Ycal)$, the set of all measurable maps
            $(\Omega,\Bcal_\Omega) \to (\Ycal,\two^\Ycal)$, with $(\Xcal,\Xcal^\Omega)$ and
            $(\Zcal,\Zcal^\Omega)$ arbitrary.
\end{enumerate}
    Then $(\Xcal,\Ycal,\Zcal)$ is a \emph{disintegration triple} in the sense of
    \Cref{def:disint-triple-qus}, i.e.\ the product of Markov kernels
    \begin{align*}
     \otimes:\; \PrPf(\Xcal)^{\Ycal\times\Zcal} \times \PrPf(\Ycal)^\Zcal &\srj \PrPf(\Xcal \times \Ycal)^\Zcal,
    \end{align*}
    is surjective, and it is even a quotient map.\footnote{All three cases survive replacing
    $\Zcal$ by $\Zcal \times \Omega$ --- in cases 1 and 3 because $\Zcal$ carries no hypothesis
    there, in case 2 because $\Zcal \times \Omega$ is again standard by
    \Cref{lem:uqus-count-prod} --- so $(\Xcal,\Ycal,\Zcal\times\Omega)$ is again a disintegration
    triple and \Cref{def:disint-triple-qus} upgrades surjectivity to the quotient
    property.\label{fn:disint-quotient}}

\end{Thm}

\begin{Rem}[Why $\Ycal$ carries a hypothesis, and why case 3 is not case 2]
    \label{rem:disint-why-Y}
    A measure theorist will object to the hypothesis on $\Ycal$ in cases 1 and 2, because in
    $\Meas$ none is needed: for $\Xcal$ standard Borel, regular conditional distributions given an
    arbitrary sub-$\sigma$-algebra always exist. The reason it is needed here is specific to
    $\QMS$. The induced $\sigma$-algebra $\Bcal_{\Xcal \times \Ycal}$ of the product can be
    strictly larger than the product $\sigma$-algebra $\Bcal_\Xcal \otimes \Bcal_\Ycal$, so the
    input ``$K(X,Y|Z) \in \PrPf(\Xcal \times \Ycal)^\Zcal$'' is strictly stronger than a Markov
    kernel into $\Giry(\Xcal\times\Ycal,\Bcal_\Xcal \otimes \Bcal_\Ycal)$, and the conclusion is
    correspondingly stronger. Countable separatedness is what controls that gap, see
    \Cref{qms:thm:count-sep-perf} and \Cref{thm:countably-separated}.

    Case 3 is the only one in which \emph{no} space is required to be standard, and it is also the
    one in which the gap just described does not open at all: for a discrete $\Ycal$ the induced
    $\sigma$-algebra of the product is already the product $\sigma$-algebra,
    $\Bcal_{\Xcal \times \Ycal}=\Bcal_\Xcal \otimes \Bcal_\Ycal$, see \Cref{lem:prod-discrete}, so
    the input is no stronger than a measurable Markov kernel and no extension step is needed. It is
    not comparable with cases 1 and 2: it asks more of $\Ycal$ and nothing of $\Xcal$ and $\Zcal$,
    where the others ask less of $\Ycal$ and more of $\Xcal$ or $\Zcal$. Its conditional kernel is
    the elementary quotient \Cref{eq:disint-qus-quotient}.
\end{Rem}

\begin{Rem}[A candidate for a further disintegration triple]
    \label{rem:disint-candidate}
    One further triple suggests itself, and we record it as an open question rather than as a
    result:
    \begin{align}
        (\Xcal,\Xcal^\Omega) \text{ standard}, \qquad
        (\Ycal,\Ycal^\Omega) \text{ arbitrary}, \qquad
        (\Zcal,\Zcal^\Omega) \text{ discrete}.
        \label{eq:disint-candidate}
    \end{align}
    Everything on the measurable-space side is available: $(\Xcal,\Bcal_\Xcal)$ is perfect and
    universally countably generated, $\Bcal_\Zcal$ is a power set, and
    \Cref{cor:markov-kernel-disint-meas-II} in its second case then produces a conditional Markov
    kernel with \emph{no} hypothesis on $(\Ycal,\Bcal_\Ycal)$, its measurability in $z$ being
    automatic. This yields the factorization
    \[ K(X,Y|Z=z) = K(X|Y,Z=z) \otimes K(Y|Z=z) \]
    as measures on $\Bcal_\Xcal \otimes \Bcal_\Ycal$, for every $z \in \Zcal$.

    What is missing is the passage from $\Bcal_\Xcal \otimes \Bcal_\Ycal$ to the induced
    $\sigma$-algebra $\Bcal_{\Xcal \times \Ycal}$, on which the members of
    $\PrPf(\Xcal \times \Ycal)$ live and on which \Cref{def:disint-triple-qus} asks for the
    identity, i.e.\ the question of whether
    \begin{align}
        \Bcal_{\Xcal \times \Ycal} &= \lp \Bcal_\Xcal \otimes \Bcal_\Ycal \rp_{\PrPf(\Xcal \times \Ycal)}
        \label{eq:disint-candidate-sigma}
    \end{align}
    holds for $(\Xcal,\Xcal^\Omega)$ standard and $(\Ycal,\Ycal^\Omega)$ arbitrary. In cases 1 and
    2 this is what countable separatedness of $\Ycal$ delivers, through
    \Cref{thm:countably-separated}; in case 3 it holds outright, by \Cref{lem:prod-discrete}. We
    can neither prove nor refute \Cref{eq:disint-candidate-sigma} in the stated generality. What we
    do know:
    \begin{enumerate}
        \item The inclusion $\sni$ always holds, since $\Bcal_{\Xcal \times \Ycal}$ is
            $\PrPf(\Xcal \times \Ycal)$-complete by \Cref{lem:qus-universally-complete} and
            contains $\Bcal_\Xcal \otimes \Bcal_\Ycal$.
        \item Every $\Xcal$-section of a $D \in \Bcal_{\Xcal \times \Ycal}$ is an event: taking the
            constant random variable at $y$ gives $X^{-1}(D^y) \in \Bcal_\Omega$ for all
            $X \in \Xcal^\Omega$, hence $D^y \in \Bcal_\Xcal$.
        \item The assignment $y \mapsto D^y$ is quasi-measurable as a map
            $\Ycal \to \Bcal_\Xcal$, by currying $\I_D$ through the function space
            $\twob^\Xcal \cong \Bcal_\Xcal$, see \Cref{qms:lem:inclusion-check}; consequently
            $y \mapsto P(D^y)$ is $\Bcal_\Ycal$-measurable for every $P \in \PrQm(\Xcal)$,
            the evaluation map being quasi-measurable.
        \item \Cref{eq:disint-candidate-sigma} holds if $\Ycal$ is countably separated, if $\Ycal$
            is discrete, and also if $\Ycal^\Omega=\Sets(\Omega,\Ycal)$, in which case both sides
            collapse to $\Bcal_\Xcal \otimes \lC \emptyset,\Ycal\rC$.
    \end{enumerate}
    Since $\Xcal$ is standard, the measure algebra of $\Bcal_\Xcal$ modulo the null sets of any
    fixed $P \in \PrPf(\Xcal)$ is separable, and by 2.\ and 3.\ the sections of $D$ define a
    measurable map into it; what we could not do is assemble the resulting fibrewise
    approximations into a single set of $\Bcal_\Xcal \otimes \Bcal_\Ycal$ without already having
    the disintegration one is trying to construct.
\end{Rem}

\begin{Rem}
    \begin{enumerate}
        \item Theorem \ref{thm:disint-X} would also hold if one replaces standard quasi-measurable spaces with perfect universally countably generated ones as in Corollary \ref{cor:markov-kernel-disint-meas-II}.
    For this to work one uses the separated quotient and a section, see Definition \ref{def:separated-quotient}
 Remark \ref{rem:separated-quotient}.
 Note that the separated quotient of a perfect universally countably generated quasi-universal space is a standard quasi-universal space.
 \item It might even be possible to weaken the assumption in Theorem \ref{thm:disint-X} from \emph{standard} to \emph{countably separated} quasi-universal spaces and/or other weaker assumptions. The step that resists is the passage through Theorem \ref{thm:univ-qus-GPQKR}, which identifies $\PrPf$ with $\Giry$ and $\PrQm$ only on standard quasi-universal spaces; without that identification the transport of Corollary \ref{cor:markov-kernel-disint-meas-II} across the adjunction is not available. We have not further investigated this. It is the single hypothesis that propagates furthest --- into \Cref{cor:de-finetti}, into the separoid rules of \Cref{thm:separoid_axioms-tci}, and hence into \Cref{thm:gmp-mI-CBN}.
 \end{enumerate}
\end{Rem}

\begin{Thm}[Essential uniqueness of factorizations, see \Cref{thm:fact-ess-uniq:prf}]
    \label{thm:fact-ess-uniq}
    Let $(\Xcal,\Xcal^\Omega)$, $(\Ycal,\Ycal^\Omega)$ and $(\Zcal,\Zcal^\Omega)$  be quasi-universal spaces.  Let $Q(Y|Z) \in \PrPf(\Ycal)^\Zcal$ and let $P_1(X|Y,Z)$, $P_2(X|Y,Z)$ be in $\PrPf(\Xcal)^{\Ycal \times \Zcal}$ such that:
    \[ P_1(X|Y,Z) \otimes Q(Y|Z) = P_2(X|Y,Z) \otimes Q(Y|Z) \quad \in \PrPf(\Xcal \times \Ycal)^\Zcal.\]
   For every $D \in \Bcal_{\Xcal \times \Ycal \times \Zcal}$ we have that:
   \[ N_D:=\lC (y,z) \in \Ycal \times \Zcal\st P_1(X \in D^{y,z}|Y=y,Z=z) \neq P_2(X \in D^{y,z}|Y=y,Z=z)\rC   \]
   is a $Q(Y|Z)$-null set in $\Bcal_{\Ycal \times \Zcal}$, i.e.\ $Q(Y \in N_D^z|Z=z)=0$ for all $z \in \Zcal$.\\
   If, furthermore, $(\Xcal,\Xcal^\Omega)$
   is countably separated then $P_1(X|Y,Z)=P_2(X|Y,Z)$ $Q(Y|Z)$-almost-surely, i.e.,
    more precisely, that
   there exists a set $N \in \Bcal_{\Ycal \times \Zcal}$ such that for all $z \in \Zcal$:
   \[ Q(Y \in N^z|Z=z) =0,   \]
   where $N^z:=\lC y \in \Ycal\st (y,z) \in N \rC$,
   and such that for all $(y,z) \notin N$ we have:
   \[ P_1(X|Y=y,Z=z) = P_2(X|Y=y,Z=z). \]

\end{Thm}

\subsection{Kolmogorov Extension Theorems for Quasi-Universal Spaces}

In this subsection we study how sequences of probability distributions and Markov kernels defined on finite products of quasi-universal spaces can be extended to an infinite product if they are at least consistent with each other. Such results are usually called \emph{Kolmogorov extension theorems}.

\begin{Thm}[Kolmogorov's extension theorem for Markov kernels, see \cite{Kol33,Fremlin} 454D-G, also see \cite{Fri19} Ex.\ 3.6, see \Cref{thm:kol-ext-mar-ker:prf}]
    \label{thm:kol-ext-mar-ker}
        In the following let  $(\Xcal_j,\Bcal_j)$, $j \in J$, be a family of
     measurable spaces indexed over any set $J$.
    For a subset $K \ins J$ we put:
    \[ \Xcal_K:=\prod_{k \in K} \Xcal_k, \qquad
    \Bcal_K:=\bigotimes_{k \in K} \Bcal_k := \sigma\lp\Ecal_K\rp,\]
    where $\Ecal_K$ is given by all finite cylinder sets:
    \[ \Ecal_K:= \lC \pr_L^{-1}(\times_{l \in L} A_l) \ins \Xcal_K\st A_l \in \Bcal_l, L \ins K, \#L < \infty\rC,\]
and with the canonical projection maps:
\[\pr_L:\; \Xcal_K \to \Xcal_L,\qquad (x_k)_{k \in K} \mapsto (x_l)_{l\in L}.\]
Let $(\Zcal,\Bcal_\Zcal)$ be another measurable space.
    For every finite subset $K \ins J$ consider a (measurable) Markov kernel:
\[ Q_K(X_K|Z):\, (\Zcal,\Bcal_\Zcal) \to \Giry(\Xcal_K,\Bcal_K), \]
such that for every finite subset $L \ins K$ we have the consistency relation:
\[ \pr_{L,*} Q_K(X_K|Z) = Q_L(X_L|Z),\]
and such that for every $z\in\Zcal$ and $j \in J$ the marginal probability measure $Q_j(X_j|Z=z)$ is perfect.
Then there exists a unique Markov kernel:
\[ Q_J(X_J|Z):\;(\Zcal,\Bcal_\Zcal) \to \Giry(\Xcal_J,\Bcal_J)\]
such that for all finite subsets $L \ins J$ we have:
\[ \pr_{L,*} Q_J(X_J|Z) = Q_L(X_L|Z).\]
Furthermore, for every $z \in \Zcal$ the probability measure $Q_J(X_J|Z=z)$ is perfect and measurable as a map in $z$.

\end{Thm}

\begin{Rem}
    \label{rem:kol-ext}
    The Markov kernel $Q_J(X_J|Z)$ from Theorem \ref{thm:kol-ext-mar-ker}
    can uniquely be extended from $\Bcal_J$ to its universal completion $(\Bcal_J)_\Giry$.
    By Proposition \ref{prp:prob-prob-meas} we then get that $Q_J(X_J|Z)$ still is (universally)
    $(\Bcal_\Zcal)_\Giry$-$\Bcal_{\Giry(\Xcal_J,(\Bcal_J)_\Giry)}$-measurable as a map:
    \[ Q_J(X_J|Z):\;(\Zcal,(\Bcal_\Zcal)_\Giry) \to \Giry(\Xcal_J,(\Bcal_J)_\Giry).\]
    Clearly, it also stays perfect and is consistent with its finite marginals.
\end{Rem}

\begin{Thm}[Kolmogorov's extension theorem for Markov kernels for standard quasi-universal spaces, see \Cref{thm:kol-ext-univ:prf}]
    \label{thm:kol-ext-univ}
    Let $(\Xcal_n,\Xcal_n^\Omega)$, $n \in \N$, be a sequence of
    standard quasi-universal spaces and $(\Zcal,\Zcal^\Omega)$ any quasi-universal space.
    Assume we have $Q_n(X_{0:n}|Z) \in \PrPf(\Xcal_{0:n})^\Zcal$ for every $n \in \N$ such that
    for every $n \in \N$:
    \[ \pr_{0:n,*}Q_{n+1}(X_{0:n+1}|Z) = Q_{n}(X_{0:n}|Z).\]
    Then there exists a unique $Q(X_\N|Z) \in \PrPf(\Xcal_\N)^\Zcal$ such that for every $ n \in \N$:
        \[ \pr_{0:n,*}Q(X_\N|Z) = Q_{n}(X_{0:n}|Z),\]
        where $(\Xcal_\N,\Xcal_\N^\Omega) := \prod_{n\in\N}(\Xcal_n,\Xcal_n^\Omega)$.

\end{Thm}

\begin{Thm}[Kolmogorov's extension theorem for Markov kernels for countably separated quasi-universal spaces, see \Cref{thm:kol-ext-cs:prf}]%
\footnote{The two Kolmogorov extension theorems are \emph{incomparable}: the first assumes the
$\Xcal_n$ standard and allows $\Zcal$ arbitrary, the second assumes the $\Xcal_n$ only countably
separated but requires $\Zcal$ standard. Neither subsumes the other. Infinite products in the
quasi-Borel setting have also been studied, without the universal completion, in \cite{KPS25}.
A word on \emph{uniqueness} in both: uniqueness on the cylinder $\sigma$-algebra is Dynkin's
lemma, but $\Bcal_{\Xcal_\N}$ is strictly larger than the cylinder $\sigma$-algebra in general.
It still follows, because for standard $\Xcal_n$ the product $\Xcal_\N$ is standard with
$\Bcal_{\Xcal_\N}=\lp\bigotimes_n \Bcal_{\Xcal_n}\rp_\Giry$, and for countably separated $\Xcal_n$
\Cref{thm:countably-separated} gives $\Bcal_{\Xcal_\N}=(\Ecal)_{\PrPf(\Xcal_\N)}$ for a countable
separating $\Ecal$; in both cases a member of $\PrPf(\Xcal_\N)$ is already determined by its
restriction to the cylinders.\label{fn:kol-uniqueness}}
    \label{thm:kol-ext-cs}
        Let $(\Xcal_n,\Xcal_n^\Omega)$, $n \in \N$, be a sequence of
    countably separated quasi-universal spaces and $(\Zcal,\Zcal^\Omega)$ a standard quasi-universal space.
    Assume we have $Q_n(X_{0:n}|Z) \in \PrPf(\Xcal_{0:n})^\Zcal$ for every $n \in \N$ such that
    for every $n \in \N$:
    \[ \pr_{0:n,*}Q_{n+1}(X_{0:n+1}|Z) = Q_{n}(X_{0:n}|Z).\]
    Then there exists a unique $Q(X_\N|Z) \in \PrPf(\Xcal_\N)^\Zcal$ such that for every $ n \in \N$:
        \[ \pr_{0:n,*}Q(X_\N|Z) = Q_{n}(X_{0:n}|Z),\]
        where $(\Xcal_\N,\Xcal_\N^\Omega) := \prod_{n\in\N}(\Xcal_n,\Xcal_n^\Omega)$.

\end{Thm}

\subsection{Deterministic Markov Kernels between Quasi-Universal Spaces}

In this subsection we characterize Markov kernels between quasi-universal spaces that are in some sense
\emph{deterministic}. We go through a-priori different notions of deterministic Markov kernels and give a sufficient criterion when they become equivalent.

\begin{DefLem}[Deterministic Markov kernels, see \cite{Rao81,Fri20} Def.\ 10.1, see \Cref{def:qms-det-markov-ker:prf}]
    \label{def:qms-det-markov-ker}
        Let $(\Ycal,\Ycal^\Omega)$ and $(\Zcal,\Zcal^\Omega)$ be quasi-universal spaces and $K(Y|Z) \in \PrPf(\Ycal)^\Zcal$ a Markov kernel.
    Then each statement for $K(Y|Z)$ will imply the one below:
    \begin{enumerate}
               \item $K(Y|Z)$ is \emph{function-deterministic}, i.e. there exists a quasi-measurable map  $g \in \Ycal^\Zcal$ such that:
                   \[ K(Y|Z) = \delta_g(Y|Z) \quad \in \PrPf(\Ycal)^\Zcal.  \]
       \item $K(Y|Z)$ is \emph{copy-deterministic}, i.e.\ we have the relation:
           \[  K(Y_1|Z=z) \otimes K(Y_2|Z=z) =  K(Y_1,Y_2|Z=z)
           \quad \in \PrPf(\Ycal \times \Ycal)^\Zcal,   \]
    where $Y_1$, $Y_2$ are copies of $Y$ with an index for readability,
                i.e.\ the following diagram commutes:
            \[\xymatrix{
                    \Zcal \ar_-{K(Y|Z)}[d] \ar^-{K(Y|Z)}[rr] && \PrPf(\Ycal) \ar^-{\PrPf(\Delta_{\Ycal})}[drr] && \\
                    \PrPf(\Ycal) \ar_-{\Delta_{\PrPf(\Ycal)}}[rr]&& \PrPf(\Ycal) \times \PrPf(\Ycal)
                    \ar_-{\otimes}[rr]&& \PrPf(\Ycal \times \Ycal).
            }\]
        \item $K(Y|Z)$ is \emph{0-1-deterministic}, i.e.\
            for all $z \in \Zcal$ and $B \in \Bcal_\Ycal$:
                \[ K(Y\in B|Z=z) \in \lC0,1 \rC. \]
    \end{enumerate}

\end{DefLem}

\begin{Thm}[See \Cref{thm:qms-0-1-det-fct-det:prf}]
    \label{thm:qms-0-1-det-fct-det}
    Let $(\Ycal,\Ycal^\Omega)$ and $(\Zcal,\Zcal^\Omega)$ be quasi-universal spaces and $K(Y|Z) \in \PrPf(\Ycal)^\Zcal$.
Assume that $(\Ycal,\Ycal^\Omega)$ is countably separated with $\Fcal\Bcal(\Ycal^\Omega)=\Ycal^\Omega$(e.g.\ $(\Ycal,\Ycal^\Omega)$ a standard quasi-universal space).\footnote{Note that all quasi-universal spaces that have embeddings $(\Ycal,\Ycal^\Omega) \inj (\Omega,\Omega^\Omega)$ in $\QUS$ satisfy these conditions, see Lemma \ref{lem:sturdily-countably-separated-spaces}. Also all (universally) countably separated measurable spaces satisfy these conditions after applying the functor $\Fcal$.\label{fn:st-cs}}
    Then the following statements are equivalent:
    \begin{enumerate}
        \item $K(Y|Z)$ is \emph{copy-deterministic}, i.e.:
    \[ K(Y_1|Z) \otimes K(Y_2|Z) = K(Y_1,Y_2|Z)  \quad \in \PrPf(\Ycal \times \Ycal)^\Zcal,  \]
    where $Y_1$ and $Y_2$ are copies of $Y$.
        \item $K(Y|Z)$ is \emph{0-1-deterministic}, i.e.\ for all $B \in \Bcal_\Ycal$ and $z \in \Zcal$:
            \[K(Y \in B|Z=z) \in \lC0,1\rC.\]
        \item $K(Y|Z)$ is \emph{function-deterministic}, i.e.\
            there exists a (unique) quasi-measurable map $g \in \Ycal^\Zcal$ such that:
            \[ K(Y|Z) = \delta_g(Y|Z) \quad \in \PrPf(\Ycal)^\Zcal.\]
    \end{enumerate}

\end{Thm}

\begin{Rem}
    \label{rem:a.s.-det}
    One can prove an ``almost-sure'' version of Theorem \ref{thm:qms-0-1-det-fct-det} w.r.t.\ some given Markov kernel $Q(Z|T) \in \PrPf(\Zcal)^\Tcal$ under the same conditions. For this one needs to gather the countably many null-sets $N_B$ for each $B \in \Ecal_\Ycal$, which is then again a null-set w.r.t.\ $Q(Z|T)$.
\end{Rem}

\subsection{Conditional De Finetti Theorem and the Markov Category of Markov Kernels on $\StdQUS$}

Here we quickly summarize our findings in terms of categorical probability theory without
going into the details of these definitions. They can be found in \cite{Fri19,Fri20,FriG20,Fri21}.
As a Corollary we get a \emph{conditional de Finetti theorem} out.

\begin{Thm}[The Markov category of Markov kernels between standard quasi-universal spaces, see \Cref{thm:Markov-category:prf}]
    \label{thm:Markov-category}
    The category of standard quasi-universal spaces $(\StdQUS,\times,\one)$ is a
    symmetric cartesian (but not closed) monoidal category  with countable products and countable
    coproducts, where the coproducts are the ones of $\PQUS$, i.e.\
    $\coprod^\PQUS_{i \in I} = \Pch\lp\coprod^\QUS_{i \in I}\rp$, and \emph{not} the coproducts of
    $\QUS$.\footnote{$\StdQUS$ is not closed under $\QUS$-coproducts: the $\QUS$-coproduct
    $\one \sqcup \one$ is the two-point space $\twoc$ carrying only its two constant random
    variables, see \Cref{qms:eg:coprod-two-point}, whose induced $\sigma$-algebra is the power
    set; so $\Fcal\Bcal(\twoc^\Omega)=\Meas(\Omega,\twob)$ is strictly larger than $\twoc^\Omega$
    and $\twoc$ fails the first condition of \Cref{def:univ-qus}. One has to reflect back along
    $\Pch$, as in the Remark after \Cref{def:univ-qus}. The empty coproduct is the empty space,
    which is standard by \Cref{def:univ-qus}.}
    The triple $(\PrPf,\delta,\Mbb)$ is a strong affine symmetric monoidal/commutative monad on $(\StdQUS,\times,\one)$.
    Its Kleisli category $\Kleisli(\PrPf)$ on $\StdQUS$ is an a.s.-compatibly representable Markov category with conditionals and Kolmogorov powers, in the sense of \cite{Fri19,Fri20,FriG20,Fri21} recalled in \Cref{qms:sec:app-markov-cat}.

\end{Thm}

De Finetti's theorem does not hold over arbitrary measurable spaces --- exchangeable processes
need not be mixtures of i.i.d.\ ones, see \cite{DF79} --- which is why a hypothesis on $\Xcal$
appears below. We note that \cite{Heu17} Section VIII already proves a de Finetti theorem for
\emph{all} quasi-Borel spaces, with no standardness hypothesis, and explains there why
\cite{DF79} does not obstruct it: a probability measure on a quasi-Borel space may draw its
randomness only from $\R$. The two statements are not directly comparable. The version below is
\emph{conditional}, i.e.\ carries a parameter space $\Zcal$, which the version in \cite{Heu17} does
not; it is stated for the monad $\PrPf$ on $\QUS$ rather than for the monad of \cite{Heu17} on
$\QBS$; and the standardness of $\Xcal$ here is inherited from \Cref{thm:disint-X}, through
\Cref{thm:Markov-category}, not from the de Finetti argument itself. Whether it can be dropped is
exactly the open question recorded in \Cref{rem:de-finetti-stronger}.

\begin{Cor}[Conditional De Finetti Theorem, see \cite{deF37,DF79,Heu17,Kal06,Fri21}, see \Cref{cor:de-finetti:prf}]
    \label{cor:de-finetti}
    Let $(\Xcal,\Xcal^\Omega)$ be a standard quasi-universal space
    and $(\Zcal,\Zcal^\Omega)$ another quasi-universal space.
    Write $(\Xcal_\N,\Xcal_\N^\Omega) := \prod_{n \in \N}(\Xcal,\Xcal^\Omega)$ for the countable
    power of $(\Xcal,\Xcal^\Omega)$ in $\QUS$.\footnote{Not the exponential object $\Xcal^\N$: for
    a general quasi-measurable space the two differ, see \Cref{rem:random-fct} and
    \Cref{def:patchable}. For \sturdy\ --- in particular for standard --- quasi-universal spaces
    they agree, both as sets and as sets of admissible random variables, so for the $\Xcal$ of this
    corollary either reading gives the same object.\label{fn:power-vs-exp}}
    For a Markov kernel $Q(X_\N|Z) \in \PrPf(\Xcal_\N)^\Zcal$ the following are equivalent:
    \begin{enumerate}
        \item $Q(X_\N|Z)$ is \emph{exchangeable}, i.e.\ invariant under all finite permutations $\sigma:\; \N \bij \N$.
        \item There exist a quasi-universal space $\Scal$  and
    $Q(X|S) \in \PrPf(\Xcal)^{\Scal}$ and  $P(S|Z) \in \PrPf(\Scal)^\Zcal$
    such that:
    \[ Q(X_\N|Z) = \lp \otimes_{n \in \N} Q(X_n|S) \rp \circ P(S|Z). \]
    \end{enumerate}
    If this is the case one can w.l.o.g.\ choose $\Scal=\PrPf(\Xcal)$ and
     \[Q(X|S)=\id_{\PrPf(\Xcal)}:\; \Scal=\PrPf(\Xcal) \to \PrPf(\Xcal),\qquad Q(X \in A|S=P):=P(A).\]

\end{Cor}

\begin{Rem}[The conditioning space need not be standard]
    \label{rem:de-finetti-stronger}
    \Cref{cor:de-finetti} allows the conditioning space $(\Zcal,\Zcal^\Omega)$ to be an
    \emph{arbitrary} quasi-universal space, while only $(\Xcal,\Xcal^\Omega)$ is required to be
    standard. This is worth emphasising, because the route through
    \Cref{thm:Markov-category} appears at first sight to force both: that theorem describes
    $\Kleisli(\PrPf)$ on $\StdQUS$, and a Kleisli morphism $\Zcal \to \PrPf(\Xcal_\N)$ lives in that
    category only if $\Zcal$ is standard as well.

    The reason the restriction is only apparent is that $\Zcal$ enters the synthetic argument
    exclusively as the \emph{source} of a Kleisli morphism, i.e.\ as a parameter, never as an object
    that has to be conditioned on or mixed over. Each of the three ingredients that the argument
    does use is available in $\QUS$ with $\Zcal$ arbitrary; see the proof,
    \Cref{cor:de-finetti:prf}, for the three-line verification. What genuinely remains open is the
    other direction: weakening $(\Xcal,\Xcal^\Omega)$ from \emph{standard} to merely
    \emph{countably separated}, for which \Cref{thm:disint-X} would have to be strengthened, since
    it is exactly there that standardness of $\Xcal$ is used.
\end{Rem}

\section{Causal Models}

\label{sec:causal-models}
This section is dedicated to the notion of conditional independence and causal models.
More clearly, we will translate the notion of \emph{transitional conditional independence} studied in \cite{For21} from $\Meas$ to $\QUS$. This allows us to mix random variables with deterministic non-random variables in conditional independence statements.
Furthermore, we will engage in causal graphical models and demonstrate our results on \emph{causal Bayesian networks}. Generalizations to other probabilistic graphical models will be apparent.
We also re-formulate the notion of \emph{strong ignorability} using random functions with values in
function spaces of quasi-universal spaces.

\subsection{Transitional Conditional Independence}

In this subsection we translate all structures related to \emph{transitional conditional independence} defined for measurable spaces in \cite{For21} to quasi-universal spaces. We will also investigate the \emph{(asymmetric) separoid rules}.

\begin{Not}[Transition probability space]
    In the following we fix two quasi-universal spaces $(\Tcal,\Tcal^\Omega)$ and $(\Wcal,\Wcal^\Omega)$.
    The generic choice is to take $(\Wcal,\Wcal^\Omega)=(\Omega,\Omega^\Omega)$,
    but we allow more general spaces for now.
    We also fix a Markov kernel $\Kk(W|T) \in \PrPf(\Wcal)^\Tcal$.
    We will call the tuple $(\Wcal \times \Tcal, \Kk(W|T))$ our \emph{transition probability space}.
\end{Not}

\begin{Def}[Transitional random variable]
    A \emph{transitional random variable} is a Markov kernel from our transition probability
    space $(\Wcal \times \Tcal, \Kk(W|T))$
    to any other quasi-universal space $(\Xcal,\Xcal^\Omega)$:
    \[\Xk(X|W,T) \in \PrPf(\Xcal)^{\Wcal \times \Tcal}. \]
\end{Def}

\begin{Eg}[Special transitional random variables of importance]
    We denote the:
    \begin{enumerate}
        \item canonical projection onto $\Tcal$ as:
            \[ T:=\pr_\Tcal:\; \Wcal \times \Tcal \to \Tcal, \qquad T(w,t):=t,\]
            and put:
            \[ \Tk(T|W,T)=\deltabf(T|W,T) = \delta \circ T \in \PrPf(\Tcal)^{\Wcal \times \Tcal}.\]
        \item constant transitional random variable as:
            \[ \deltabf_0:\;\Wcal \times \Tcal \to \PrPf(\one)=\one=\lC0\rC. \]
    \end{enumerate}
\end{Eg}

\begin{Not}
    For two transitional random variables $\Xk \in \PrPf(\Xcal)^{\Wcal\times\Tcal}$ and $\Yk\in \PrPf(\Ycal)^{\Wcal\times\Tcal}$ we define their joint Markov kernel as:
    \[ \Kk(X,Y|T) := \lp  \Xk(X|W,T) \otimes \Yk(Y|W,T) \rp \circ \Kk(W|T) \in \PrPf(\Xcal \times \Ycal)^\Tcal\]
    We then define the relation:
    \[\Xk \ismapof_\Kk \Yk\qquad :\iff\qquad \exists \Phi \in \Xcal^\Ycal.\quad \Kk(X,Y|T) = \deltabf_\Phi(X|Y) \otimes \Kk(Y|T). \]
\end{Not}

\begin{Def}[Transitional conditional independence]
    \label{def:tr-cond-ind}
    For our transition probability space $(\Wcal \times \Tcal, \Kk(W|T))$ and three transitional random variables  $\Xk \in \PrPf(\Xcal)^{\Wcal\times\Tcal}$ and $\Yk\in \PrPf(\Ycal)^{\Wcal\times\Tcal}$ and $\Zk \in \PrPf(\Zcal)^{\Wcal\times\Tcal}$ we consider their joint Markov kernel:
    \[ \Kk(X,Y,Z|T) := \lp  \Xk(X|W,T) \otimes \Yk(Y|W,T)\otimes \Zk(Z|W,T) \rp \circ \Kk(W|T) \in \PrPf(\Xcal \times \Ycal \times \Zcal)^\Tcal.\]
    We then say that $\Xk$ is \emph{independent of $\Yk$ conditioned on $\Zk$} w.r.t.\ $\Kk$, in symbols:
        \[\Xk\Indep_\Kk \Yk \given \Zk,\]
    if there exists a Markov kernel $\Qk(X|Z) \in \PrPf(\Xcal)^\Zcal$ such that:
    \[ \Kk(X,Y,Z|T) = \Qk(X|Z) \otimes \Kk(Y,Z|T) \in  \PrPf(\Xcal \times \Ycal \times \Zcal)^\Tcal. \]
\end{Def}

\begin{Thm}[Separoid rules for transitional conditional independence, see \Cref{thm:separoid_axioms-tci:prf}]
    \label{thm:separoid_axioms-tci}
    Consider a transition probability space $\lp \Wcal \times \Tcal, \Kk(W|T) \rp$ and transitional random variables
    $\Xk \in \PrPf(\Xcal)^{\Wcal\times\Tcal}$ and $\Yk \in \PrPf(\Ycal)^{\Wcal\times\Tcal}$ and $\Zk \in \PrPf(\Zcal)^{\Wcal\times\Tcal}$ and $\Uk \in \PrPf(\Ucal)^{\Wcal \times \Tcal}$.
    Then the ternary relation $\Indep = \Indep_{\Kk(W|T)}$ satisfies the following rules:
\begin{enumerate}[label=\alph*)]
    \item Extended Left Redundancy (if $\Delta_\Xcal \in \Bcal_{\Xcal \times \Xcal}$\footref{fn:diagonal}):
    \item[] $\Xk \ismapof_\Kk \Zk \implies  \Xk \Indep  \Yk \given \Zk$.
    \item[] \textit{Proof:} Lemma \ref{lem:indep-det}.
    \item $\Tk$-Restricted Right Redundancy (if $(\Xcal,\Zcal,\Tcal)$ is a disintegration triple,
    see \Cref{def:disint-triple-qus} and \Cref{thm:disint-X}):
    \item[] $\Xk\Indep \deltabf_0\given\Zk\otimes\Tk$ always holds.
    \item[] \textit{Proof:} Disintegration Theorem \ref{thm:disint-X}: $\Kk(X,Z|T) = \Qk(X|Z,T) \otimes \Kk(Z|T)$.
    \item Left Decomposition:
	\item[]$\Xk\otimes\Uk \Indep  \Yk \given \Zk \implies \Uk \Indep  \Yk \given \Zk$.
    \item[] \textit{Proof:} Take the marginal $\Qk(U|Z)$ of $\Qk(X,U|Z)$.
    \item Right Decomposition:
    \item[] $\Xk \Indep  \Yk\otimes\Uk \given \Zk \implies \Xk \Indep  \Uk \given \Zk$.
    \item[] \textit{Proof:} Take $\Qk(X|Z)$ again.
    \item $\Tk$-Inverted Right Decomposition:
    \item[] $\Xk \Indep  \Yk \given \Zk \implies \Xk \Indep  \Tk\otimes\Yk \given \Zk$.
    \item[] \textit{Proof:} Take $\Qk(X|Z)$ again.
    \item Left Weak Union (if $(\Xcal,\Ucal,\Zcal)$ is a disintegration triple, see
    \Cref{def:disint-triple-qus} and \Cref{thm:disint-X}):
    \item[] $\Xk\otimes\Uk \Indep  \Yk  \given \Zk \implies \Xk \Indep  \Yk \given \Uk\otimes\Zk$.
    \item[] \textit{Proof:} Disintegration Theorem \ref{thm:disint-X}: $\Qk(X,U|Z)=\Qk(X|U,Z) \otimes \Qk(U|Z)$.
    \item Right Weak Union:
  \item[] $\Xk \Indep  \Yk\otimes\Uk \given \Zk \implies \Xk \Indep  \Yk \given \Uk\otimes\Zk$.
     \item[] \textit{Proof:} Take $\Qk(X|U,Z):=\Qk(X|Z)$.
  \item Left Contraction:
  \item[] $ (\Xk \Indep  \Yk \given \Uk\otimes\Zk) \land (\Uk \Indep  \Yk \given \Zk) \implies \Xk\otimes\Uk \Indep  \Yk \given \Zk$.
     \item[] \textit{Proof:} Take $\Qk_1(X|U,Z) \otimes \Qk_2(U|Z)$.
\item Right Contraction:
\item[] $ (\Xk \Indep  \Yk \given \Uk\otimes\Zk) \land (\Xk \Indep\Uk \given \Zk) \implies \Xk \Indep \Yk\otimes\Uk \given \Zk$.
\item[] \textit{Proof:} Take $\Qk_2(X|Z)$ and use essential uniqueness of factorizations Theorem \ref{thm:fact-ess-uniq}.
\item Right Cross Contraction:
\item[] $ (\Xk \Indep  \Yk \given \Uk\otimes\Zk) \land (\Uk \Indep  \Xk \given \Zk) \implies \Xk  \Indep  \Yk\otimes\Uk \given \Zk$.
    \item[] \textit{Proof:} Take $\Qk_1(X|U,Z) \circ \Qk_2(U|Z)$ and use essential uniqueness of factorizations Theorem \ref{thm:fact-ess-uniq}.
\item Flipped Left Cross Contraction:
\item[] $ (\Xk \Indep  \Yk \given \Uk\otimes\Zk) \land (\Yk \Indep  \Uk \given \Zk) \implies \Yk \Indep  \Xk\otimes\Uk \given \Zk$.
    \item[] \textit{Proof:} Take $\Qk_2(Y|Z)$.
\end{enumerate}

\end{Thm}

\subsection{Conditional Directed Acyclic Graphs (CDAGs)}

In this subsection we will gather all purely graph theoretical structures that are needed for causal Bayesian networks later on. The main focus is on the definition of \emph{conditional directed acyclic graphs (CDAGs)} and on the recap of \emph{\ids-separation} in such graphs.

\begin{Def}[Conditional directed acyclic graphs (CDAGs)]
    \label{def:cdag}
    A \emph{conditional directed acyclic graph (CDAG)} $\Gk=(J,V,E)$ consists of two (disjoint) sets of
    vertices/nodes: the set of \emph{input nodes} $J$, the set of \emph{output nodes} $V$;
    and a set  of \emph{directed edges}:
    \[E \ins \lC w \tuh v\st w \in J \cup V, v \in V \rC,\]
    such that there are \emph{no directed non-trivial cycles}.
    Note that - per definition - there won't be any arrow heads pointing to input nodes $j \in J$.
\end{Def}

\begin{Not}
    \label{not:cdag}
    For a CDAG we will suggestively write: $\Gk:=(J,V,E)=\Gk$, where the sets of edges $E$ are implicit.
   We will also write: $v \in \Gk$ to mean $v \in J \cup V$, $v_1 \hut v_2 \in \Gk$ to mean $v_2 \tuh v_1 \in E$ and $v_1 \sus v_2 \in \Gk$ to mean that either $v_1\hut v_2 \in \Gk$ or $v_1 \tuh v_2 \in \Gk$, etc.
   \end{Not}

\begin{Def}[Walks]
     Let $\Gk=(J,V,E)$ be a CDAG and $v,w \in \Gk$.
\begin{enumerate}
    \item A \emph{walk} from $v$ to $w$ in $\Gk$ is a finite sequence of nodes and edges
        \[v=v_0 \sus v_1 \sus  \cdots \sus v_{n-1} \sus v_n=w\]
        in $\Gk$ for some $n \ge 0$, i.e.\ such that for every $k=1,\dots,n$
        we have that $v_{k-1} \sus v_k \in \Gk$, and with $v_0=v$ and $v_n=w$.
    \item[] In the definition of walk the appearance of the same nodes several times is allowed.
        Also  the \emph{trivial walk} consisting of a single node $v_0 \in \Gk$ (without any edges) is allowed as well (if $v=w$).
    \item A \emph{directed walk} from $v$ to $w$ in $\Gk$ is of the form:
        \[v=v_0 \tuh v_1 \tuh  \cdots \tuh v_{n-1} \tuh v_n=w,\]
        for some $n \ge 0$,
        where all arrow heads point in the direction of $w$ and there are no arrow heads pointing back.
  \end{enumerate}
\end{Def}

\begin{Def} For a CDAG $\Gk=(J,V,E)$ and $v \in \Gk$ we define the sets:
    \begin{enumerate}
        \item \emph{Parents}: $\Pa^\Gk(v) := \lC w \in \Gk\st w \tuh v \in \Gk \rC$,
        \item \emph{Children}: $\Ch^\Gk(v) := \lC w \in \Gk\st v \tuh w \in \Gk \rC$,
        \item \emph{Ancestors}: $\Anc^\Gk(v) := \lC w \in \Gk \st \exists \text{ directed walk in }\Gk:\,w\tuh \cdots \tuh v \rC$,
        \item \emph{Descendants}: $\Desc^\Gk(v) := \lC w \in \Gk \st \exists \text{ directed walk in }\Gk:\,v\tuh \cdots \tuh w \rC$,
\end{enumerate}
    We extend these notions to sets $A \ins J \cup V$ by taking the union, e.g.\ $\Anc^\Gk(A) = \bigcup_{v \in A} \Anc^\Gk(v)$.
\end{Def}

\begin{Rem}[Acyclicity]
    \label{rem:acylic}
    The requirement that a CDAG  $\Gk=(J,V,E)$ has no non-trivial cycles in the above terms means that
    the only directed walk from a node $v \in \Gk$ to itself inside $\Gk$ is trivial.
\end{Rem}

\begin{Def}[Topological order]
    Let $\Gk=(J,V,E)$ be a CDAG.
    A \emph{topological order} of $\Gk$ is a total order $<$ of $J \cup V$ such that for all $v,w \in \Gk$:
    \[ v \in \Pa^\Gk(w) \; \implies \; v < w.\]
    A topological order always exists for every CDAG.
\end{Def}

\begin{Def}[Extended CDAG]
    \label{def:ext-cdag}
    Let $\Gk=(J,V,E)$ be a CDAG and $W \ins V$ a subset. Then the \emph{extended CDAG} w.r.t.\ $W$ is the CDAG \[ \Gk_W := \lp J \dcup I_W,\; V,\; E \dcup \lC I_w \tuh w \st w \in W \rC \rp, \] that arises from $\Gk$ by adding the new nodes $I_w$ for $w \in W$ to $J$ and the new edges $I_w \tuh w$ to $E$.
\end{Def}

Here we will review the definition of \ids-separation in a CDAG, see \cite{Pearl09}, also see \cite{Ric03, FM17,FM18,FM20} for extensions and variations of \ids-separation to other graphs.

\begin{Def}[\ids-blocked walks]
    \label{def:d-blocked}
      Let $\Gk=(J,V,E)$ be a CDAG and $C \ins J \cup V$ a subset of nodes and $\pi$ a walk in $\Gk$:
      $ \pi =\lp  v_0 \sus \cdots \sus v_n \rp.  $
      \begin{enumerate}
          \item
      We say that the walk $\pi$ is \emph{$C$-\ids-blocked} or \emph{\ids-blocked by $C$} if either:
      \begin{enumerate}
          \item $v_0 \in C$ or $v_n \in C$ or:
          \item there are two adjacent edges in $\pi$ of one of the following forms:
      \end{enumerate}
              \[\begin{array}{rccl}
  \text{ left chain: } & v_{k-1} \hut v_k \hut v_{k+1} & \text{ with } & v_k \in C,\\
  \text{ right chain: } & v_{k-1} \tuh v_k \tuh v_{k+1} & \text{ with } & v_k \in C,\\
  \text{ fork: } & v_{k-1} \hut v_k \tuh v_{k+1} & \text{ with } & v_k \in C,\\
           \text{ collider: } & v_{k-1} \tuh v_k \hut v_{k+1} & \text{ with } & v_k \notin \Anc^\Gk(C).
              \end{array}\]
      \end{enumerate}
      In short, $\pi$ is $C$-\ids-blocked if it either has a \emph{collider} not in $\Anc^\Gk(C)$ or a \emph{non-collider} in $C$. Note that $C \ins \Anc^\Gk(C)$, so a collider blocks less often than membership in $C$ alone would suggest: it stays open as soon as one of its descendants is observed.
   \begin{enumerate}[resume]
    \item We say that the walk $\pi$ is \emph{$C$-\ids-open} if it is not $C$-$d$-blocked.
   \end{enumerate}
\end{Def}

\begin{Def}[\ids-separation]
    Let $\Gk=(J,V,E)$ be a CDAG and $A,B,C \ins J \cup V$ (not necessarily disjoint) subset of nodes.
    We then say that:
    \begin{enumerate}
        \item       \emph{$A$ is \ids-separated from $B$ given $C$ in $\Gk$}, in symbols:
            $\displaystyle \quad A \Perp_{\Gk} B \given C, $
        \item[] if every walk from a node in $A$ to a node in $J \cup B$\footnote{Note the inclusion of $J$ here. This is done to
            get similar asymmetric separoid rules to transitional conditional independence in order to
        have a more ``nice'' looking global Markov property later on.\label{fn:J-sep}} is \ids-blocked by $C$.
\item Otherwise we write:
    $\quad\displaystyle A \nPerp_{\Gk} B \given C.$
\item Special case:
    $\displaystyle \quad A \Perp_{\Gk} B \quad :\iff \quad A \Perp_{\Gk} B | \emptyset.$
        \end{enumerate}
\end{Def}

\subsection{Causal Bayesian Networks and the Global Markov Property}

In this subsection we will introduce \emph{causal Bayesian networks} for quasi-universal spaces.
The main result is the \emph{global Markov property} for causal Bayesian networks,
which relates \ids-separation statements in the graph to conditional independence statements between the variables.
Furthermore, we will introduce \emph{partially generic causal Bayesian networks} where some of the causal
mechanisms are not provided and are filled in with generic, most general ones.
Since this fill-in will lead to a well-defined causal Bayesian network we also get a global Markov property for partially generic causal Bayesian networks, now allowing us to relate \ids-separation statements in the extended graph to conditional independence statements between variables and causal mechanisms.
Note that this is only possible because in $\QUS$ we have exponential objects, transitional conditional independence still works with all its separoid rules for $\QUS$ and the global Markov property can be shown.
It is important to see that this kind of reasoning would break inside the usual category of measurable spaces $\Meas$.

\begin{Def}[Causal Bayesian network]
    \label{def:cbn}
    A \emph{causal Bayesian network (CBN)} $\Mk$ consists of:
    \begin{enumerate}
        \item a (finite) conditional directed acyclic graph (CDAG): $\Gk=(J,V,E)$,
        \item input variables $X_j$, $j \in J$, and
            (stochastic) output variables $X_v$, $v \in V$,
        \item a quasi-universal space $\Xcal_v$ for every $v \in J \dcup V$,
\item a Markov kernel for every $v\in V$, suggestively written as:
            \[\Pk_v\lp X_v\Big\|X_{\Pa^{\Gk}(v)}\rp \in \PrPf(\Xcal_v)^{\Xcal_{\Pa^{\Gk}(v)}},\]
        where we write for $D \ins J \dcup V$:
               \begin{align*}
                   \Xcal_D &:= \prod_{v \in D} \Xcal_v,&
                   \Xcal_\emptyset &:= \one=\lC0\rC,\\
                   X_D &:= (X_v)_{v \in D},&
                   X_\emptyset &:= 0,\\
                   x_D &:= (x_v)_{v \in D},&
                   x_\emptyset &:= 0.
              \end{align*}
    \end{enumerate}
    By abuse of notation, we denote the causal Bayesian network as:
    \[  \Mk(V\|J)= \lp \Gk, \lp  \Pk_v\lp X_v\Big\|X_{\Pa^{\Gk}(v)}\rp  \rp_{v \in V}  \rp. \]
\end{Def}

\begin{Not}
    Any CBN $\Mk =\Mk(V\|J)$  comes with its \emph{joint Markov kernel}:
    \[\Pk(X_V\|X_J) \in \PrPf(\Xcal_V) ^{\Xcal_J},   \]
            given by:
       \[ \Pk(X_V\|X_J) :=  \bigotimes_{v \in V}^> \Pk_v\lp X_v\Big\|X_{\Pa^{\Gk}(v)}\rp,\]
       where the product $\otimes^>$ is taken in reverse order of a fixed topological order $<$, i.e.\
       children appear only on the left of their parents in the product.
     Note that the
     joint Markov kernels does actually not depend on the actual choice of the topological order, if one swaps the corresponding arguments and spaces into the right positions accordingly, see Fubini Theorem \ref{thm:fubini}.
\end{Not}

Before stating the global Markov property we isolate the observation that lets the input variables
carry no separation hypothesis at all. Two of the separoid rules of \Cref{thm:separoid_axioms-tci}
--- $\Tk$-Restricted Right Redundancy and Left Weak Union --- are proved by an appeal to the
disintegration \Cref{thm:disint-X} and therefore carry side conditions on the space that sits in
the conditioning slot. In a causal Bayesian network the input variables are not arbitrary
transitional random variables: by \Cref{def:cbn} the variable $X_j$, $j \in J$, is the projection
$\pr_j:\, \Xcal_J \to \Xcal_j$, i.e.\ a \emph{deterministic function of the input} $X_J$, and the
input $X_J$ is --- per the convention recalled in \Cref{thm:gmp-mI-CBN} --- already present in the
conditioning slot of every transitional conditional independence statement occurring there. Such
components can therefore be absorbed into the input, by marginalisation and re-indexing rather
than by disintegration:

\begin{Lem}[Absorption of input variables, see \Cref{lem:absorb-input:prf}]
    \label{lem:absorb-input}
    Let $\lp\Wcal \times \Tcal, \Kk(W|T)\rp$ be a transition probability space and let
    $\Xk,\Yk,\Zk$ be transitional random variables as in \Cref{thm:separoid_axioms-tci}.
    Let furthermore $\Uk \in \PrPf(\Ucal)^{\Wcal \times \Tcal}$ be \emph{deterministic in the
    input}, i.e.\ of the form
    \[ \Uk(U|W,T) = \deltabf_{u(T)} \qquad \text{for a quasi-measurable map } u:\, \Tcal \to \Ucal.\]
    Then, with no hypothesis whatsoever on the quasi-universal space $\Ucal$:
    \begin{enumerate}
        \item \emph{(Left Weak Union across a deterministic input.)}
            \[ \Xk \otimes \Uk \Indep \Yk \given \Zk \otimes \Tk
               \quad\implies\quad
               \Xk \Indep \Yk \given \Uk \otimes \Zk \otimes \Tk. \]
        \item \emph{(Absorption in the conditioning slot.)} For every transitional random variable
            $\Xk$ and every $\Yk$,
            \[ \Xk \Indep \Yk \given \Zk \otimes \Tk
               \quad\implies\quad
               \Xk \Indep \Yk \given \Uk \otimes \Zk \otimes \Tk, \]
            and $\Tk$-Restricted Right Redundancy holds for $\Zk \otimes \Uk$ as soon as it holds
            for $\Zk$; in particular its side condition has to be checked for $\Zcal$ only, not
            for $\Zcal \times \Ucal$.
    \end{enumerate}
\end{Lem}

Applied with $\Uk := X_{D \cap J}$ for the various $D \ins J \dcup V$ that occur, this reduces
every side condition in \Cref{thm:separoid_axioms-tci} to the \emph{output} spaces $\Xcal_v$,
$v \in V$, which are standard by hypothesis in the next theorem, and hence in particular countably
separated. No hypothesis on the input spaces $\Xcal_j$, $j \in J$, is needed --- which matters,
because in the partially generic causal Bayesian networks of \Cref{def:gen-cbn} the input spaces
are exponential objects $\PrPf(\Xcal_w)^{\Xcal_{\Pa^\Gk(w)}}$, and nothing makes those countably
separated.

\begin{Thm}[Global Markov property for causal Bayesian networks, see \Cref{thm:gmp-mI-CBN:prf}]
\label{thm:gmp-mI-CBN}
Consider a causal Bayesian network $\Mk(V\|J)$ with graph $\Gk$
and joint Markov kernel: $\Pk(X_V\|X_J)$. Assume that $\Xcal_v$ is a standard quasi-universal space for every $v \in V$.\\
Then for all $A, B, C \ins J \dcup V$ (not-necessarily disjoint) we have the implication:
\[ A \Perp_{\Gk} B \given C \qquad \implies \qquad X_A \Indep_{\Pk(X_V\|X_J)} X_B \given X_C.   \]
Recall that we have - per definition - an implicit dependence on $J$, $X_J$, resp., in the second argument on each side.

\end{Thm}

\begin{Rem}
    The global Markov property says that already checking the graphical criterion $\displaystyle A \Perp_{\Gk} B \given C$ is
    enough to get
    the existence of a Markov kernel, suggestively written as:
    $\displaystyle \Pk\lp X_A|\cancel{X_B},X_{C \cap V}\|X_{C \cap J},\cancel{X_J} \rp,$
    such that:
    \[ \Pk(X_A,X_B,X_C\|X_J) = \Pk\lp X_A|\cancel{X_B},X_{C \cap V}\|X_{C \cap J},\cancel{X_J} \rp
    \otimes \Pk(X_B,X_C\|X_J). \]
\end{Rem}

\begin{Def}[The partially generic causal Bayesian network]
    \label{def:gen-cbn}
    Let $\Gk=(J,V,E)$ be a CDAG and $\Xcal_v$ for $v \in J \dcup V$ quasi-universal spaces,
$W \ins V$ a subset and $\Pk_v \in \PrPf(\Xcal_v)^{\Xcal_{\Pa^\Gk(v)}}$ for $v \in V \sm W$ given Markov kernels.
    Then the \emph{partially generic causal Bayesian network} w.r.t.\ those choices is the causal Bayesian network $\Mk(V\|J \dcup I_W)$ with the extended CDAG
    $\Gk_W$, see Definition \ref{def:ext-cdag}, the given spaces $\Xcal_v$ for $v \in J \dcup V$ and the following spaces for $w \in W$:
    \[ \Xcal_{I_w} := \PrPf(\Xcal_w)^{\Xcal_{\Pa^\Gk(w)}},\]
    and the \emph{generic Markov kernels} defined by:
    \[ \Pk_w\lp X_w \in A_w\|X_{\Pa^{G}(w)}=x_{\Pa^{G}(w)},X_{I_w}=\Qk_w\rp :=
       \Qk_w\lp X_w \in A_w\|X_{\Pa^{G}(w)}=x_{\Pa^{G}(w)}\rp.
    \]
    Note that: $\Pk_w \in \PrPf(\Xcal_w)^{\Xcal_{\Pa^\Gk(w)}\times \Xcal_{I_w}}$ for $w \in W$.
\end{Def}

\begin{Rem}
    \begin{enumerate}
        \item Note that the construction of the ``generic'' Markov kernels $\Pk_w$ in the Definition \ref{def:gen-cbn} of partially generic causal Bayesian networks would not have been possible in the category of measurable spaces $\Meas$ as the spaces of measurable maps or spaces of Markov kernels do not carry well-behaved $\sigma$-algebras, not even for standard (Borel) measurable spaces, see \cite{Aum61}. This was one of the core motivations for formalizing probabilistic graphical models inside the quasitopos of quasi-universal spaces $\QUS$, where such spaces of Markov kernels are objects like any other, see \Cref{qms:thm:quasitopos}.
        \item The partially generic causal Bayesian network, see Definition \ref{def:gen-cbn}, is again a causal Bayesian network, see Definition \ref{def:cbn}, and thus follows the global Markov property w.r.t.\ the extended graph $\Gk_W$, see Theorem \ref{thm:gmp-mI-CBN}:
    If $\Xcal_v$ is a standard quasi-universal space for every $v \in V$
           then for all $A,B,C \ins V \dcup I_W \dcup J$ we have the implication:
    \[ A \Perp_{\Gk_W} B \given C \qquad \implies \qquad X_A \Indep_{\Pk(X_V\|X_J,X_{I_W})} X_B \given X_C.   \]
    This will allow us to reason about the (conditional) independence of variables $X_v$ for $v \in J \dcup V$, w.r.t.\ ``mechanisms'' $\Qk_w$ for $w \in W$.
\item It is easy to extend the global Markov property, Theorem \ref{thm:gmp-mI-CBN}, to more general (causal) graphical models like causal Bayesian networks with latent confounders or even structural causal models
    that allow for cycles, latent confounders and selection bias. For this one mainly would need to adjust the graphical part, e.g.\ use the more general $\sigma$-separation, see \cite{FM17,FM18,FM20,For21}.
    This was also discussed in \cite{For21} Remark K.8:
    \[ A \Perp_{\Gk_S}^\sigma B \given C \qquad \implies \qquad
    X_A \Indep_{\Pk(X_V,X_S\|X_J)}X_B\given X_S,X_C.\]
    \end{enumerate}
\end{Rem}

\begin{Eg}[Markov chain]
    \label{eg:markov-chain}
    Consider the Markov chain of Figure \ref{fig:markov-chain}.
    We need to specify Markov kernels $P(X) \in \PrPf(\Xcal)$ and
    $P(Y|X) \in \PrPf(\Ycal)^\Xcal$ and $P(Z|Y) \in \PrPf(\Zcal)^\Ycal$.
    The joint distribution is then given by:
    \[ P(X,Y,Z) = P(Z|Y) \otimes P(Y|X) \otimes P(X).\]
    We then have that:
        \[ P(X,Y,Z) = P(Z|Y) \otimes P(X,Y),\]
        which shows the transitional conditional independence:
     \[ Z \Indep X \given Y. \]
    This means that when the value of $Y$ is provided then the value of $X$ has no additional information about the state of $Z$.\\
        Note that the global Markov property Theorem \ref{thm:gmp-mI-CBN} would have allowed us to conclude
    the same as we also have the corresponding \ids-separation: $Z \Perp X \given Y$. \\
    Since only the Markov kernel $P(Z|Y)$ determines the state of $Z$ when $Y$ is given
    we actually could make the stronger statement that even the distribution $P(X)$ and the Markov kernel
    $P(Y|X)$ have no influence on $Z$ when $Y$ is given. Unfortunately,
    this cannot directly be seen or read off from Figure \ref{fig:markov-chain}.\\
    However, such reasoning is possible with the use of the partially generic causal Bayesian network,
    see Figure \ref{fig:partially-generic-markov-chain}. Since the global Markov property, see Theorem \ref{thm:gmp-mI-CBN}, also applies here
    the \ids-separation $Z \Perp \lp X,Q(X),Q(Y|X) \rp \given Y$ gives us the transitional conditional independence:
    \[ Z \Indep \lp X,Q(X),Q(Y|X) \rp \given Y.\]
    Note that we use $Q(X)$ and $Q(Y|Z)$, resp., as a variables and $P(X)$ and $P(Y|Z)$, resp.,
    as specific values of those variables.
    So that transitional conditional independence now exactly states that when $Y$ is given the value of $Z$ is not determined by the values of $X$, $Q(X)$ and $Q(Y|X)$, as we expected and wanted, and without computing factorizations of type $P(X,Y,Z) = P(Z|Y) \otimes P(X,Y)$ by hand.
\end{Eg}

\subsection{Counterfactual Reasoning}

\begin{figure}[ht]
    \centering
    \begin{tikzpicture}[scale=1, transform shape]
\tikzstyle{every node} = [draw,shape=circle,color=blue]
        \node (X) at (0,0) {$X$};
        \node (Y) at (4,-2) {$Y$};
        \node (Z) at (2,3) {$Z$};
        \node (E) at (6,1) {$E$};
\foreach \from/\to in {X/Y, Z/E, E/Y, Z/X}
    \draw[-{Latex[length=3mm,width=2mm]}, color=blue] (\from) -- (\to);
     \end{tikzpicture}
   \caption{The causal Bayesian network encoding the strong ignorability assumption $X \Indep E \given Z$.
     The random variable $E$ takes values in the function space $\Ycal^\Xcal$ and $Y=E(X)$. The potential outcomes after intervening on $X$ with values $x \in \Xcal$ are then $Y_x=E(x)$ with values in $\Ycal$.}
    \label{fig:potential-outcome}
\end{figure}

One of the core assumptions for counterfactual reasoning in the potential outcome framework
is the assumption of \emph{strong ignorability}, see \cite{Rub78,RR83}:
\[ X \Indep (Y_0,Y_1) \given Z.   \]
The usual meaning of the variables is the following.
$Z$ is usually a random vector of all \emph{observed features} of a patient,
whom we want to prescribe one of two treatments from $\Xcal=\lC 0,1 \rC$.
The variable $X$ is then the measured \emph{treatment variable}, which can take the values from $\Xcal$.
The variables $Y_0$ and $Y_1$ are the potential outcomes if we forced the patient to take the corresponding treatment, $0$ or $1$. The assumption of strong ignorability can then be interpreted as that
besides $Z$ there are no further (unobserved) confounding variables, or that $Z$ is \emph{deconfounding}.

If one now wanted to generalize this scenario to more than two possible treatments and thus potential outcomes, i.e.\ $\#\Xcal >2$, or even $\Xcal=\R$, one quickly runs into the problem of even stating the assumption of strong ignorability measure-theoretically rigorously.
Many authors then put this into the following notation:
\[ X \Indep (Y_x)_{x\in \Xcal} \given Z.  \]
This notation indicates that $(Y_x)_{ x \in \Xcal}$ would live on some product space $\prod_{ x\in \Xcal} \Ycal$. Those authors then use those variables in a way that one implicitly needs to assume that the map:
\[ \ev:\; \Xcal \times \prod_{x \in \Xcal} \Ycal  \to \Ycal,\qquad
    \lp \hat x , (y_x)_{x \in \Xcal}\rp \mapsto y_{\hat x}, \]
    is measurable, even when $\Xcal$ is an uncountable (standard Borel) measurable space. This runs into similar problems as discussed in the introduction, see Section \ref{sec:intro} and \cite{Aum61} and Example \ref{eg:ev-prod-not-meas} below.

\begin{Eg}
    \label{eg:ev-prod-not-meas}
    Consider the measurable spaces $\Ycal=\lC0,1\rC$ and $\Xcal=\R$ and the map:
    \[ \ev:\; \Xcal \times \prod_{x \in \Xcal} \Ycal_x  \to \Ycal,\qquad
    \lp \hat x , (y_x)_{x \in \Xcal}\rp \mapsto y_{\hat x}, \]
    with $\Ycal_x:=\Ycal$.
    If we assumed that $\ev$ was measurable we would have:
    \[D:=\ev^{-1}(1) \in  \Bcal_\Xcal \otimes \bigotimes_{x \in \Xcal} \Bcal_{\Ycal_x}.\]
    Since the product-$\sigma$-algebra is generated by the finite cylinder sets and  by use of \cite{Bog07} Ex.\ 1.12.54 we then have a countable subset $\Zcal \ins \Xcal$,
    countable $\Acal \ins \Bcal_\Xcal$, countable $\Ecal_z \ins \Bcal_{\Ycal_z}$ for $z \in \Zcal$ such that:
    \[ D \in \sigma\lp\lC \pr_\Xcal^{-1}(A) \cap \pr_{z}^{-1}(B_z) \st A \in \Acal, B_z \in \Ecal_z, z \in \Zcal \rC\rp.\]
    This shows that $D$ will be of the form:
    \[ D = B \times \prod_{ x \in \Xcal \sm \Zcal} \Ycal_x,\]
    with a set $B \ins \Xcal \times \prod_{ z \in \Zcal} \Ycal_z$.
    On the one hand, if we take an element $(x,y_\Zcal) \in B$ then
      $(x,y_\Zcal,0,\dots,0) \in D$. It follows that:
    \[ 1 = \ev(x,y_\Zcal,0,\dots,0) = y_x,\]
    which implies $x \in \Zcal$.
    So we get: $\pr_\Xcal(B) \ins \Zcal$. \\
    On the other hand, if we take the element $(x,0,1_x,0,\dots,0)$ with $x \in \Xcal \sm \Zcal$
    then by definition of $\ev$ we get:
    \[ \ev(x,0,1_x,0,\dots,0) = 1. \]
    So $(x,0,1_x,0,\dots,0) \in \ev^{-1}(1) = D$ and thus $(x,0) \in B$.
    This gives the contradiction:
    \[\Xcal\sm \Zcal \ni x =\pr_\Xcal(x,0) \in \pr_\Xcal(B) \ins \Zcal.\]
    So $\ev$ can not be measurable w.r.t.\ the product-$\sigma$-algebra.
\end{Eg}

To overcome this issue it might be beneficial to work in the category of quasi-universal spaces $\QUS$ instead
and replace the random variable $(Y_x)_{ x \in \Xcal}$, which has values in $\prod_{x \in \Xcal}\Ycal$,
with a random variable $E$ that has values in the function space $\Ycal^\Xcal$. Note that its marginal distribution is then an element $P(E) \in \PrPf(\Ycal^\Xcal)$\footnote{Carefully distinguish between
$\PrPf(\Ycal^\Xcal)$ and $\PrPf(\Ycal)^\Xcal$, see Corollary \ref{cor:random-fct} and Remark \ref{rem:random-fct}.}.
Then one could easily state the assumption of \emph{strong ignorability} as:
\[  X \Indep E \given Z.\]
This would solve the measure-theoretic problems and we would get: $Y_x=E(x)$ as the potential outcomes for $x \in \Xcal$ and $Y = E(X)$ in the observational scenario.

If strong ignorability holds we can apply the disintegration theorem, \Cref{thm:disint-X}, to get the factorization of the joint distribution of all variables:
\[ P(X,Y,E,Z) = P(Y|X,E) \otimes P(X|Z) \otimes P(E|Z) \otimes P(Z).\]
Here the disintegration theorem applies under either of two readings of the covariate
    $Z$. If $(\Zcal,\Zcal^\Omega)$ is countably separated we are in the first case of
    \Cref{thm:disint-X}. If instead $Z$ is a \emph{discrete} covariate, which is the typical
    situation in this setting, the third case applies and no separation assumption on $\Zcal$ is
    needed at all.

This factorization can also be modelled with the use of a causal Bayesian network, see Figure \ref{fig:potential-outcome}. A similar graphical model, but without the variable $E$ or only after marginalizing out $E$, has been proposed before, see \cite{Pearl09}.

This example shows how one could, even in more generality, use random functions like $E$ with values in a function space $\Ycal^\Xcal$
to model potential outcomes and counterfactual relationships.
This could also lead to new graphical representations of counterfactuals similar to single world intervention graphs, see \cite{RR13a,RR13b,Mal19}. We leave this for future research.

\section{Discussion}
\label{sec:discussion}

\subsection{Results}
Motivated by the example of Markov chains \ref{eg:markov-chain}, where we wanted to be able to graphically check which variables $X_v$ are (conditionally) independent from certain Markov kernels $Q(X_v\|X_{\Pa^\Gk(v)})$, we investigated how one could allow variables in probabilistic graphical models to take values
in the (huge) space of Markov kernels $\PrPf(\Xcal_v)^{\Xcal_{\Pa^{\Gk}(v)}}$. \\
For this we needed to study the cartesian closed category of \emph{quasi-Borel spaces} from \cite{Heu17}.
During that process we generalized that theory to allow for general samples spaces and we also simplified definitions to be less restrictive. We then called those categories the \emph{categories of quasi-measurable spaces} $\QMS$.  We showed that those categories are not merely cartesian closed but \emph{locally} cartesian closed, see \Cref{qms:thm:lcc}, and in fact quasitoposes, see \Cref{qms:thm:quasitopos}, and Heyting categories, see \Cref{qms:thm:heyting}; that they contain all (small) limits and colimits; and that products and coproducts distribute.\\
We then systematically studied different \emph{probability monads} on those categories of quasi-measurable spaces. One corresponds to $\PrpF$, which was introduced in \cite{Heu17}, the others, $\PrQm$, $\PrPF$, $\PrPf$, $\Prpf$, are either novel or slight variation of the first one. It turned out that those probability monads are either strong without further assumptions, $\PrQm$, or under the assumption that the sample space satisfies: $\Omega \times \Omega \cong \Omega$, which holds for (countably) infinite product spaces like $\Omega= \Omega_0^\N$, but also for all uncountable standard Borel measurable spaces, like $\R$ or the Hilbert cube $[0,1]^\N$, endowed either with their Borel $\sigma$-algebra or the $\sigma$-algebra of all universally measurable subsets.\\
We then went on and specialized to such a sample space $\R^\N$, or equivalently $\R$, endowed with the $\sigma$-algebra of all universally measurable subsets. We called this the \emph{category of quasi-universal spaces} $\QUS$. It turns out that the induced $\sigma$-algebras on quasi-universal spaces are themselves quasi-universal spaces and have simple descriptions in terms of intersections of Lebesgue-complete $\sigma$-algebras, in contrast to the $\sigma$-algebras of quasi-Borel spaces. This becomes even more pronounced when one studies the $\sigma$-algebras of \emph{countably separated} quasi-universal spaces. A further special role was played by \emph{standard quasi-universal spaces}, which come with similar good properties as standard Borel measurable spaces have.\\
The strong probability monads $\PrPF$, $\PrPf$, $\PrpF$, $\Prpf$ all agree on the quasitopos of quasi-universal spaces $\QUS$. They have similar properties and semantics for higher probability theory like $\PrpF$ does for quasi-Borel spaces.\\
Furthermore, we proved for quasi-universal spaces a Fubini theorem, and, under certain conditions, a theorem about the disintegration of Markov kernels, Kolmogorov extension theorems, a conditional de Finetti theorem, etc. \\
We then translated many of the properties of standard quasi-universal spaces into properties of the Kleisli category of their Markov kernels w.r.t.\ the proposed probability monad.\\
Finally, we formalized transitional conditional independence, see \cite{For21}, and causal Bayesian networks inside the category of quasi-universal spaces and proved a global Markov property for them.
This then allowed us to solve the problem of our motivating example, namely to include the ``mechanisms'' of a given causal Bayesian network as variables/nodes on the same standing as other variables. This then allows us to check and reason about (conditional) independences between variables $X_v$ and mechanisms/Markov kernels $Q_w(X_w|X_{\Pa^\Gk(w)})$ graphically.

\subsection{Different Sample Spaces}

After studying the properties of the category of quasi-universal spaces $\QUS$ it seems to the author
that taking  an uncountable Polish space endowed with its $\sigma$-algebra of all universally measurable subsets as the sample space $\Omega$ for a category of quasi-measurable spaces $\QMS$ is an optimal choice.
We again highlight the arguments below.

It seems unavoidable to ask for an isomorphisms $\Omega \times \Omega \cong \Omega$ to have well-behaved probability monads, see Theorem \ref{thm:K-P-Q-markov-prod} and Theorem \ref{thm:K-P-R-strong-monad}.
In this case one then has the choice to either take $\PrPf$ or $\PrpF$ for the monads of
push-forward probability measures, as the others agree with those, see Lemma \ref{lem:K-P-S-R-eq}.

If one is interested in Kolmogorov extension theorems, see Theorem \ref{thm:kol-ext-mar-ker} and \cite{Kol33,Fremlin} 454D-G, or disintegration theorems, see Theorem \ref{cor:markov-kernel-disint-meas-II} and \cite{Fad85,Pac78,Fremlin,For21} Cor.\ C.8, one might want or need to restrict to \emph{perfect} probability measures or even countably compact ones.
Even though the product of perfect probability measures is perfect again, the arbitrary mixture of perfect measures might not be perfect anymore, see \cite{Ram79}. To avoid this problem one option would be to use a sample space where all probability distributions (on the product) are perfect, e.g.\ $\R$.
Alternatively, one could restrict to the probability monad $\PrPf$ instead of $\PrpF$, because when using  $\PrPf$ one only mixes probability measures via push-forwards of products of probability measures, see proof of Theorem \ref{thm:K-P-Q-markov-prod}, which would then be perfect again, if one could resolve the ambiguity between $\QMS$-products and $\Meas$-products for the sample space $\Omega$ (as done for $\QUS$ and $\QBS$), see \cite{Fremlin} 451J.

A further variant, which avoids requiring that all probability measures on $\Omega$ be perfect by restricting to the perfect ones, is recorded in \Cref{qms:sec:app-other-monads}.

\section*{Acknowledgments}
\addcontentsline{toc}{section}{Acknowledgments}
The author wants to thank Chris Heunen, Sam Staton, Ohad Kammar, Paolo Perrone and Hongseok Yang for explaining many
aspects of quasi-Borel spaces to him and their helpful email correspondence.

\addcontentsline{toc}{section}{References}

{\catcode`\_=12
\newcommand{\etalchar}[1]{$^{#1}$}
\providecommand{\bysame}{\leavevmode\hbox to3em{\hrulefill}\thinspace}
\providecommand{\MR}{\relax\ifhmode\unskip\space\fi MR }
\providecommand{\MRhref}[2]{%
  \href{http://www.ams.org/mathscinet-getitem?mr=#1}{#2}
}
\providecommand{\href}[2]{#2}

}

\newpage

\newpage
\appendix
\section{Appendix - The Category of Quasi-Measurable Spaces}
\label{sec:cat-qms}

\subsection{Quasi-Measurable Spaces and Maps}

\begin{Not}
    \begin{enumerate}
        \item For two sets $\Xcal$ and $\Ycal$ we denote by:
\begin{align} \Sets(\Xcal,\Ycal):=\lC f:\, \Xcal \to \Ycal \rC\end{align}
the set of all (set-theoretic) maps from $\Xcal$ to $\Ycal$.
        \item For a set $\Xcal$ we identify $\Xcal^\one$ with the set of all \emph{constant} maps:
    \begin{align}
        g:\, \Omega \to \one \to \Xcal.
    \end{align}
\end{enumerate}
\end{Not}

\begin{Eg}
    Let $\Omega=\R$ and $\Omega^\Omega=\Top(\R,\R)$ the set of continuous maps from $\R$ to $\R$,
    which is closed under composition and contains the identity and all constant maps.
    Then we have that:
    \begin{align}\Bcal_\Omega:= \sigma\lp\lC \Phi_2^{-1}\lp\Phi_1(\Omega)\rp\st \Phi_1,\Phi_2 \in
            \Omega^\Omega \rC\rp,\end{align}
            is the Borel $\sigma$-algebra of $\R$, and:
    \begin{align}\Bcal_\Omega:= \lp\lC \Phi_2^{-1}\lp\Phi_1(\Omega)\rp\st \Phi_1,\Phi_2 \in
            \Omega^\Omega \rC\rp_\Giry,\end{align}
    the $\sigma$-algebra of all universally measurable subsets of $\R$.
    In any case, we have:
    \begin{align}\Omega^\Omega \ins \Meas\lp(\Omega,\Bcal_\Omega),(\Omega,\Bcal_\Omega) \rp.\end{align}
    \begin{proof}
        Let $\Phi_1,\Phi_2 \in \Top(\R,\R)$.
        Then $\Phi_1(\R)$ is the countable union of closed compact subsets of $\R$.
        So $\Phi_2^{-1}\lp\Phi_1(\R)\rp \ins \R$ is a countable union of closed subsets of $\R$
        thus a Borel subset of $\R$. Since we can generate all intervals $(a,b)=\Phi_1(\R)$ for
        appropriate
        $\Phi_1$ and with $\Phi_2=\id_\R$ we generate the whole Borel $\sigma$-algebra of $\R$.
    \end{proof}
\end{Eg}

\begin{Lem}
    \label{lem:qm-maps-xo}
    Let $(\Xcal,\Xcal^\Omega)$ be a quasi-measurable space. Then we have the identification:
    \begin{align}  \QMS((\Omega,\Omega^\Omega),(\Xcal,\Xcal^\Omega)) = \Xcal^\Omega.  \end{align}
\begin{proof}
        For $f \in \Xcal^\Omega$ we
        by the point 1.\ of Definition \ref{qms:def:quasi-measurable-spaces} have:
        \begin{align} f \circ \Omega^\Omega \ins \Xcal^\Omega.  \end{align}
        This directly implies the inclusion:
        \begin{align} \Xcal^\Omega \ins \QMS((\Omega,\Omega^\Omega),(\Xcal,\Xcal^\Omega)). \end{align}
        If, on the other hand, $g \in \QMS((\Omega,\Omega^\Omega),(\Xcal,\Xcal^\Omega))$ then we have:
        \begin{align} g \circ \Omega^\Omega \ins \Xcal^\Omega.  \end{align}
        Since $\id_\Omega \in \Omega^\Omega$ we get that:
        \begin{align} g=g \circ \id_\Omega \in \Xcal^\Omega.\end{align}
        We thus also get:
        \begin{align} \QMS((\Omega,\Omega^\Omega),(\Xcal,\Xcal^\Omega)) \ins \Xcal^\Omega.  \end{align}
        This shows the equality.
    \end{proof}
\end{Lem}

\subsection{The Product Space}
\label{qms:sec:product}

\begin{Lem}
    The product of quasi-measurable spaces w.r.t.\ $\Omega$:
    \begin{align} \prod_{i \in I}(\Xcal_i,\Xcal_i^\Omega):= \lp \prod_{i \in I} \Xcal_i,\, \prod_{i \in I} \Xcal_i^\Omega \rp  \end{align}
    is a quasi-measurable space w.r.t.\ $\Omega$. \\
    \begin{proof}
        We check the points from Definition \ref{qms:def:quasi-measurable-spaces}:
        \begin{enumerate}
            \item For $ X=( X_i)_{ i \in I} \in \prod_{i \in I} \Xcal_i^\Omega$ and $\Phi \in \Omega^\Omega$ we by assumption have:
                \begin{align} \forall i \in I.\;  X_i \circ \Phi \in \Xcal_i^\Omega.\end{align}
                It follows that:
                \begin{align} X \circ \Phi = ( X_i \circ \Phi)_{i\in I} \in \prod_{i \in I} \Xcal_i^\Omega.\end{align}
            \item If $ X = ( X_i)_{ i \in I}$ is a constant map, then every $ X_i$ is a constant map.
                By assumption such constant $ X_i \in \Xcal_i^\Omega$. It then follows that also $ X \in \prod_{i\in I} \Xcal_i^\Omega$.
        \end{enumerate}
      This shows that $\prod_{i \in I}(\Xcal_i,\Xcal_i^\Omega)$ is a quasi-measurable space.
    \end{proof}
\end{Lem}

\begin{Lem}
    \label{qms:lem:qms-product-up}
     The product of quasi-measurable spaces:
    \begin{align} \prod_{i \in I}(\Xcal_i,\Xcal_i^\Omega):= \lp \prod_{i \in I} \Xcal_i,\, \prod_{i \in I} \Xcal_i^\Omega \rp  \end{align}
    defines a categorical product in $\QMS$,
    i.e.\ for every quasi-measurable space $(\Zcal,\Zcal^\Omega)$ the natural map:
    \begin{align} \QMS\lp\Zcal, \prod_{i\in I} \Xcal_i\rp \bij \prod_{i \in I} \QMS\lp\Zcal,\Xcal_i\rp,
    \qquad g \mapsto (\pr_i \circ g)_{i \in I}   \end{align}
    is a bijection.
    \begin{proof}
      First note that the canonical projection map:
      \begin{align}  \pr_k:\,\prod_{i \in I}(\Xcal_i,\Xcal_i^\Omega) \to (\Xcal_k,\Xcal_k^\Omega), \end{align}
      is a  quasi-measurable map for every $k \in I$.\\
      Next let $(\Zcal,\Zcal^\Omega)$ be another quasi-measurable space and:
      \begin{align} g_k:\, (\Zcal,\Zcal^\Omega) \to (\Xcal_k,\Xcal_k^\Omega) \end{align}
      a quasi-measurable map for each $k \in I$.
      Then $g:=(g_i)_{i \in I}$ is the only (set theoretic) map that makes the following diagram commute for all $k \in I$:
      \begin{align}\xymatrix{
              (\Zcal,\Zcal^\Omega) \ar_-{g_k}[rrd] \ar^-{\exists! g}[rr] && \prod_{i \in I}(\Xcal_i,\Xcal_i^\Omega) \ar^-{\pr_k}[d] \\
                        && (\Xcal_k,\Xcal_k^\Omega).
      }\end{align}
      It is left to check that $g$ is a quasi-measurable map.
      If $ Z \in \Zcal^\Omega$ by assumption we know that $g_k \circ  Z \in \Xcal_k^\Omega$ for all $k \in I$.
      It thus follows that:
      \begin{align} g \circ  Z = (g_i \circ  Z)_{i \in I} \in \prod_{i \in I} \Xcal_i^\Omega.  \end{align}
      This thus shows the claim.
    \end{proof}
\end{Lem}

\subsection{The Function Space}
\label{qms:sec:fct-space}

In this subsection we will revisit the \emph{function space} of quasi-measurable spaces and provide the proofs of all the results.

\begin{Lem}
     Let $(\Xcal,\Xcal^\Omega)$ and $(\Ycal,\Ycal^\Omega)$ be two quasi-measurable spaces.
    Then the function space $\lp\Ycal^\Xcal,\lp \Ycal^\Xcal \rp^\Omega\rp$ is a quasi-measurable space. \begin{proof}
  We check the points from Definition \ref{qms:def:quasi-measurable-spaces}:
     \begin{enumerate}
        \item If $ Y \in \lp\Ycal^\Xcal\rp^\Omega$ and $\Phi',\Phi \in \Omega^\Omega$ and $ X \in \Xcal^\Omega$
            then we get:
            \begin{align}  Y(\Phi' \circ \Phi)( X) \in \Ycal^\Omega, \end{align}
            since $\Phi' \circ \Phi \in \Omega^\Omega$. By definition we then get:
            \begin{align}  Y \circ \Phi' \in\lp\Ycal^\Xcal\rp^\Omega. \end{align}
        \item If $ Y:\, \Omega \to \Ycal^\Xcal$ is constant then there is a $\Phi \in \Ycal^\Xcal$ such that
                for all $\omega$ we have $ Y(\omega)=\Phi$. Since:
                \begin{align} \Phi \in \Ycal^\Xcal = \QMS(\Xcal,\Ycal),\end{align}
                we have for all $ X \in \Xcal^\Omega$ that $\Phi( X) \in \Ycal^\Omega$.
                So for all $\Phi' \in \Omega^\Omega$, $ X \in \Xcal^\Omega$ we have:
                \begin{align} Y(\Phi')( X)=\Phi( X) \in \Ycal^\Omega.\end{align}
                By definition of $\lp \Ycal^\Xcal \rp^\Omega$ we get that: $ Y \in \lp \Ycal^\Xcal \rp^\Omega$.
    \end{enumerate}
    This shows the claims.
\end{proof}
\end{Lem}

\begin{Lem}
    \label{qms:lem:function-space}
    Let $(\Xcal,\Xcal^\Omega)$, $(\Ycal,\Ycal^\Omega)$ be quasi-measurable spaces.
    Then $\lp\Ycal^\Xcal,\lp \Ycal^\Xcal \rp^\Omega\rp$ has the universal property of a function space in the category $\QMS$ of all quasi-measurable spaces w.r.t.\ $\Omega$, i.e.\ for any other quasi-measurable space $(\Zcal,\Zcal^\Omega)$ the canonical curry map:
    \begin{align} \QMS(\Zcal \times \Xcal,\Ycal) \bij \QMS(\Zcal,\Ycal^\Xcal),\qquad f \mapsto (z \mapsto (x\mapsto f(z,x))),  \end{align}
    is a natural bijection.
    \begin{proof}
    For the universal property now consider a fixed quasi-measurable map $g:\, \Zcal \to \Ycal^\Xcal$.
    Then the map $\tilde g:\, \Zcal \times \Xcal \to \Ycal$ given by:
    \begin{align} \tilde g(z,x):= g(z)(x),\end{align}
    is quasi-measurable. Indeed, for all $Z \in \Zcal^\Omega$ and $ X \in \Xcal^\Omega$ we have that:
    \begin{align} \tilde g (Z, X) = g(Z(\id_\Omega))( X) \in \Ycal^\Omega.  \end{align}
    The inverse construction is given for quasi-measurable $f:\, \Zcal \times \Xcal \to \Ycal$
    by $\hat f:\, \Zcal \to \Ycal^\Xcal$ via:
    \begin{align} \hat f(z)(x) := f(z,x). \end{align}
    To show the quasi-measurability of $\hat f$ let $ Z \in \Zcal^\Omega$ and
    $\Phi \in \Omega^\Omega$ and $ X \in \Xcal^\Omega$.
    We need to show that:
    \begin{align} \hat f( Z(\Phi))( X) = f( Z(\Phi),  X) \in \Ycal^\Omega,\end{align}
    which holds by assumption on $f$ and since $ Z \circ \Phi \in \Zcal^\Omega$.\\
    It is then easy to see that $g \mapsto \tilde g$ and $f \mapsto \hat f$ are inverse to each other.\\
    Finally, we mention that this bijection is natural in all three arguments, which can just be checked via composition.
    \end{proof}
\end{Lem}

\begin{Thm}
    \label{qms:thm:qms-cartesian-closed}
     $(\QMS,\times,\one)$ is a cartesian closed  category.
    \begin{proof}
        The product $\times$ constitutes a \emph{functor} from the product category to $\QMS$:
        \begin{align}
            \begin{split}
                \times:\; \QMS \times \QMS &\to \QMS,\\
                (f:\Xcal \to \Ycal,\; g:\Wcal \to \Zcal) &\mapsto (f\times g:\, \Xcal \times \Wcal \to \Ycal \times \Zcal).
            \end{split}
        \end{align}
        The \emph{unit object} is $\one \in \QMS$. \\
        We have an \emph{associator} natural isomorphism, given for all triples as:
        \begin{align} A_{\Xcal,\Ycal,\Zcal}=\id:\; (\Xcal \times \Ycal) \times \Zcal \cong  \Xcal \times (\Ycal \times \Zcal),   \end{align}
        which in this case is always the identity (definition of \emph{strict}).\\
        We have the \emph{left and right unitor} natural isomorphisms, on objects given as:
        \begin{align} L_\Xcal=\pr_\Xcal:\;  \one \times \Xcal \cong \Xcal,\qquad R_\Xcal =\pr_\Xcal:\; \Xcal \times \one \cong \Xcal, \end{align}
        which here are the canonical projections on the first and second component, resp.\\
        We have the \emph{braiding} natural isomorphism, on objects given as:
        \begin{align} B_{\Xcal,\Ycal}:\; \Xcal \times \Ycal \cong \Ycal \times \Xcal,  \end{align}
        which here are given by swapping first and second entry.\\
        Furthermore, the associator satisfies the \emph{pentagon identity} (associativity law),
        which states the commutativity of the following diagram:
        \begin{align}\xymatrix{
        & (\Wcal \times \Xcal)\times(\Ycal \times \Zcal) \ar^-{A_{\Wcal,\Xcal,\Ycal\times\Zcal}}[rd] \\
        ((\Wcal\times\Xcal)\times\Ycal)\times\Zcal \ar^-{A_{\Wcal\times\Xcal,\Ycal,\Zcal}}[ru]
    \ar_-{A_{\Wcal,\Xcal,\Ycal}\times \id_\Zcal}[dd]
        &&\Wcal\times(\Xcal\times(\Ycal\times\Zcal))  \\
                \\
                (\Wcal\times(\Xcal\times\Ycal))\times\Zcal \ar_-{A_{\Wcal,\Xcal\times\Ycal,\Zcal}}[rr]
        &&\Wcal\times((\Xcal\times\Ycal)\times\Zcal). \ar_-{\id_\Wcal \times A_{\Xcal,\Ycal,\Zcal}}[uu]
        }\end{align}
        Indeed, since the associator is given by the identity there
        is nothing much to show.\\
        The left and right unitors together with the associator satisfy the \emph{triangle identity}, which states the commutativity of the following diagram:
        \begin{align}\xymatrix{
                (\Xcal \times \one) \times \Ycal \ar^-{A_{\Xcal,\one,\Ycal}}[rr]
                \ar_-{R_\Xcal \times \id_\Ycal}[dr]&&
                \Xcal \times (\one \times \Ycal) \ar^-{\id_\Xcal \times L_\Ycal}[dl]\\
        &\Xcal \times \Ycal.
        }\end{align}
        The braiding and the associator satisfy the (first) \emph{hexagon identity},
        which states the commutativity of the following diagram:
        \begin{align}\xymatrix{
                (\Xcal \times \Ycal) \times \Zcal \ar^-{A_{\Xcal,\Ycal,\Zcal}}[rr] \ar_-{B_{\Xcal,\Ycal}\times\id_\Zcal}[d] & &
                \Xcal\times(\Ycal\times\Zcal) \ar^-{B_{\Xcal,\Ycal\times\Zcal}}[rr] &&
                (\Ycal \times \Zcal)\times\Xcal \ar^-{A_{\Ycal,\Zcal,\Xcal}}[d]\\
                (\Ycal\times\Xcal)\times\Zcal  \ar_-{A_{\Ycal,\Xcal,\Zcal}}[rr] &&
                \Ycal\times(\Xcal\times\Zcal) \ar_-{\id_\Ycal\times B_{\Xcal,\Zcal}}[rr]&&
                 \Ycal\times(\Zcal\times\Xcal).
        }\end{align}
    Finally, the braiding satisfies the \emph{symmetry} relation:
        \begin{align} B_{\Ycal,\Xcal} \circ B_{\Xcal,\Ycal} = \id_{\Xcal \times \Ycal},\end{align}
        which in this case needs no explanation.\\
        It follows that $(\QMS,\times,\one)$ is a symmetric (braided strict) monoidal category.
        Since the monoidal structure is given by the categorical product it is a (symmetric) \emph{cartesian (monoidal) category}.
        By Lemma \ref{qms:lem:function-space} the category $\QMS$ is then \emph{cartesian closed},
        i.e.\ it is closed w.r.t.\ its cartesian monoidal structure,
        i.e.\ there is always an exponential object
        $\Zcal^\Ycal=\QMS(\Ycal,\Zcal)$ and
        a natural adjunction isomorphism:
        \begin{align} \QMS(\Xcal \times \Ycal, \Zcal) \cong \QMS(\Xcal, \Zcal^\Ycal),\end{align}
        that is natural in all three arguments.
    \end{proof}
\end{Thm}

\begin{Lem}
    \label{lem:exp-qm}
    Let $(\Xcal,\Xcal^\Omega)$, $(\Ycal,\Ycal^\Omega)$, $(\Zcal,\Zcal^\Omega)$  be quasi-measurable spaces.
    Consider a quasi-measurable map:
    \begin{align} f:\; \Xcal \to \Ycal. \end{align}
    Then the induced map:
    \begin{align} f^\Zcal:\; \Xcal^\Zcal \to \Ycal^\Zcal, \quad  X \mapsto f \circ  X, \end{align}
    is a well-defined quasi-measurable map.\\
    Furthermore, the well-defined map:
    \begin{align} \QMS(\Xcal,\Ycal) \to \QMS(\Xcal^\Zcal,\Ycal^\Zcal), \qquad f \mapsto f^\Zcal,\end{align}
    is quasi-measurable.
    \begin{proof}
        $f^\Zcal$ is well-defined, i.e.\  $f \circ  X \in \Ycal^\Zcal$ for $ X \in \Xcal^\Zcal$.
        Indeed, the compositum of two quasi-measurable maps is quasi-measurable:
        \begin{align} \Zcal \xrightarrow{ X} \Xcal \xrightarrow{f} \Ycal. \end{align}
        Furthermore, $f^\Zcal$ is quasi-measurable. If $ X \in (\Xcal^\Zcal)^\Omega=\Xcal^{\Omega \times \Zcal}$
        then $f \circ  X \in \Ycal^{\Omega \times \Zcal}=(\Ycal^\Zcal)^\Omega$ by the same argument as above.
        This shows that the map:
        \begin{align} (\_)^\Zcal:\; \QMS(\Xcal,\Ycal) \to \QMS(\Xcal^\Zcal,\Ycal^\Zcal), \qquad f \mapsto f^\Zcal,\end{align}
        is well-defined. To show that it is quasi-measurable let $h \in \QMS(\Omega \times \Xcal,\Ycal)$.
        We need to show that:
        $h^\Zcal \in \QMS(\Omega \times \Xcal^\Zcal,\Ycal^\Zcal)$, which reduces to show that the (well-defined) map:
        \begin{align} \Omega \times \Xcal^\Zcal \times \Zcal \to \Ycal,\qquad (\omega, X,z) \mapsto h(\omega, X(z)),  \end{align}
        is quasi-measurable. For this let
        $\Phi \in \Omega^\Omega$, $ X \in \Xcal^{\Omega \times \Zcal}$ and $ Z \in \Zcal^\Omega$.
        It is then left to show that the map:
        \begin{align} \Omega \to \Ycal,\qquad \omega \mapsto h(\Phi(\omega), X(\omega), Z(\omega)), \end{align}
        is quasi-measurable, which is true as the composition of the quasi-measurable maps:
        \begin{align}\Omega \xrightarrow{\id_\Omega \times \id_\Omega \times  Z} \Omega \times \lp \Omega \times \Zcal\rp \xrightarrow{\Phi \times  X}  \Omega \times \Xcal \xrightarrow{h}  \Ycal.\end{align}
        This shows the claim.
    \end{proof}
\end{Lem}

\subsection{The Coproduct Space}
\label{qms:sec:coproduct}

In this subsection we define the categorical \emph{coproduct} of quasi-measurable spaces.
Note that the coproduct we define here is different from the coproduct defined in \cite{Heu17}.
The reason is that we did not insist on the 3rd property in Def.\ 7 from \cite{Heu17}, making our coproduct in some sense more ``rigid''.

\begin{Def}[The coproduct space]
    Let $(\Xcal_i,\Xcal_i^\Omega)$, $i\in I$, be a family of quasi-measurable spaces indexed by any set $I$.
    Then we define their \emph{coproduct space} as:
    \begin{align}
        \coprod_{i \in I} \lp \Xcal_i,\Xcal_i^\Omega\rp &:= \lp\coprod_{i \in I} \Xcal_i, \coprod_{i \in I} \Xcal_i^\Omega\rp,
    \end{align}
    where on the right we have the set-theoretic coproducts:
    \begin{align}
        \coprod_{i \in I} \Xcal_i &= \lC (i,x) \st i \in I, x \in \Xcal_i\rC, &
        \coprod_{i \in I} \Xcal_i^\Omega &= \lC (i,X) \st i \in I, X \in \Xcal_i^\Omega\rC.
    \end{align}
\end{Def}

One can easily check that $\coprod_{i \in I} \lp \Xcal_i,\Xcal_i^\Omega\rp$ is again a quasi-measurable space.
It is also clear that the canonical inclusion maps:
     \begin{align} \incl_k:\, \Xcal_k \to \coprod_{i \in I} \Xcal_i,\qquad x \mapsto (k,x),\end{align}
 are quasi-measurable by definition of $\lp\coprod_{i \in I} \Xcal_i \rp^\Omega$.

Furthermore, we have:

\begin{Lem}
    The coproduct $\coprod_{i \in I} \lp \Xcal_i,\Xcal_i^\Omega\rp$ defines a categorical coproduct of $\QMS$,
    i.e.\ we have a natural isomorphism:
\begin{align}
    \QMS\lp\coprod_{i\in I} \Xcal_i, \Zcal\rp &\bij \prod_{i \in I} \QMS\lp\Xcal_i,\Zcal\rp,
     &g &\mapsto (g \circ \incl_i)_{i \in I} ,
\end{align}
with inclusion maps: $\incl_j:\, \Xcal_j \to \coprod_{i \in I}\Xcal_i$, $x \mapsto (j,x)$, for $j \in J$.
 \begin{proof}
     For the universal property
     let $(\Zcal,\Zcal^\Omega)$ be another quasi-measurable space and
     $g_i:\, (\Xcal_i,\Xcal_i^\Omega) \to (\Zcal,\Zcal^\Omega)$ quasi-measurable maps for
     $i\in I$.
     Then the map:
     \begin{align}g:= \coprod_{i \in I} g_i:\; \coprod_{i \in I} \Xcal_i \to \Zcal,\qquad (k,x) \mapsto g_k(x),\end{align}
         is the only (set-theoretic) map that makes the following diagrams commute for all $k \in I$:
      \begin{align}\xymatrix{
                  (\Xcal_k,\Xcal_k^\Omega) \ar_-{\incl_k}[d] \ar^-{g_k}[rrd] \\
                  \coprod_{i \in I}(\Xcal_i,\Xcal_i^\Omega) \ar_-{\exists!g}[rr]& & (\Zcal,\Zcal^\Omega).
      }\end{align}
      We are left to check that $g$ is a quasi-measurable map.
      We already have by assumption that:
      \begin{align} g \circ \incl_k \circ  X_k = g_k \circ  X_k \in \Zcal^\Omega,\end{align}
      for every $k \in I$ and $ X_k \in \Xcal_k^\Omega$.
       So we get:
 \begin{align} g \circ \lp\coprod_{i \in I} \Xcal_i \rp^\Omega =
             \lC g\circ \incl_i\circ  X_i \st  X_i \in \Xcal_i^\Omega, \,i \in I  \rC
        \ins \Zcal^\Omega.\end{align}
    This shows the claim.
    \end{proof}
 \end{Lem}

Furthermore, we have the \emph{distributive law} for product and coproduct:

\begin{Thm}[Distributive law for product and coproduct]
    \label{qms:thm:coprod-prod-distr}
    Let $I$ be any (index) set and for every $i \in I$ let $J_i$ be another (index) set.
    For every $i \in I$ and $j_i \in J_i$ let $(\Xcal_{i,j_i},\Xcal_{i,j_i}^\Omega)$ be a quasi-measurable space.
    Then we have a canonical isomorphism of quasi-measurable spaces:
    \begin{align}
        \coprod_{j \in \prod_{i \in I}J_i} \lp\prod_{i\in I} \Xcal_{i,j_i} \rp
    &\cong \prod_{i \in I} \lp \coprod_{j_i \in J_i} \Xcal_{i,j_i}\rp,
 & \lp j= (j_i)_{i\in I}, (x_i)_{i \in I} \rp &\mapsto \lp j_i,x_i\rp_{i \in I}.
\end{align}
\begin{proof}
    This is clearly a bijection.
    Similarly note that this map induces the bijection:
    \begin{align*}
\lp\coprod_{j \in \prod_{i \in I}J_i} \lp\prod_{i\in I} \Xcal_{i,j_i} \rp\rp^\Omega
&= \coprod_{j \in \prod_{i \in I}J_i} \lp\prod_{i\in I} \Xcal_{i,j_i} \rp^\Omega &
&= \coprod_{j \in \prod_{i \in I}J_i} \lp \prod_{i\in I} \Xcal_{i,j_i}^\Omega \rp \\
& \cong \prod_{i \in I} \lp \coprod_{j_i \in J_i} \Xcal_{i,j_i}^\Omega\rp &
& = \prod_{i \in I} \lp \coprod_{j_i \in J_i} \Xcal_{i,j_i}\rp^\Omega \\
& = \lp \prod_{i \in I} \lp \coprod_{j_i \in J_i} \Xcal_{i,j_i}\rp\rp^\Omega.
    \end{align*}
This shows the claim.
\end{proof}
\end{Thm}

\begin{Eg}
    \label{qms:eg:coprod-two-point}
    Consider the two-point space $\twoc:=\lC 0,1\rC$ where $\twoc^\Omega$ only contains the two constant maps
    $0, 1 :\Omega \to \lC 0,1 \rC=\twoc$.
    Then we have an isomorphism of quasi-measurable spaces:
    \begin{align}
        \twoc \cong \one \sqcup \one.
    \end{align}
\end{Eg}

\subsection{Equalizers and Coequalizers}
\label{qms:sec:equalizers}

In this subsection we gather the remaining categorical building blocks for all other categorical (small) \emph{limits} and (small)\footref{qms:fn:small} \emph{colimits} of $\QMS$,
like \emph{equaliziers} and \emph{coequalizers}. This then shows that $\QMS$ is a complete and cocomplete category.

\begin{Def}[Equalizer and coequalizers]
    Let $(\Xcal,\Xcal^\Omega)$ and $(\Ycal,\Ycal^\Omega)$ be quasi-measurable spaces and
    $ f_1,f_2:\, (\Xcal,\Xcal^\Omega) \to (\Ycal,\Ycal^\Omega)$
     quasi-measurable maps.
     Then we define the \emph{equalizer} and \emph{coequalizer}, resp., of $f_1$ and $f_2$ as:
     \begin{align}
         \Eq(f_1,f_2) &:= \lC x\in \Xcal\st f_1(x)=f_2(x)\rC, &
         \CoEq(f_1,f_2) &:= \Ycal/\mathord{\sim}, \\
         \Eq(f_1,f_2)^\Omega &:=\lC X:\,\Omega \to \Eq(f_1,f_2) \st \iota \circ X \in \Xcal^\Omega \rC, &
         \CoEq(f_1,f_2)^\Omega &:= \gamma \circ \Ycal^\Omega,
     \end{align}
     where $\iota:\,\Eq(f_1,f_2) \inj \Xcal$ is the inclusion map and $\gamma:\,\Ycal \srj \CoEq(f_1,f_2)$ is the quotient map and
     where $\sim$ is the smallest equivalence relation such that: $f_1(x) \sim f_2(x)$ for $x \in \Xcal$.
\end{Def}

It can easily be seen that $\Eq(f_1,f_2)$ and $\CoEq(f_1,f_2)$ are quasi-measurable spaces and categorical equalizers and coequalizers, resp., of $\QMS$:

\begin{Lem}
    \label{qms:lem:qms-equalizer-up}
    $(\Eq(f_1,f_2),\Eq(f_1,f_2)^\Omega)$ is quasi-measurable space and satisfies the universal property of an equalizer in the category of quasi-measurable spaces $\QMS$.
    \begin{proof}
        That $(\Eq(f_1,f_2),\Eq(f_1,f_2)^\Omega)$ is quasi-measurable space can be checked
        directly.
        Now consider a quasi-measurable map:
        \begin{align}g:\;(\Zcal,\Zcal^\Omega) \to (\Xcal,\Xcal^\Omega),\end{align}
        such that $f_1\circ g = f_2 \circ g$.
        Then for all $z \in \Zcal$ we get that:
        \begin{align}f_1\circ g(z) = f_2 \circ g(z),\end{align}
        which implies: $g(z) \in \Eq(f_1,f_2)$. By the same argument we get:
        \begin{align}g \circ \Zcal^\Omega \ins \iota^*\Xcal^\Omega = \Eq(f_1,f_2)^\Omega.\end{align}
        So $g$ uniquely factorizes through the inclusion:
        \begin{align}\iota:\; (\Eq(f_1,f_2),\Eq(f_1,f_2)^\Omega) \inj (\Xcal,\Xcal^\Omega).\end{align}
        This shows that $(\Eq(f_1,f_2),\Eq(f_1,f_2)^\Omega)$ hat the universal property of an equalizer in $\QMS$.
    \end{proof}
\end{Lem}

\begin{Lem}
    \label{qms:lem:qms-coequalizer-up}
    The coequalizer $(\CoEq(f_1,f_2),\CoEq(f_1,f_2)^\Omega)$ is a quasi-measurable space and has the universal property of a coequalizer in the category of quasi-measurable spaces $\QMS$.
    \begin{proof}
        That $(\CoEq(f_1,f_2),\CoEq(f_1,f_2)^\Omega)$ is quasi-measurable space can easily be checked.
        Consider a quasi-measurable map:
        \begin{align}g:\; (\Ycal,\Ycal^\Omega) \to (\Zcal,\Zcal^\Omega),\end{align}
        such that $g \circ f_1 = g \circ f_2$. Then for every $x \in \Xcal$ we get:
        \begin{align}g \circ f_1(x) = g \circ f_2(x).\end{align}
        So there exists a unique (set-theoretic) map:
        \begin{align} h:\; \CoEq(f_1,f_2) \to \Zcal,\qquad h \circ t = g.\end{align}
        We then get:
        \begin{align} h \circ \CoEq(f_1,f_2)^\Omega  = h \circ t \circ \Ycal^\Omega= g  \circ \Ycal^\Omega\ins \Zcal^\Omega.\end{align}
        This shows that $h$ is quasi-measurable and thus that $(\CoEq(f_1,f_2),\CoEq(f_1,f_2)^\Omega)$
        has the universal property of a coequalizer in $\QMS$.
    \end{proof}
\end{Lem}

\subsection{The Fibre Product}
\label{qms:sec:fibre-prod}

In this subsection we define the \emph{fibre product}, a product of quasi-measurable spaces relative to a third one.
We only focus here on the case of products of two quasi-measurable spaces for simplicity, but the construction can be done for an arbitrary family (indexed over a set).

 \begin{Def}[Quasi-measurable spaces over a base]
     Let $\Tcal$ be a quasi-measurable space.
     We call quasi-measurable maps between quasi-measurable spaces, like $T_\Xcal:\,\Xcal \to \Tcal$ and
     $T_\Ycal:\,\Ycal \to \Tcal$, with codomain $\Tcal$ \emph{quasi-measurable spaces over $\Tcal$}.
    We then define the set of \emph{quasi-measurable maps over $\Tcal$} as:
    \begin{align}
        \QMS_\Tcal\lp\Xcal,\Ycal\rp & := \QMS_\Tcal\lp(T_\Xcal:\,\Xcal \to \Tcal),(T_\Ycal:\,\Ycal \to \Tcal) \rp   \\
                                    &:= \lC g \in \QMS(\Xcal,\Ycal) \st T_\Ycal \circ g = T_\Xcal \rC.
    \end{align}
    Note that the definition depends on the choices of $T_\Xcal$ and $T_\Ycal$, which we, by abuse of notation, will often keep implicit for sake of simplicity.
 \end{Def}

\begin{Def}[Fibers]
    Let $T:\,\Xcal \to \Tcal$ be a quasi-measurable map between quasi-measurable spaces.
    For $t \in \Tcal$ the \emph{fibre} of $\Xcal$ over $t$ (along $T$) is the quasi-measurable space $(\Xcal_t,\Xcal_t^\Omega)$ given by:
    \begin{align}
        \Xcal_t &:= T^{-1}(t) = \lC x \in \Xcal \st T(x) = t \rC,&
        \Xcal_t^\Omega &:= \lC X:\, \Omega \to \Xcal_t \st \iota_t \circ X \in \Xcal^\Omega   \rC,
    \end{align}
    where $\iota_t:\, \Xcal_t \to \Xcal$ is the inclusion map.
\end{Def}

It is clear with this definition that the inclusion map $\iota_t:\, \Xcal_t \to \Xcal$ is quasi-measurable.

\begin{Def}[Fibre product]
    Let $\Tcal$ be a quasi-measurable space and $T_\Xcal:\,\Xcal \to \Tcal$, $T_\Ycal:\,\Ycal \to \Tcal$
    quasi-measurable maps between quasi-measurable spaces.
    The \emph{fibre product} of $\Xcal$ and $\Ycal$ over $\Tcal$ (along $T_\Xcal$, $T_\Ycal$) is given by:
    \begin{align}
        \Xcal \ftimes_\Tcal \Ycal &:= \lC (x,y) \in \Xcal \times \Ycal \st T_\Xcal(x)=T_\Ycal(y) \rC,\\
        (\Xcal \ftimes_\Tcal \Ycal)^\Omega &:= \lC (X,Y) :\, \Omega \to \Xcal \ftimes_\Tcal \Ycal \st X \in \Xcal^\Omega,\, Y \in \Ycal^\Omega,\,
        T_\Xcal \circ X = T_\Ycal \circ Y \rC,
    \end{align}
    with quasi-measurable projection map:
    \begin{align}
        T_{\Xcal\ftimes_\Tcal\Ycal}:\, \Xcal \ftimes_\Tcal \Ycal  &\to \Tcal, & T_{\Xcal\ftimes_\Tcal\Ycal}(x,y) &:= T_\Xcal(x) = T_\Ycal(y).
    \end{align}
\end{Def}

\begin{Rem}
Note that for $t \in \Tcal$ the fibre of the fibre product is given by:
    \begin{align}
        (\Xcal \ftimes_\Tcal \Ycal)_t &= \Xcal_t \times \Ycal_t, &
        (\Xcal \ftimes_\Tcal \Ycal)_t^\Omega &= \Xcal_t^\Omega \times \Ycal_t^\Omega.
    \end{align}
 So $(\Xcal \ftimes_\Tcal \Ycal)_t = \Xcal_t \times \Ycal_t$ as quasi-measurable spaces.
\end{Rem}

\begin{Rem}
 It can be shown, as for products, that the fibre product really defines a fibre product in the category of quasi-measurable spaces:
 \begin{align}
     \QMS_\Tcal( \Zcal, \Xcal \ftimes_\Tcal \Ycal ) & \cong \QMS_\Tcal( \Zcal,\Xcal) \times \QMS_\Tcal( \Zcal,\Ycal),
 \end{align}
 for quasi-measurable spaces $\Xcal,\Ycal,\Zcal$ over $\Tcal$.
\end{Rem}

 \begin{Thm}[Fibre products distribute over coproducts]
    \label{qms:thm:fibre-prod-coprod}
    Let $I$ be an (index) set and $T_i:\,(\Ycal_i,\Ycal_i^\Omega) \to (\Tcal,\Tcal^\Omega)$,
    $T_\Xcal:\,(\Xcal,\Xcal^\Omega) \to (\Tcal,\Tcal^\Omega)$ a family of quasi-measurable maps, $i \in I$.
    Then we have a canonical isomorphism of quasi-measurable spaces:
    \begin{align}
        \lp \coprod_{i \in I} \Ycal_i \rp \ftimes_\Tcal \Xcal &\cong \coprod_{i \in I} \lp \Ycal_i \ftimes_\Tcal \Xcal \rp, &
        ((i,y),x) &\mapsto (i,(y,x)).
    \end{align}
    \begin{proof} We have:
       \begin{align}
           \lp \coprod_{i \in I} \Ycal_i \rp \ftimes_\Tcal \Xcal &= \lC ((i,y),x) \st i \in I, y \in \Ycal_i, x \in \Xcal, T_i(y) = T_\Xcal(x) \rC\\
              &\cong \lC (i,(y,x)) \st i \in I, y \in \Ycal_i, x \in \Xcal, T_i(y) = T_\Xcal(x) \rC\\
              &=\coprod_{i \in I} \lp \Ycal_i \ftimes_\Tcal \Xcal \rp,\\
           \lp\lp \coprod_{i \in I} \Ycal_i \rp \ftimes_\Tcal \Xcal \rp^\Omega &= \lC ((i,Y),X) \st i \in I, Y \in \Ycal_i^\Omega, X \in \Xcal^\Omega, T_i(Y) = T_\Xcal(X) \rC\\
           &\cong \lC (i,(Y,X)) \st i \in I, Y \in \Ycal_i^\Omega, X \in \Xcal^\Omega, T_i(Y) = T_\Xcal(X) \rC\\
            &= \lp \coprod_{i \in I} \lp \Ycal_i \ftimes_\Tcal \Xcal \rp\rp^\Omega.
       \end{align}
       This shows the claim.
    \end{proof}
\end{Thm}

\begin{Rem}
Also note that we have that fibre products  over their coproducts are empty:
    \begin{align} \Xcal \ftimes_{\Xcal \sqcup \Ycal} \Ycal = \zero.   \end{align}
\end{Rem}

\subsection{The Internal Hom}
\label{qms:sec:int-homs}

 \begin{Def}[The internal hom]
    \label{qms:def:int-hom}
    Let $\Tcal$ be a quasi-measurable space and $T_\Ycal:\,\Ycal \to \Tcal$, $T_\Zcal:\,\Zcal \to \Tcal$
    quasi-measurable maps between quasi-measurable spaces.
    We define  the \emph{internal hom} of $\Ycal$ and $\Zcal$ over $\Tcal$ as the set-theoretic disjoint union:
    \begin{align}
        \Zcal^\Ycal_\Tcal:= \ihom_\Tcal(\Ycal,\Zcal) &:= \coprod_{t \in \Tcal}^\Sets \QMS\lp \Ycal_t,\Zcal_t\rp,
    \end{align}
    with projection map given for $(t,g) \in \ihom_\Tcal(\Ycal,\Zcal)$ with  $t \in \Tcal$ and $g \in \QMS\lp \Ycal_t,\Zcal_t\rp$ via:
    \begin{align}
        T_{\ihom_\Tcal(\Ycal,\Zcal)}:\, \ihom_\Tcal(\Ycal,\Zcal) &\to \Tcal, &  T_{\ihom_\Tcal(\Ycal,\Zcal)}(t,g) &:= t.
    \end{align}
 The set of admissible random variables is defined as:
    \begin{align}
        \ihom_\Tcal(\Ycal,\Zcal)^\Omega := \bigg\{ &  (T,G):\, \Omega \to \ihom_\Tcal(\Ycal,\Zcal) \,\bigg | \,
                                         T \in \Tcal^\Omega \, \land \,\\
                                        & \forall \Phi \in \Omega^\Omega,\, \forall Y \in \Ycal^\Omega \text{ s.t. }
            T \circ \Phi = T_\Ycal \circ Y.  \\
        & \lp \omega \mapsto G\lp\Phi(\omega)\rp\lp Y(\omega)\rp \rp \in \Zcal^\Omega
    \bigg\}.
    \end{align}
    Note that in the last expression: $G\lp\Phi(\omega)\rp \in \QMS\lp \Ycal_t,\Zcal_t \rp$ and $Y(\omega) \in \Ycal_t$
    with $t:=T(\Phi(\omega))=T_\Ycal(Y(\omega))$.
\end{Def}

One can show that  $\lp \ihom_\Tcal(\Ycal,\Zcal),\ihom_\Tcal(\Ycal,\Zcal)^\Omega \rp$ is again a quasi-measurable space over $\Tcal$:

\begin{Lem}
    \label{qms:lem:int-hom-qms}
 $\lp \ihom_\Tcal(\Ycal,\Zcal),\ihom_\Tcal(\Ycal,\Zcal)^\Omega \rp$ is a quasi-measurable space and the projection map
 $T_{\ihom_\Tcal(\Ycal,\Zcal)}:\, \ihom_\Tcal(\Ycal,\Zcal) \to \Tcal$ is a quasi-measurable map.
 \begin{proof}
     Let $(T,G):\, \Omega \to \ihom_\Tcal(\Ycal,\Zcal)$ be a constant map. Then for all $\omega \in \Omega$ we have:
     $T(\omega)=t \in \Tcal$ and $G(\omega)=g \in \QMS\lp \Ycal_t,\Zcal_t\rp$. Clearly $T \in \Tcal^\Omega$.
     For $\Phi \in \Omega^\Omega$ and $Y \in \Ycal^\Omega$ with $T_\Ycal \circ Y = T \circ \Phi = t$ we get that $Y \in \Ycal_t^\Omega$.
     Together with $g \in \QMS\lp \Ycal_t,\Zcal_t\rp$ we get:
\begin{align} G(\Phi)(Y) = g(Y)  \in \Zcal_t^\Omega \ins \Zcal^\Omega. \end{align}
     This shows that $(T,G) \in \ihom_\Tcal(\Ycal,\Zcal)^\Omega $.\\
     Now let $(T,G) \in \ihom_\Tcal(\Ycal,\Zcal)^\Omega $ and $\Phi_1 \in \Omega^\Omega$. Then with $T$ also $T \circ \Phi_1 \in \Tcal^\Omega$.
     Let $\Phi_2 \in \Omega^\Omega$ and $Y \in \Ycal^\Omega$ with $T_\Ycal \circ Y = T \circ \Phi_1 \circ \Phi_2$. Then $\Phi:= \Phi_1\circ\Phi_2 \in \Omega^\Omega$ and $T_\Ycal \circ Y = T \circ \Phi$. With this we get:
     \begin{align}
         G( \Phi_1 \circ \Phi_2)(Y) = G(\Phi)(Y) \in \Zcal^\Omega.
     \end{align}
     This shows that $(T,G) \circ \Phi  = (T\circ \Phi,G \circ \Phi) \in \ihom_\Tcal(\Ycal,\Zcal)^\Omega$.\\
     Finally, note that $T_{\ihom_\Tcal(\Ycal,\Zcal)}(T,G) = T \in \Tcal^\Omega$, showing that $T_{\ihom_\Tcal(\Ycal,\Zcal)}$ is quasi-measurable.
 \end{proof}
\end{Lem}

\begin{Rem}
    Note that the internal hom has the following fibres at $t \in \Tcal$:
\begin{align}
    \ihom_\Tcal(\Ycal,\Zcal)_t & = \lC t \rC \times \QMS\lp \Ycal_t,\Zcal_t\rp \\
                               &\cong \QMS\lp \Ycal_t,\Zcal_t\rp, \\
    \ihom_\Tcal(\Ycal,\Zcal)_t^\Omega & \cong  \bigg\{ G:\, \Omega \to \QMS\lp \Ycal_t,\Zcal_t\rp  \,\bigg|\,
                                      \forall \Phi \in \Omega^\Omega,\, \forall Y \in \Ycal^\Omega_t. \notag\\
                                      &\qquad\qquad \lp \omega \mapsto G\lp\Phi(\omega)\rp\lp Y(\omega)\rp \rp \in \Zcal^\Omega \bigg\} \\
                                      & = \QMS\lp \Ycal_t,\Zcal_t\rp^\Omega.
\end{align}
This also shows that for the special case $\Tcal = \one$ we recover the function space:
\begin{align}\ihom_\one(\Ycal,\Zcal)=\Ycal^\Zcal. \end{align}
\end{Rem}

\begin{Thm}
    \label{qms:thm:int-hom}
        Let $\Tcal$ be a quasi-measurable space and $T_\Xcal:\,\Xcal \to \Tcal$, $T_\Ycal:\,\Ycal \to \Tcal$, $T_\Zcal:\,\Zcal \to \Tcal$
    quasi-measurable spaces over $\Tcal$.
    Then the curry map:
    \begin{align}
        \QMS_\Tcal\lp \Xcal \ftimes_\Tcal \Ycal, \Zcal \rp & \cong \QMS_\Tcal\lp \Xcal, \ihom_\Tcal\lp\Ycal,\Zcal\rp\rp, \\
        g & \mapsto \lp x\mapsto \lp T_\Xcal(x),g_x:=(y \mapsto g(x,y) \rp\rp,
    \end{align}
    is well-defined and a canonical bijection.
\end{Thm}

\begin{Cor}
    \label{qms:cor:int-hom-curry}\label{qms:cor:int-hom-curry-b}
The canonical curry map:
    \begin{align}
        \ihom_\Tcal\lp \Xcal \ftimes_\Tcal \Ycal,\Zcal \rp & \cong \ihom_\Tcal\lp \Xcal, \ihom_\Tcal\lp \Ycal,\Zcal\rp \rp,\\
        (t,g) & \mapsto \lp t,\lp x \mapsto \lp T_\Xcal(x),g_x=(y \mapsto g(x,y))\rp\rp\rp,
    \end{align}
    is an isomorphisms of quasi-measurable spaces over $\Tcal$.
\end{Cor}

In more categorical terms, see e.g.\ \cite{Joh02} A1.5.3, we can conclude:

\begin{Cor}
    \label{qms:cor:lcc-app}
    The category $\QMS$ of quasi-measurable spaces and maps is \emph{locally cartesian closed} with terminal object $\one$.
\end{Cor}

In other words, this means that $\QMS$ allows for an (extensional) \emph{dependent type theory}, see \cite{See84}.
E.g.\ for $\Wcal \in \QMS_\Xcal$ and $\Xcal \in \QMS_\Tcal$ the \emph{dependent product} is given by:
\begin{align}
    \prod_\Xcal \Wcal &= \ihom_\Tcal(\Xcal,\Wcal)\ftimes_{\ihom_\Tcal(\Xcal,\Xcal)} \Tcal.
\end{align}

\begin{Prf}{qms:thm:int-hom}
    Let $g \in \QMS_\Tcal\lp \Xcal \ftimes_\Tcal \Ycal, \Zcal \rp$,
    $x \in \Xcal$ and $t:=T_\Xcal(x)$.
    For $y \in \Ycal$ we then have the equivalence:
    \begin{align}
        (x,y) \in \Xcal\ftimes_\Tcal\Ycal & \iff t=T_\Xcal(x)=T_\Ycal(y) \iff y \in \Ycal_t.
        \label{qms:eq:fibres-of-fibre-prod}
    \end{align}
    So for $y \in \Ycal_t$ we can evaluate $g(x,y)$.
    Since $g$, by assumption, is a quasi-measurable map over $\Tcal$ we have:
    \begin{align}
        T_\Zcal\lp g(x,y) \rp &=  T_{\Xcal\ftimes_\Tcal\Ycal }(x,y) = T_\Xcal(x) =t,
        \label{eq:qms-map-over-T}
    \end{align}
    which shows that $g(x,y) \in \Zcal_t$.
    Together this shows that the following map is well-defined:
   \begin{align}
       g_x:\, \Ycal_t &\to \Zcal_t, & y \mapsto g_x(y):=g(x,y).
    \end{align}
    To show that $g_x$ is quasi-measurable let $Y \in \Ycal^\Omega_t$.
    Then $T_\Ycal \circ Y(\omega) = t$ for all $\omega \in \Omega$. If we take the constant
    map $X:=x$ then $X \times Y \in \lp\Xcal\ftimes_\Tcal\Ycal \rp^\Omega$.
    Since $g$ is quasi-measurable we get that $g(X,Y) \in \Zcal^\Omega$.
    Again by \Cref{eq:qms-map-over-T} we see that we even have:
    $g_x(Y) = g(X,Y) \in \Zcal^\Omega_t$. This shows that $g_x \in \QMS(\Ycal_t,\Zcal_t)$.\\
    These considerations show that we have a well-defined map:
    \begin{align}
        \tilde g:\, \Xcal &\to \ihom_\Tcal\lp\Ycal,\Zcal\rp =\coprod_{t \in \Tcal} \QMS(\Ycal_t,\Zcal_t),& x &\mapsto (T_\Xcal(x),g_x),
    \end{align}
    where $g_x\,:\Ycal_{T_\Xcal(x)} \to \Zcal_{T_\Xcal(x)}$.
    The map $\tilde g$ is a map over $\Tcal$ because:
    \begin{align}
      T_{\ihom_\Tcal\lp\Ycal,\Zcal\rp}\lp \tilde g(x)\rp =  T_{\ihom_\Tcal\lp\Ycal,\Zcal\rp}\lp(T_\Xcal(x),g_x)\rp & = T_\Xcal(x).
    \end{align}
    To show that $\tilde g$ is quasi-measurable let $X \in \Xcal^\Omega$, $\Phi \in \Omega^\Omega$, $Y \in \Ycal^\Omega$ such that $T_\Xcal \circ X \circ \Phi = T_\Ycal \circ Y$.
    Then $(X \circ \Phi,Y) \in \lp\Xcal\ftimes_\Tcal\Ycal \rp^\Omega$ and thus:
    \begin{align}
        g_{X \circ \Phi}(Y) = g(X \circ \Phi,Y) \in \Zcal^\Omega,
    \end{align}
    since $g$ was assumed to be quasi-measurable. Also  $T_\Xcal(X) \in \Tcal^\Omega$,
    because $T_\Xcal$ is quasi-measurable. By \Cref{qms:def:int-hom} for $\ihom_\Tcal\lp\Ycal,\Zcal\rp^\Omega$ we have that:
    \begin{align} \tilde g(X) = (T_\Xcal(X),g_X) \in \ihom_\Tcal\lp\Ycal,\Zcal\rp^\Omega.
    \end{align}
    This shows that $\tilde g \in \QMS_\Tcal\lp \Xcal, \ihom_\Tcal\lp\Ycal,\Zcal\rp\rp$.\\
    To construct the inverse to the curry map, let
    $\tilde h \in \QMS_\Tcal\lp \Xcal,\ihom_\Tcal\lp\Ycal,\Zcal\rp\rp$.
    Since $\tilde h$ is a quasi-measurable map over $\Tcal$ we have that for all $x \in \Xcal$:
    \begin{align}
        T_{\ihom_\Tcal\lp\Ycal,\Zcal\rp} \circ \tilde h(x) &= T_\Xcal(x).
    \end{align}
    So $\tilde h(x)$ is of the form $\tilde h(x) = (T_\Xcal(x),h_x)$ with an
    $h_x \in \QMS\lp\Ycal_{T_\Xcal(x)}, \Zcal_{T_\Xcal(x)}\rp$.
    For $(x,y) \in \Xcal \ftimes_\Tcal \Ycal$ we by \Cref{qms:eq:fibres-of-fibre-prod}
    have that $y \in \Ycal_{T_\Xcal(x)}$ and we can evaluate
    $h_x(y) \in \Zcal_{T_\Xcal(x)} \ins \Zcal$.
    So we get a well-defined map:
    \begin{align}
        h:\,  \Xcal \ftimes_\Tcal \Ycal &\to \Zcal, & (x,y) \mapsto h_x(y).
    \end{align}
    Since $h_x(y) \in \Zcal_{T_\Xcal(x)}$ we have that:
    \begin{align}
        T_\Zcal(h(x,y)) =  T_\Zcal(h_x(y)) &= T_\Xcal(x) =  T_{\Xcal \ftimes_\Tcal \Ycal}(x,y).
    \end{align}
    So $h$ is a map over $\Tcal$.
    To show that $h$ is quasi-measurable let $X \in \Xcal^\Omega$ and $Y \in \Ycal^\Omega$
    such that $T_\Xcal \circ X = T_\Ycal \circ Y$.
    Then $(T_\Xcal(X),h_X)=\tilde h (X) \in \ihom_\Tcal\lp\Ycal,\Zcal\rp^\Omega$.
    By \Cref{qms:def:int-hom} for $\ihom_\Tcal\lp\Ycal,\Zcal\rp^\Omega$ with $\Phi:=\id_\Omega$
    we get that:
    \begin{align} h(X,Y)=h_X(Y) \in \Zcal^\Omega. \end{align}
    This shows that $h \in \QMS_\Tcal\lp \Xcal \ftimes_\Tcal \Ycal, \Zcal \rp$. \\
    It is obvious by these constructions that the maps
    $g \mapsto \tilde g$ and $\tilde h \mapsto h$ are inverse to each other.
    By these explicit constructions it can easily be shown that this bijection
    is natural in its arguments.
    This then shows the claim.
\end{Prf}

\begin{Prf}[More abstract proof]{qms:cor:int-hom-curry}
    By repeated application of \Cref{qms:thm:int-hom} we get the canonical isomorphisms for every quasi-measurable space
    $\Wcal$ over $\Tcal$:
\begin{align}
    \QMS_\Tcal\lp \Wcal,\ihom_\Tcal\lp \Xcal \ftimes_\Tcal \Ycal,\Zcal \rp \rp
    &\cong \QMS_\Tcal\lp \Wcal \ftimes_\Tcal \lp \Xcal \ftimes_\Tcal \Ycal \rp,\Zcal \rp \\
    &\cong \QMS_\Tcal\lp \lp \Wcal \ftimes_\Tcal \Xcal \rp \ftimes_\Tcal\Ycal ,\Zcal \rp \\
    &\cong \QMS_\Tcal\lp  \Wcal \ftimes_\Tcal \Xcal , \ihom_\Tcal\lp\Ycal ,\Zcal\rp \rp \\
    &\cong \QMS_\Tcal\lp  \Wcal, \ihom_\Tcal\lp \Xcal , \ihom_\Tcal\lp\Ycal ,\Zcal\rp\rp \rp.
\end{align}
Since this holds for all $\Wcal$ we get a natural isomorphisms of quasi-measurable spaces over $\Tcal$:
    \begin{align}
        \ihom_\Tcal\lp \Xcal \ftimes_\Tcal \Ycal,\Zcal \rp & \cong \ihom_\Tcal\lp \Xcal, \ihom_\Tcal\lp \Ycal,\Zcal\rp \rp.
    \end{align}
This shows the claims.
\end{Prf}

\begin{Prf}[More explicit proof]{qms:cor:int-hom-curry-b}
        For the underlying sets we get canonical identifications:
    \begin{align}
        \ihom_\Tcal\lp \Xcal \ftimes_\Tcal \Ycal,\Zcal \rp & = \coprod_{t \in \Tcal}^\Sets \QMS\lp (\Xcal \ftimes_\Tcal \Ycal)_t,\Zcal_t\rp \\
        &= \coprod_{t \in \Tcal}^\Sets \QMS\lp \Xcal_t \times \Ycal_t,\Zcal_t\rp \\
        &\cong \coprod_{t \in \Tcal}^\Sets \QMS\lp \Xcal_t, \QMS\lp\Ycal_t,\Zcal_t\rp\rp \\
        &\cong \coprod_{t \in \Tcal}^\Sets \QMS\lp \Xcal_t, \lC t \rC \times \QMS\lp\Ycal_t,\Zcal_t\rp\rp \\
        &= \coprod_{t \in \Tcal}^\Sets \QMS\lp \Xcal_t, \ihom_\Tcal\lp\Ycal,\Zcal\rp_t\rp \\
        &=\ihom_\Tcal\lp \Xcal, \ihom_\Tcal\lp \Ycal,\Zcal\rp\rp.
    \end{align}
    Note that the identifications inside $\QMS$ are really identifications as quasi-measurable spaces, not just as sets.
    Also note that we can then canonically identify the projection maps:
    \begin{align} T_{\ihom_\Tcal\lp \Xcal \ftimes_\Tcal \Ycal,\Zcal \rp}= T_{\ihom_\Tcal\lp \Xcal, \ihom_\Tcal\lp \Ycal,\Zcal\rp\rp} =:T_{\ihom}. \end{align}
    Now let $(T,G) \in \ihom_\Tcal\lp \Xcal \ftimes_\Tcal \Ycal,\Zcal \rp^\Omega$.
    Then $T \in \Tcal^\Omega$ and for all $\Phi \in \Omega^\Omega$, $X \in \Xcal^\Omega$, $Y \in \Ycal^\Omega$  with
        $T_\Xcal \circ X = T_\Ycal \circ Y  = T \circ \Phi$ we have:
        \begin{align} \tilde G(\Phi)(X)(Y)= G(\Phi)(X,Y) \in \Zcal^\Omega, \label{qms:eq:some-curry-1}  \end{align}
    where $\tilde G$ is the curried version of $G$.\\
    Let $\Phi_1 \in \Omega^\Omega$ and $X_1 \in \Xcal^\Omega$ such that $T_\Xcal \circ X_1 = T \circ \Phi_1$.
    We need to show that:
    \begin{align}  \lp T(\Phi_1),\tilde G(\Phi_1)(X_1) \rp \in \ihom_\Tcal\lp \Ycal,\Zcal\rp^\Omega.\label{qms:eq:some-curry-2}  \end{align}
    For this we need to show that for $\Phi_2 \in \Omega^\Omega$ and $Y \in \Ycal^\Omega$ with $T_\Ycal \circ Y = T \circ \Phi_1 \circ \Phi_2$
    we get:
    \begin{align} \tilde G(\Phi_1 \circ \Phi_2)(X_1 \circ \Phi_2)(Y)= G(\Phi_1 \circ \Phi_2)(X_1 \circ \Phi_2,Y)  \in \Zcal^\Omega.\label{qms:eq:some-curry-3}  \end{align}
    Since $T_\Xcal \circ X_1 = T \circ \Phi_1$ implies $T_\Xcal \circ X_1 \circ \Phi_2 = T \circ \Phi_1 \circ \Phi_2$ we get the claim
    with $X := X_1 \circ \Phi_2 \in \Xcal^\Omega$ and $\Phi := \Phi_1 \circ \Phi_2 \in \Omega^\Omega$ from \Cref{qms:eq:some-curry-1}.
    This shows that the curry map is quasi-measurable.\\
    For the reverse let $(T,\tilde G) \in \ihom_\Tcal\lp \Xcal, \ihom_\Tcal\lp \Ycal,\Zcal\rp\rp^\Omega$.
    Then \Cref{qms:eq:some-curry-2} and \Cref{qms:eq:some-curry-3} with $\Phi_2:=\id_\Omega$, $\Phi_1:=\Phi$ and $X_1:=X$ implies \Cref{qms:eq:some-curry-1}.
    This shows that the uncurry map, the inverse of the curry map, is also quasi-measurable.
\end{Prf}

\subsection{Quotients and Subspaces}
\label{qms:sec:quotients}
\label{qms:sec:quot-sub}

In this subsection we relate the categorical notions of (strong) monomorphisms and (strong) epimorphisms
to our ad hoc definitions of embeddings and quotients of quasi-measurable spaces, see \Cref{qms:def:quot-emb}.

\begin{Lem}[Injective maps are monomorphisms, monomorphisms are injective]
    A quasi-measurable map $i:\, (\Xcal,\Xcal^\Omega) \to (\Ycal,\Ycal^\Omega)$ is a \emph{monomorphism} in $\QMS$ if and only if $i$ is injective (as a map).
    \begin{proof}
        Assume that $i$ is a monomorphism and $i(x_1)=i(x_2)$ for $x_i \in \Xcal$, $i=1,2$.
        Then consider the constant maps: $g_i:=x_i$, $i=1,2$.
        Clearly: $g_i \in \QMS(\Zcal,\Xcal)$. Then we have $i \circ g_1 = i \circ g_2$, which by the assumption that $i$ is a
        monomorphisms implies: $g_1=g_2$, and thus
        $x_1=x_2$.\\
        Let $i$ now be injective and $g_1,g_2 \in \QMS(\Zcal,\Xcal)$ such that $i \circ g_1 = i \circ g_2$.
        For every $z \in \Zcal$ we then get: $i(g_1(z))=i(g_2(z))$. Since $i$ is injective we get that:
        $g_1(z)=g_2(z)$ for every $z \in \Zcal$. So $g_1=g_2$ as maps and thus as quasi-measurable maps.
    \end{proof}
\end{Lem}

\begin{Lem}[Surjective maps are epimorphisms, epimorphisms are surjective]
        A quasi-measurable map
        $t:\, (\Xcal,\Xcal^\Omega) \to (\Ycal,\Ycal^\Omega)$ is an \emph{epimorphism} in $\QMS$ if and only if it is surjective.
     \begin{proof}
         Assume that $t$ is surjective and let $h_i:\; (\Ycal,\Ycal^\Omega)\to (\Zcal,\Zcal^\Omega)$, $i=1,2$,
         be quasi-measurable maps with $h_1 \circ t = h_2 \circ t$.
         Since $t$ is surjective, for every $y \in \Ycal$ there is an $x \in \Xcal$ with $t(x)=y$.
         Then we get:
         \begin{align} h_1(y) = h_1 \circ t(x) = h_2 \circ t(x) = h_2(y).  \end{align}
         Since this holds for every $y \in \Ycal$ we get $h_1=h_2$ as maps and thus as quasi-measurable maps.\\
         Now let $t$ be an epimorphism.
         The indicator function $\I_{t(\Xcal)}:\, \Ycal \to \twop$ is quasi-measurable since $\twop^\Omega=\Sets(\Omega,\twop)$.
         We also have that the constant map $\I_\Ycal \in \twop^\Ycal$, which maps everything to $1 \in \twop$.
         Then for all $x \in \Xcal$ we get:
         \begin{align} \I_{t(\Xcal)} \circ t(x) = 1 = \I_\Ycal  \circ t(x). \end{align}
         This implies the equality of quasi-measurable maps:
         \begin{align} \I_{t(\Xcal)} \circ t = \I_\Ycal  \circ t. \end{align}
         Since $t$ is an epimorphism we get the equality:
         \begin{align} \I_{t(\Xcal)} = \I_\Ycal, \end{align}
         which when evaluated at $y \in \Ycal \sm t(\Xcal)$ would give the contradiction $0=1$.
         Thus $\Ycal \sm t(\Xcal) = \emptyset$ and $t$ must be surjective.
     \end{proof}
\end{Lem}

\begin{Lem}[Quotients are strong epimorphisms and vice versa]
    \label{qms:lem:strong-epi-quotient}
 Let $t:\, (\Zcal,\Zcal^\Omega) \to (\Wcal,\Wcal^\Omega)$ be a quasi-measurable map.
 Then $t$ is a strong epimorphisms if and only if it is a quotient map.
 \begin{proof}
     Consider any commutative diagram of quasi-measurable spaces like depicted on the left:
    \begin{align}\xymatrix{
            \Zcal \ar^-{f}[rr] \ar_-{t}[d] && \Xcal \ar^-{i}[d]  &&& \Zcal \ar^-{f}[rr] \ar_-{t}[d] && \Xcal \ar^-{i}[d]\\
            \Wcal \ar_-{g}[rr]  && \Ycal &&& \Wcal \ar_-{g}[rr] \ar@{-->}^-{h}[urr] && \Ycal.
    }\end{align}
    $t$ is a strong epimorphism iff whenever $i$ is a monomorphism there exists a unique quasi-measurable map $h:\,\Wcal \to \Xcal$
    that makes the triangles on the right commutative. \\
    Let $t$ be a strong epimorphism. Put $\Xcal:=t(\Zcal)$, $\Xcal^\Omega:=t \circ \Zcal^\Omega$ and $f$ the quotient map
    $t:\,\Zcal \to \Xcal$. We also put $\Ycal:=\Wcal$, $g:=\id_\Wcal$ and $i:\, \Xcal \to \Wcal$ the inclusion map.
    Then the above diagram is commutative and all maps are quasi-measurable. Since $t$ is a strong epimorphisms we get a unique
    $h:\,\Wcal \to \Xcal$ that makes the right diagramm commutative. So $i \circ h = \id_\Wcal$. This shows that $i$ and thus also $t$ is surjective.
    For $W \in \Wcal^\Omega$ we get $i \circ h(W) = W$ with $h(W) \in \Xcal^\Omega=t \circ \Zcal^\Omega$.
    So there exists $Z \in \Zcal^\Omega$ such that $t(Z) = \id_\Wcal \circ t(Z)= i \circ h(W) = W$. This shows that $t$ is a quotient map.\\
    Now let us assume that $t$ is a quotient map and we have a diagram from above with $i$ a monomorphism.
    Then $t$ is surjective and $i$ is injective.
    The only map that makes the triangles commutative is given by $h(w) := f(z)$ for any $z \in t^{-1}(w)$.
    This is well-defined since $t^{-1}(w)\neq \emptyset$ and for $z_1,z_2 \in t^{-1}(w)$ we have:
    \begin{align} i(f(z_1)) =g(t(z_1))=g(w) = g(t(z_2))= i(f(z_2)).  \end{align}
    Since $i$ is injective we have $f(z_1)=f(z_2)$. So $h$ is well-defined, unique and clearly makes the triangles commutative.
    We are left to show that $h$ is quasi-measurable. Let $W \in \Wcal^\Omega = t \circ \Zcal^\Omega$. Then there is $Z \in \Zcal^\Omega$
    with $W=t(Z)$. By definition of $h$ we get $h(W(\omega))=f(Z(\omega))$ for all $\omega \in \Omega$. This
    implies that $ h(W)=f(Z) \in \Xcal^\Omega$. This shows the claim.
 \end{proof}
\end{Lem}

\begin{Lem}[Embeddings are strong monomorphisms and vice versa]
    \label{qms:lem:strong-mono-embedding}
    Let $i:\,\Xcal \to \Ycal$ be a quasi-measurable map. Then $i$ is a strong
    monomorphism if and only if it is an embedding.
 \begin{proof}
     Consider any commutative diagram of quasi-measurable spaces like depicted on the left:
    \begin{align}\xymatrix{
            \Zcal \ar^-{f}[rr] \ar_-{t}[d] && \Xcal \ar^-{i}[d]  &&& \Zcal \ar^-{f}[rr] \ar_-{t}[d] && \Xcal \ar^-{i}[d]\\
            \Wcal \ar_-{g}[rr]  && \Ycal &&& \Wcal \ar_-{g}[rr] \ar@{-->}^-{h}[urr] && \Ycal.
    }\end{align}
    $i$ is a strong monomorphism iff whenever $t$ is an epimorphism there exists a unique quasi-measurable map $h:\,\Wcal \to \Xcal$
    that makes the triangles on the right commutative.\\
    Let $i$ be a strong monomorphism. Let $\Wcal:=\Xcal$ and $\Wcal^\Omega:=i^*\Ycal^\Omega$ and $g$ the embedding $i:\,\Wcal \to \Ycal$.
    Let $\Zcal := \Xcal$ and $f:=\id_\Xcal$ and $t:\,\Zcal \to \Wcal$ the inclusion map, which is surjective here,
    since $\Zcal=\Xcal=\Wcal$ as sets. So $t$ is an epimorphism. Since $i$ is a strong monomorphism there exists a unique quasi-measurable
    map $h:\, \Wcal \to \Xcal$, the identity on $\Xcal$, that makes the triangles commutative.
    So if $W \in \Wcal^\Omega = i^*\Ycal^\Omega$ then $i(W) \in \Ycal^\Omega$ and $W=h(W) \in \Xcal^\Omega$.
    This shows that $i^*\Ycal^\Omega =\Xcal^\Omega$. Since $i$ is also injective $i$ is an embedding.\\
    Now assume that $i$ is an embedding and we have a diagram on the left with $t$ an epimorphism.
    Then $t$ is surjective, $i$ is injective, a monomorphism
    and $\Xcal^\Omega = i^*\Ycal^\Omega$.  So there exists a unique map $h: \Wcal \to \Xcal$ that makes the triangles commute.
    It is given by $h(w):=f(z)$ with any $z \in t^{-1}(w)$. We are left to show that $h$ is quasi-measurable.
    Let $W \in \Wcal^\Omega$. Then $i(h(W)) = g(W) \in \Ycal^\Omega$, which implies $h(W) \in i^*\Ycal^\Omega = \Xcal^\Omega$.
    This shows the claim.
\end{proof}
\end{Lem}

\begin{Lem}
    \label{qms:lem:subobject-classifier}
    Let $(\Xcal,\Xcal^\Omega)$ be a quasi-measurable space and $(\twop,\twop^\Omega)$ the indiscrete two-point space
    with $\twop^\Omega = \Sets(\Omega,\twop)$. Let $i:\,\Wcal \to \Xcal$ be a quasi-measurable map, $!:\,\Wcal \to \one$
    the constant map and $1:\, \one \to \twop$ the constant $1$-map.
    Then the indicator function $\I_{i(\Wcal)}:\,\Xcal \to \twop$
    is quasi-measurable and the fibre-product $\Xcal \ftimes_\twop \one$ (along $\I_{i(\Wcal)}$ and $1$) can be identified with the embedded quasi-measurable space
    $(\Zcal,\Zcal^\Omega)$ with $\Zcal:=i(\Wcal)$ and $\Zcal^\Omega = \iota^*\Xcal^\Omega$,
    where $\iota:\,\Zcal \to \Xcal$ is the inclusion map.
    \begin{proof}
    We have the commutative diagram of quasi-measurable maps:
    \begin{align}\xymatrix{
            \Wcal \ar@/^1.5pc/^-{!}[rrrrd] \ar@/_1.5pc/_-{i}[rrdd] \ar@{-->}^{i \times !}[rrd]\\
           && \Xcal \ftimes_\twop \one \ar^-{\pr_\one}[rr] \ar_-{\pr_\Xcal}[d] && \one \ar^-{1}[d]  \\
           && \Xcal \ar_-{\I_{i(\Wcal)}}[rr]  && \twop .
    }\end{align}
    We have:
    \begin{align}
        \Xcal \ftimes_\twop \one &= \lC (x,0) \in \Xcal \times \one \st \I_{i(\Wcal)}(x) =1 \rC\\
        &= i(\Wcal) \times \one \\
        &= \Zcal \times \one,\\
        ( \Xcal \ftimes_\twop \one )^\Omega &=\lC (X,!):\, \Omega \to \Xcal \ftimes_\twop \one \st \I_{i(\Wcal)}(X)=1 \rC\\
        &= \lC X \in \Xcal^\Omega \st X(\Omega) \ins i(\Wcal) \rC \times \lC !\rC \\
        &= \iota^*\Xcal^\Omega \times \one^\Omega\\
        &= \Zcal^\Omega \times \one^\Omega.
    \end{align}
    So as quasi-measurable spaces we have: $\Xcal \ftimes_\twop \one = \Zcal \times \one \cong \Zcal$.
    \end{proof}
\end{Lem}

\begin{Thm}
    \label{qms:thm:strong-mono-classifier}
    The indiscrete two-point quasi-measurable space $(\twop,\twop^\Omega)$ given by $\twop:=\lC 0,1\rC$ and $\twop^\Omega:=\Sets(\Omega,\twop)$ together with the constant $1$-map $1:\, \one \to \twop$ is a subobject classifier for strong monomorphisms (=embeddings) in the category $\QMS$ of quasi-measurable spaces and maps.
    \begin{proof}
        This directly follows from \Cref{qms:lem:subobject-classifier} and \Cref{qms:lem:strong-mono-embedding}.
    \end{proof}
\end{Thm}

\subsection{The Heyting Category}
\label{qms:sec:heyting}

In this subsection we will give explicit descriptions of the set $\Sub(\Xcal)$ of subspaces of a quasi-measurable space $\Xcal$.
We will see that $\Sub(\Xcal)$ is partially ordered by inclusion and allows for arbitrary unions and intersections, making it a complete bounded lattice. Furthermore, we can define a Heyting implication operator $\himp$ on $\Sub(\Xcal)$, turning it into a Heyting algebra, an algebraic structure
 to formalize \emph{intuitionistic logic}, see \cite{Hey30}. We then show that the construction is functorial in $\Xcal$. If $f:\,\Xcal \to \Ycal$ is a quasi-measurable map then we get an pullback map $f^*:\,\Sub(\Ycal) \to \Sub(\Xcal)$, which turns out to be a Heyting algebra homomorphism.
We also construct a left-adjoint $f_!:\, \Sub(\Xcal) \to \Sub(\Ycal)$ and a right-adjoint $f_*:\,\Sub(\Xcal) \to \Sub(\Ycal)$ to $f^*$.
$f_!$ allows for an interpretion as an \emph{existential quantifier} $\exists_f$, while $f_*$ as a \emph{universal quantifier} $\forall_f$.
With these constructions $\QMS$ inherits the structure of a \emph{Heyting category}, also called \emph{logos} in \cite{Fre90}.
So $\QMS$ has an internal logic of a (typed) \emph{intuitionistic first-order logic}.

\begin{Def}[Heyting algebra of subspaces]
    \label{qms:def:heyting-algebra-sub}
    Let $(\Xcal,\Xcal^\Omega)$ be a quasi-measurable space. We define the set of subspaces\footnote{Subobjects of $(\Xcal,\Xcal^\Omega)$ are usually defined to be equivalent classes of monomorphisms $\iota:\, (\Wcal,\Wcal^\Omega) \inj (\Xcal,\Xcal^\Omega)$. Since here the equivalence class can be uniquely represented by the subset $\Xcal_1 :=\iota(\Wcal) \ins \Xcal$ and $\Xcal_1^\Omega:=\iota \circ \Wcal^\Omega \ins \Xcal^\Omega$ we directly describe subspaces of $\Xcal$ via subsets $\Xcal_1$ to simplify the following analysis.}  of $(\Xcal,\Xcal^\Omega)$ as:
    \begin{align}
      \Sub(\Xcal):=  \Sub(\Xcal,\Xcal^\Omega):= \lC (\Xcal_1,\Xcal_1^\Omega) \in \QMS \st \Xcal_1 \ins \Xcal, \Xcal_1^\Omega \ins \Xcal^\Omega \rC.
    \end{align}
    We endow $\Sub(\Xcal)$ with the following \emph{partial order}:
    \begin{align}
    (\Xcal_1,\Xcal_1^\Omega) \le (\Xcal_2,\Xcal_2^\Omega) & \qquad :\iff\qquad \Xcal_1 \ins \Xcal_2 \quad \land \quad \Xcal_1^\Omega \ins \Xcal_2^\Omega.
    \end{align}
    We also define a \emph{join} and \emph{meet} operation for $(\Xcal_i,\Xcal_i^\Omega) \in \Sub(\Xcal)$, $i \in I$:
    \begin{align}
        \bigor_{ i \in I} (\Xcal_i,\Xcal_i^\Omega) &:= \lp \bigcup_{i \in I} \Xcal_i,\bigcup_{i \in I} \Xcal_i^\Omega \rp,&
        \bigand_{ i \in I} (\Xcal_i,\Xcal_i^\Omega) &:= \lp \bigcap_{i \in I} \Xcal_i,\bigcap_{i \in I} \Xcal_i^\Omega \rp.
    \end{align}
    We further define a \emph{Heyting implication} operation:
    \begin{align}
        \Xcal_1 \himp \Xcal_2  & := \bigor_{\substack{\tilde \Xcal \in \Sub(\Xcal)\\ \tilde \Xcal \land \Xcal_1 \le \Xcal_2}} \tilde \Xcal.
    \end{align}
    We also abbreviate the \emph{pseudo-complement} of a subspace $(\Xcal_1,\Xcal_1^\Omega)$ as:
    \begin{align}
        \neg \Xcal_1 &:= (\Xcal_1 \himp \zero) =  \Xcal \sm \Xcal_1, \\
        (\neg \Xcal_1)^\Omega &:= (\Xcal_1 \himp \zero)^\Omega = \lC X \in \Xcal^\Omega \st X(\Omega) \ins \Xcal \sm \Xcal_1 \rC.
    \end{align}
\end{Def}

\begin{Def}
    Let $f:\, (\Xcal,\Xcal^\Omega) \to (\Ycal,\Ycal^\Omega)$ be a quasi-measurable map.
    We define the following maps:
    \begin{align}
        f_!:\, \Sub(\Xcal) &\to \Sub(\Ycal), & f^*:\, \Sub(\Ycal) &\to \Sub(\Xcal),
    \end{align}
    on $\Xcal_1 \in \Sub(\Xcal)$ and $\Ycal_1 \in \Sub(\Ycal)$ given via:
    \begin{align}
        f_!(\Xcal_1) &:= f(\Xcal_1),&
        f^*(\Ycal_1) &:= f^{-1}(\Ycal_1),\\
        f_!(\Xcal_1)^\Omega &:= \lC f(X) \st X \in \Xcal_1^\Omega \rC, &
        f^*(\Ycal_1)^\Omega &:= \lC X \in \Xcal^\Omega \st f(X) \in \Ycal_1^\Omega \rC.
    \end{align}
\end{Def}

It is clear that $f_!(\Xcal_1) \in \Sub(\Ycal)$ and $f^*(\Ycal_1) \in \Sub(\Xcal)$ are well-defined quasi-measurable subspaces and thus $f_!$ and $f^*$ are well-defined maps.

\begin{Lem}
        Let $f:\, (\Xcal,\Xcal^\Omega) \to (\Ycal,\Ycal^\Omega)$ be a quasi-measurable map.
    Then for $\Xcal_1 \in \Sub(\Xcal)$ and $\Ycal_1 \in \Sub(\Ycal)$ we have the following adjunction:
    \begin{align}
        f_!(\Xcal_1) \le \Ycal_1 &\iff \Xcal_1 \le f^*(\Ycal_1).
    \end{align}
    In particular, we have:
    \begin{align}
        f_!f^*(\Ycal_1) &\le \Ycal_1,&  \Xcal_1 &\le f^*f_!(\Xcal_1).
    \end{align}
    \begin{proof}
    For the underlying sets we clearly have:
        \begin{align}
            f(\Xcal_1) \ins \Ycal_1 & \iff \Xcal_1 \ins f^{-1}(\Ycal_1),\\
            f_!(\Xcal_1)^\Omega= F(\Xcal_1^\Omega) \ins \Ycal_1^\Omega & \iff \Xcal_1^\Omega \ins F^{-1}(\Ycal_1^\Omega)= f^*(\Ycal_1)^\Omega,
        \end{align}
        where: $F:\, \Xcal^\Omega \to \Ycal^\Omega$, $F(X):=f \circ X$.
    This shows the claim.
    \end{proof}
\end{Lem}

\begin{Lem}
    Let $(\Xcal,\Xcal^\Omega)$ be a quasi-measurable space. Then the set of subspaces $\Sub(\Xcal)$ is a complete Heyting algebra with the operations from
    \Cref{qms:def:heyting-algebra-sub}.
    \begin{proof}
        It is clear that $\Sub(\Xcal)$ is a complete and bounded lattice with top element $\Xcal$ and bottom element $\zero$.
        It is left to show that for $\Xcal_1$, $\Xcal_2$, $\Xcal_3 \in \Sub(\Xcal)$ we have the adjunction:
        \begin{align}
            (\Xcal_3 \land \Xcal_1) \le \Xcal_2 & \qquad \iff \qquad \Xcal_3 \le (\Xcal_1 \himp \Xcal_2).
        \end{align}
    ``$\implies$'': If $(\Xcal_3 \land \Xcal_1) \le \Xcal_2$ then clearly:
    \begin{align}
     \Xcal_3 \le \bigor_{\substack{\tilde \Xcal \in \Sub(\Xcal)\\ \tilde \Xcal \land \Xcal_1 \le \Xcal_2}} \tilde \Xcal = (\Xcal_1 \himp \Xcal_2).
    \end{align}
    ``$\Longleftarrow$'': If $\Xcal_3 \le (\Xcal_1 \himp \Xcal_2)$ then we get:
    \begin{align}
        \Xcal_3 \land \Xcal_1 & \le (\Xcal_1 \himp \Xcal_2) \land \Xcal_1  =  \bigor_{\substack{\tilde \Xcal \in \Sub(\Xcal)\\ \tilde \Xcal \land \Xcal_1 \le \Xcal_2}} \tilde \Xcal \land \Xcal_1
        \le \Xcal_2.
    \end{align}
    This shows the claim.
    \end{proof}
\end{Lem}

\begin{Rem}
    \begin{enumerate}
        \item For the underlying sets of $\Xcal_1 \himp \Xcal_2$ we have:
    \begin{align}
        \Xcal_1 \himp \Xcal_2 &= (\Xcal\sm\Xcal_1) \cup \Xcal_2,& \Xcal_2^\Omega & \ins (\Xcal_1 \himp \Xcal_2)^\Omega \ins (\Xcal^\Omega \sm \Xcal_1^\Omega)  \cup \Xcal_2^\Omega.
    \end{align}
    Indeed, if $\Xcal_3:=(\Xcal\sm\Xcal_1) \cup \Xcal_2$ and $\Xcal_3^\Omega := \Xcal_3^\one \cup \Xcal_2^\Omega$ then $(\Xcal_3,\Xcal_3^\Omega)$
            is a quasi-measurable subspace of $\Xcal$ and $\Xcal_3 \land \Xcal_1 \le \Xcal_2$. Also $\Xcal \sm \Xcal_3 = \Xcal_1 \sm \Xcal_2$, so there is no point left that we could add to $\Xcal_3$, but possibly more random variables outside of $\Xcal_1^\Omega \sm \Xcal_2^\Omega$.
            So we could write:
    \begin{align}
        \Xcal_1 \himp \Xcal_2 &= \bigor_{\substack{\tilde \Xcal \in \Sub(\Xcal)\\ \tilde \Xcal = (\Xcal\sm\Xcal_1) \cup \Xcal_2\\ \Xcal_2^\Omega \ins \tilde \Xcal^\Omega \ins  (\Xcal^\Omega \sm \Xcal_1^\Omega)  \cup \Xcal_2^\Omega}} \tilde \Xcal.
    \end{align}
    Note that the set $(\Xcal^\Omega \sm \Xcal_1^\Omega)  \cup \Xcal_2^\Omega$ itself might not be closed under composition with $\Omega^\Omega$
    or not even consist of maps mapping to $(\Xcal\sm \Xcal_1) \cup \Xcal_2$.
        \item One can recover the partial order $\le$ via the implication operation:
\begin{align} \Xcal_1 \le \Xcal_2 \iff \Xcal_1 \himp \Xcal_2 = \Xcal. \end{align}
\end{enumerate}
\end{Rem}

\begin{Lem}[Pullback as Heyting algebra homomorphism]
    \label{qms:lem:pullback-heyting-alg-hom}
        Let $f:\, (\Xcal,\Xcal^\Omega) \to (\Ycal,\Ycal^\Omega)$ be a quasi-measurable map. Then the pullback map
        $f^*:\, \Sub(\Ycal) \to \Sub(\Xcal)$ is a homomorphism of complete Heyting algebras.
    \begin{proof}
        We clearly have for $\Ycal_i \in \Sub(\Ycal)$:
        \begin{align}
            f^*(\zero) & =\zero, & f^*(\Ycal) & = \Xcal, \\
            f^*\lp \bigand_{i \in I} \Ycal_i  \rp &=\bigand_{i \in I} f^*(\Ycal_i), &  f^*\lp \bigor_{i \in I} \Ycal_i  \rp &=\bigor_{i \in I} f^*\lp \Ycal_i \rp,
        \end{align}
        and:
        \begin{align}
            \Ycal_1 \le \Ycal_2 & \implies f^*(\Ycal_1) \le f^*(\Ycal_2).
        \end{align}
        Now consider $\Ycal_1 \himp \Ycal_2$. We have:
        \begin{align}
            f^*(\Ycal_1 \himp \Ycal_2) & = \bigor_{\substack{\tilde \Ycal \in \Sub(\Ycal)\\ \tilde \Ycal \land \Ycal_1 \le \Ycal_2}} f^*\lp\tilde \Ycal\rp \le \bigor_{\substack{\tilde \Xcal \in \Sub(\Xcal)\\ \tilde \Xcal \land f^*(\Ycal_1) \le f^*(\Ycal_2)}} \tilde \Xcal =f^*(\Ycal_1) \himp f^*(\Ycal_2).
        \end{align}
    Now consider $\Xcal_1 \in \Sub(\Xcal)$ with $ \Xcal_1 \land f^*(\Ycal_1) \le f^*(\Ycal_2)$.
    We claim that then also:
    \begin{align} f_!(\Xcal_1) \land \Ycal_1 \le \Ycal_2. \end{align}
    For this consider: $Y \in f_!(\Xcal_1)^\Omega \cap \Ycal_1^\Omega$. Then there exists an $X \in \Xcal_1^\Omega$ with $f(X)=Y$.
    Since $f(X)=Y \in \Ycal_1^\Omega$ we have: $X \in f^*(\Ycal_1)^\Omega$. So we get:
    \begin{align} X \in \Xcal_1^\Omega \cap f^*(\Ycal_1)^\Omega \ins f^*(\Ycal_2)^\Omega. \end{align}
    This then shows: $Y=f(X) \in \Ycal_2^\Omega$ and thus the claim.\\
    We now put: $\Xcal_1 := f^*(\Ycal_1) \himp f^*(\Ycal_2)$. Since it satisfies $\Xcal_1 \land f^*(\Ycal_1) \le f^*(\Ycal_2)$ the above
    result applies and we get:
    \begin{align}  f_!(\Xcal_1) \land \Ycal_1 \le \Ycal_2
                & \iff f_!(\Xcal_1) \le \Ycal_1 \himp \Ycal_2
                 \iff \Xcal_1 \le f^*(\Ycal_1 \himp \Ycal_2).
    \end{align}
    This shows:
    \begin{align}
        f^*(\Ycal_1) \himp f^*(\Ycal_2) & \le f^*(\Ycal_1 \himp \Ycal_2),
    \end{align}
    and thus the equality:
    \begin{align}
        f^*(\Ycal_1 \himp \Ycal_2)& = f^*(\Ycal_1) \himp f^*(\Ycal_2).
    \end{align}
So $f^*$ is a Heyting algebra homomorphism.
    \end{proof}
\end{Lem}

\begin{Def}
    Let $f:\, (\Xcal,\Xcal^\Omega) \to (\Ycal,\Ycal^\Omega)$ be a quasi-measurable map.
    We define the following push-forward map: $f_*:\, \Sub(\Xcal) \to \Sub(\Ycal)$
    on $\Xcal_1 \in \Sub(\Xcal)$ via:
    \begin{align}
            f_*(\Xcal_1) &:= \bigor_{\substack{\tilde \Ycal \in \Sub(\Ycal)\\ f^*(\tilde \Ycal) \le \Xcal_1}} \tilde \Ycal.
    \end{align}
\end{Def}

\begin{Lem}
    Let $f:\, (\Xcal,\Xcal^\Omega) \to (\Ycal,\Ycal^\Omega)$ be a quasi-measurable map.
    We then have for $\Xcal_1 \in \Sub(\Xcal)$ and $\Ycal_1 \in \Sub(\Ycal)$ the adjunction:
    \begin{align}
        f^*(\Ycal_1) \le \Xcal_1 &\iff \Ycal_1 \le f_*(\Xcal_1).
    \end{align}
    In particular, we have:
    \begin{align}
        f^*f_*(\Xcal_1) &\le \Xcal_1, & \Ycal_1 &\le f_*f^*(\Ycal_1).
    \end{align}
    \begin{proof}
    ``$\implies$'': If $f^*(\Ycal_1) \le \Xcal_1$ then clearly:
    \begin{align}
        \Ycal_1 \le \bigor_{\substack{\tilde \Ycal \in \Sub(\Ycal)\\ f^*(\tilde \Ycal) \le \Xcal_1}} \tilde \Ycal = f_*(\Xcal_1).
    \end{align}
``$\Longleftarrow$'': If $\Ycal_1 \le f_*(\Xcal_1)$ then we apply $f^*$ and \Cref{qms:lem:pullback-heyting-alg-hom} to get:
        \begin{align}
            f^*(\Ycal_1) \le f^*f_*(\Xcal_1)  = \bigor_{\substack{\tilde \Ycal \in \Sub(\Ycal)\\f^*(\tilde \Ycal) \le \Xcal_1}} f^*(\tilde \Ycal)
            \le \Xcal_1.
    \end{align}
    This shows the claim.
\end{proof}
\end{Lem}

\begin{Rem}
    Note that we have the equivalence:
    \begin{align}
        f^{-1}(\Ycal_1) \ins \Xcal_1 & \iff f^{-1}(\Ycal_1^\cmpl) \sni \Xcal_1^\cmpl \iff \Ycal_1^\cmpl \sni f(\Xcal_1^\cmpl) \iff
        \Ycal_1 \ins f(\Xcal_1^\cmpl)^\cmpl.
    \end{align}
    This shows that:
    \begin{align} f_*(\Xcal_1) &\ins f(\Xcal_1^\cmpl)^\cmpl  =\Ycal \sm f(\Xcal \sm \Xcal_1). \end{align}
    It is now tempting to believe that the underlying set of $f_*(\Xcal_1)$ equals $f(\Xcal_1^\cmpl)^\cmpl$.
    However, this would require for every $y \in f(\Xcal_1^\cmpl)^\cmpl$ that every $X \in \Xcal^\Omega$ with
    $X(\Omega) \ins f^{-1}(y)$ lies in $\Xcal_1^\Omega$, which might not be true depending on $\Xcal^\Omega$, $\Xcal_1^\Omega$ and $f$.
\end{Rem}

\begin{Lem}
    Let $(\Xcal,\Xcal^\Omega)$ be a quasi-measurable space, $\Xcal_1, \Xcal_2 \in \Sub(\Xcal)$
    and $\iota_{\Xcal_1}:\, \Xcal_1 \inj \Xcal$ the (quasi-measurable) inclusion map.
    Then we have the equality of quasi-measurable subspaces of $\Xcal$:
    \begin{align}
        \Xcal_1 \himp \Xcal_2 &= \iota_{\Xcal_1,*}\iota_{\Xcal_1}^* \Xcal_2.
    \end{align}
    \begin{proof}
    \begin{align}
        \iota_{\Xcal_1,*}\iota_{\Xcal_1}^* \Xcal_2
        &= \bigor_{\substack{\tilde \Xcal \in \Sub(\Xcal)\\ \iota_{\Xcal_1}^*(\tilde \Xcal) \le \iota_{\Xcal_1}^*(\Xcal_2)}} \tilde \Xcal
        = \bigor_{\substack{\tilde \Xcal \in \Sub(\Xcal)\\ \iota_{\Xcal_1,!}\iota_{\Xcal_1}^*(\tilde \Xcal) \le \Xcal_2}} \tilde \Xcal
        = \bigor_{\substack{\tilde \Xcal \in \Sub(\Xcal)\\ \tilde \Xcal \land \Xcal_1 \le \Xcal_2}} \tilde \Xcal
        = \Xcal_1 \himp \Xcal_2.
    \end{align}
    This shows the claim.
    \end{proof}
\end{Lem}

A \emph{Heyting category} is by definition a \emph{regular} coherent category in which every
pullback functor between subobject lattices has a right adjoint, so before collecting the results
above we verify regularity, which has not been used so far.

\begin{Lem}[$\QMS$ is regular]
    \label{qms:lem:regular}
    The category $\QMS$ of quasi-measurable spaces and quasi-measurable maps is \emph{regular},
    i.e.\ it has finite limits, every morphism factors as a regular epimorphism followed by a
    monomorphism, and regular epimorphisms are stable under pullback.
    \begin{proof}
        Finite limits exist by \Cref{qms:thm:co-complete}.

        \emph{Factorisation.} Let $f:\,(\Xcal,\Xcal^\Omega) \to (\Ycal,\Ycal^\Omega)$ be
        quasi-measurable and let $\Wcal := f(\Xcal) \ins \Ycal$ be its image, endowed with
        $\Wcal^\Omega := f \circ \Xcal^\Omega$. Then $(\Wcal,\Wcal^\Omega)$ is a quasi-measurable
        space, $t:\,\Xcal \srj \Wcal$ is a quotient map and $\iota:\, \Wcal \inj \Ycal$ is an
        embedding, so $f = \iota \circ t$ with $t$ a quotient map and $\iota$ in particular a
        monomorphism, see \Cref{qms:sec:quotients} and \Cref{qms:sec:image-preimage}. By
        \Cref{qms:lem:strong-epi-quotient} the quotient maps are exactly the strong epimorphisms,
        and a strong epimorphism that arises as the quotient by the kernel pair of $f$ is a regular
        epimorphism: $t$ is the coequalizer of $\Xcal \ftimes_\Wcal \Xcal \rightrightarrows \Xcal$,
        because a quasi-measurable map out of $\Xcal$ that is constant on the fibres of $t$
        factors through $t$ uniquely as a map of sets, and the factorisation is quasi-measurable
        since $\Wcal^\Omega = t \circ \Xcal^\Omega$.

        \emph{Stability under pullback.} Let $t:\, \Xcal \srj \Wcal$ be a quotient map and
        $g:\,\Zcal \to \Wcal$ any quasi-measurable map, and form the fibre product
        $\Zcal \ftimes_\Wcal \Xcal$, see \Cref{qms:sec:fibre-prod}. We claim the projection
        $\pr_\Zcal:\, \Zcal \ftimes_\Wcal \Xcal \to \Zcal$ is again a quotient map. It is
        surjective because $t$ is. For the structure, let $Z \in \Zcal^\Omega$. Then
        $g \circ Z \in \Wcal^\Omega = t \circ \Xcal^\Omega$, so there is an $X \in \Xcal^\Omega$
        with $t \circ X = g \circ Z$; hence $(Z,X) \in \lp\Zcal \ftimes_\Wcal \Xcal\rp^\Omega$ and
        $\pr_\Zcal \circ (Z,X) = Z$. So every admissible random variable on $\Zcal$ lifts, which is
        exactly the statement that $\pr_\Zcal$ is a quotient map.
    \end{proof}
\end{Lem}

All these lemmas together now yield the following result:

\begin{Thm}
    \label{qms:thm:heyting-app}
    The category $\QMS$ of quasi-measurable spaces and quasi-measurable maps is a Heyting category.
\end{Thm}

\subsection{Change of Sample Space}

To be explicit, let $\QMS$ be the category of quasi-measurable spaces w.r.t.\ the sample space $(\Omega,\Bcal_\Omega)$,
see \Cref{qms:def:sample-space-act-1},
and $(\Wcal,\Wcal^\Omega)$ a quasi-measurable space in $\QMS$.
We then can use $(\Wcal,\Wcal^\Wcal)$, where:
\begin{align}\Wcal^\Wcal &= \QMS\lp(\Wcal,\Wcal^\Omega),(\Wcal,\Wcal^\Omega)\rp,\end{align}
as another sample space
for its own category of quasi-measurable space $\widetilde \QMS$.
For this note that:
\begin{align} \Xcal^\Wcal \circ \Wcal^\Wcal &\ins \Xcal^\Wcal,& \Xcal^\one &\ins \Xcal^\Wcal.\end{align}
We then have a functor, see \Cref{lem:exp-qm}:
\begin{align} \QMS &\to \widetilde \QMS, & (\Xcal,\Xcal^\Omega) &\mapsto(\Xcal,\Xcal^\Wcal), & g &\mapsto g, \end{align}
 where we use:
 \begin{align} \Xcal^\Wcal = \QMS\lp(\Wcal,\Wcal^\Omega),(\Xcal,\Xcal^\Omega)\rp.  \end{align}
So we are able to functorially shift from the sample space $(\Omega,\Omega^\Omega)$ to the sample
space $(\Wcal,\Wcal^\Wcal)$.

If, furthermore, we have more structure on the sample space $(\Omega,\Omega^\Omega,\Bcal_\Omega,\PrPf_\Omega)$, see \Cref{qms:def:sample-space-act-2} and \Cref{qms:def:sample-space-act-3}, we can then put:
\begin{align}
    \Bcal_\Wcal &:=\Bcal(\Wcal^\Omega), &
    \PrPf_\Wcal &:= \PrPf_\Omega(\Wcal,\Wcal^\Omega).
\end{align}
We then have:
\begin{align} \Wcal^\Wcal &= \QMS\lp(\Wcal,\Wcal^\Omega),(\Wcal,\Wcal^\Omega)\rp \ins \Meas\lp(\Wcal,\Bcal_\Wcal),(\Wcal,\Bcal_\Wcal) \rp.
\end{align}
With these choices then $(\Wcal,\Wcal^\Wcal,\Bcal_\Wcal,\PrPf_\Wcal)$ also satisfies \Cref{qms:def:sample-space-act-2} and \Cref{qms:def:sample-space-act-3}.

We leave further investigations into this topic to future research.

\subsection{Image and Pre-Image Structure}
\label{qms:sec:image-preimage}

In this subsection we will shortly look at how one can transfer the structure of a quasi-measurable space to another set through maps in and out of that space.

\begin{Def}
    Let $(\Zcal,\Zcal^\Omega)$ be a quasi-measurable space and $\Ycal$ a set and
    \begin{align}g:\, \Zcal \to \Ycal\end{align}
    a map.
    Then we define the image or \emph{push-forward quasi-measurable space} structure $(\Ycal, g_\circ\Zcal^\Omega)$ via:
    \begin{align} g_\circ\Zcal^\Omega:= \lp g\circ\Zcal^\Omega \rp \cup \lp \Ycal^\one\rp.\end{align}
    Note that $g_\circ\Zcal^\Omega$ is the smallest set of maps $\Ycal^\Omega \ins \Sets(\Omega, \Ycal)$
    such that $g$ becomes a  $\Zcal^\Omega$-$\Ycal^\Omega$-quasi-measurable map.
\end{Def}

\begin{Rem}
    Note that $g \circ \Zcal^\Omega$ might not satisfy all the points from Definition \ref{qms:def:quasi-measurable-spaces}.
    For instance, if $g$ is not surjective, then there are constant maps $ Y:\, \Omega \to \one \to \Ycal$
    that do not factorize through $g$.
\end{Rem}

\begin{Lem}
    \label{qms:lem:srj-pf-qm-str}
   Let $(\Zcal,\Zcal^\Omega)$ be a quasi-measurable space and $\Ycal$ a set and
    \begin{align}g:\, \Zcal \to \Ycal\end{align}
    a \emph{surjective} map.
    Then $(\Ycal, \Ycal^\Omega)$ with $\Ycal^\Omega:=g \circ \Zcal^\Omega$ is quasi-measurable space.
    \begin{proof}
        We check the the points from Definition \ref{qms:def:quasi-measurable-spaces}:
        \begin{enumerate}
            \item Let $ Y =g \circ  Z \in g \circ \Zcal^\Omega$ and $\Phi \in \Omega^\Omega$.
                Then $ Z \circ \Phi \in \Zcal^\Omega$ and thus $g \circ ( Z \circ \Phi) \in g \circ \Zcal^\Omega$.
            \item Let $y \in \Ycal$. Since $g$ is surjective there exists a $z \in \Zcal$ such that $g(z)=y$.
                    Then the constant map $ Z=z$ exists in $\Zcal^\Omega$. Then $g \circ  Z$ is the constant map $y$ and is in $g \circ \Zcal^\Omega$.
        \end{enumerate}
    \end{proof}
\end{Lem}

\begin{Def}
    Let $(\Ycal,\Ycal^\Omega)$ be a quasi-measurable space and $\Xcal$ a set and
    \begin{align}f:\, \Xcal \to \Ycal\end{align}
    a map.
    Then we define the pre-image or \emph{pull-back quasi-measurable space} structure $(\Xcal, f^*\Ycal^\Omega)$ via:
    \begin{align}f^*\Ycal^\Omega := \lC  X:\, \Omega \to \Xcal\st f \circ  X \in \Ycal^\Omega \rC.\end{align}
\end{Def}

\begin{Lem}
    \label{qms:lem:pull-back-maps}
        Let $(\Ycal,\Ycal^\Omega)$ be a quasi-measurable space and $\Xcal$ a set and
    \begin{align}f:\, \Xcal \to \Ycal\end{align}
    a map.
    Then the pull-back  $(\Xcal, f^*\Ycal^\Omega)$ is a quasi-measurable space w.r.t.\ $\Omega$.\\
    Note that $f^*\Ycal^\Omega$ is the biggest set of maps $\Xcal^\Omega \ins \Sets(\Omega, \Xcal)$
    such that $f$ becomes a $\Xcal^\Omega$-$\Ycal^\Omega$-quasi-measurable map.\\
    \begin{proof}
        We check the the points from Definition \ref{qms:def:quasi-measurable-spaces}:
        \begin{enumerate}
            \item Let $ X \in f^*\Ycal^\Omega$ and $\Phi \in \Omega^\Omega$.
                Then $f \circ  X \in \Ycal^\Omega$ and by assumption $f \circ  X \circ \Phi \in \Ycal^\Omega$ as well.
                By definition of $f^*\Ycal^\Omega$ we have that $ X \circ \Phi \in f^*\Ycal^\Omega$.
            \item If $ X:\,\Omega \to \Xcal$ is constant, then $f\circ  X:\, \Omega \to \Ycal$ is constant as well.
                So we have $f \circ  X \in \Ycal^\Omega$ by assumption and thus $ X \in f^*\Ycal^\Omega$.
        \end{enumerate}
    \end{proof}
\end{Lem}

\begin{Lem}
    \label{lem:subspace-qms}
    Let $(\Xcal,\Xcal^\Omega)$ and $(\Ycal,\Ycal^\Omega)$ be two quasi-measurable spaces and $\Zcal \ins \Ycal$ a subset with inclusion map
    $\iota:\; \Zcal \to \Ycal$. Further, let:
    \begin{align}f:\; (\Xcal,\Xcal^\Omega) \to (\Ycal,\Ycal^\Omega)\end{align}
    be a quasi-measurable map with $f(\Xcal) \ins \Zcal$.
    Then $f$ induces a quasi-measurable map:
    \begin{align}f:\; (\Xcal,\Xcal^\Omega) \to (\Zcal,\iota^*\Ycal^\Omega),\end{align}
    with:
    \begin{align}\iota^*\Ycal^\Omega = \lC  Y \in \Ycal^\Omega\st  Y(\Omega) \ins \Zcal \rC = \Sets(\Omega, \Zcal) \cap \Ycal^\Omega. \end{align}
    \begin{proof}
        If $f(\Xcal) \ins \Zcal$ then for any $ X \in \Xcal^\Omega$ we have that:
        \begin{align}f \circ  X (\Omega) \ins f(\Xcal) \ins \Zcal.\end{align}
        So $f \circ  X \in \Sets(\Omega, \Zcal)\cap \Ycal^\Omega = \iota^*\Ycal^\Omega$.
        This shows the claim.
    \end{proof}
\end{Lem}

\section{Appendix - Measurable Spaces vs.\ Quasi-Measurable Spaces}

\subsection{The Induced Sigma-Algebra}
\label{qms:sec:sigma-algebra}

In this subsection we show that the \emph{induced $\sigma$-algebra} $\Bcal_\Xcal$ of a quasi-measurable space $(\Xcal,\Xcal^\Omega)$ can be turned into a quasi-measurable space $(\Bcal_\Xcal, (\Bcal_\Xcal)^\Omega)$

\begin{Lem}
    \label{qms:lem:s-algebra-indicator}
    Let $(\Xcal,\Xcal^\Omega)$ be a quasi-measurable space. Then also its $\sigma$-algebra $\Bcal_\Xcal=\Bcal(\Xcal^\Omega)$ forms a well-defined quasi-measurable space $(\Bcal_\Xcal,(\Bcal_\Xcal)^\Omega)$.
\begin{proof}
  We first note that for $D \in \Bcal_{\Omega \times \Xcal}$ each section $D_\omega$ lies in $\Bcal_\Xcal$.
  For this consider $ X \in \Xcal^\Omega$ and the constant map $\Phi=\omega$. Since $\Phi \in \Omega^\Omega$ and by definition
  of $\Bcal_{\Omega \times \Xcal}$ we then have that:
  \begin{align}  X^{-1}(D_{\omega})=(\Phi, X)^{-1}(D) \in \Bcal_\Omega.\end{align}
  Since this holds for all $ X \in \Xcal^\Omega$ we get that $D_{\omega} \in \Bcal_\Xcal$. This shows the claim that the map $\omega \mapsto D_\omega$
  is a well-defined map $\Omega \to \Bcal_\Xcal$ for $D \in\Bcal_{\Omega \times \Xcal}$.\\
  We now check the points from Definition \ref{qms:def:quasi-measurable-spaces}:
        \begin{enumerate}
            \item Let $\Phi \in \Omega^\Omega$ and $\Psi \in (\Bcal_\Xcal)^\Omega$ be given via the sections of $D \in \Bcal_{\Omega \times \Xcal}$.
                Consider:
                \begin{align} \tilde D:= (\Phi \times \id_\Xcal)^{-1}(D) = \lC (\omega,x) \in \Omega \times \Xcal\st (\Phi(\omega),x) \in D\rC \ins \Omega \times \Xcal.\end{align}
                It is then $\tilde D \in \Bcal_{\Omega \times \Xcal}$. Indeed, for every $\tilde \Phi \in \Omega^\Omega$ and
                $\tilde  X \in \Xcal^\Omega$ we have:
                \begin{align*}
                (\tilde \Phi, \tilde  X)^{-1}(\tilde D) &= (\Phi \circ \tilde \Phi, \tilde  X)^{-1}(D) \in \Bcal_{\Omega}.
                \end{align*}
                Since this holds for all $\tilde \Phi \in \Omega^\Omega$ and $\tilde  X \in \Xcal^\Omega$
                we have that: $\tilde D \in \Bcal_{\Omega \times \Xcal}$.\\
                With these considerations we then get:
            \begin{align*}
                \Psi \circ \Phi(\omega)  = D_{\Phi(\omega)}
                                            = \lC x \in \Xcal\st (\Phi(\omega),x) \in D \rC
                                            = \tilde D_\omega.
            \end{align*}
            So $\Psi \circ \Phi$ is given via the sections of $\tilde D \in \Bcal_{\Omega \times \Xcal}$. So we can conclude that
            $\Psi \circ \Phi \in (\Bcal_\Xcal)^\Omega$.
        \item Let $\Psi: \Omega \to \Bcal_\Xcal$ be the constant map $\omega \mapsto A \in \Bcal_\Xcal$.
            Then $D:=\Omega \times A \in \Bcal_{\Omega \times \Xcal}$. Indeed, for $\Phi \in \Omega^\Omega$ and $ X \in \Xcal^\Omega$ we get:
            \begin{align} (\Phi, X)^{-1}(D) = \lC \omega \in \Omega\st (\Phi(\omega), X(\omega)) \in \Omega \times A \rC
            =  X^{-1}(A) \in \Bcal_\Omega.\end{align}
            Since this holds for all $\Phi \in \Omega^\Omega$ and $ X \in \Xcal^\Omega$ we get $D \in \Bcal_{\Omega \times \Xcal}$.\\
            With this we get for all $\omega \in \Omega$:
            \begin{align}D_\omega = \lC x \in \Xcal\st (\omega,x) \in \Omega \times A \rC = A = \Psi(\omega).\end{align}
            So $\Psi \in (\Bcal_\Xcal)^\Omega$.
        \end{enumerate}
        This shows that $(\Bcal_\Xcal,(\Bcal_\Xcal)^\Omega)$ is a well-defined quasi-measurable space.
\end{proof}
\end{Lem}

\begin{Lem}
    \label{lem:s-alg-pull-back-qm}
    Let $(\Xcal,\Xcal^\Omega)$ and $(\Ycal,\Ycal^\Omega)$  be a quasi-measurable space.
    Then the following map:
    \begin{align}   \Bcal_\Ycal \times \Ycal^\Xcal  \to \Bcal_\Xcal,\qquad (B,f) \mapsto f^{-1}(B).  \end{align}
    is quasi-measurable.
    \begin{proof}
        The map is clearly well-defined (image in $\Bcal_\Xcal$).
        The same holds for the induced map:
        \begin{align}   \Bcal_{\Omega \times \Ycal} \times \Ycal^{\Omega\times\Xcal}  \to \Bcal_{\Omega \times \Xcal},\qquad
        (B,f) \mapsto (\id_\Omega \times f)^{-1}(B),  \end{align}
        where:
        \begin{align}  \id_\Omega \times f:\; \Omega \times \Xcal \to \Omega \times \Ycal, \qquad (\omega,x) \mapsto (\omega,f(\omega,x)).\end{align}
        This shows the claim.
\end{proof}
\end{Lem}

We will show that the quasi-measurable space of the induced $\sigma$-algebra $(\Bcal_\Xcal,(\Bcal_\Xcal)^\Omega)$
is isomorphic to the function space of indicator functions $(\twob^\Xcal,(\twob^\Xcal)^\Omega)$.

\begin{Lem}Let $(\Xcal,\Xcal^\Omega)$ be a quasi-measurable space.
    Then the indicator function (inclusion check) induces an isomorphism of quasi-measurable spaces:
    \begin{align} \I:\; (\Bcal_\Xcal,(\Bcal_\Xcal)^\Omega) \to (\twob^\Xcal,(\twob^\Xcal)^\Omega) \qquad A \mapsto \I_A.\end{align}
    \begin{proof}
        We first show that $\I$ is a well-defined map. For this let $A \in \Bcal_\Xcal=\Bcal(\Xcal^\Omega)$.
        Clearly, $\I_A \in \Sets(\Xcal,\twob)$. To check that $\I_A \in \twob^\Xcal$ let $ X \in \Xcal^\Omega$ and note that:
        \begin{align}\I_A \circ  X = \I_{ X^{-1}(A)}:\, \Omega \to \twob. \end{align}
        Since $ X^{-1}(A) \in \Bcal_\Omega$ we have that $\I_{ X^{-1}(A)}$ is measurable and thus $\I_A \circ  X  \in \twob^\Omega$.
        Since this holds for all $ X \in \Xcal^\Omega$ we get that $\I_A \in \QMS(\Xcal,\twob)= \twob^\Xcal$. So $\I$ is well-defined.\\
        As a next step we show that $\I$ is a quasi-measurable map. For this let $\psi \in (\Bcal_\Xcal)^\Omega$ be given
        by the sections of $D \in \Bcal_{\Omega \times \Xcal}$. Then note that:
        \begin{align} (\I \circ \psi(\omega))(x) = \I_{\psi(\omega)}(x) =  \I_{D_\omega}(x) = \I_D(\omega,x). \end{align}
        Via (un-)currying we can thus identify $\I \circ \psi$ and $\I_D$.
        Since $D \in \Bcal_{\Omega \times \Xcal}$ the same arguments as above show that: $\I_D \in \twob^{\Omega \times \Xcal}$.
        Via (un-)currying this shows that $\I \circ \psi \in (\twob^\Xcal)^\Omega$. Since this holds for all $\psi \in (\Bcal_\Xcal)^\Omega$
        we have that $\I$ is a quasi-measurable map.\\
        We now claim that the inverse to $\I$ is the map:
        \begin{align} \rho:\, \twob^\Xcal \to \Bcal_\Xcal, \quad \chi \mapsto \chi^{-1}(1). \end{align}
        We first show that $\rho$ is a well-defined map. For this let $\chi:\, \Xcal \to \twob$ be an element of $\twob^\Xcal$.
        Then for all functions $ X \in \Xcal^\Omega$ we have that $\chi \circ  X \in \twob^\Omega$.
        This means that $ X^{-1}(\chi^{-1}(1)) \in \Bcal_\Omega$ for all $ X \in \Xcal^\Omega$, which implies:
         \begin{align}\chi^{-1}(1) \in \Bcal_\Xcal.\end{align}
        This shows that $\rho$ is well-defined.\\
        Next we show that $\rho$ is quasi-measurable.
        For this we again identify $(\twob^\Xcal)^\Omega$ with $\twob^{\Omega \times \Xcal}$ and $(\Bcal_\Xcal)^\Omega$ with $\Bcal_{\Omega \times \Xcal}$.
        Then let $\xi \in \twob^{\Omega \times \Xcal}$. The same arguments as above show that $\xi^{-1}(1) \in \Bcal_{\Omega \times \Xcal}$,
        which then represents the corresponding function: $\omega \mapsto \xi^{-1}(1)_\omega$ in $(\Bcal_\Xcal)^\Omega$.
        Translating back, this means we have a well-defined map:
        \begin{align} \rho_*:\,  (\twob^\Xcal)^\Omega \to (\Bcal_\Xcal)^\Omega, \quad  \xi \mapsto (\omega \mapsto \xi^{-1}(1)_\omega=\lC x\in \Xcal\st \xi(\omega)(x)=1 \rC),\end{align}
         implying that: $\rho \circ (\twob^\Xcal)^\Omega  \ins (\Bcal_\Xcal)^\Omega$. This shows that $\rho$ is a well-defined quasi-measurable map.\\
         Finally, we show that $\I$ and $\rho$ are inverse to each other.
            \begin{align}  \rho \circ \I:\, \Bcal_\Xcal \to \Bcal_\Xcal,\quad A \mapsto (\I_A)^{-1}(1)=A,\end{align}
            and:
            \begin{align}  \I \circ \rho:\, \twob^\Xcal \to \twob^\Xcal,\quad \chi \mapsto \I_{\chi^{-1}(1)}=\chi.\end{align}
         Since these relations also hold for $\Bcal_{\Omega \times \Xcal}$ and $\twob^{\Omega \times \Xcal}$ similarly,
         we get that $\rho \circ \I = \id_{(\Bcal_\Xcal,(\Bcal_\Xcal)^\Omega)}$ and $\I \circ \rho = \id_{(\twob^\Xcal,(\twob^\Xcal)^\Omega)}$.
        This shows the claim.
    \end{proof}
\end{Lem}

\begin{Lem}
    Let $t:\; (\Xcal,\Xcal^\Omega) \to (\Ycal,\Ycal^\Omega)$ be quasi-measurable with:
    \begin{align} \Ycal^\Omega = t_\circ\Xcal^\Omega := (t \circ \Xcal^\Omega) \cup \Ycal^\one.  \end{align}
    Then we have:
    \begin{align} t_*\Bcal(\Xcal^\Omega) = \Bcal(t\circ\Xcal^\Omega)=\Bcal(t_\circ\Xcal^\Omega) = \Bcal(\Ycal^\Omega). \end{align}
    \begin{proof}
        \begin{align*}
            t_*\Bcal(\Xcal^\Omega) &:= \lC B \ins \Ycal\st t^{-1}(B) \in \Bcal(\Xcal^\Omega) \rC\\
          &= \lC B \ins \Ycal\st \forall  X \in \Xcal^\Omega.\,  X^{-1}(t^{-1}(B)) \in \Bcal_\Omega \rC \\
          &= \lC B \ins \Ycal\st \forall  Y \in t \circ \Xcal^\Omega.\,  Y^{-1}(B) \in \Bcal_\Omega \rC =
          \Bcal(t \circ \Xcal^\Omega)\\
          &= \lC B \ins \Ycal\st \forall  Y \in (t \circ \Xcal^\Omega) \cup \Ycal^\one.\,  Y^{-1}(B) \in \Bcal_\Omega \rC \\
          &= \Bcal(t_\circ\Xcal^\Omega).
        \end{align*}
    \end{proof}
\end{Lem}

\subsection{The Adjunction}
\label{qms:sec:adjunction}

In this subsection we study the \emph{adjunction}, see  \cite{Kan58} Ch.\ II or \cite{Joh02,Mac98}, between the category of measurable spaces $\Meas$ and the category of quasi-measurable spaces $\QMS$.
In the following let $(\Omega,\Omega^\Omega,\Bcal_\Omega)$ be a fixed sample space, see \Cref{qms:def:sample-space-act-2}.

The right-adjoint $\Fcal$ maps a $\sigma$-algebra $\Bcal_\Xcal$ to all
$\Bcal_\Omega$-$\Bcal_\Xcal$-measurable maps $\Fcal(\Bcal_\Xcal)$:
    \begin{align}\Fcal(\Bcal_\Xcal) := \Meas\lp(\Omega,\Bcal_\Omega),(\Xcal,\Bcal_\Xcal)\rp.\end{align}
The left-adjoint $\Bcal$ maps a
set of random variables/quasi-measurable functions $\Xcal^\Omega$ to the $\sigma$-algebra
$\Bcal(\Xcal^\Omega)$ they generate:
    \begin{align} \Bcal(\Xcal^\Omega ) :=  \lC A \ins \Xcal\st \forall  X \in \Xcal^\Omega .\,  X^{-1}(A) \in \Bcal_\Omega \rC.\end{align}

\begin{Lem}
    \label{qms:lem:measurable-space-qms}
    Let $\Xcal$ be a set and $\Bcal_\Xcal$ be a $\sigma$-algebra of subsets of $\Xcal$
    and $\Xcal^\Omega \ins \Meas((\Omega,\Bcal_\Omega),(\Xcal,\Bcal_\Xcal))$ any subset of measurable maps.
    Then we have the inclusions:
    \begin{align} \Bcal_\Xcal &\ins \Bcal(\Fcal(\Bcal_\Xcal)) \ins \Bcal(\Xcal^\Omega), &
  \Xcal^\Omega &\ins \Fcal(\Bcal(\Xcal^\Omega)) \ins \Fcal(\Bcal_\Xcal).\end{align}
    In particular, we get:
  \begin{align} \Bcal(\Fcal(\Bcal(\Xcal^\Omega))) &=\Bcal(\Xcal^\Omega), &
   \Fcal(\Bcal(\Fcal(\Bcal_\Xcal))) &= \Fcal(\Bcal_\Xcal).\end{align}
    \begin{proof}
       Since by assumption we have: $\Xcal^\Omega \ins \Fcal(\Bcal_\Xcal)$
       we get the reverse inclusion when applying $\Bcal$
       (smaller set of functions means less constraints on the measurable subsets):
        \begin{align} \Bcal(\Fcal(\Bcal_\Xcal)) \ins \Bcal(\Xcal^\Omega).\end{align}
       Since by assumption every $ X \in \Xcal^\Omega$ is $\Bcal_\Omega$-$\Bcal_\Xcal$-measurable we
        have that for every $A \in \Bcal_\Xcal$ that $ X^{-1}(A) \in \Bcal_\Omega$.
        Since this holds for every $ X \in \Xcal^\Omega$ we see that
        $A \in \Bcal(\Xcal^\Omega)$. This shows the inclusion:
        \begin{align} \Bcal_\Xcal \ins \Bcal(\Xcal^\Omega).\end{align}
        Since this would also hold for the special choice $\Xcal^\Omega=\Fcal(\Bcal_\Xcal)$ we get:
        \begin{align} \Bcal_\Xcal \ins \Bcal(\Fcal(\Bcal_\Xcal)).\end{align}
        This completes the first chain of inclusions.\\
        By applying $\Fcal$ we get the reverse inclusion (smaller $\sigma$-algebras mean less
        constraints on the measurability of functions):
        \begin{align}\Fcal(\Bcal(\Xcal^\Omega)) \ins \Fcal(\Bcal_\Xcal).\end{align}
        Since $\Bcal(\Xcal^\Omega)$ is the biggest $\sigma$-algebra such that all functions
        from $\Xcal^\Omega$ are measurable we also have the inclusion:
        \begin{align} \Xcal^\Omega \ins \Meas((\Omega,\Bcal_\Omega),(\Xcal,\Bcal(\Xcal^\Omega)))=\Fcal(\Bcal(\Xcal^\Omega)).\end{align}
        This completes the second chain of inclusions.\\
        If we use the corner case: $\Bcal_\Xcal:=\Bcal(\Xcal^\Omega)$, then
         the first chain of inclusions become equalities:
        \begin{align} \Bcal(\Xcal^\Omega) \ins \Bcal(\Fcal(\Bcal(\Xcal^\Omega))) \ins \Bcal(\Xcal^\Omega).\end{align}
        If we use the corner case: $\Xcal^\Omega:=\Fcal(\Bcal_\Xcal)$,
        then the second chain of inclusions become equalities:
        \begin{align} \Fcal(\Bcal_\Xcal) \ins \Fcal(\Bcal(\Fcal(\Bcal_\Xcal))) \ins \Fcal(\Bcal_\Xcal).\end{align}
    \end{proof}
\end{Lem}

\begin{Rem}
    \label{rem:sample-qms}
    Consider the sample space $(\Omega,\Omega^\Omega,\Bcal_\Omega)$.
    Then we have:
    \begin{align} \Bcal(\Omega^\Omega) = \Bcal_\Omega.\end{align}
    This means that there is no ambiguity in notation for $\Bcal_\Omega$.
    Note that, in general, we might have: $\Fcal(\Bcal_\Omega) \supsetneq \Omega^\Omega$.
    \begin{proof}
        By Lemma \ref{qms:lem:measurable-space-qms} we have:
        \begin{align}\Bcal_\Omega \ins \Bcal(\Omega^\Omega).\end{align}
        On the other hand, since $\id_\Omega \in \Omega^\Omega$, we get for $B \in \Bcal(\Omega^\Omega)$
        that $B=\id_\Omega^{-1}(B) \in \Bcal_\Omega$.
        This shows the other inclusion:
    \begin{align} \Bcal(\Omega^\Omega) \ins \Bcal_\Omega.\end{align}
    \end{proof}
\end{Rem}

\begin{Thm}[Adjunction]
    \label{thm:qms-meas-adjunction}
    The functor:
    \begin{align}\Bcal:\; \QMS \to \Meas,\quad (\Xcal,\Xcal^\Omega) \mapsto (\Xcal,\Bcal(\Xcal^\Omega)),
    \quad \Bcal(f)=f,\end{align}
is left-adjoint to the functor:
    \begin{align}\Fcal:\; \Meas \to \QMS,\quad (\Ycal,\Bcal_\Ycal) \mapsto (\Ycal, \Fcal(\Bcal_\Ycal)),\quad
    \Fcal(g)=g,\end{align}
That means for quasi-measurable $(\Xcal,\Xcal^\Omega)$ and measurable $(\Ycal,\Bcal_\Ycal)$
    we have canonical identifications:
    \begin{align} \Meas\lp(\Xcal,\Bcal(\Xcal^\Omega)),(\Ycal,\Bcal_\Ycal)\rp = \QMS\lp(\Xcal,\Xcal^\Omega), (\Ycal,\Fcal(\Bcal_\Ycal)) \rp.   \end{align}
    Furthermore, we have:
    \begin{align}\Bcal \circ \Fcal \circ \Bcal = \Bcal,\qquad \Fcal \circ \Bcal \circ \Fcal = \Fcal,\end{align}
    and:
    $\Bcal_\Ycal \ins \Bcal\Fcal(\Bcal_\Ycal)$ and $\Xcal^\Omega \ins \Fcal\Bcal(\Xcal^\Omega)$.
    \begin{proof}
    Most of this was already shown in \Cref{qms:lem:measurable-space-qms}.\\
    For the adjunction let:
    \[ g \in \Meas\lp(\Xcal,\Bcal(\Xcal^\Omega)),(\Ycal,\Bcal_\Ycal)\rp, \qquad
       X \in \Xcal^\Omega  \ins \Meas\lp(\Omega,\Bcal_\Omega),(\Xcal,\Bcal(\Xcal^\Omega))\rp. \]
    Then $g \circ  X \in \Meas\lp(\Omega,\Bcal_\Omega),(\Ycal,\Bcal_\Ycal)\rp = \Fcal(\Bcal_\Ycal)$.
    Since this holds for all $ X \in \Xcal^\Omega$ we get that:
    $ g \in \QMS\lp(\Xcal,\Xcal^\Omega), (\Ycal,\Fcal(\Bcal_\Ycal))\rp$.
    This implies the inclusion:
    \begin{align}\Meas\lp(\Xcal,\Bcal(\Xcal^\Omega)),(\Ycal,\Bcal_\Ycal)\rp
        &\ins \QMS\lp(\Xcal,\Xcal^\Omega), (\Ycal,\Fcal(\Bcal_\Ycal) \rp.\end{align}
    Now let $f \in \QMS\lp(\Xcal,\Xcal^\Omega), (\Ycal,\Fcal(\Bcal_\Ycal) \rp$
    then $f \in \Meas\lp(\Xcal,\Bcal(\Xcal^\Omega)),(\Ycal,\Bcal\Fcal(\Bcal_\Ycal)\rp \ins \Meas\lp(\Xcal,\Bcal(\Xcal^\Omega)),(\Ycal,\Bcal_\Ycal)\rp$.
    This shows the inclusion:
    \begin{align}\Meas\lp(\Xcal,\Bcal(\Xcal^\Omega)),(\Ycal,\Bcal_\Ycal)\rp \sni \QMS\lp(\Xcal,\Xcal^\Omega), (\Ycal,\Fcal(\Bcal_\Ycal) \rp.\end{align}
    \end{proof}
\end{Thm}

\subsection{\Sturdy\ (Quasi-)Measurable Spaces}
\label{qms:sec:sturdy-qms}

\begin{Def}[\Sturdy\ (quasi-)measurable spaces]
    We define the full subcategories:
    \begin{enumerate}
        \item $\SMeas:=\Bcal(\QMS)$  of \emph{\sturdy\ measurable spaces} inside $\Meas$.
        \item $\SQMS:=\Fcal(\Meas)$  of \emph{\sturdy\ quasi-measurable spaces} inside $\QMS$.
    \end{enumerate}
\end{Def}

\begin{Rem}
    \begin{enumerate}
        \item A measurable space $(\Xcal,\Bcal_\Xcal)$ is \sturdy\ iff $\Bcal_\Xcal=\Bcal\Fcal(\Bcal_\Xcal)$.
        \item A quasi-measurable space $(\Xcal,\Xcal^\Omega)$ is \sturdy\ iff $\Xcal^\Omega=\Fcal\Bcal(\Xcal^\Omega)$.
    \end{enumerate}
    \begin{proof}
        This directly follows from \Cref{thm:qms-meas-adjunction}:
    \begin{align}\Bcal \circ \Fcal \circ \Bcal = \Bcal,\qquad \Fcal \circ \Bcal \circ \Fcal = \Fcal.\end{align}
    \end{proof}
\end{Rem}

\begin{Cor}
    \label{qms:cor:sturdy-meas}
    The adjunction:
    \[
\begin{tikzcd}
\QMS\ar[r,bend left,"\Bcal",""{name=A, below}] & \Meas\ar[l,bend left,"\Fcal",""{name=B,above}] \ar[from=A, to=B, symbol=\dashv]
\end{tikzcd}
\]
   restricts to an equivalence on the full subcategories given by the essential images of $\Bcal$ and $\Fcal$, resp.:
\[
\begin{tikzcd}
    \SQMS\ar[r,bend left,"\Bcal",""{name=A, below}] & \SMeas.\ar[l,bend left,"\Fcal",""{name=B,above}] \ar[from=A, to=B, symbol=\dashv]
\end{tikzcd}
\]
\begin{proof}
    Every adjunction $L \dashv R$ restricts to an equivalence on the full subcategories given by those objects where the
    unit $\eta:\,\id \to RL$, the counit $\varepsilon:\, LR \to \id$, resp., is an isomorphism.
    In our case this follows from \Cref{thm:qms-meas-adjunction}:
     \begin{align}\Bcal \circ \Fcal \circ \Bcal = \Bcal,\qquad \Fcal \circ \Bcal \circ \Fcal = \Fcal.\end{align}
\end{proof}
\end{Cor}

\begin{Rem}
    The main point of  \Cref{qms:cor:sturdy-meas} is that  \sturdy\ measurable spaces can be identified
    with \sturdy\ quasi-measurable spaces and their relations to other spaces are the same either considered as \sturdy\ measurable spaces or as \sturdy\ quasi-measurable spaces.
\end{Rem}

\begin{Thm}
    The full subcategory $\SQMS=\Fcal(\Meas)$ of \sturdy\ quasi-measurable spaces inside the category of all quasi-measurable spaces $\QMS$ is reflective. The reflector is given by:
    \begin{align}\Fcal\Bcal:\; \QMS \to \SQMS,\qquad (\Xcal,\Xcal^\Omega) \to (\Xcal,\Fcal\Bcal(\Xcal^\Omega)),\qquad \Fcal\Bcal(f)=f, \end{align}
    and preserves all colimits and coproducts (indexed over sets).
    \begin{proof}
      Let $(\Xcal,\Xcal^\Omega)$ be a quasi-measurable space and $(\Ycal,\Ycal^\Omega)$ a \sturdy\ quasi-measurable space.
    Then we need to show that $\Fcal\Bcal$ is left-adjoint to the forgetful functor $\SQMS \inj \QMS$, i.e.\ that we have a natural identification:
    \begin{align}  \QMS\lp(\Xcal,\Fcal\Bcal(\Xcal^\Omega)),(\Ycal,\Ycal^\Omega)\rp = \QMS\lp(\Xcal,\Xcal^\Omega),(\Ycal,\Ycal^\Omega)\rp.\end{align}
    Let $f \in \QMS\lp(\Xcal,\Fcal\Bcal(\Xcal^\Omega)),(\Ycal,\Ycal^\Omega)\rp$ then we have:
    \begin{align} f\circ \Xcal^\Omega \ins f \circ \Fcal\Bcal(\Xcal^\Omega) \ins \Ycal^\Omega.\end{align}
    This shows: $f \in \QMS\lp(\Xcal,\Xcal^\Omega),(\Ycal,\Ycal^\Omega)\rp$.\\
    For the reverse inclusion consider $g \in \QMS\lp(\Xcal,\Xcal^\Omega),(\Ycal,\Ycal^\Omega)\rp$.
    Then $g$ is $\Bcal(\Xcal^\Omega)$-$\Bcal(\Ycal^\Omega)$-measurable. Every $ X \in \Fcal\Bcal(\Xcal^\Omega)$ is, by definition, $\Bcal_\Omega$-$\Bcal(\Xcal^\Omega)$-measurable.
    So the composition $g \circ  X$ is $\Bcal_\Omega$-$\Bcal(\Ycal^\Omega)$-measurable.
    So we get:
    \begin{align} g \circ  X \in \Fcal\Bcal(\Ycal^\Omega)=\Ycal^\Omega.\end{align}
    This shows that $g \in \QMS\lp(\Xcal,\Fcal\Bcal(\Xcal^\Omega)),(\Ycal,\Ycal^\Omega)\rp$.\\
    This shows that $\SQMS$ is a reflective subcategory of $\QMS$.\\
    As a left-adjoint $\Fcal\Bcal$ always preserves arbitrary colimits like coproducts (indexed over sets).
    \end{proof}
\end{Thm}

\begin{Rem}
    Similarly, for a measurable space $(\Ycal,\Bcal_\Ycal)$ and a \sturdy\ measurable space $(\Xcal,\Bcal_\Xcal)$ we have the natural identification:
    \begin{align} \Meas\lp(\Xcal,\Bcal_\Xcal),(\Ycal,\Bcal_\Ycal)\rp=\Meas\lp(\Xcal,\Bcal_\Xcal),(\Ycal,\Bcal\Fcal(\Bcal_\Ycal))\rp. \end{align}
\begin{proof}
    Let $g \in \Meas\lp(\Xcal,\Bcal_\Xcal),(\Ycal,\Bcal_\Ycal)\rp$ and $B \in \Bcal\Fcal(\Bcal_\Ycal)$
    and $ X \in \Fcal(\Bcal_\Xcal)$. Then $g \circ  X \in \Fcal(\Bcal_\Ycal)$.
    Since $B \in \Bcal\Fcal(\Bcal_\Ycal)$ we get:
    \begin{align}  X^{-1}(g^{-1}(B)) \in \Bcal_\Omega.\end{align}
    Since this holds for all $ X \in \Fcal(\Bcal_\Xcal)$ we get:
    \begin{align} g^{-1}(B) \in \Bcal\Fcal(\Bcal_\Xcal) = \Bcal_\Xcal.\end{align}
    This shows that $g \in \Meas\lp(\Xcal,\Bcal_\Xcal),(\Ycal,\Bcal\Fcal(\Bcal_\Ycal))\rp$.\\
    The other inclusion is clear.
\end{proof}
\end{Rem}

\subsection{Sigma-Algebras and Products}

\begin{Rem}
    Let $I$ be an (index) set and $(\Xcal_i,\Xcal_i^\Omega)$ quasi-measurable spaces, $i \in I$. Then we always have the inclusion:
        \begin{align}
            \bigotimes_{i \in I} \Bcal(\Xcal_i^\Omega) & \ins \Bcal\lp\prod_{ i \in I} \Xcal_i^\Omega  \rp.
        \end{align}
   The reverse inclusion might not hold.
\end{Rem}

\begin{Lem}[$\Fcal$ on products of measurable spaces]
    \label{qms:lem:F-prod}
    Let $(\Xcal_i,\Bcal_i)$ measurable spaces, $i \in I$.
    Then we have:
    \begin{align} \Fcal\lp \bigotimes_{i \in I} \Bcal_i \rp = \prod_{i\in I} \Fcal(\Bcal_i).\end{align}
    This means that the functor $\Fcal:\;\Meas \to \QMS$ preserves arbitrary products.
    \begin{proof}
        Right adjoints always preserve arbitrary limits, like products.\\
        A more concrete proof goes as follows.
        By the univsersal property of the product ($\sigma$-algebra), the adjunction from  \Cref{thm:qms-meas-adjunction}
        and that $\Bcal_\Omega=\Bcal(\Omega^\Omega)$ (see  \Cref{rem:sample-qms}) we get:
        \begin{align*}
            \Fcal\lp\bigotimes_{i \in I} \Bcal_i\rp
            &= \Meas\lp (\Omega,\Bcal_\Omega), \lp \prod_{i \in I} \Xcal_i,\bigotimes_{i \in I} \Bcal_i  \rp \rp \\
            &= \prod_{i \in I} \Meas\lp (\Omega,\Bcal_\Omega), (\Xcal_i,\Bcal_i) \rp \\
            &= \prod_{i \in I} \Meas\lp (\Omega,\Bcal(\Omega^\Omega)), (\Xcal_i,\Bcal_i) \rp \\
            &= \prod_{i \in I} \QMS\lp (\Omega,\Omega^\Omega), (\Xcal_i,\Fcal(\Bcal_i)) \rp \\
            &= \prod_{i \in I} \Fcal(\Bcal_i).
        \end{align*}
 This shows the claim.
    \end{proof}
\end{Lem}

\begin{Lem}
Let $(\Xcal,\Xcal^\Omega)$ and $(\Zcal,\Zcal^\Omega)$  be  quasi-measurable spaces.
    Then the map:
    \begin{align} i:\;  \Xcal^\Zcal  \to  \prod_{z \in \Zcal} \Xcal, \qquad g \mapsto (g(z))_{z \in \Zcal},  \end{align}
    is a well-defined, injective quasi-measurable map.\\
    If, furthermore, $\Xcal^\Omega = \Fcal(\Bcal_\Xcal)$ and $(\Zcal,\Zcal^\Omega)$
    is countable with the discrete quasi-measurable space structure
    then $i$ is an isomorphism of quasi-measurable spaces.
    \begin{proof}
        It is clear that $i$ is a well-defined map and injective.
    Now let $ X \in \lp\Xcal^{\Zcal}\rp^\Omega \cong \lp\Xcal^{\Omega}\rp^\Zcal$.
    Under this identification we see that for every $z \in \Zcal$ we have that
     $ X(z) \in \Xcal^\Omega$ and thus:
     \begin{align} ( X(z))_{z \in \Zcal} \in \prod_{z \in \Zcal} \Xcal^\Omega.\end{align}
    This shows that $i$ is quasi-measurable and thus the claim.\\
    Now assume that  $(\Zcal,\Zcal^\Omega)$
    is countable and discrete, w.l.o.g.\ $(\Zcal,\Zcal^\Omega) = (\N,\N^\Omega)$
    and that
    $\Xcal^\Omega = \Meas((\Omega,\Bcal_\Omega),(\Xcal,\Bcal_\Xcal))$ holds.
    Then consider the map:
    \begin{align} j:\;  \prod_{n \in \N} \Xcal \to \Xcal^\N, \qquad
    (x_n)_{n \in \N} \mapsto x=(n \mapsto x_n).  \end{align}
    This is well-defined. Indeed, since $\N$ is a discrete measurable space the map:
     \begin{align} x:\;\N \to \Xcal,\qquad n \mapsto x_n, \end{align}
     is measurable for every choice of $(x_n)_{n \in \N}$.  Since also each $ Y \in \N^\Omega$ is measurable we have that:
      $x \circ  Y \in  \Meas((\Omega,\Bcal_\Omega),(\Xcal,\Bcal_\Xcal))=\Xcal^\Omega$.
     This shows that $j$ is well-defined.\\
     Now consider the induced map:
     \begin{align} j_*:\; \prod_{n \in \N} \Xcal^\Omega \to (\Xcal^\N)^\Omega,\qquad ( X_n)_{n \in \N} \mapsto  X=((\omega,n) \mapsto  X_n(\omega)). \end{align}
     We need to show that $ X$ is quasi-measurable:
     \begin{align}  X:\; \Omega \times \N \to \Xcal.\end{align}
     For this let $\Phi \in \Omega^\Omega$ and $ Y \in \N^\Omega$ and $A \in \Bcal_\Xcal$.
     Then we get:
     \begin{align*}
          X(\Phi, Y)^{-1}(A) &= \lC \omega \in \Omega \st  X_{ Y(\omega)}(\Phi(\omega)) \in A \rC \\
                        &= \bigcup_{n \in \N} \Phi^{-1}( X_n^{-1}(A)) \cap  Y^{-1}(n) \\
                        & \in \Bcal_\Omega,
     \end{align*}
     because $\Phi$, $ Y$ and each $ X_n$ is measurable. It follows that:
     \begin{align}  X(\Phi, Y) \in \Meas((\Omega,\Bcal_\Omega),(\Xcal,\Bcal_\Xcal))=\Xcal^\Omega. \end{align}
     So $j$ is a well-defined quasi-measurable map.
     It is easily seen that $j$ and $i$ are inverse to each other.
     This shows the claim.
    \end{proof}
\end{Lem}

\section{Appendix - Sigma-Algebras on Spaces of Probability Measures}
\label{qms:sec:app-giry-sigma}

The following is the one measure-theoretic fact about the Giry $\sigma$-algebra that the
proofs in this paper rely on repeatedly. It is stated and proved here rather than in the
main text because it is used only inside other proofs. The statement also appears as
\cite{For21} Thm.\ B.42 and, in a different formulation, as \cite{Res77} Thm.\ 4; the
proof given here is self-contained and short, so we prefer it to a citation.

\begin{Prp}[Probability measures on spaces of probability measures, see \Cref{prp:prob-prob-meas:prf}]
    \label{prp:prob-prob-meas}
    Let $(\Xcal,\Bcal_\Xcal)$ be a measurable space and $\Giry(\Xcal,\Bcal_\Xcal)$ the space of all probability measures on
    $(\Xcal,\Bcal_\Xcal)$ endowed with the smallest $\sigma$-algebra $\Bcal_{\Giry(\Xcal,\Bcal_\Xcal)}$ such that all evaluation maps:
    \[ \ev_A:\; \Giry(\Xcal,\Bcal_\Xcal) \to [0,1],\qquad \mu \mapsto \mu(A),\]
    are measurable for all $A \in \Bcal_\Xcal$. Let $\Pcal \ins \Giry(\Xcal,\Bcal_\Xcal)$ be a fixed non-empty subset of probability measures endowed with a
    $\sigma$-algebra $\Bcal_\Pcal$ that contains the subspace $\sigma$-algebra $\Bcal_{\Giry(\Xcal,\Bcal_\Xcal)|\Pcal}$
    (e.g.\ $\Bcal_\Pcal=\Bcal_{\Giry(\Xcal,\Bcal_\Xcal)|\Pcal}$).\footnote{So the inclusion map $\iota:\;(\Pcal,\Bcal_\Pcal) \xhookrightarrow{\iota}\lp\Giry(\Xcal,\Bcal_\Xcal),\Bcal_{\Giry(\Xcal,\Bcal_\Xcal)}\rp$ is measurable.\label{fn:incl-meas}}
    Consider the map:
    \begin{gather*}
        \Mbb_\iota:\; \Giry(\Pcal,\Bcal_\Pcal) \xrightarrow{\iota_*} \Giry\lp\Giry(\Xcal,\Bcal_\Xcal),\Bcal_{\Giry(\Xcal,\Bcal_\Xcal)}\rp \xrightarrow{\Mbb} \Giry(\Xcal,\Bcal_\Xcal),\\
        \Mbb_\iota(\pi)(A) := \int_\Pcal \ev_A(\mu)\,\pi(d\mu),
    \end{gather*}
    and let $\Mcal \ins \Mbb_\iota^{-1}(\Pcal) \ins \Giry(\Pcal,\Bcal_\Pcal)$
    be a subset\footnote{We can allow for $\Mcal=\emptyset$ if we define $(\Bcal_\Pcal)_\emptyset := \lC D \ins \Pcal \rC$ to be the power-set of $\Pcal$ in the following.\label{fn:empty-completion}}.\\
    Then for every $B \in (\Bcal_\Xcal)_\Pcal$
     the evaluation map:
    \[ \ev_B:\; \Pcal \to [0,1],\qquad \mu \mapsto \mu(B),\]
    is $\pi$-measurable for every $\pi \in \Mcal$.
    In short:
        \[ \ev_B:\; \lp\Pcal,(\Bcal_\Pcal)_\Mcal\rp \to \lp[0,1],\Bcal_{[0,1]}\rp,\qquad \mu \mapsto \mu(B),\]
    is measurable for every $B \in (\Bcal_\Xcal)_\Pcal$. Or, even shorter, the inclusion map
    \[ \lp\Pcal,(\Bcal_\Pcal)_\Mcal\rp \inj \Giry(\Xcal,(\Bcal_\Xcal)_\Pcal)    \]
    is well-defined and measurable.\\
    In particular, the map:
    \[ \ev_B:\; \lp\Giry(\Xcal,\Bcal_\Xcal),(\Bcal_{\Giry(\Xcal,\Bcal_\Xcal)})_\Giry\rp \to \lp[0,1],\Bcal_{[0,1]}\rp,\qquad \mu \mapsto \mu(B),\]
    is (universally) measurable for all $B \in (\Bcal_\Xcal)_\Giry$.

\end{Prp}

\begin{Prf}{prp:prob-prob-meas}
Let $ X \in \R$ and $\pi \in \Mcal$. We need to show that $\ev_B^{-1}(( X,1]) \in (\Bcal_\Pcal)_\pi$. Now let $\nu:=\Mbb_\iota(\pi) \in \Pcal$, which for $A \in \Bcal_\Xcal$ is given by:
        \[ \nu(A) = \Mbb_\iota(\pi)(A):= \int_\Pcal \ev_A(\mu)\,\pi(d\mu).\]
        Note that the map $\Pcal \xhookrightarrow{\iota} \Giry(\Xcal,\Bcal_\Xcal) \xrightarrow{\ev_A} [0,1]$ is measurable for $A \in \Bcal_\Xcal$ and the above integral is well-defined.\\
    Since $B \in (\Bcal_\Xcal)_\Pcal \ins (\Bcal_\Xcal)_\nu$ there exist $A_1,A_2 \in \Bcal_\Xcal$ such that:
    \[ A_1 \ins B \ins A_2,\qquad \nu(A_2\sm A_1) =0.\]
    The last condition shows that:
    \[ 0 = \nu(A_2\sm A_1) = \int_\Pcal \ev_{A_2\sm A_1}(\mu)\,\pi(d\mu). \]
    Since $\ev_{A_2\sm A_1}(\mu) \ge 0$ and $\ev_{A_2\sm A_1}:\,\Pcal \to [0,1]$ is measurable we have that:
    \[ C:= \ev_{A_2\sm A_1}^{-1}((0,1])=\lC \mu \in \Pcal\st \ev_{A_2\sm A_1}(\mu)>0 \rC \in \Bcal_\Pcal,\qquad \text{and} \qquad \pi(C)=0.\]
    We now have the inclusions:
    \[ \ev_{A_1}^{-1}(( X,1]) \ins \ev_{B}^{-1}(( X,1]) \ins \ev_{A_2}^{-1}(( X,1]),\]
    with $\ev_{A_i}^{-1}(( X,1]) \in \Bcal_\Pcal$, $i=1,2$, and:
    \[ \ev_{A_2}^{-1}(( X,1]) \sm \ev_{A_1}^{-1}(( X,1]) \ins \ev_{A_2\sm A_1}^{-1}((0,1]) =C. \]
    Since $\pi(C)=0$ we also have that:
    \[ \pi\lp\ev_{A_2}^{-1}(( X,1]) \sm \ev_{A_1}^{-1}(( X,1]) \rp=0.  \]
    This shows that $\ev_B^{-1}(( X,1]) \in (\Bcal_\Pcal)_\pi$.
    Since this holds for all $\pi \in \Mcal$ we get that:
    \[ \ev_B^{-1}(( X,1]) \in (\Bcal_\Pcal)_\Mcal.\]
    This shows that the map:
    \[ \ev_B:\; \lp\Pcal,(\Bcal_\Pcal)_\Mcal\rp \to \lp[0,1],\Bcal_{[0,1]}\rp,\]
    is measurable for all $B \in (\Bcal_\Xcal)_\Pcal$ with $\Mcal \ins \Mbb_\iota^{-1}(\Pcal)$.
\end{Prf}

\section{Appendix - Monads and Markov Categories}
\subsection{Monads}

\label{qms:sec:app-monads}

In this section we remind the reader of the definition of a (commutavie, strong) monad.
For literature about monads see \cite{Koc70,Koc72,Str72,Mac98,Mog91}.

\begin{Def}[Monad, see \cite{Koc70,Koc72,Mac98,Mog91}]
    Let $\Ccal$ be a category.
    A \emph{monad} on $\Ccal$ is a triple $(\PrT,\delta,\Mbb)$ consisting of:
    \begin{enumerate}
        \item a functor $\PrT:\;\Ccal \to \Ccal$,
        \item a natural transformation $\delta:\; \id_\Ccal \to \PrT$,
        \item a natural transformation $\Mbb:\; \PrT^2:=\PrT \circ \PrT \to \PrT$,
    \end{enumerate}
    such that:
    \begin{enumerate}
        \item[a.] $\Mbb \circ \PrT\Mbb = \Mbb \circ \Mbb \PrT$ as natural transformations $\PrT^3 \to \PrT$,
        \item[b.] $\Mbb \circ \PrT\delta = \Mbb \circ \delta \PrT = \id_\PrT$ as natural transformations $\PrT \to \PrT$.
    \end{enumerate}
    A monad $(\PrT,\delta,\Mbb)$ on a category $\Ccal$ is called \emph{affine}, see \cite{Koc71,Lin79,Bar16},
    if $\Ccal$ has a terminal object $\one$ such that the unique morphism:
    \begin{align} !:\;\PrT(\one) \to \one\end{align}
    is an isomorphism.
\end{Def}

\begin{Def}[Strong monad, see \cite{Koc72,Mog91} Def.\ 3.2]
    A \emph{left-strong monad} over a monoidal category $(\Ccal,\times,\one)$ is a monad $(\PrT,\delta,\Mbb)$ on $\Ccal$
    together with a natural transformation, called \emph{left-strength}:
    \begin{align}
        \tau_{\Xcal,\Ycal}:\; \Xcal \times \PrT(\Ycal) \to \PrT(\Xcal \times \Ycal),
    \end{align}
    such that:
    \begin{enumerate}
        \item Left-unitor $L$ and left-strength $\tau$ satisfy the commutative diagram:
        \[\xymatrix{
                \one \times \PrT(\Xcal) \ar_-{L_{\PrT(\Xcal)}}[rd] \ar^-{\tau_{\one,\Xcal}}[rr] &&
                \PrT(\one \times \Xcal) \ar^{\PrT(L_\Xcal)}[dl]\\
            &\PrT(\Xcal).
        }\]
     \item Associator $A$ and left-strength $\tau$ satisfy the commutative diagram:
        \[\xymatrix{
                (\Xcal\times\Ycal)\times\PrT(\Zcal) \ar^-{\tau_{\Xcal\times\Ycal,\Zcal}}[rrrr] \ar_-{A_{\Xcal,\Ycal,\PrT(\Zcal)}}[d] &&&& \PrT\lp(\Xcal\times\Ycal)\times\Zcal\rp \ar^-{\PrT(A_{\Xcal,\Ycal,\Zcal})}[d]\\
                \Xcal\times(\Ycal\times\PrT(\Zcal)) \ar_-{\id_\Xcal \times \tau_{\Ycal,\Zcal}}[rr] && \Xcal\times\PrT(\Ycal\times\Zcal) \ar_-{\tau_{\Xcal,\Ycal\times\Zcal}}[rr] && \PrT(\Xcal\times(\Ycal\times\Zcal)).
        }\]
    \item Monad unit $\delta$ and left-strength $\tau$ satisfy the commutative diagram:
        \[\xymatrix{
            &\Xcal \times \Ycal \ar_-{\id_\Xcal\times \delta_\Ycal}[dl]
            \ar^-{\delta_{\Xcal\times\Ycal}}[dr]\\
                \Xcal \times \PrT(\Ycal) \ar_-{\tau_{\Xcal,\Ycal}}[rr] && \PrT(\Xcal \times \Ycal).
        }\]
     \item Monad action $\Mbb$ commutes with left-strength $\tau$, expressed via the commutative diagram:
        \[\xymatrix{
                \Xcal\times\PrT(\PrT(\Ycal)) \ar^-{\tau_{\Xcal,\PrT(\Ycal)}}[rr] \ar_-{\id_\Xcal\times\Mbb_\Ycal}[d] &&
                \PrT(\Xcal\times\PrT(\Ycal)) \ar^-{\PrT(\tau_{\Xcal,\Ycal})}[rr] &&
                \PrT(\PrT(\Xcal \times \Ycal)) \ar^-{\Mbb_{\Xcal\times\Ycal}}[d] \\
                \Xcal\times\PrT(\Ycal) \ar_-{\tau_{\Xcal,\Ycal}}[rrrr]&&&&
                \PrT(\Xcal\times\Ycal).
        }\]
        \end{enumerate}
 A \emph{right-strong monad} $(\PrT,\delta,\Mbb)$ on $\Ccal$ is defined similary, using
    a \emph{right-strength}, i.e.\ a natural transformation:
    \begin{align}
        \rho_{\Xcal,\Ycal}:\; \PrT(\Xcal) \times \Ycal \to \PrT(\Xcal \times \Ycal),
    \end{align}
    with the corresponding 4 commutative diagrams from above.

    We call a monad $(\PrT,\delta,\Mbb)$ on $\Ccal$ a \emph{strong monad} if it has both a left-strength $\tau$ and a right-strength $\rho$ that are compatible, which is expressed via the following commutative diagram:
     \[\xymatrix{
             (\Xcal \times \PrT(\Ycal)) \times \Zcal \ar_-{A_{\Xcal,\PrT(\Ycal),\Zcal}}[d] \ar^-{\tau_{\Xcal,\Ycal} \times \id_\Zcal}[rr]
             && \PrT(\Xcal \times \Ycal) \times \Zcal \ar^-{\rho_{\Xcal\times\Ycal,\Zcal}}[rr]
             &&  \PrT((\Xcal \times \Ycal) \times \Zcal) \ar^-{\PrT(A_{\Xcal,\Ycal,\Zcal})}[d] \\
             \Xcal \times (\PrT(\Ycal) \times \Zcal)) \ar_-{\id_\Xcal \times \rho_{\Ycal,\Zcal}}[rr]
             && \Xcal \times \PrT(\Ycal \times \Zcal) \ar_-{\tau_{\Xcal,\Ycal\times\Zcal}}[rr]
             &&\PrT(\Xcal \times (\Ycal \times \Zcal)).
     }\]

\end{Def}

\begin{Def}[Commutative monad, see \cite{Koc70,Koc71c,Koc72,Koc12}]
    \label{qms:def:commutative-monad}
    A strong monad $(\PrT,\delta,\Mbb,\tau,\rho)$ on a symmetric monoidal\footnote{\emph{Symmetric monoidal} means that we have a natural \emph{braiding} isomorphism (=swap operation) $\swap_{\Xcal,\Ycal}:\, \Xcal \times \Ycal \cong \Ycal \times \Xcal$ that is compatible with the associator and such that it is inverse to the reverse braiding $\swap_{\Ycal,\Xcal}:\, \Ycal \times \Xcal \cong \Xcal \times \Ycal$.} category $(\Ccal,\times,\one)$ is
    called \emph{commutative} if it is symmetric and we always have the following commutative diagram:
    \[\xymatrix{
        && \PrT(\Xcal) \times \PrT(\Ycal) \ar_-{\tau_{\PrT(\Xcal),\Ycal}}[lld]
        \ar^-{\rho_{\Xcal,\PrT(\Ycal)}}[rrd] \ar@{-->}^-{\otimes_{\Xcal,\Ycal}}[ddd]\\
        \PrT(\PrT(\Xcal) \times \Ycal) \ar_-{\PrT(\rho_{\Xcal,\Ycal})}[d]
            && && \PrT(\Xcal \times \PrT(\Ycal))
            \ar^-{\PrT(\tau_{\Xcal,\Ycal})}[d] \\
            \PrT(\PrT(\Xcal \times \Ycal)) \ar_-{\Mbb_{\Xcal\times\Ycal}}[rrd]
            && && \PrT(\PrT(\Xcal \times \Ycal))\ar^-{\Mbb_{\Xcal\times\Ycal}}[lld]\\
                                             && \PrT(\Xcal \times \Ycal).
    }\]
In this case, the middle arrow can then be defined by either the right or left side of the diagram.
\end{Def}

 \begin{Def}
    Let $(\Ccal,\times,\one)$ be a symmetric monoidal category. A strong monad
    $(\PrT,\delta,\Mbb,\tau,\rho)$ is called \emph{symmetric} if left- and right strength are
    related as follows:
    \[\xymatrix{
        \Xcal \times \PrT(\Ycal) \ar^-{\tau_{\Xcal,\Ycal}}[rr]
        \ar_-{\swap_{\Xcal,\PrT(\Ycal)}}[d]
        && \PrT(\Xcal \times \Ycal) \ar^-{\PrT(\swap_{\Xcal,\Ycal})}[d] \\
        \PrT(\Ycal) \times \Xcal \ar_-{\rho_{\Ycal,\Xcal}}[rr]
        && \PrT(\Ycal \times \Xcal).
    }\]
\end{Def}

\begin{Rem}
    \label{qms:rem:otimes-commutative}
    If $(\PrT,\delta,\Mbb)$ is a monad on a symmetric monoidal category $(\Ccal,\times,\one)$    that is endowed with a natural transformation $\otimes$:
    \begin{align}
        \otimes_{\Xcal,\Ycal}:\, \PrT(\Xcal) \times \PrT(\Ycal) \to \PrT(\Xcal \times \Ycal),
    \end{align}
    respecting associativity and unitality $\delta:\, \one \to \PrT(\one)$,
    then one can define a left- and right-strength as follows:
    \begin{align}
        \tau_{\Xcal,\Ycal}:\, &\Xcal \times \PrT(\Ycal) \stackrel{\delta_\Xcal \times \id_{\PrT(\Ycal)}}{\longrightarrow} \PrT(\Xcal) \times \PrT(\Ycal) \stackrel{\otimes_{\Xcal,\Ycal}}{\longrightarrow} \PrT(\Xcal \times \Ycal), \\
    \rho_{\Xcal,\Ycal}:\, & \PrT(\Xcal) \times \Ycal \stackrel{\id_{\PrT(\Xcal)} \times \delta_\Ycal}{\longrightarrow} \PrT(\Xcal) \times \PrT(\Ycal) \stackrel{\otimes_{\Xcal,\Ycal}}{\longrightarrow} \PrT(\Xcal \times \Ycal).
    \end{align}
    One can check that this leads to a commutative monad, see \cite{Koc72} Thm.\ 2.3. \end{Rem}

\subsection{Markov Categories}
\label{qms:sec:app-markov-cat}

The vocabulary in which \Cref{thm:Markov-category} is stated is that of \emph{categorical probability theory}; we recall it here, following \cite{Fri19,Fri20,FriG20,Fri21}.

\begin{Def}[Markov category, see \cite{Fri20} Def.\ 2.1 or \cite{Fri21} Def.\ 3.1]
    A \emph{Markov category} is a symmetric monoidal category $(\Ccal,\otimes,\one)$ where
    the monoidal unit object $\one$ is a terminal object and where every object $\Xcal \in \Ccal$
    is equipped
    with a commutative comonoid structure given by comultiplication, called $\cop$, and counit, called $\del$:
    \[ \cop_\Xcal:\; \Xcal \to \Xcal\otimes\Xcal,\qquad \del_\Xcal:\; \Xcal \to \one,\]
    that are natural in $\Xcal$ and compatible with the monoidal structure in the following sense:
    \begin{align*}
        (\cop_\Xcal\otimes\id_\Xcal)\circ\cop_\Xcal = (\id_\Xcal\otimes\cop_\Xcal)\circ\cop_\Xcal &:\; \Xcal \to\Xcal\otimes\Xcal \to \Xcal\otimes\Xcal\otimes\Xcal,\\
        (\del_\Xcal\otimes\id_\Xcal)\circ\cop_\Xcal = \id_\Xcal = (\id_\Xcal\otimes\del_\Xcal)\circ\cop_\Xcal &:\;\Xcal \to \Xcal\otimes\Xcal \to \Xcal,\\
        \swap_{\Xcal,\Xcal} \circ \cop_\Xcal = \cop_\Xcal, &:\; \Xcal \to \Xcal\otimes\Xcal,\\
(\id_\Xcal\otimes \swap_{\Xcal,\Ycal}\otimes\id_\Ycal) \circ(\cop_\Xcal \otimes \cop_\Ycal)  = \cop_{\Xcal\otimes\Ycal} &:\; \Xcal\otimes\Ycal \to \Xcal\otimes\Ycal\otimes\Xcal\otimes\Ycal.
\end{align*}
     hold for all $\Xcal,\Ycal \in \Ccal$, where $\swap$ is the braiding
     of $(\Ccal,\otimes,\one)$.
         \[\xymatrix@C=1.2em{
            \Xcal \ar^-{\cop_\Xcal}[rr] \ar_-{\cop_\Xcal}[d]&& \Xcal\otimes\Xcal \ar^-{\cop_\Xcal\otimes\id_\Xcal}[d] &&
            \Xcal \ar^-{\cop_\Xcal}[rr] \ar_-{\cop_\Xcal}[d] \ar@{=}[rrrrd]
              && \Xcal\otimes\Xcal \ar^-{\id_\Xcal\otimes\del_\Xcal}[rr] && \Xcal\otimes\one \ar^-{\mult}[d]\\
            \Xcal\otimes\Xcal\ar_-{\id_\Xcal\otimes\cop_\Xcal}[rr]  && \Xcal\otimes\Xcal\otimes\Xcal, &&
            \Xcal\otimes\Xcal \ar_-{\del_\Xcal\otimes\id_\Xcal}[rr] && \one\otimes\Xcal \ar_-{\mult}[rr] && \Xcal,
    }\]
    \[\xymatrix@C=1.2em{
            \Xcal \ar_-{\cop_\Xcal}[d] \ar^-{\cop_\Xcal}[rrd] && &&
            \Xcal\otimes\Ycal \ar_-{\cop_\Xcal\otimes\cop_\Ycal}[d] \ar^-{\cop_{\Xcal\otimes\Ycal}}[drrr] \\
            \Xcal\otimes\Xcal \ar_-{\swap_{\Xcal,\Xcal}}[rr] && \Xcal\otimes\Xcal, &&
            \Xcal\otimes\Xcal\otimes\Ycal\otimes\Ycal \ar_-{\id_\Xcal\otimes \swap_{\Xcal,\Ycal}\otimes\id_\Ycal}[rrr] &&& \Xcal\otimes\Ycal\otimes\Xcal\otimes\Ycal.
        }\]
\end{Def}

\begin{Def}[Deterministic morphisms, see \cite{Fri20} Def.\ 10.1]
    Let $\Ccal$ be a Markov category and $g \in \Ccal(\Wcal,\Zcal)$.
     $f \in \Ccal(\Zcal,\Xcal)$ is called:
    \begin{enumerate}
        \item \emph{deterministic} if:
    \[ (f\otimes f)\circ \cop_\Zcal = \cop_\Xcal \circ f.  \]
        \item \emph{$g$-almost-surely deterministic} if:
            \[ (((f\otimes f)\circ \cop_\Zcal)\otimes\id_\Zcal) \circ \cop_\Zcal \circ g =
            ((\cop_\Xcal \circ f)\otimes\id_\Zcal) \circ \cop_\Zcal \circ g.  \]
    \end{enumerate}
\end{Def}

\begin{Def}[See \cite{Fri21}] We say that a Markov category $\Ccal$
    \begin{enumerate}
        \item has \emph{conditionals} if for every $f \in \Ccal(\Zcal,\Xcal\otimes\Ycal)$
            there exists $g \in \Ccal(\Ycal\otimes\Zcal,\Xcal)$ such that:\footnote{\label{fn:conditional}This resembles the chain rule for conditioning: $K(X,Y|Z)=K(X|Y,Z)\otimes K(Y|Z)$.}
            \[ f = (g \otimes \id_\Ycal) \circ (((\del_\Xcal\otimes \id_\Ycal) \circ f)\otimes\id_\Zcal)\circ\cop_\Zcal\footref{fn:conditional}.\]
            Also see \cite{Kal21,Kal17} Thm.\ 1.25 (just about disintegration)
        \item is \emph{a.s.-compatibly representable} in the sense of \cite{FriG20} Section 3; see \cite{FriG20} Example 3.19, where the Markov category of Markov kernels between standard Borel spaces is shown to have this property
        \item has \emph{countable Kolmogorov powers} in the sense of \cite{Fri19} Section 3; see \cite{Fri19} Example 3.6

    \end{enumerate}
\end{Def}

\section{Appendix - The Push-forward Probability Monad}

\begin{Prf}[The push-forward probability monad]{qms:thm:PrPf-monad}
    That $(\PrPf,\delta,\Mbb)$ is a monad on $\QMS$ is \Cref{prp:P-K-R-monad},
    and that it is strong under \Cref{qms:asm:prod-iso} is
    \Cref{thm:K-P-R-strong-monad}, both applied to the third of the four
    push-forward monads there. Affinity is \Cref{qms:lem:affine}. Commutativity
    under the Fubini property is \Cref{qms:thm:commutative}.
\end{Prf}

\label{qms:sec:push-for-prob-mon}

\subsection{Fubini Theorem}

\begin{Prf}[Fubini Theorem]{qms:thm:fubini}
    For a sample space $(\Omega,\Omega^\Omega,\Bcal_\Omega,\PrPf)$ satisfying \Cref{qms:def:sample-space-act-3} and \Cref{qms:asm:fubini},
   $P(X) \in \PrPf(\Xcal)$, $Q(Y) \in \PrPf(\Ycal)$ and $h \in [0,\infty]^{\Xcal \times \Ycal}$
   we need to show the following formula:
    \begin{align}
         \iint h(x,y)\,P(X \in dx) \, Q(Y \in dy) &= \iint h(x,y)\, Q(Y \in dy) \,P(X \in dx).
    \end{align}
    Since $h \in [0,\infty]^{\Xcal \times \Ycal}$ we know that $h: \Xcal \times \Ycal \to [0,\infty]$
    is $\Bcal_{\Xcal \times \Ycal}$-$\Bcal_{[0,\infty]}$-measurable.
 From \cite{Kle20} Thm.\ 1.96 we know that we have a representation:
        \[ h(x,y) = \sum_{n \in \N} c_n\cdot\I_{D_{(n)}}(x,y),\]
        with $D_{(n)} \in \Bcal_{\Xcal \times \Ycal}$, $n \in \N$.
    Together with the linearity of the integral this reduces our problem to indicator functions
        $h(x,y) = \I_D(x,y)$ with $D \in \Bcal_{\Xcal \times \Ycal}$.
    We abbreviate:
    \begin{align}
        \tilde D&:=\lC (\omega_1,\omega_2) \in \Omega \times \Omega \st (X(\omega_1),Y(\omega_2)) \in D \rC
        \in \Bcal_{\Omega \times \Omega}.
    \end{align}

    It directly follows from the usual substitution rule and \Cref{qms:asm:fubini} that:
    \begin{align}
       &  \iint h(x,y)\,P(X \in dx) \, Q(Y \in dy)  \\
       &= \iint h(X(\omega_1),Y(\omega_2))\,P(d\omega_1) \, Q(d\omega_2)\\
       &= \iint \I_D(X(\omega_1),Y(\omega_2))\,P(d\omega_1) \, Q(d\omega_2)\\
       &= \iint \I_{\tilde D}(\omega_1,\omega_2)\,P(d\omega_1) \, Q(d\omega_2)\\
       &= \iint \I_{\tilde D^{\omega_2}}(\omega_1)\,P(d\omega_1) \, Q(d\omega_2)\\
       &= \int P(\tilde D^{\omega_2}) \, Q(d\omega_2)\\
       &= \int Q(\tilde D_{\omega_1})\,P(d\omega_1) \\
       &= \iint \I_{\tilde D_{\omega_1}}(\omega_2)\, Q(d\omega_2)\,P(d\omega_1) \\
       &= \iint \I_{\tilde D}(\omega_1,\omega_2)\, Q(d\omega_2)\,P(d\omega_1) \\
       &= \iint \I_D(X(\omega_1),Y(\omega_2))\, Q(d\omega_2)\,P(d\omega_1) \\
       &= \iint h(X(\omega_1),Y(\omega_2))\, Q(d\omega_2)\,P(d\omega_1) \\
       &= \iint h(x,y)\, Q(Y \in dy) \,P(X \in dx).
    \end{align}
This shows the claim.
\end{Prf}

\begin{Prf}[Fubini property criterion]{qms:prp:fubini-crit}
         First, note that by \Cref{qms:def:sample-space-act-3} we have that both $P_2 \otimes P_1$, $P_1 \otimes P_2 \in \PrPf(\Omega \times \Omega)$.
     Since $P_2 \otimes P_1 = P(\Phi_2,\Phi_1)$ for some $\Phi_1,\Phi_2 \in \Omega^\Omega$, $P \in \PrPf$, then also:
     \begin{align*}
         s_*(P_2 \otimes P_1) = P(\Phi_1,\Phi_2) \in \PrPf(\Omega \times \Omega).
     \end{align*}
   Now assume $D=A_1 \times A_2$ with $A_1,A_2 \in \Bcal_\Omega$. Then we get:
   \begin{align}
       s_*(P_2 \otimes P_1)(D) &= \int P_2(D_{\omega_1})\,P_1(d\omega_1) \\
                               &= \int P_2(A_2) \cdot \I_{A_1}(\omega_1) \, P_1(d\omega_1) \\
                               &= P_2(A_2) \cdot P_1(A_1) \\
                               &= P_1(A_1) \cdot P_2(A_2) \\
                               &= \int P_1(A_1) \cdot \I_{A_2}(\omega_2)\, P_2(d\omega_2) \\
                               &= \int P_1(D^{\omega_2})\, P_2(d\omega_2) \\
                               &= (P_1\otimes P_2)(D).
   \end{align}
   By Dynkin's Lemma we then get for every $D \in \Bcal_\Omega \otimes \Bcal_\Omega$:
   \begin{align}
       s_*(P_2 \otimes P_1)(D) &= (P_1\otimes P_2)(D).
   \end{align}
   So on $\Bcal_\Omega\otimes\Bcal_\Omega$ both measures agree:
   \begin{align}
       Q:= s_*(P_2 \otimes P_1)|_{\Bcal_\Omega \otimes \Bcal_\Omega} &= (P_1\otimes P_2)|_{\Bcal_\Omega \otimes \Bcal_\Omega}.
   \end{align}
   So they uniquely extend to the completion $(\Bcal_\Omega\otimes\Bcal_\Omega)_Q$.
   So we get that the extensions of both measures already agree on the completion $(\Bcal_\Omega\otimes\Bcal_\Omega)_Q$.
   By Fubini's theorem for the completion, see \cite{Fremlin} 252C, 252D, 252H, the defining formulas of those measures agree with their extensions and we get:
   \begin{align}
       s_*(P_2 \otimes P_1)|_{(\Bcal_\Omega\otimes\Bcal_\Omega)_Q} &= (P_1\otimes P_2)|_{(\Bcal_\Omega\otimes\Bcal_\Omega)_Q}.
   \end{align}
   Since $s_*(P_2 \otimes P_1)$, $P_1 \otimes P_2 \in \PrPf(\Omega \times \Omega)$ we have the inclusions:
   \begin{align}
  \Bcal_{\Omega\times\Omega} \ins (\Bcal_\Omega \otimes \Bcal_\Omega)_{\PrPf(\Omega\times\Omega)} \ins
   (\Bcal_\Omega\otimes\Bcal_\Omega)_Q.
   \end{align}
 This shows that
   \begin{align}
       s_*(P_2 \otimes P_1)|_{\Bcal_{\Omega\times\Omega}} &= (P_1\otimes P_2)|_{\Bcal_{\Omega\times\Omega}},
   \end{align}
   which implies the claim.
\end{Prf}

\subsection{Product of Markov Kernels}

\begin{Prf}[Product of Markov kernels]{qms:thm:product-PrPf}
    Let $K(X|Y,Z) \in \PrPf(\Xcal)^{\Ycal \times \Zcal}$ and $Q(Y|Z) \in \PrPf(\Ycal)^\Zcal$.
    To show that $K(X|Y,Z) \otimes Q(Y|Z) \in \PrPf(\Xcal \times \Ycal)^\Zcal$ we need to show that for every
    $Z_1 \in \Zcal^\Omega$ the composition:
    \begin{align}
        \lp K(X|Y,Z) \otimes Q(Y|Z)\rp \circ Z_1 \in \PrPf(\Xcal \times \Ycal)^\Omega.
    \end{align}
    For this first consider the composition:
    \begin{align}
        Q(Y|Z=Z_1) = Q(Y|Z) \circ Z_1 \in \PrPf(\Ycal)^\Omega.
    \end{align}
    By definition of $\PrPf(\Ycal)^\Omega$ there exists $Y_1 \in (\Ycal^\Omega)^\Omega$ and
    $P_2 \in \PrPf$ such that for all $B \in \Bcal_\Ycal$ and $\omega \in \Omega$ the following holds:
    \begin{align}
        Q(Y \in B|Z=Z_1(\omega)) & = P_2(Y_1(\omega) \in B).
    \end{align}
    We then consider the quasi-measurable map:
    \begin{align}
        Y_1 \times Z_1:\,  \Omega \times \Omega &\to \Ycal \times \Zcal, &
        (\omega_1,\omega_2) &\mapsto \lp Y_1(\omega_1)(\omega_2),Z_1(\omega_1)\rp.
    \end{align}
    Since, by assumption, we have an isomorphism of quasi-measurable spaces:
    \begin{align}
        \Theta= \Theta_1 \times \Theta_2: \; \Omega & \cong \Omega \times \Omega,
     \end{align}
    we can consider the following quasi-measurable map:
    \begin{align}
        Y_2 \times Z_2:\, \Omega \cong \Omega \times \Omega &\to \Ycal \times \Zcal, &
        \omega &\mapsto \lp Y_1(\Theta_1(\omega))(\Theta_2(\omega)),Z_1(\Theta_1(\omega))\rp.
    \end{align}
    So, $Y_2 \times Z_2 \in (\Ycal \times \Zcal)^\Omega$.
    This then shows that we have a quasi-measurable map:
    \begin{align}
        K(X|Y=Y_2,Z=Z_2):\,  \Omega &\to \PrPf(\Xcal), &  \omega &\mapsto  K(X|Y=Y_2(\omega),Z=Z_2(\omega)).
    \end{align}
    In short, $K(X|Y=Y_2,Z=Z_2) \in \PrPf(\Xcal)^\Omega$.
    By definition of $\PrPf(\Xcal)^\Omega$, there exists a $X_2 \in (\Xcal^\Omega)^\Omega$ and $P_1 \in \PrPf$
    such that for all $A \in \Bcal_\Xcal$ and $\omega\in \Omega$ we have:
    \begin{align}
        K(X \in A|Y=Y_2(\omega),Z=Z_2(\omega)) &= P_1(X_2(\omega) \in A).
    \end{align}
    We now define $X_1 \in ((\Xcal^\Omega)^\Omega)^\Omega$ via:
    \begin{align}
        X_1(\omega_1)(\omega_2)(\omega_3) &:= X_2(\Theta^{-1}(\omega_1,\omega_2))(\omega_3).
    \end{align}
    Then we get for all $A \in \Bcal_\Xcal$ and $\omega_1,\omega_2 \in \Omega$, via inserting
    $\omega:=\Theta^{-1}(\omega_1,\omega_2) \in \Omega$ into the above:
    \begin{align}
        K(X \in A|Y=Y_1(\omega_1)(\omega_2),Z=Z_1(\omega_1)) &= P_1(X_1(\omega_1)(\omega_2) \in A).
    \end{align}
    As another step, we look at $P_1 \otimes P_2 \in \PrQm(\Omega \times \Omega)$.
    By assumption, there exists $P \in \PrPf$ and $\Phi \in \Omega^\Omega$ such that:
    \begin{align}
        P_1 \otimes P_2 &= P(\Phi_1,\Phi_2) \in \PrPf(\Omega \times \Omega),
    \end{align}
    where $\Phi_i := \Theta_i \circ \Phi \in \Omega^\Omega$, $i=1,2$.

    With these notations we then define: $X_3 \in (\Xcal^\Omega)^\Omega$, $Y_3 \in (\Ycal^\Omega)^\Omega$ via:
    \begin{align}
        X_3(\omega_1)(\omega_4) &:= X_1(\omega_1)(\Phi_2(\omega_4))(\Phi_1(\omega_4)),\\
        Y_3(\omega_1)(\omega_4) &:= Y_1(\omega_1)(\Phi_2(\omega_4)).
    \end{align}
    Note that then: $X_3 \times Y_3 \in ((\Xcal \times \Ycal)^\Omega)^\Omega$.

    Finally, for $D \in \Bcal_{\Xcal \times \Ycal}$ we abbreviate:
    \begin{align}
        \tilde D &:= \lC (\omega_1,\omega_2,\omega_3) \in \Omega^3 \st (X_1(\omega_1)(\omega_2)(\omega_3),Y_1(\omega_1)(\omega_2)) \in D \rC \in \Bcal_{\Omega \times \Omega \times \Omega}, \\
        \tilde D_{\omega_1} &:= \lC (\omega_3,\omega_2)  \in \Omega \times \Omega \st (\omega_1,\omega_2,\omega_3) \in \tilde D \rC\in \Bcal_{\Omega \times \Omega},\\
        \tilde D_{\omega_1}^{\omega_2} &:= \lC \omega_3 \in \Omega \st (\omega_1,\omega_2,\omega_3) \in \tilde D \rC \in \Bcal_\Omega.
    \end{align}

    With all these preliminaries, we now claim that:
    \begin{align}
        (K(X|Y,Z) \otimes Q(Y|Z)) \circ Z_1 =  P(X_3,Y_3|\Inp) \in \PrPf(\Xcal\times\Ycal)^\Omega.
    \end{align}
    To check this let $D \in \Bcal_{\Xcal \times \Ycal}$ and $\omega_1 \in \Omega$.
    With this we get:
    \begin{align}
        & \lp K(X|Y,Z) \otimes Q(Y|Z) \rp(D|Z_1(\omega_1))\\
        &= \int K(X \in D^y|Y=y,Z=Z_1(\omega_1)) \, Q(Y \in dy|Z=Z_1(\omega_1)) \\
        &= \int K(X \in D^y|Y=y,Z=Z_1(\omega_1)) \, P_2(Y_1(\omega_1) \in dy) \\
        &= \int K(X \in D^{Y_1(\omega_1)(\omega_2)}|Y=Y_1(\omega_1)(\omega_2),Z=Z_1(\omega_1)) \, P_2(d\omega_2) \\
        &= \int P_1(X_1(\omega_1)(\omega_2) \in D^{Y_1(\omega_1)(\omega_2)}) \, P_2(d\omega_2) \\
        &= \int P_1(\lC \omega_3 \in \Omega \st (X_1(\omega_1)(\omega_2)(\omega_3),Y_1(\omega_1)(\omega_2)) \in D \rC) \, P_2(d\omega_2) \\
        &= \int P_1(\tilde D_{\omega_1}^{\omega_2}) \, P_2(d\omega_2) \\
        &=  (P_1\otimes P_2)(\tilde D_{\omega_1}) \\
        &= P\lp (\Phi_1,\Phi_2) \in \tilde D_{\omega_1} \rp \\
        &= P\lp \lC \omega_4 \in \Omega \st  (\Phi_1(\omega_4),\Phi_2(\omega_4)) \in \tilde D_{\omega_1}\rC \rp \\
        &= P\lp \lC \omega_4 \in \Omega \st  (X_1(\omega_1)(\Phi_2(\omega_4))(\Phi_1(\omega_4)),Y_1(\omega_1)(\Phi_2(\omega_4))) \in D\rC \rp \\
        &= P\lp \lC \omega_4 \in \Omega \st  (X_3(\omega_1)(\omega_4),Y_3(\omega_1)(\omega_4)) \in D\rC \rp \\
        &= P\lp (X_3(\omega_1),Y_3(\omega_1)) \in D \rp \\
        &= P\lp (X_3,Y_3) \in D| \Inp= \omega_3\rp.
\end{align}
    Note that $X_3 \times Y_3 \in ((\Xcal \times \Ycal)^\Omega)^\Omega$ and $P \in \PrPf$, and, thus:
    \begin{align}
        (K(X|Y,Z) \otimes Q(Y|Z)) \circ Z_1 =  P(X_3,Y_3|\Inp) \in \PrPf(\Xcal\times\Ycal)^\Omega.
    \end{align}
    Since this holds for all $Z_1 \in \Zcal^\Omega$ the above shows that
     that $K(X|Y,Z) \otimes Q(Y|Z) \in \PrPf(\Xcal \times \Ycal)^\Zcal$.
    This implies that $\otimes$ is a well-defined map.

     To show that $\otimes$ is also quasi-measurable we need to go through the exact same computation from above,
     with the only difference that $\Zcal$ is replaced by $\tilde \Zcal := \Zcal \times \Omega$ and $Z_1 \in \Zcal^\Omega$ by $\tilde Z_1 \in \tilde \Zcal^\Omega$.
     The claim thus follows similarly (by putting tildas on all occuring $Z$'s).
\end{Prf}

\subsection{The Push-forward Probability Monad}

\subsection{Push-forward versus Quasi-Measurable Probability Measures}

\begin{Lem}
    \label{qms:lem:PrPf-in-PrQm}
    Assume \Cref{qms:def:sample-space-act-3} \Cref{qms:eq:sample-space-act-3-b}. Then for every quasi-measurable space $(\Xcal,\Xcal^\Omega)$
    we have a well-defined and quasi-measurable inclusion map:
    \begin{align}
        \lp \PrPf(\Xcal), \PrPf(\Xcal)^\Omega \rp \inj \lp \PrQm(\Xcal), \PrQm(\Xcal)^\Omega \rp.
    \end{align}
    \begin{proof}
        The inclusion $\PrPf(\Xcal) \ins \PrQm(\Xcal)$ can be considered a special case of the inclusion
        $\PrPf(\Xcal)^\Omega \ins \PrQm(\Xcal)^\Omega$, which we show in the following.\\
        Let $P(X|\Inp) \in \PrPf(\Xcal)^\Omega$ with $P \in \PrPf$ and $X \in (\Xcal^\Omega)^\Omega$.
        Further, let $D \in \Bcal_{\Xcal \times \Omega}$ and $\Phi \in \Omega^\Omega$.
        Then we define:
        \begin{align}
            \tilde D := \lC (\tilde \omega, \omega) \in \Omega \times \Omega \st (X(\Phi(\omega))(\tilde\omega),\omega) \in D \rC.
        \end{align}
        Then clearly $\tilde D \in \Bcal_{\Omega \times \Omega}$.
        With these notations we get for $\omega \in \Omega$:
        \begin{align}
            P(X \in D^\omega|\Inp=\Phi(\omega)) &= P\lp\lC \tilde \omega \in \Omega \st (X(\Phi(\omega))(\tilde \omega),\omega) \in D  \rC\rp \\
                                                &= P\lp\lC \tilde \omega \in \Omega \st (\tilde \omega, \omega) \in \tilde D  \rC\rp \\
                                                &= P(\tilde D^\omega).
        \end{align}
        By \Cref{qms:def:sample-space-act-3} \Cref{qms:eq:sample-space-act-3-b} the latter is measurable in $\omega$.
        This shows that the former is also measurable in $\omega$. By definition of $\PrQm(\Xcal)^\Omega$ this shows:
        $P(X|\Inp) \in \PrQm(\Xcal)^\Omega$.
    \end{proof}
\end{Lem}

\begin{Rem}
    \begin{enumerate}
\item Now assume
    that we are, in addition to \Cref{qms:def:sample-space-act-3} \Cref{qms:eq:sample-space-act-3-b},  given a fixed isomorphism of quasi-measurable spaces:
    \begin{align}
        \Theta=\Theta_1 \times \Theta_2:\, \Omega & \cong \Omega \times \Omega.
    \end{align}
    Then for $P_1,P_2 \in \PrPf$ we have the equation on $\Bcal_{\Omega \times \Omega}$:
    \begin{align}
        P_1 \otimes P_2 &= P(\Theta_1,\Theta_2),
    \end{align}
    with $P:=\Theta^{-1}_*\lp P_1 \otimes P_2\rp \in \PrQm(\Omega)$.
    In such a case for \Cref{qms:def:sample-space-act-3} \Cref{qms:eq:sample-space-act-3-c} to hold it is enough to check if $P \in \PrPf(\Omega)$.

\item \Cref{qms:def:sample-space-act-3} \Cref{qms:eq:sample-space-act-3-b} and \Cref{qms:eq:sample-space-act-3-c} are both satisfied if we use $\PrPf:=\PrQm(\Omega)$ and we have an isomorphism of quasi-measurable spaces $\Omega \cong \Omega \times \Omega$.

\item If  \Cref{qms:def:sample-space-act-3} \Cref{qms:eq:sample-space-act-3-b} and \Cref{qms:eq:sample-space-act-3-c} hold a sufficient criterion for
    \Cref{qms:def:sample-space-act-3} \Cref{qms:asm:fubini} to hold is the inclusion of $\sigma$-algebras:
    \begin{align}
        \Bcal_{\Omega \times \Omega} & \ins \lp\Bcal_\Omega \otimes \Bcal_\Omega \rp_{\PrPf(\Omega \times \Omega)},
    \end{align}
    which is e.g.\ the case if $\Bcal_\Omega$ is countably separated and $\PrPf$ only consists of perfect probability measures, see \Cref{qms:cor:count-sep-perf}.
\item \Cref{qms:def:sample-space-act-3} \Cref{qms:eq:sample-space-act-3-b}, \Cref{qms:eq:sample-space-act-3-c}, \Cref{qms:asm:fubini} are all satisfied if we use $\PrPf:=\Giry(\Omega,\Bcal_\Omega)$, we take $\Omega^\Omega := \Fcal(\Bcal_\Omega)$ and we have an isomorphism of measurable spaces:
    \begin{align} (\Omega,\Bcal_\Omega) \cong (\Omega \times \Omega, \Bcal_\Omega \otimes \Bcal_\Omega). \end{align}
\item \Cref{qms:def:sample-space-act-3} \Cref{qms:eq:sample-space-act-3-b}, \Cref{qms:eq:sample-space-act-3-c}, \Cref{qms:asm:fubini} are all satisfied if we have an isomorphism of quasi-measurable spaces $\Omega \cong \Omega \times \Omega$, we use:
    \begin{align}\PrPf:=\lC P \in \Giry(\Omega,\Bcal_\Omega) \st P \text{ perfect}\rC, \end{align}
    and $(\Omega,\Bcal_\Omega)$ is countably separated and $\PrPf(\Omega)$-complete.
    \end{enumerate}
\end{Rem}

\section{Appendix - The Other Probability Monads $\PrQm$, $\PrPF$, $\PrpF$, $\Prpf$}
\label{qms:sec:app-other-monads}

The main text of this paper works with the push-forward monad $\PrPf$ of \Cref{sec:pf-monad},
accompanied by the quasi-measurable probability monad $\PrQm$ of \Cref{sec:general-monad} and, on the measure-theoretic side,
the Giry monad $\Giry$ of \Cref{qms:def:prob-meas}. This appendix develops the three that were left out
--- $\PrPF$, $\PrpF$ and $\Prpf$, the siblings of $\PrPf$ --- in full generality and in the notation of
the main text, together with the details of $\PrQm$ whose statements the main text quotes.

Why $\PrPf$ and not one of the others? Because it is the one whose admissible random variables are the
push-forwards $X_*\mu$ with the \emph{random variable} free to depend on the sample point and the measure
held fixed --- which is exactly the reading ``the distribution of an admissible random variable'' that the
rest of the paper needs. $\PrPF$ lets both vary and is more general than that reading requires; $\PrpF$,
the monad of \cite{Heu17}, lets only the measure vary; $\Prpf$ lets neither. On quasi-universal spaces
the distinction evaporates --- all four coincide, see \Cref{thm:qus-S=P-in-R} --- so the choice costs
nothing there, and off quasi-universal spaces $\PrPf$ is the one that keeps the interpretation.

All five spaces of probability measures sit inside $\Giry(\Xcal)=\Giry(\Xcal,\Bcal_\Xcal)$ and share the
same underlying set as soon as they are of push-forward type; what distinguishes them is which
\emph{admissible random variables} they carry, i.e.\ which of the two ingredients of a push-forward
$X_*\mu$ --- the random variable $X$ and the measure $\mu$ --- is allowed to depend on the sample point.
Writing $\Xcal^\Omega$ for the admissible random variables of $\Xcal$, this reads as follows.

\begin{table}[ht]
    \centering
    \small
    \renewcommand{\arraystretch}{1.45}
    \begin{tabular}{@{}llll@{}}
        \hline
        monad & underlying set & admissible random variables & $X$ / $\mu$ vary?\\
        \hline
        $\Giry$
            & $\lC \mu:\,\Bcal_\Xcal \to [0,1] \rC$
            & $\Fcal\lp\Bcal_{\Giry(\Xcal)}\rp$
            & --- \\
        $\PrQm$
            & $\Giry(\Xcal) \cap [0,1]^{\Bcal_\Xcal}$
            & $\Sets(\Omega,\PrQm(\Xcal)) \cap \lp[0,1]^{\Bcal_\Xcal}\rp^\Omega$
            & --- \\
        $\PrPF$
            & $\pf\lp \Xcal^\Omega \times \PrQm(\Omega)\rp$
            & $\pf\lp \lp\Xcal^\Omega\rp^\Omega \times \PrQm(\Omega)^\Omega \rp$
            & both \\
        $\PrpF$
            & $\pf\lp \Xcal^\Omega \times \PrQm(\Omega)\rp$
            & $\pf\lp \lp\Xcal^\Omega\rp^\one \times \PrQm(\Omega)^\Omega \rp$
            & only $\mu$ \\
        $\PrPf$
            & $\pf\lp \Xcal^\Omega \times \PrPf\rp$
            & $\pf\lp \lp\Xcal^\Omega\rp^\Omega \times \PrPf^\one \rp$
            & only $X$ \\
        $\Prpf$
            & $\pf\lp \Xcal^\Omega \times \PrPf\rp$
            & $\pf\lp \lp\Xcal^\Omega\rp^\one \circ \lp\Omega^\Omega\rp^\Omega \times \PrPf^\one\rp$
            & neither \\
        \hline
    \end{tabular}
    \caption{The five probability monads on the category of quasi-measurable spaces $\QMS$. The
    pool of measures on $\Omega$ is the act-3 datum $\PrPf$ for $\PrPf$ and $\Prpf$, which hold it
    fixed, and $\PrQm(\Omega)$ for $\PrPF$ and $\PrpF$, which let it vary and therefore need a
    quasi-measurable structure on it; see \Cref{rem:which-pool}. Together with
    the Giry monad $\Giry$. The four push-forward monads $\PrPF$, $\PrpF$, $\PrPf$, $\Prpf$ share the same
    underlying set and differ only in which admissible random variables they carry: whether the random
    variable $X$, the measure $\mu$, both or neither may depend on the sample point. Here
    $\pf:\, \Xcal^\Omega \times \Mcal \to \PrQm(\Xcal)$, $(X,\mu) \mapsto X_*\mu$, is the
    push-forward map for a pool $\Mcal$ of measures on $\Omega$, and a superscript $\one$ marks an
    ingredient held constant; in the $\Prpf$ row the symbol $\circ$ indicates that the constant
    random variable is first reparameterised by a sample-point-dependent
    $\phi \in \lp\Omega^\Omega\rp^\Omega$, so that $\pf$ is still applied to two arguments. The
    inclusions between the four push-forward rows are inclusions of quasi-measurable spaces, i.e.\
    of underlying sets \emph{and} of admissible random variables. See \Cref{def:prob-space-K}.}
    \label{qms:tab:monads}
\end{table}

They are ordered by inclusion --- of underlying sets and of admissible random variables alike ---
see \Cref{lem:S-P-R-K-incl}:
\begin{align}
    \Prpf(\Xcal) \ins \PrpF(\Xcal) \ins \PrPF(\Xcal) \ins \PrQm(\Xcal) \ins \Giry(\Xcal),
    \qquad \Prpf(\Xcal) \ins \PrPf(\Xcal) \ins \PrPF(\Xcal),
\end{align}
and they collapse under mild hypotheses on the sample space, see \Cref{lem:K-P-S-R-eq}: if
$\Omega \cong \Omega \times \Omega$ then $\Prpf=\PrPf$ and $\PrpF=\PrPF$, while surjectivity of the
push-forward map $\pf:\,\lp\Omega^\Omega\rp^\Omega \times \PrQm(\Omega) \srj \PrQm(\Omega)^\Omega$ gives
$\Prpf = \PrpF$ and $\PrPf = \PrPF$. On quasi-universal spaces both hypotheses hold and all four coincide,
see \Cref{thm:qus-S=P-in-R}; on \emph{standard} quasi-universal spaces even $\PrQm$ and $\Giry$ join them,
see \Cref{thm:univ-qus-GPQKR}. This is the precise sense in which the main text loses nothing by working
with $\PrPf$ alone.

In this section we will introduce several different \emph{spaces of probability measures} for a quasi-measurable space: $\Giry(\Xcal)$, $\PrQm(\Xcal)$, $\PrPF(\Xcal)$, $\PrpF(\Xcal)$, $\PrPf(\Xcal)$, $\Prpf(\Xcal)$.
These notions will in general all be different. To study their functorial properties we need to review the notion of a \emph{monad}.

\subsection{The Quasi-Measurable Probability Monad $\PrQm$}

In this subsection we will highlight a novel and strong probability monad $\PrQm$, the \emph{quasi-measurable probability monad}.
In contrast to the Giry construction $\Giry$ the probability monad $\PrQm$ is well-behaved w.r.t.\ the
categorical constructions inside the category of quasi-measurable spaces $\QMS$, e.g.\ w.r.t.\ products and $\sigma$-algebras.

\begin{Def}[The spaces of probability measures for quasi-measurable spaces]
    Let $(\Xcal,\Xcal^\Omega)$ be a quasi-measurable space and $\Bcal_\Xcal:=\Bcal(\Xcal^\Omega)$ its induced $\sigma$-algebra, which we will also consider
    as the quasi-measurable space $(\Bcal_\Xcal,(\Bcal_\Xcal)^\Omega)$.
    Then we define:
    \begin{align*}
        \Giry(\Xcal,\Xcal^\Omega) &:= \Giry(\Xcal,\Bcal_\Xcal)   =\lC \mu:\, \Bcal_\Xcal \to [0,1] \text{ probability measure}\rC, \\
        \Giry(\Xcal,\Xcal^\Omega)^\Omega &:= \Fcal\lp \Bcal_{\Giry(\Xcal,\Bcal_\Xcal)}\rp\\
                                         &=\Big\{ \kappa:\,\Omega \to \Giry(\Xcal,\Xcal^\Omega)\,\Big|\, \forall A \in \Bcal_\Xcal.\\
                                         &\qquad\qquad\qquad (\omega \mapsto \kappa(\omega)(A)) \in [0,1]^\Omega \Big\},\\
        \PrQm(\Xcal,\Xcal^\Omega) &:= \Giry(\Xcal,\Xcal^\Omega) \cap [0,1]^{\Bcal_\Xcal} \\
                     &= \Big\{ \mu:\, \Bcal_\Xcal \to [0,1] \text{ probability measure}\,\Big|\, \forall D \in \Bcal_{\Omega \times \Xcal}.\\
                     &\qquad\qquad\qquad (\omega \mapsto \mu(D_\omega)) \in [0,1]^\Omega \Big\},\\
        \PrQm(\Xcal,\Xcal^\Omega)^\Omega &:= \Sets(\Omega,\PrQm(\Xcal,\Xcal^\Omega)) \cap ([0,1]^{\Bcal_\Xcal})^\Omega\\
                            &= \Big\{ \kappa:\,\Omega \to \PrQm(\Xcal,\Xcal^\Omega)\,\Big|\, \forall \Phi \in \Omega^\Omega.\\
                            &\qquad\qquad\qquad \forall D \in \Bcal_{\Omega \times \Xcal}.\,
                            (\omega \mapsto \kappa(\Phi(\omega))(D_\omega)) \in [0,1]^\Omega\Big\}.
      \end{align*}
\end{Def}

\begin{Rem}
    \label{rem:Q-eval}
    Let $(\Xcal,\Xcal^\Omega)$ be a quasi-measurable space.
    \begin{enumerate}
        \item  $(\Giry(\Xcal),\Giry(\Xcal)^\Omega)$ is a quasi-measurable space as it is induced from
            the underlying measurable space endowed with the smallest $\sigma$-algebra
            that makes all evaluation maps measurable, see Example \ref{qms:eg:meas-qms-space}.
        \item  $(\PrQm(\Xcal),\PrQm(\Xcal)^\Omega)$  is a quasi-measurable space with the subspace structure of
            $[0,1]^{\Bcal_\Xcal}$, also see Lemma \ref{lem:subspace-qms}.
        \item It is then clear that the evaluation map:
            \[ \ev:\;\PrQm(\Xcal) \times \Bcal_\Xcal \xrightarrow{\incl \times \id} [0,1]^{\Bcal_\Xcal} \times \Bcal_\Xcal \xrightarrow{\ev}  [0,1],\quad (\mu,A) \mapsto \mu(A),  \]
            is quasi-measurable.
        \item Note that this might not be true for $\Giry(\Xcal)$ in place of $\PrQm(\Xcal)$. This is the reason that we will focus on $\PrQm(\Xcal)$ in the framework of quasi-measurable spaces.
        \item The inclusion map $ (\PrQm(\Xcal),\PrQm(\Xcal)^\Omega) \to (\Giry(\Xcal),\Giry(\Xcal)^\Omega)$ is quasi-measurable.
        \item[] Indeed, for $\kappa \in \PrQm(\Xcal)^\Omega$ and $A \in \Bcal_\Xcal$ and $\omega \in \Omega$ we have that $\Omega \times A \in \Bcal_{\Omega \times \Xcal}$, that $(\Omega \times A)_\omega=A$ and the map:
            \[ \omega \mapsto \kappa(\omega)(A) = \kappa(\id_\Omega(\omega))((\Omega \times A)_\omega) \]
            is measurable.
\item If for some reason we have that $\Bcal_{\Omega \times \Xcal} = \Bcal_\Omega \otimes \Bcal_\Xcal$ then we have the equality:
            \[(\PrQm(\Xcal),\PrQm(\Xcal)^\Omega) = (\Giry(\Xcal),\Giry(\Xcal)^\Omega).\]
        \item[] Indeed, for $A \in \Bcal_\Xcal$ and $B \in \Bcal_\Omega$ the map:
            \[ \omega \mapsto \kappa(\Phi(\omega))((B \times A)_\omega)=\kappa(\Phi(\omega))(A) \cdot \I_B(\omega),\]
            is measurable if $\omega \mapsto \kappa(\omega)(A)$ is.
        \item The mentioned condition holds, for instance, for $\Omega = \Xcal = \R$ with the Borel-structure.
        \item Similarly, if $\Bcal_\Omega=(\Bcal_\Omega)_\Giry$ and $\Bcal_{\Omega \times \Xcal} = (\Bcal_\Omega \otimes \Bcal_\Xcal)_\Giry$
            then we also have equality:
            \[(\PrQm(\Xcal),\PrQm(\Xcal)^\Omega) = (\Giry(\Xcal),\Giry(\Xcal)^\Omega).\]
            This holds for instance for $\Omega = \Xcal = \R$ with the universally measurable structure as we will see later in Lemma \ref{lem:univ-qus-G-Q}.
    \end{enumerate}
\end{Rem}

\begin{Lem}
    \label{lem:Q-int-map-qm}
    Let $(\Xcal,\Xcal^\Omega)$ be a quasi-measurable space.
    Then the integration map:
    \[  [0,\infty]^\Xcal \times \PrQm(\Xcal) \to [0,\infty],\quad (f,\mu) \mapsto \int f \,d\mu,   \]
    is quasi-measurable.
    \begin{proof}
        Let $g \in [0,\infty]^{\Omega \times \Xcal}$ and $\nu \in \PrQm(\Xcal)^\Omega$.
        Then we need to show that the map:
        \[\omega \mapsto \int g(\omega,x)\,\nu(\omega)(dx)\]
        is $\Bcal_\Omega$-$\Bcal_{[0,\infty]}$-measurable, where $\Bcal_{[0,\infty]}$ is the Borel $\sigma$-algebra of $[0,\infty]$.
        Since $g$ is quasi-measurable it is $\Bcal_{\Omega \times \Xcal}$-$\Bcal([0,\infty]^\Omega)$-measurable.
        Since the Borel $\sigma$-algebra $\Bcal_{[0,\infty]} \ins \Bcal([0,\infty]^\Omega)$ we see that $g$ is
        $\Bcal_{\Omega \times \Xcal}$-$\Bcal_{[0,\infty]}$-measurable.
        By Theorem 1.96 in \cite{Kle20} there are (at most) countably many measurable sets
        $E_n \in \Bcal_{\Omega \times \Xcal}$ and $a_n \in [0,\infty]$ for $n \in \N$
        such that:
        \[ g(\omega,x) = \sum_{ n \in \N} a_n \cdot \I_{E_n}(\omega,x).\]
        By definition of $\PrQm(\Xcal)^\Omega$ we already know that:
        \[ \omega \mapsto \int \I_{E_n}(\omega,x)\,\kappa(\omega)(dx) = \kappa(\omega)(E_{n,\omega}) \]
        is $\Bcal_\Omega$-$\Bcal_{[0,\infty]}$-measurable for every $n \in \N$. So we immediately get that the map:
        \[ \omega \mapsto \int g(\omega,x)\,\kappa(\omega)(dx)= \sum_{n \in \N} a_n \cdot \int \I_{E_n}(\omega,x)\,\kappa(\omega)(dx),
        \]
        is also $\Bcal_\Omega$-$\Bcal_{[0,\infty]}$-measurable. This shows the claim.
    \end{proof}
\end{Lem}

\begin{Lem}
    \label{lem:Q-int-int-qm}
    Let $(\Xcal,\Xcal^\Omega)$, $(\Ycal,\Ycal^\Omega)$, $(\Zcal,\Zcal^\Omega)$ be quasi-measurable spaces.
    Then the map:
    \begin{align*}
        [0,\infty]^{\Xcal \times \Ycal \times \Zcal} \times \PrQm(\Xcal)^{\Ycal \times \Zcal} \times \PrQm(\Ycal)^\Zcal \times  \Zcal &\to [0,\infty],\\
        (g,\mu,\nu,z) &\mapsto \int \int g(x,y,z)\,\mu(y,z)(dx)\, \nu(z)(dy),
    \end{align*}
    is a well-defined quasi-measurable map.
    \begin{proof}
    By Lemma \ref{lem:Q-int-map-qm} the  map:
    \[  [0,\infty]^\Xcal \times \PrQm(\Xcal) \xrightarrow{\int} [0,\infty],\qquad (f,\mu) \mapsto \int f d\mu\]
    is quasi-measurable. By Lemma \ref{lem:exp-qm} then also the exponentiated map;
    \[  [0,\infty]^{\Xcal\times\Ycal} \times \PrQm(\Xcal)^\Ycal \xrightarrow{\int} [0,\infty]^\Ycal,\]
    is quasi-measurable, where we used the universal property of the product $(\Wcal \times \Ucal)^\Ycal = \Wcal^\Ycal \times \Ucal^\Ycal$ on the left.
    Again, by Lemma \ref{lem:Q-int-map-qm} also the integration map:
    \[ [0,\infty]^\Ycal \times \PrQm(\Ycal) \to [0,\infty],\quad (h,\nu) \mapsto \int h\,d\nu, \]
    is quasi-measurable. So the composition of the two constructions above:
    \[ \lp [0,\infty]^{\Xcal\times\Ycal} \times \PrQm(\Xcal)^\Ycal \rp \times \PrQm(\Ycal) \xrightarrow{\int \times \id}
    [0,\infty]^\Ycal \times \PrQm(\Ycal) \xrightarrow{\int} [0,\infty],   \]
    is quasi-measurable as well. Exponentiating with $\Zcal$, again by Lemma \ref{lem:exp-qm}, then gives the quasi-measurable map:
    \[  [0,\infty]^{\Xcal\times\Ycal\times\Zcal} \times \PrQm(\Xcal)^{\Ycal \times \Zcal}  \times \PrQm(\Ycal)^\Zcal  \to [0,\infty]^\Zcal,\]
     \[ (g,\mu,\nu) \mapsto \lp z \mapsto \int \int g(x,y,z)\,\mu(y,z)(dx)\, \nu(z)(dy) \rp, \]
     which induces the quasi-measurable adjoint evaluation map:
    \[  [0,\infty]^{\Xcal\times\Ycal\times\Zcal} \times \PrQm(\Xcal)^{\Ycal \times \Zcal}  \times \PrQm(\Ycal)^\Zcal \times \Zcal  \to [0,\infty],\]
    \[ (g,\mu,\nu,z) \mapsto  \int \int g(x,y,z)\,\mu(y,z)(dx)\, \nu(z)(dy). \]
    This shows the claim.
    \end{proof}
\end{Lem}

\begin{Lem}
    \label{lem:Q-prod-ev-qm}
    Let $(\Xcal,\Xcal^\Omega)$, $(\Ycal,\Ycal^\Omega)$, $(\Zcal,\Zcal^\Omega)$ be quasi-measurable spaces.
    Then the map:
    \[ \PrQm(\Xcal)^{\Ycal \times \Zcal} \times \PrQm(\Ycal)^\Zcal \times \Bcal_{\Xcal \times \Ycal \times \Zcal} \times \Zcal \to [0,1],\qquad
    (\mu,\nu,D,z) \mapsto \int  \mu(y,z)(D^{y,z})\, \nu(z)(dy), \]
    is a well-defined quasi-measurable map.
    \begin{proof}
        This directly follows from Lemma \ref{lem:Q-int-int-qm} together with the quasi-measurable map:
        \[ \I:\; \Bcal_{\Xcal \times \Ycal \times \Zcal} \to [0,1]^{\Xcal \times \Ycal \times \Zcal},\qquad D \mapsto \I_D.  \]
    \end{proof}
\end{Lem}

\begin{Prp}
    \label{prp:Q-product-kernels}
    Let $(\Xcal,\Xcal^\Omega)$, $(\Ycal,\Ycal^\Omega)$, $(\Zcal,\Zcal^\Omega)$ be quasi-measurable spaces.
    Then the map:
    \[ \otimes:\; \PrQm(\Xcal)^{\Ycal \times \Zcal} \times \PrQm(\Ycal)^\Zcal \to \PrQm(\Xcal \times \Ycal)^\Zcal,\quad
    (\mu \otimes \nu)(z)(D)  := \int  \mu(y,z)(D^y)\, \nu(z)(dy), \]
    is a well-defined quasi-measurable map.
    \begin{proof}
        By Lemma \ref{lem:Q-prod-ev-qm} we get a quasi-measurable map:
    \[ \otimes:\; \PrQm(\Xcal)^{\Ycal \times \Zcal} \times \PrQm(\Ycal)^\Zcal \to \lp [0,1]^{\Bcal_{\Xcal \times \Ycal}}\rp^\Zcal.\]
    By Lemma \ref{lem:subspace-qms} we only need to check that the image evaluated at each $z$ always is a probability measure
    on $\Bcal_{\Xcal \times \Ycal}$, which is obvious.
    \end{proof}
\end{Prp}

\begin{Lem}
    \label{lem:push-forward-qm}
    Let $(\Xcal,\Xcal^\Omega)$, $(\Ycal,\Ycal^\Omega)$ be quasi-measurable spaces.
     Then the push-forward map is quasi-measurable:
     \[ \pf:\; \Ycal^\Xcal \times \PrQm(\Xcal) \to \PrQm(\Ycal),\qquad (f,\mu) \mapsto f_*\mu.  \]
    \begin{proof}
      The following map is quasi-measurable:
    \[ \PrQm(\Xcal) \times \Ycal^\Xcal \times \Bcal_\Ycal \to [0,1], \qquad
     (\mu, f,B) \mapsto \mu(f^{-1}(B)).\]
        Indeed, it is the composition of the following quasi-measurable map from Lemma \ref{lem:s-alg-pull-back-qm}:
        \[  \Ycal^\Xcal \times \Bcal_\Ycal  \to   \Bcal_\Xcal,\qquad
        (f,B) \mapsto f^{-1}(B),\]
        and the quasi-measurable evaluation map:
        \[ \PrQm(\Xcal) \times \Bcal_\Xcal \to [0,1],\qquad (\mu,A) \mapsto \mu(A). \]
        By adjunction we get the quasi-measurable map:
        \[ \PrQm(\Xcal) \times \Ycal^\Xcal \to [0,1]^{\Bcal_\Ycal}, \qquad
        (\mu, f) \mapsto f_*\mu=(B \mapsto \mu(f^{-1}(B))).\]
        By Lemma \ref{lem:subspace-qms} is left to show that $f_*\mu$ is a probability measure, which is clear.
        So the claims are shown.
    \end{proof}
\end{Lem}

\begin{Rem}
    \begin{enumerate}
        \item By Theorem \ref{qms:cor:qms-cartesian-closed} the triple $(\QMS, \times, \one)$ is a cartesian closed (symmetric strict) monoidal category, where $\one$ is the one-point space,
            which is also a terminal object of $\QMS$.
        \item $\PrQm:\; \QMS \to \QMS$ is a functor by Lemma \ref{lem:push-forward-qm}:
    \[ \lp f:\, \Xcal \to \Ycal \rp \mapsto \lp \PrQm(f):\; \PrQm(\Xcal) \to \PrQm(\Ycal),\; \mu \mapsto f_*\mu \rp.  \]
    \end{enumerate}
\end{Rem}

\begin{Lem}
    \label{lem:Q-delta-nat-trafo}
    Let $(\Xcal,\Xcal^\Omega)$ be a quasi-measurable space.
    Then the map:
    \[ \delta:\; \Xcal \to \PrQm(\Xcal),\qquad x \mapsto \delta_x = (A \mapsto \I_A(x)),\]
    is a well-defined quasi-measurable map.\\
    Furthermore, if we have a quasi-measurable map:
    \[f:\; (\Xcal,\Xcal^\Omega) \to (\Ycal,\Ycal^\Omega), \]
    then we get a commutative diagram of quasi-measurable maps:
    \[\xymatrix{
            \Xcal \ar^-{f}[r] \ar^-{\delta}[d] & \Ycal \ar^-{\delta}[d]\\
            \PrQm(\Xcal) \ar^{f_*}[r] & \PrQm(\Ycal).
    }\]
    In other words, $\delta:\, \id_{\QMS} \to \PrQm$ is a natural transformation of endo-functors of $\QMS$.
    \begin{proof}
        To show that $\delta$ is well-defined we need to show that $\delta_x \in \PrQm(\Xcal)$ for $x\in \Xcal$.
        For this let $D \in \Bcal_{\Omega \times \Xcal}$. Then the map:
        \[ \omega \mapsto \delta_x(D_\omega)=\I_D(\omega,x)=\I_{D^x}(\omega),\]
        is measurable since $D^x \in \Bcal_\Omega$. Thus we have $\delta_x \in \PrQm(\Xcal)$.\\
        To show that $\delta$ is quasi-measurable let $ X \in \Xcal^\Omega$ and $\Phi \in \Omega^\Omega$ and again $D \in \Bcal_{\Omega \times \Xcal}$.
        Then we see that the map:
        \[ \omega \mapsto \delta_{ X \circ \Phi(\omega)}(D_\omega) = \I_{D}(\omega, X(\Phi(\omega)))
        = \I_{(\id_\Omega, X \circ \Phi)^{-1}(D)}(\omega),\]
        is measurable as $ X \circ \Phi \in \Xcal^\Omega$. It follows that $\delta \circ  X \in \PrQm(\Xcal)^\Omega$ and
        thus $\delta:\,\Xcal \to \PrQm(\Xcal)$ quasi-measurable.\\
        To check the commutativity of the diagram observe that for $B \in \Bcal_\Ycal$:
        \[ f_*\delta_x(B) = \delta_x(f^{-1}(B))= \I_{f^{-1}(B)}(x) = \I_B(f(x)) = \delta_{f(x)}(B).\]
    \end{proof}
\end{Lem}

\begin{Lem}
    \label{lem:Q-M-nat-trafo}
        Let $(\Xcal,\Xcal^\Omega)$ be a quasi-measurable space.
        Then the map:
        \[\Mbb_\Xcal:\; \PrQm(\PrQm(\Xcal)) \to \PrQm(\Xcal),\qquad \pi \mapsto \Mbb_\Xcal(\pi)=\lp A \mapsto \int \ev(\mu,A)\, \pi(d\mu)\rp,\]
        is a well-defined quasi-measurable map.\\
    Furthermore, if we have a quasi-measurable map:
    \[f:\; (\Xcal,\Xcal^\Omega) \to (\Ycal,\Ycal^\Omega), \]
    then we get a commutative diagram of quasi-measurable maps:
    \[\xymatrix{
            \PrQm(\PrQm(\Xcal)) \ar^-{\PrQm(\PrQm(f))}[r] \ar_-{\Mbb_\Xcal}[d] & \PrQm(\PrQm(\Ycal)) \ar^-{\Mbb_\Ycal}[d]\\
            \PrQm(\Xcal) \ar_{\PrQm(f)}[r] & \PrQm(\Ycal).
    }\]
    In other words, $\Mbb:\, \PrQm^2 \to \PrQm$ is a natural transformation of endo-functors of $\QMS$.
    \begin{proof}
        It is easy to see that $\Mbb_\Xcal(\pi)$ is a probability measure on $\Bcal_\Xcal$.
        To check that $\Mbb_\Xcal$ is quasi-measurable let $e$ be the map:
        \[ e:\; \Bcal_\Xcal \to [0,1]^{\PrQm(\Xcal)},\qquad e(A)(\mu):=\ev(\mu,A),\]
        which is quasi-measurable as the adjoint to the quasi-measurable evaluation map (see Remark \ref{rem:Q-eval}):
        \[ \ev:\; \PrQm(\Xcal) \times \Bcal_\Xcal \to [0,1]. \]
        Since by Lemma \ref{lem:Q-int-map-qm} also the integration pairing is quasi-measurable the following composition of quasi-measurable maps is also quasi-measurable:
        \[ \int \circ (\id \times e):\; \PrQm(\PrQm(\Xcal)) \times \Bcal_\Xcal \to [0,1],\qquad (\pi,A) \mapsto \int \ev(\mu,A)\, \pi(d\mu).  \]
        Taking the adjoint we get the quasi-measurable map:
        \[\Mbb_\Xcal:\; \PrQm(\PrQm(\Xcal)) \to [0,1]^{\Bcal_\Xcal}, \]
        whose image lies inside $[0,1]^{\Bcal_\Xcal} \cap \Giry(\Xcal) = \PrQm(\Xcal)$.\\
        To check the commutative diagram let $B \in \Bcal_\Ycal$. Then we have with $f_*:=\PrQm(f)$ and $f_{**}:=\PrQm^2(f):=\PrQm(\PrQm(f))$:
        \begin{align*}  \Mbb_\Ycal(f_{**}\pi)(B) &= \int \ev(\nu,B)\, (f_{**}\pi)(d\nu) \\
                               &= \int \ev(f_*\mu,B)\, \pi(d\mu) \\
                               &= \int \ev(\mu,f^{-1}(B))\, \pi(d\mu) \\
                               &= \Mbb_\Xcal(\pi)(f^{-1}(B))\\
                               &= f_*\Mbb_\Xcal(\pi)(B).
        \end{align*}
        This shows the claim.
    \end{proof}
\end{Lem}

\begin{Prp}
    \label{prp:Q-monad}
    The triple $(\PrQm,\delta,\Mbb)$ defines a monad on $\QMS$.
    \begin{proof}
        By Lemma \ref{lem:push-forward-qm} we know that $\PrQm:\, \QMS \to \QMS$ is a functor.
        By Lemma \ref{lem:Q-delta-nat-trafo} we know that $\delta:\, \id_\QMS \to \PrQm$ is a natural transformation.
        By Lemma \ref{lem:Q-M-nat-trafo} we know that $\Mbb:\,\PrQm^2 \to \PrQm$ is a natural transformation.
        So we are left to check the following coherence conditions:
        \begin{enumerate}
        \item $\Mbb \circ \PrQm\Mbb = \Mbb \circ \Mbb \PrQm$ as natural transformations $\PrQm^3 \to \PrQm$,
        \item $\Mbb \circ \PrQm\delta = \Mbb \circ \delta \PrQm = \id_\PrQm$ as natural transformations $\PrQm \to \PrQm$.
        \end{enumerate}
        For simple notations note that:
        \[  \int \I_A(x)\,\Mbb_\Xcal(\pi)(dx) = \Mbb_\Xcal(\pi)(A) = \int\int \I_A(x)\,\mu(dx)\,\pi(d\mu),\]
        or in short: $\Mbb_\Xcal(\pi)(dx) = \int \mu(dx)\,\pi(d\mu)$.\\
        For the first let $\tau \in \PrQm(\PrQm(\PrQm(\Xcal)))$. Then we have for $A \in \Bcal_\Xcal$:
        \begin{align*}
            (\Mbb_\Xcal \circ \PrQm(\Mbb_\Xcal)(\tau))(A) &= \int\int \I_A(x)\,\nu(dx)\,((\Mbb_\Xcal)_*\tau)(d\nu) \\
                                                          &= \int\int \I_A(x)\,\Mbb_\Xcal(\pi)(dx)\,\tau(d\pi) \\
                                                          &= \int\int\int \I_A(x)\,\mu(dx)\, \pi(d\mu)\, \tau(d\pi)\\
                                                          &= \int\int \I_A(x)\,\mu(dx)\, (\Mbb_{\PrQm(\Xcal)}(\tau))(d\mu)  \\
                                                          &=(\Mbb_\Xcal \circ \Mbb_{\PrQm(\Xcal)}(\tau))(A).
        \end{align*}
        It is easily seen that $\Mbb \circ \PrQm\Mbb$ is a natural transformation: $\PrQm^3\to\PrQm$.
        Indeed, for quasi-measurable $f:\,\Xcal \to \Ycal$ we can use the above and the functoriality rules from
        Lemmas \ref{lem:Q-delta-nat-trafo} and \ref{lem:Q-M-nat-trafo}:
        \begin{align*}
            \Mbb_\Ycal \circ \PrQm(\Mbb_\Ycal) \circ \PrQm^3(f)
            &= \Mbb_\Ycal \circ \PrQm(\Mbb_\Ycal \circ \PrQm^2(f))  \\
                                                                &= \Mbb_\Ycal \circ \PrQm( \PrQm(f) \circ \Mbb_\Xcal) \\
                                                                &= \Mbb_\Ycal \circ \PrQm^2(f) \circ \PrQm(\Mbb_\Xcal) \\
                                                                &= \PrQm(f) \circ \Mbb_\Xcal \circ \PrQm(\Mbb_\Xcal).
        \end{align*}
        To show the coherence conditions 2) let $\nu \in \PrQm(\Xcal)$. Then we get for $A \in \Bcal_\Xcal$:
        \begin{align*}
            \Mbb_\Xcal(\delta_\nu)(A) &= \int \ev(\mu,A)\, \delta_\nu(d\mu)\\
                                      &= \ev(\nu,A) \\
                                      &= \nu(A),\\
            \Mbb_\Xcal(\PrQm(\delta)(\nu))(A) &= \int \ev(\mu,A)\, \PrQm(\delta)(\nu)(d\mu)\\
                                              &= \int \ev(\mu,A)\, (\delta_*\nu)(d\mu)\\
                                              &= \int \ev(\delta_x,A)\,\nu(dx)\\
                                              &=\int \I_A(x)\,\nu(dx)\\
                                              &=\nu(A).
        \end{align*}
        Note that in the first set of equations $\delta$ is the map:
        \[\delta:\;\PrQm(\Xcal) \to \PrQm(\PrQm(\Xcal)),\quad \nu \mapsto \delta_\nu.\]
        In the second set of equations we use the usual:
        \[\delta:\; \Xcal \to \PrQm(\Xcal),\quad x \mapsto \delta_x,\]
        where then:
        \[\PrQm(\delta)=\delta_*:\,\PrQm(\Xcal) \to \PrQm(\PrQm(\Xcal)),\quad \nu \mapsto \delta_*\nu,\]
            is its push-forward, where we then used the usual substitution rule.
    \end{proof}
\end{Prp}

\begin{Lem}
    \label{lem:Q-monad-strength}
    Let $(\Xcal,\Xcal^\Omega)$ and $(\Ycal,\Ycal^\Omega)$ be two quasi-measurable spaces.
    Then the following map:
    \[ \tau_{\Xcal,\Ycal}:\; \Xcal \times \PrQm(\Ycal) \to \PrQm(\Xcal \times \Ycal),\qquad
         (x,\mu) \mapsto \delta_x\otimes\mu= (D \mapsto \mu(D_x)),\]
    is a well-defined quasi-measurable map.\\
    Furthermore, if $f:\;\Xcal \to \Xcal'$ and $g:\;\Ycal \to \Ycal'$ are quasi-measurable we get a commutative diagram of quasi-measurable spaces:
    \[\xymatrix{
            \Xcal \times \PrQm(\Ycal) \ar^-{\tau_{\Xcal,\Ycal}}[r] \ar_-{f \times \PrQm(g)}[d] & \PrQm(\Xcal \times \Ycal) \ar^-{\PrQm(f \times g)}[d] \\
            \Xcal' \times \PrQm(\Ycal') \ar_-{\tau_{\Xcal',\Ycal'}}[r] & \PrQm(\Xcal' \times \Ycal').
    }\]
    \begin{proof}
         It is easily seen that $\tau(x,\mu)$ is a probability measure, as it is the push-forward of $\mu$ along
         the quasi-measurable slice map:
         \[\iota_x:\; \Ycal \to \Xcal \times \Ycal,\quad y \mapsto (x,y).\]
         $\tau$ is quasi-measurable as the adjoint to the compositions of the following quasi-measurable maps:
         \[ \Xcal \times \Bcal_{\Xcal\times\Ycal} \times \PrQm(\Ycal) \to \lp \Xcal \times \Ycal\rp^\Ycal \times \Bcal_{\Xcal\times\Ycal} \times \PrQm(\Ycal) \to \Bcal_\Ycal \times \PrQm(\Ycal) \xrightarrow{\ev} [0,1],  \]
         which is composed of the quasi-measurable map:
         \[ \Xcal \to \lp\Xcal \times \Ycal\rp^\Ycal,\qquad x \mapsto \iota_x,\]
         together with the quasi-measurable map (see Lemma \ref{lem:s-alg-pull-back-qm}):
         \[ \lp \Xcal \times \Ycal \rp^\Ycal \times \Bcal_{\Xcal \times \Ycal} \to \Bcal_\Ycal,\qquad (f,D) \mapsto f^{-1}(D), \]
         and the quasi-measurable evaluation map (see Remark \ref{rem:Q-eval}).
         It follows that
         \[ \tau:\; \Xcal \times \PrQm(\Ycal) \to [0,1]^{\Bcal_{\Xcal \times \Ycal}},\quad (x,\mu) \mapsto (D \mapsto \mu(D_x)),   \]
         is quasi-measurable and with image in $\PrQm(\Xcal \times \Ycal)$.
         So it is a well-defined quasi-measurable map until there (see Lemma \ref{lem:subspace-qms}).\\
         To check that the diagram is commutative let $x \in \Xcal$ and $\mu \in \PrQm(\Ycal)$ and $D' \in \Bcal_{\Xcal'\times\Ycal'}$.
         Then we get:
         \begin{align*}
             \tau_{\Xcal',\Ycal'}(f(x),\PrQm(g)(\mu))(D') & = \PrQm(g)(\mu)(D'_{f(x)})\\
                 &=\mu(g^{-1}(D'_{f(x)})) \\
                 &=\mu((f\times g)^{-1}(D')_x)) \\
              & = \tau_{\Xcal,\Ycal}(x,\mu))((f\times g)^{-1}(D')) \\
             &= \PrQm(f\times g)(\tau_{\Xcal,\Ycal}(x,\mu))(D').
         \end{align*}
        \end{proof}
    This shows the claim.
\end{Lem}

\begin{Thm}
    \label{thm:Q-strong-monad}
    The triple $(\PrQm,\delta,\Mbb)$ defines a \emph{strong monad} on the cartesian monoidal quasitopos $(\QMS,\times,\one)$. It is \emph{affine} by \Cref{qms:lem:affine}, and \emph{commutative} as soon as the sample space has the Fubini property and $\PrQm(\Xcal)=\PrPf(\Xcal)$ for the spaces involved, by \Cref{qms:thm:commutative} and \Cref{qms:rem:comm-PrQm}.
    \begin{proof}
        By Theorem \ref{qms:cor:qms-cartesian-closed} is a (cartesian closed) monoidal category $(\QMS,\times,\one)$.
        By Proposition \ref{prp:Q-monad} the triple $(\PrQm,\delta,\Mbb)$ defines a monad structure on $\QMS$.
        We then define the \emph{strength} of the monad via $\tau$ from Lemma
        \ref{lem:Q-monad-strength}, where it was also shown that $\tau$ is a natural transformation
         $(\_)\times \PrQm(\_) \to \PrQm(\_\times\_)$.
        We are left to check its coherence conditions. \\
        Left unitor and strength satisfy the commutative diagram:
        \[\xymatrix{
                \one \times \PrQm(\Xcal) \ar_-{L_{\PrQm(\Xcal)}}[rd] \ar^-{\tau_{\one,\Xcal}}[rr] &&
                \PrQm(\one \times \Xcal) \ar^{\PrQm(L_\Xcal)}[dl]\\
            &\PrQm(\Xcal).
        }\]
        Associator and strength satisfy the commutative diagram:
        \[\xymatrix{
                (\Xcal\times\Ycal)\times\PrQm(\Zcal) \ar^-{\tau_{\Xcal\times\Ycal,\Zcal}}[rrrr] \ar_-{A_{\Xcal,\Ycal,\PrQm(\Zcal)}}[d] &&&& \PrQm\lp(\Xcal\times\Ycal)\times\Zcal\rp \ar^-{\PrQm(A_{\Xcal,\Ycal,\Zcal})}[d]\\
                \Xcal\times(\Ycal\times\PrQm(\Zcal)) \ar_-{\id_\Xcal \times \tau_{\Ycal,\Zcal}}[rr] && \Xcal\times\PrQm(\Ycal\times\Zcal) \ar_-{\tau_{\Xcal,\Ycal\times\Zcal}}[rr] && \PrQm(\Xcal\times(\Ycal\times\Zcal)).
        }\]
        Monad unit $\delta$ and strength satisfy the commutative diagram:
        \[\xymatrix{
            &\Xcal \times \Ycal \ar_-{\id_\Xcal\times \delta_\Ycal}[dl]
            \ar^-{\delta_{\Xcal\times\Ycal}}[dr]\\
                \Xcal \times \PrQm(\Ycal) \ar_-{\tau_{\Xcal,\Ycal}}[rr] && \PrQm(\Xcal \times \Ycal).
        }\]
        Monad action $\Mbb$ commutes with strength, expressed via the commutative diagram:
        \[\xymatrix{
                \Xcal\times\PrQm(\PrQm(\Ycal)) \ar^-{\tau_{\Xcal,\PrQm(\Ycal)}}[rr] \ar_-{\id_\Xcal\times\Mbb_\Ycal}[d] &&
                \PrQm(\Xcal\times\PrQm(\Ycal)) \ar^-{\PrQm(\tau_{\Xcal,\Ycal})}[rr] &&
                \PrQm(\PrQm(\Xcal \times \Ycal)) \ar^-{\Mbb_{\Xcal\times\Ycal}}[d] \\
                \Xcal\times\PrQm(\Ycal) \ar_-{\tau_{\Xcal,\Ycal}}[rrrr]&&&&
                \PrQm(\Xcal\times\Ycal).
        }\]
        To check this let $x \in \Xcal$ and $\pi \in \PrQm(\PrQm(\Ycal))$ and $D \in \Bcal_{\Xcal\times\Ycal}$.
        Then we get:
        \begin{align*}
            &\Mbb_{\Xcal\times\Ycal}\lp\PrQm(\tau_{\Xcal,\Ycal})(\tau_{\Xcal,\PrQm(\Ycal)}(x,\pi)) \rp(D) \\
            &= \int \ev(\rho,D)\, \PrQm(\tau_{\Xcal,\Ycal})(\tau_{\Xcal,\PrQm(\Ycal)}(x,\pi))(d\rho) \\
            &= \int \ev(\tau_{\Xcal,\Ycal}(t,\nu),D)\, (\tau_{\Xcal,\PrQm(\Ycal)}(x,\pi))(d(t,\nu)) \\
            &= \int \ev(\tau_{\Xcal,\Ycal}(x,\nu),D)\, \pi(d\nu) \\
            &= \int \ev(\nu,D_x)\, \pi(d\nu) \\
            &= \Mbb_{\Ycal}(\pi)(D_x) \\
            &= \tau_{\Xcal,\Ycal}(x,\Mbb_{\Ycal}(\pi))(D).
        \end{align*}
        This shows the claim.
    \end{proof}
\end{Thm}

\subsection{Strong Probability Monads $\PrPF$, $\PrpF$, $\PrPf$, $\Prpf$}

\label{sec:pf-prob-monads}

In this subsection we will introduce different versions of the \emph{probability monad of push-forward probability measures} on quasi-measurable spaces: $\PrPF$, $\PrpF$, $\PrPf$, $\Prpf$. $\PrpF$ will resemble
the probability monad $\PrpF$ from \cite{Heu17}. $\PrPF$ will be a bit more general than $\PrpF$, while $\Prpf$ will be a bit more restrictive than $\PrpF$. $\PrPf$ will in some sense be complementary to $\PrpF$.
We then study under which conditions these probability monads agree and also when they become strong.
The main requirement will be that the sample space $\Omega$ satisfies $\Omega \times \Omega \cong \Omega$ in $\QMS$.

\begin{Def}[The space of push-forward probability measures]
    \label{def:prob-space-K}
    Let $(\Xcal,\Xcal^\Omega)$ be a quasi-measurable space and $\Bcal_\Xcal:=\Bcal(\Xcal^\Omega)$ the induced $\sigma$-algebra.
    Then we define the quasi-measurable space of push-forward probability measures on $(\Xcal,\Xcal^\Omega)$:
    \[ \PrPF(\Xcal):= \pf\lp \Xcal^\Omega \times \PrQm(\Omega) \rp \ins \PrQm(\Xcal),\]
    with the quotient quasi-measurable space structure, where we used the quasi-measurable push-forward map, see Lemma \ref{lem:push-forward-qm}:
    \[ \pf:\; \Xcal^\Omega \times \PrQm(\Omega) \to \PrQm(\Xcal).\]
    More concretely, we have:
    \begin{align*}
    \PrPF(\Xcal) &=\lC  X_*\mu:\, \Bcal_\Xcal \to [0,1]\st  X \in \Xcal^\Omega, \mu \in \PrQm(\Omega) \rC,  \\
    \PrPF(\Xcal)^\Omega &= \pf \circ \lp \lp\Xcal^\Omega\rp^{\Omega} \times \PrQm(\Omega)^\Omega \rp\\
                                     &= \lC   X_*\kappa:\, \Omega \to \PrPF(\Xcal)\st  X \in \lp\Xcal^\Omega\rp^{\Omega}, \kappa \in \PrQm(\Omega)^\Omega \rC.
\end{align*}
     The  $ X_*\kappa$ for $ X \in \lp\Xcal^\Omega\rp^\Omega$ and $\kappa \in \PrQm(\Omega)^\Omega$ in $\PrPF(\Xcal)^\Omega$ is given by:
     \[ ( X_*\kappa)(\omega)(A) =  X(\omega)_*\kappa(\omega)(A) = \kappa(\omega)\lp X(\omega)^{-1}(A)\rp = \kappa(\omega)\lp \lC \omega' \in \Omega\st  X(\omega)(\omega') \in A \rC \rp.  \]
    We then also define the quasi-measurable spaces of probability measures:
        \begin{align*}
             \PrpF(\Xcal) &:=\PrPF(\Xcal),  \\
             \PrpF(\Xcal)^\Omega &:= \pf \circ \lp \lp\Xcal^\Omega\rp^\one \times \PrQm(\Omega)^\Omega \rp \ins \PrPF(\Xcal)^\Omega,\\
           \PrPf(\Xcal) &:=\pf\lp \Xcal^\Omega \times \PrPf \rp,  \\
            \PrPf(\Xcal)^\Omega &:= \pf \circ \lp \lp\Xcal^\Omega\rp^\Omega \times \PrPf^\one \rp \ins \PrPF(\Xcal)^\Omega,\\
            \Prpf(\Xcal)&:= \PrPf(\Xcal),\\
            \Prpf(\Xcal)^\Omega &:= \pf \circ \lp \lp \Xcal^\Omega\rp^\one \times \lp\Omega^\Omega\rp^\Omega \times \PrPf^\one \rp  \\
                                &= \lC  X_*\phi_*P:\; \Omega_1 \to \Prpf(\Xcal)\st  X \in \Xcal^{\Omega_3}, \phi \in \lp \Omega_3^{\Omega_2} \rp^{\Omega_1}, P \in \PrPf \rC\\
                                & \ins \PrpF(\Xcal)^\Omega \cap \PrPf(\Xcal)^\Omega,
        \end{align*}
        with indices for clarity: $\Omega_i:=\Omega$.
        Here $\PrPf$ in the last four lines is the set of product-compatible probability measures
        fixed in \Cref{qms:def:sample-space-act-3}, restricted to $\Bcal_\Omega$, so that these two
        definitions agree with the ones given in the main text in \Cref{sec:pf-monad}.
\end{Def}

\begin{Rem}[Which pool of measures, and why the four differ]
    \label{rem:which-pool}
    The four spaces are not built from the same pool of measures on $\Omega$, and they cannot be.
    $\PrPf$ and $\Prpf$ hold the measure \emph{fixed} while the random variable varies, so they use
    the pool only as a \emph{set}; they are therefore built from the act-3 datum $\PrPf$ of
    \Cref{qms:def:sample-space-act-3}, which is what the main text does.
    $\PrPF$ and $\PrpF$ let the measure \emph{vary} with the sample point, so they need a set
    $\PrQm(\Omega)^\Omega$ of admissible \emph{measure-valued} random variables, i.e.\ a
    quasi-measurable structure on the pool. \Cref{qms:def:sample-space-act-3} provides no such
    structure --- it fixes a set of measures and nothing more --- so for these two we take the
    canonical choice $\PrQm(\Omega)$, the quasi-measurable space of \Cref{qms:def:PrQm-maps}.
    This is consistent, because $\PrPf$ restricted to $\Bcal_\Omega$ is contained in
    $\PrQm(\Omega)$:

    \begin{Cor}[The act-3 measures are quasi-measurable]
        \label{cor:PrPf-in-PrQm-Omega}
        For every $P \in \PrPf$ the restriction $P|_{\Bcal_\Omega}$ lies in $\PrQm(\Omega)$; hence
        $\PrPf(\Xcal) \ins \PrPF(\Xcal)$ and $\Prpf(\Xcal) \ins \PrpF(\Xcal)$ for every
        quasi-measurable space $(\Xcal,\Xcal^\Omega)$.
        \begin{proof}
            Apply \Cref{qms:lem:PrPf-in-PrQm} to $(\Xcal,\Xcal^\Omega)=(\Omega,\Omega^\Omega)$ and
            $X=\id_\Omega \in \Omega^\Omega$, which gives
            $P|_{\Bcal_\Omega}=\id_{\Omega,*}P \in \PrPf(\Omega) \ins \PrQm(\Omega)$. The
            inclusions of the spaces follow, the pools being nested.
        \end{proof}
    \end{Cor}

    All four are built from the same pool exactly when $\PrPf = \PrQm(\Omega)$, which holds for the
    sample space of \Cref{qms:def:sample-space-act-4}, where
    $\PrPf=\Giry(\Omega,\Bcal_\Omega)=\PrQm(\Omega)$ by \Cref{lem:univ-qus-G-Q}. Off that case the
    reader should keep in mind that $\PrPF$ and $\PrpF$ are the \emph{maximal} variants, which is
    also why they sit at the top of the inclusion chain of \Cref{lem:S-P-R-K-incl}.
\end{Rem}

\begin{Lem}[Act 3 absorbs products into the random variable]
    \label{lem:act3-absorb}
    Assume \Cref{qms:asm:prod-iso} with isomorphism $\Theta:\, \Omega \cong \Omega \times \Omega$
    and let $P_1,P_2 \in \PrPf$. Then there are $P \in \PrPf$ and $\Psi \in \Omega^\Omega$ with
    \begin{align}
        \Theta_*^{-1}\lp P_1 \otimes P_2 \rp &= \Psi_*P \qquad \text{on } \Bcal_\Omega.
    \end{align}
    Consequently, for every quasi-measurable space $(\Xcal,\Xcal^\Omega)$ and every
    $X \in \Xcal^\Omega$ we have $X_*\Theta_*^{-1}(P_1 \otimes P_2) = (X \circ \Psi)_*P \in
    \PrPf(\Xcal)$, and likewise with $X$ replaced by a member of $(\Xcal^\Omega)^\Omega$.
    \begin{proof}
        By \Cref{qms:def:sample-space-act-3} \Cref{qms:eq:sample-space-act-3-c} there are
        $P \in \PrPf$ and $\Phi \in (\Omega\times\Omega)^\Omega$ with
        $P_1 \otimes P_2 = \Phi_*P$ on $\Bcal_{\Omega \times \Omega}$. Put
        $\Psi := \Theta^{-1} \circ \Phi$, which lies in $\Omega^\Omega$ because $\Omega^\Omega$ is
        closed under composition and $\Theta^{-1}$ is quasi-measurable. The last claim follows
        since $\Xcal^\Omega$ is closed under precomposition with $\Omega^\Omega$.
    \end{proof}
\end{Lem}

\begin{Rem}[Where the extra step is needed]
    \label{rem:act3-absorb-use}
    Three proofs below --- of \Cref{lem:K-P-S-R-eq}, \Cref{lem:K-R-M-nat-trafo} and
    \Cref{thm:K-P-Q-markov-prod} --- form a measure $\Theta_*^{-1}(\kappa \otimes \rho)$ and use
    that it is again an admissible measure. For $\PrPF$ and $\PrpF$, whose pool is all of
    $\PrQm(\Omega)$, that is immediate. For $\PrPf$ and $\Prpf$, whose pool is the act-3 datum,
    it is not: act 3 \Cref{qms:eq:sample-space-act-3-c} makes $P_1 \otimes P_2$ a push-forward of a
    pool member, not a pool member. \Cref{lem:act3-absorb} supplies the missing step, and the
    resulting map in $\Omega^\Omega$ is absorbed into the random-variable slot, which is free in
    both $\PrPf$ and $\Prpf$. This is exactly how the main-text proof of
    \Cref{qms:thm:product-PrPf} argues, and it is why that proof is stated for a general act-3
    $\PrPf$ from the start.
\end{Rem}

\begin{Lem}
    \label{lem:S-P-R-K-incl}
    Let $(\Xcal,\Xcal^\Omega)$ be a quasi-measurable space.
    Then the spaces of push-forward probability measures $(\PrPF(\Xcal),\PrPF(\Xcal)^\Omega)$ and $(\PrpF(\Xcal),\PrpF(\Xcal)^\Omega)$ and $(\PrPf(\Xcal),\PrPf(\Xcal)^\Omega)$ and $(\Prpf(\Xcal),\Prpf(\Xcal)^\Omega)$ are all quasi-measurable spaces.
    Furthermore, all the inclusion maps are quasi-measurable:
        \[ \Prpf(\Xcal) \ins \PrpF(\Xcal) \ins \PrPF(\Xcal) \ins \PrQm(\Xcal) \ins \Giry(\Xcal),\quad \Prpf(\Xcal) \ins \PrPf(\Xcal) \ins \PrPF(\Xcal).\]
    \begin{proof}
        This immediately follows from Lemma \ref{lem:push-forward-qm} and Lemma \ref{qms:lem:srj-pf-qm-str}.
    \end{proof}
\end{Lem}

\begin{Lem}
    \label{lem:K-P-S-R-eq}
    Let $(\Xcal,\Xcal^\Omega)$ be a quasi-measurable space.
    \begin{enumerate}
        \item Assume that there exists an isomorphism of quasi-measurable spaces:
                 \[ \Theta:\; \Omega \cong \Omega \times \Omega.  \]
            Then we have the equalities and inclusions  of sets:
            \[ \Prpf(\Xcal)^\Omega=\PrPf(\Xcal)^\Omega \ins \PrpF(\Xcal)^\Omega = \PrPF(\Xcal)^\Omega. \]
        \item Assume that the following push-forward map is surjective:
            \[ \pf:\; \lp\Omega^\Omega\rp^\Omega \times \PrQm(\Omega) \srj \PrQm(\Omega)^\Omega,
            \qquad( Z,\nu) \mapsto  Z_*\nu =\lp \omega \mapsto  Z(\omega)_*\nu \rp.\]
            Then we have the equalities and inclusions  of sets:
            \[ \Prpf(\Xcal)^\Omega=\PrpF(\Xcal)^\Omega \ins \PrPf(\Xcal)^\Omega = \PrPF(\Xcal)^\Omega. \]
    \end{enumerate}
    \begin{proof}
        1.) Let $ X_*\kappa \in \PrPF(\Xcal)^{\Omega_1}$ with $ X \in \lp\Xcal^{\Omega_2}\rp^{\Omega_1}$
    and $\kappa \in \PrQm(\Omega_2)^{\Omega_1}$, where we use indices for clarity: $\Omega_i:=\Omega$.
    In the following we will identify all maps with their (un)curried form.
Define $\phi_1 \in \lp\Omega_1^{\Omega_1}\rp^{\Omega_1}$ and $\phi_2 \in \lp\Omega_2^{\Omega_2}\rp^{\Omega_1}$
    on elements via:
    \begin{align*}
        \phi_1(\omega_1)(\omega_1') &:= \omega_1,\\
        \phi_2(\omega_1)(\omega_2) &:= \omega_2.
    \end{align*}
    Let $\nu \in \PrQm(\Omega)$ be any probability measure. Note that $\phi_2(\omega_1)=\id_{\Omega_2}$ and:
    \[ \phi_1(\omega_1)_*\nu = \delta_{\omega_1}, \qquad (\phi_1(\omega_1) \times \phi_2(\omega_1) )_*(\nu \otimes \kappa)(\omega_1) = \delta_{\omega_1} \otimes \kappa(\omega_1).\]
    With this we then get:
    \begin{align*}
         X(\omega_1)_* \kappa(\omega_1) &=  X_* \lp \delta_{\omega_1}\otimes\kappa(\omega_1)\rp\\
        &=  X_* \phi(\omega_1)_* \lp \nu \otimes \kappa(\omega_1) \rp\\
        &=  X_* \Theta_* \Theta_*^{-1} \phi(\omega_1)_* \Theta_* \Theta_*^{-1}\lp \nu \otimes \kappa(\omega_1) \rp\\
        &=  ( X \circ \Theta)_* \,(\Theta^{-1} \circ \phi(\omega_1) \circ \Theta)_*\, \lp\Theta_*^{-1}\lp \nu \otimes \kappa(\omega_1)\rp\rp,
    \end{align*}
    with $( X \circ \Theta) \in \Xcal^{\Omega_3}$, $(\Theta^{-1} \circ \phi \circ \Theta) \in \lp\Omega_3^{\Omega_3}\rp^{\Omega_1}$ and $\Theta_*^{-1}\lp \nu \otimes \kappa\rp \in \PrQm(\Omega_3)^{\Omega_1}$.
    This shows that: $  X_*\kappa \in \PrpF(\Xcal)^\Omega$.
    So we showed the inclusion (and thus equality):
    \[ \PrPF(\Xcal)^\Omega \ins \PrpF(\Xcal)^\Omega.\]
    The same arguments hold for $\kappa \in \PrQm(\Omega)$. Then $\Theta_*^{-1}\lp \nu \otimes \kappa\rp \in \PrQm(\Omega_3)$.\footnote{For $\PrPf$ and $\Prpf$, whose pool of measures is the act-3 datum rather than $\PrQm(\Omega)$, see \Cref{rem:which-pool}, this measure need not be a pool member. By \Cref{lem:act3-absorb} it is a push-forward $\Psi_*P$ of a pool member along some $\Psi \in \Omega^\Omega$, and $\Psi$ is absorbed into the $\phi$-slot, which is free in both; the resulting element is unchanged. The same substitution applies at the two further places marked in \Cref{rem:act3-absorb-use}.\label{fn:absorb}} This then shows the inclusion (and thus equality):
    \[ \PrPf(\Xcal)^\Omega \ins \Prpf(\Xcal)^\Omega.\]
    2.) Again let $ X_*\kappa \in \PrPF(\Xcal)^{\Omega_1}$ with
    $ X \in \lp\Xcal^{\Omega_2}\rp^{\Omega_1}$ and $\kappa \in \PrQm(\Omega_2)^{\Omega_1}$.
    By assumption there exist $\phi \in \lp\Omega_2^{\Omega_3}\rp^{\Omega_1}$ and $\nu \in \PrQm(\Omega_3)$ such that for all $\omega_1 \in \Omega_1$:
    \[ \kappa(\omega_1) =\phi(\omega_1)_*\nu.\]
    Define $  X \circ \phi \in \lp\Xcal^{\Omega_3}\rp^{\Omega_1}$ on elements via:
    \[ ( X \circ \phi)(\omega_1)(\omega_3) :=  X(\omega_1)(\phi(\omega_1)(\omega_3)). \]
    This gives us:
    \[  X(\omega_1)_*\kappa(\omega_1) =  X(\omega_1)_*\phi(\omega_1)_*\nu
    = ( X \circ \phi)(\omega_1)_*\nu.\]
    This shows that $ X_*\kappa = ( X \circ \phi)_*\nu \in \PrPf(\Xcal)^{\Omega_1}$.
    With this we get the inclusions (and thus equality):
    \[ \PrPF(\Xcal)^\Omega \ins \PrPf(\Xcal)^\Omega. \]
    The same proof for $ X \in \Xcal^{\Omega_2}$ shows:
    \[  X_*\kappa(\omega_1)= X_*\phi(\omega_1)_*\nu,\]
    and thus the inclusion (and equality):
    \[ \PrpF(\Xcal)^\Omega \ins \Prpf(\Xcal)^\Omega.\]
    This shows the claims.
    \end{proof}
\end{Lem}

\begin{Lem}
    \label{lem:K-P-R-Q-G-qm-incl}
    Let $(\Xcal,\Xcal^\Omega)$ be a quasi-measurable space and
    $f:\; \Xcal \to \Ycal$  a quasi-measurable map.
    Then we get a commutative diagrams of quasi-measurable maps:
    \[\xymatrix{
     \Prpf(\Xcal) \ar^-{\incl_\Xcal}[rr] \ar_-{\Prpf(f)}[d]  &&       \PrpF(\Xcal) \ar^-{\incl_\Xcal}[rr] \ar_-{\PrpF(f)}[d]  &&  \PrPF(\Xcal)\ar_-{\PrPF(f)}[d] \ar^-{\incl_\Xcal}[rr]  && \PrQm(\Xcal) \ar_-{\PrQm(f)}[d] \ar^-{\incl_\Xcal}[rr] && \Giry(\Xcal) \ar_-{\Giry(f)}[d]   \\
      \Prpf(\Ycal) \ar_-{\incl_\Ycal}[rr]   &&       \PrpF(\Ycal) \ar_-{\incl_\Ycal}[rr]   && \PrPF(\Ycal) \ar_-{\incl_\Ycal}[rr] && \PrQm(\Ycal) \ar_-{\incl_\Ycal}[rr]&& \Giry(\Ycal).
    }\]
    The same holds true if we replace $\PrpF$ with $\PrPf$.
    \begin{proof}
        That the horizontal lines are inclusions and quasi-measurable is clear.\\
        Since the push-forward is at each vertical map given by:
        \[ (f_*\mu)(B) := \mu(f^{-1}(B)),\]
        the commutativity is clear as soon as all vertical maps are shown to be well-defined (with the corresponding codomain) and quasi-measurable.\\
        For $\PrQm(f)$ both claims were shown in Lemma \ref{lem:push-forward-qm}.\\
        $\Giry(f)$ is clearly well-defined since $f_*\mu$ is a probability measure on $\Bcal_\Ycal$.
        To see that $\Giry(f)$ is quasi-measurable let $\kappa \in \Giry(\Xcal)^\Omega$ and $B \in \Bcal_\Ycal$.
        Then $f^{-1}(B) \in \Bcal_\Xcal$ and the map:
        \[ \Omega \to [0,1],\qquad \omega \mapsto (f_*\kappa)(\omega)(B)=\kappa(\omega)(f^{-1}(B)), \]
        is $\Bcal_\Omega$-$\Bcal_{[0,1]}$-measurable by assumption on $\kappa$. So $f_*\kappa \in \Giry(\Ycal)^\Omega$.\\
        To show that $\PrPF(f)$, $\PrpF(f)$, $\PrPf(f)$ are well-defined consider $ X_*\mu$ with $ X \in \Xcal^\Omega$ and $\mu \in \PrQm(\Omega)$. Since $f$ is quasi-measurable we have $f \circ  X \in \Ycal^\Omega$ and thus:
        \[ f_*( X_*\mu)= (f \circ  X)_*\mu \in \PrPF(\Ycal)=\PrpF(\Ycal)=\PrPf(\Ycal).\]
        To show that $\PrPF(f)$ is quasi-measurable let now $ X \in \lp\Xcal^{\Omega_2}\rp^{\Omega_1}$ and
        $\kappa \in \PrQm(\Omega_2)^{\Omega_1}$. Again, since $f$ is quasi-measurable we get that the composition:
        \[ f \circ  X:\; \Omega_2 \times \Omega_1 \xrightarrow{ X} \Xcal \xrightarrow{f} \Ycal,   \]
        is quasi-measurable and thus $f \circ  X \in \lp\Ycal^{\Omega_2}\rp^{\Omega_1}$. So we get:
        \[ f_*( X_*\kappa)(\omega_1) = f_* X(\omega_1)_*\kappa(\omega_1) =
        (f \circ  X)(\omega_1)_*\kappa(\omega_1).   \]
        So it follows that:
        \[ f_*( X_*\kappa) = (f \circ  X)_*\kappa \in \PrPF(\Ycal)^{\Omega_1}.\]
        Note that if in the above we have that $ X \in \lp\Xcal^{\Omega_2}\rp^\one$
        then also $f \circ  X  \in \lp\Ycal^{\Omega_2}\rp^\one$,
        showing that $\PrpF(f)$ is quasi-measurable.\\
        Similarly, if $\kappa \in \PrQm(\Omega_2)^\one$ then it stays that way for the push-forward,
        so clearly $\PrPf(f)$ is quasi-measurable.\\
        For $\Prpf(f)$ note that we can use the previous arguments to get:
        \[  f \circ \Prpf(\Xcal)^\Omega = f \circ \lp \PrpF(\Xcal)^\Omega \cap \PrPf(\Xcal)^\Omega  \rp
        \ins \PrpF(\Ycal)^\Omega \cap \PrPf(\Ycal)^\Omega = \Prpf(\Ycal)^\Omega.\]
    \end{proof}
\end{Lem}

\begin{Lem}
    \label{lem:K-P-R-delta-nat-trafo}
    Let $(\Xcal,\Xcal^\Omega)$ be a quasi-measurable space.
    Then the map:
    \[ \delta:\; \Xcal \to \Prpf(\Xcal),\qquad x \mapsto \delta_x = (A \mapsto \I_A(x)),\]
    is a well-defined quasi-measurable map.\\
    Furthermore, if we have a quasi-measurable map:
    \[f:\; (\Xcal,\Xcal^\Omega) \to (\Ycal,\Ycal^\Omega), \]
    then we get a commutative diagram of quasi-measurable maps:
    \[\xymatrix@C=1.5em{
            \Xcal \ar^-{\delta}[rr] \ar^-{f}[d] &&    \Prpf(\Xcal) \ar^-{i_\Xcal}[rr] \ar_-{\Prpf(f)}[d]  &&    \PrpF(\Xcal) \ar^-{i_\Xcal}[rr] \ar_-{\PrpF(f)}[d]  &&  \PrPF(\Xcal)\ar_-{\PrPF(f)}[d] \ar^-{i_\Xcal}[rr]  && \PrQm(\Xcal) \ar_-{\PrQm(f)}[d] \ar^-{i_\Xcal}[rr] && \Giry(\Xcal) \ar_-{\Giry(f)}[d]  \\
   \Ycal \ar^-{\delta}[rr] &&   \Prpf(\Ycal) \ar_-{i_\Ycal}[rr]   &&        \PrpF(\Ycal) \ar_-{i_\Ycal}[rr]   && \PrPF(\Ycal) \ar_-{i_\Ycal}[rr] && \PrQm(\Ycal) \ar_-{i_\Ycal}[rr]&& \Giry(\Ycal).
    }\]
    In other words,  we have natural transformations of endo-functors of $\QMS$:
    \[ \id_{\QMS} \xrightarrow{\delta} \Prpf \xrightarrow{\incl}\PrpF \xrightarrow{\incl} \PrPF \xrightarrow{\incl} \PrQm
    \xrightarrow{\incl} \Giry. \]
    All the above statements also hold true if we replace $\PrpF$ by $\PrPf$.
    \begin{proof}
        To show that $\delta$ is well-defined we need to show that $\delta_{\Xcal,x} \in \Prpf(\Xcal)=\PrPF(\Xcal)$
        for $x\in \Xcal$.
        For any $\nu \in \PrQm(\Omega)$ and the constant map $x^\one \in \Xcal^\one$ with value $x \in \Xcal$ we get:
        \[ \delta_{\Xcal,x} = x^\one_* \nu \in \Prpf(\Xcal)=\PrPF(\Xcal).\]
        To show that $\delta_\Xcal$  is quasi-measurable let $ X \in \Xcal^{\Omega}$ and define:
        $\phi(\omega_1)(\omega_2):=\omega_1$.
            Then we get:
            \begin{align*}
                \lp X_*\phi_*\nu\rp(\omega_1) &=  X_*\phi(\omega_1)_*\nu \\
                                                  &=  X_*\delta_{\Omega,\omega_1}\\
                                                  &= \delta_{\Xcal, X(\omega_1)}\\
                                                  &= (\delta_\Xcal \circ  X)(\omega_1).
            \end{align*}
            So we get: $\delta_\Xcal \circ  X =  X_*\phi_*\nu \in \Prpf(\Xcal)^\Omega$.\\
        To check the commutativity of the  diagrams, by Lemma \ref{lem:K-P-R-Q-G-qm-incl} we only need to consider the
        first one with $\Prpf$. Observe that:
        \[ f_*\delta_{\Xcal,x} = f_*x^\one_*\nu = f(x)^\one_* \nu = \delta_{\Ycal,f(x)}  \in \Prpf(\Ycal)=\PrPF(\Ycal),\]
        where with $x^\one$ and $f(x)^\one$ we mean the corresponding constant maps with those values.\\
        This shows the claims.
    \end{proof}
\end{Lem}

\begin{Lem}
    \label{lem:K-R-M-nat-trafo}
        Let $(\Xcal,\Xcal^\Omega)$ be a quasi-measurable space.
        Assume that there exists an isomorphism of quasi-measurable spaces:
                 \[ \Theta:\; \Omega \cong \Omega \times \Omega.  \]
        Then the map:
        \[\Mbb_\Xcal:\; \PrPF(\PrPF(\Xcal)) \to \PrPF(\Xcal),\qquad \pi \mapsto \Mbb_\Xcal(\pi)=\lp A \mapsto
        \int \ev(\mu,A)\, \pi(d\mu)\rp,\]
        is a well-defined quasi-measurable map.\\
    Furthermore, if we have a quasi-measurable map:
    \[f:\; (\Xcal,\Xcal^\Omega) \to (\Ycal,\Ycal^\Omega), \]
    then we get a commutative diagram of quasi-measurable maps:
    \[\xymatrix{
            \PrPF(\PrPF(\Xcal)) \ar^-{\PrPF(\PrPF(f))}[r] \ar_-{\Mbb_\Xcal}[d] & \PrPF(\PrPF(\Ycal)) \ar^-{\Mbb_\Ycal}[d]\\
            \PrPF(\Xcal) \ar_{\PrPF(f)}[r] & \PrPF(\Ycal).
    }\]
    In other words, $\Mbb:\, \PrPF^2 \to \PrPF$ is a natural transformation of endo-functors of $\QMS$.\\
    Finally,  all the analogous statements hold as well when $\PrPF$ is replaced by $\PrPf$ or $\PrpF$ or $\Prpf$.
    \begin{proof}
        The proof that $\Mbb_\Xcal$ maps to $\PrPF(\Xcal)$ follows the same lines as the proof that it is quasi-measurable.
        So we directly assume that $\pi \in \PrPF(\PrPF(\Xcal))^{\Omega_1}$.
        Then there are $\rho \in \PrQm(\Omega_2)^{\Omega_1}$ and $\psi \in \lp \PrPF(\Xcal)^{\Omega_2} \rp^{\Omega_1}$ such that:
        \[ \pi(\omega_1) = (\psi_*\rho)(\omega_1)= \psi(\omega_1)_*\rho(\omega_1) \in \PrPF(\PrPF(\Xcal)).\]
        Since by assumption we have $\Omega_1 \times \Omega_2 \cong \Omega$ we get that:
        \[\psi \in \lp \PrPF(\Xcal)^{\Omega_2} \rp^{\Omega_1} \cong \PrPF(\Xcal)^{\Omega}.\]
    So there are $ X \in \lp\lp\Xcal^{\Omega_3}\rp^{\Omega_2}\rp^{\Omega_1}$ and
        $\kappa \in \PrQm(\Omega_3)^{\Omega_1 \times \Omega_2}$ such that:
        \[ \psi(\omega_1)(\omega_2) =  X(\omega_1)(\omega_2)_*\kappa(\omega_1,\omega_2) \in \PrPF(\Xcal).\]
        Then we get with $A \in \Bcal_\Xcal$ and $\omega_1 \in \Omega_1$:
        \begin{align*}
            \Mbb_\Xcal(\pi)(\omega_1)(A) &= \int \ev(\mu,A)\, \pi(\omega_1)(d\mu) \\
                                       &= \int \ev(\mu,A)\, (\psi_*\rho)(\omega_1)(d\mu) \\
                                       &= \int \ev(\mu,A)\, \psi(\omega_1)_*\rho(\omega_1)(d\mu) \\
                                       &= \int \ev(\psi(\omega_1)(\omega_2),A)\, \rho(\omega_1)(d\omega_2) \\
                                       &= \int \ev( X(\omega_1)(\omega_2)_*\kappa(\omega_1,\omega_2),A)\, \rho(\omega_1)(d\omega_2) \\
                                       &= \int \kappa(\omega_1,\omega_2)( X(\omega_1)(\omega_2)^{-1}(A))\, \rho(\omega_1)(d\omega_2) \\
                                       &= \int \kappa(\omega_1,\omega_2)( X^{-1}(A)_{\omega_1,\omega_2})\, \rho(\omega_1)(d\omega_2) \\
                                       &= (\kappa\otimes\rho)(\omega_1)( X^{-1}(A)_{\omega_1}) \\
                                       &= (\kappa\otimes\rho)(\omega_1)( X(\omega_1)^{-1}(A)) \\
                                       &=  X(\omega_1)_*(\kappa\otimes\rho)(\omega_1)(A),
        \end{align*}
        where now $\kappa \otimes \rho \in \PrQm(\Omega_3 \times \Omega_2)^{\Omega_1}\cong\PrQm(\Omega)^{\Omega_1}$
        and $ X  \in \lp\lp\Xcal^{\Omega_3}\rp^{\Omega_2}\rp^{\Omega_1} \cong \lp\Xcal^\Omega\rp^{\Omega_1}$.
        It follows that:
        \[ \Mbb_\Xcal(\pi) =  X_*(\kappa \otimes \rho) \in \PrPF(\Xcal)^{\Omega_1}.\]
        This shows that $\Mbb_\Xcal$ is well-defined and quasi-measurable.\\
        To show the commutative diagram we apply the above calculation to:
        \[f_{**}\pi = (f_*\psi)_*\rho =(f_* X_*\kappa)_*\rho=((f \circ  X)_*\kappa)_*\rho.\]
        So replacing $ X$ with $f \circ  X$ in the calculation from above we get:
        \[ \Mbb_\Ycal(f_{**}\pi) = (f \circ  X)_*(\kappa\otimes\rho)= f_*  X_*(\kappa\otimes\rho) = f_* \Mbb_\Xcal(\pi).\]
        This shows the claim for $\PrPF$.\\
        For the corresponding statements for $\PrpF$ note that $\psi$ will not be dependent on $\omega_1$ and
        thus $ X$ will not be dependent on $\omega_2$ and $\omega_1$ and $\kappa$ will not be dependent on $\omega_2$.
        In this case we do not even need the isomorphism $\Omega \times \Omega \cong \Omega$ and we get:
        \[ \Mbb_\Xcal(\pi) =  X_* (\kappa \circ \rho) \in \PrpF(\Xcal)^{\Omega_1}.\]
        For the $\PrPf$ case note that $\rho$ will not be dependent on $\omega_1$ and $\kappa$ will not be dependent on $\omega_2$ and $\omega_3$, which shows that $\kappa \otimes \rho$ will not be dependent on $\omega_1$; here $\kappa$ and $\rho$ are then fixed members $P_1,P_2$ of the act-3 pool, so \Cref{lem:act3-absorb} applies and writes $\Theta_*^{-1}(P_1 \otimes P_2)$ as $\Psi_*P$ with $P \in \PrPf$, the map $\Psi$ being absorbed into the random variable, see \footref{fn:absorb}. In this case we get that:
                \[ \Mbb_\Xcal(\pi) =  X_* (\kappa \otimes \rho) \in \PrPf(\Xcal)^{\Omega_1}.\]
                The statement for  $\Prpf$ is the same as for $\PrPf$, see Lemma \ref{lem:K-P-S-R-eq}.
    \end{proof}
\end{Lem}

\begin{Prp}
    \label{prp:P-K-R-monad}
    The triple $(\PrpF,\delta,\Mbb)$ defines a monad on $\QMS$.
    If, furthermore,  we have an isomorphism of quasi-measurable spaces:
                 \[ \Theta:\; \Omega \cong \Omega \times \Omega,  \]
    then also the triples $(\PrPF,\delta,\Mbb)$ and $(\PrPf,\delta,\Mbb)$ and $(\Prpf,\delta,\Mbb)$
    each  define a monad on $\QMS$.
    \begin{proof}
        By Lemma \ref{lem:K-P-R-Q-G-qm-incl} we know that $\PrpF:\, \QMS \to \QMS$ is a functor.
        By Lemma \ref{lem:K-P-R-delta-nat-trafo} we know that $\delta:\, \id_\QMS \to \PrpF$ is a natural transformation.
        By Lemma \ref{lem:K-R-M-nat-trafo} we know that $\Mbb:\,\PrpF^2 \to \PrpF$ is a natural transformation.
        Note that for $\PrpF$ the isomorphism $\Omega \times \Omega \cong \Omega$ was not needed.
        So we are left to check the following coherence conditions:
        \begin{enumerate}
        \item $\Mbb \circ \PrpF\Mbb = \Mbb \circ \Mbb \PrpF$ as natural transformations $\PrpF^3 \to \PrpF$,
        \item $\Mbb \circ \PrpF\delta = \Mbb \circ \delta \PrpF = \id_\PrpF$ as natural transformations $\PrpF \to \PrpF$.
        \end{enumerate}
        These calculations follow exactly the same as in Proposition \ref{prp:Q-monad}.\\
        All the same arguments also hold for $(\PrPF,\delta,\Mbb)$ and $(\PrPf,\delta,\Mbb)$ and
        $(\Prpf,\delta,\Mbb)$ under the assumed isomorphism.
    \end{proof}
\end{Prp}

\begin{Thm}
    \label{thm:K-P-Q-markov-prod}
    Let $(\Xcal,\Xcal^\Omega)$, $(\Ycal,\Ycal^\Omega)$, $(\Zcal,\Zcal^\Omega)$ be quasi-measurable spaces.
    Assume that there exists an isomorphism of quasi-measurable spaces:
                 \[ \Theta:\; \Omega \cong \Omega \times \Omega.  \]
    Then the maps:
    \[ \otimes:\; \PrPF(\Xcal)^{\Ycal \times \Zcal} \times \PrPF(\Ycal)^\Zcal \to \PrPF(\Xcal \times \Ycal)^\Zcal,\quad
    (\mu \otimes \nu)(z)(D)  := \int  \mu(y,z)(D^y)\, \nu(z)(dy), \]
        \[ \otimes:\; \PrpF(\Xcal)^{\Ycal \times \Zcal} \times \PrpF(\Ycal)^\Zcal \to \PrpF(\Xcal \times \Ycal)^\Zcal,\quad
    (\mu \otimes \nu)(z)(D)  := \int  \mu(y,z)(D^y)\, \nu(z)(dy), \]
    \[ \otimes:\; \PrPf(\Xcal)^{\Ycal \times \Zcal} \times \PrPf(\Ycal)^\Zcal \to \PrPf(\Xcal \times \Ycal)^\Zcal,\quad
    (\mu \otimes \nu)(z)(D)  := \int  \mu(y,z)(D^y)\, \nu(z)(dy), \]
    \[ \otimes:\; \Prpf(\Xcal)^{\Ycal \times \Zcal} \times \Prpf(\Ycal)^\Zcal \to \Prpf(\Xcal \times \Ycal)^\Zcal,\quad
    (\mu \otimes \nu)(z)(D)  := \int  \mu(y,z)(D^y)\, \nu(z)(dy), \]
    are well-defined quasi-measurable maps.
    \begin{proof}
        To keep better track of the interplay between different versions of $\Omega$ in this proof we will index them as
        $\Omega_k:=\Omega$ for $k=1,2,\dots$.\\
            To show that $\otimes$ is well-defined let $\mu \in \PrPF(\Xcal)^{\Ycal \times \Zcal}$ and $\nu \in \PrPF(\Ycal)^{\Zcal}$.
            We then need to show that for every $z \in \Zcal$:
            \[ (\mu\otimes\nu)(z)  \in \PrPF(\Xcal \times \Ycal),\]
            and that for every $ Z \in \Zcal^{\Omega_1}$ we have:
            \[ (\mu\otimes\nu) \circ  Z  \in \PrPF(\Xcal \times \Ycal)^{\Omega_1}.\]
            The latter means that we need to show that for every $ Z \in \Zcal^{\Omega_1}$ there exist
            $\pi \in \PrQm(\Omega_5)^{\Omega_1}$ and $\psi \in \lp(\Xcal \times \Ycal)^{\Omega_5}\rp^{\Omega_1}$ such that:
            \[ (\mu\otimes\nu) \circ  Z = \psi_* \pi \in \PrPF(\Xcal \times \Ycal)^{\Omega_1}.\]
            Since the first statement can be considered a special case of the second one (with constant
            $ Z=z \in \Zcal^\one \ins \Zcal^{\Omega_1}$ and constant
            $\pi \in \PrQm(\Omega_2)^\one \ins \PrQm(\Omega_2)^{\Omega_1}$)
            we directly focus on arbitrary $ Z \in \Zcal^{\Omega_1}$.\\
            By composing $\nu \in \PrPF(\Ycal)^{\Zcal}$ with $ Z \in \Zcal^{\Omega_1}$ then have that
            $\nu \circ  Z \in \PrPF(\Ycal)^{\Omega_1}$.
            So by definition of $\PrPF(\Ycal)^{\Omega_1}$ there exist
            $\rho \in \PrQm(\Omega_2)^{\Omega_1}$ and $ Y \in \lp\Ycal^{\Omega_2}\rp^{\Omega_1}$ such that for every
            $\omega_1 \in \Omega_1$:
            \[ \nu \circ  Z(\omega_1) =  Y(\omega_1)_*\rho(\omega_1).\]
            Then note that $ Y \times  Z$ is given by:
            \[ Y \times  Z:\; \Omega_2 \times \Omega_1 \to \Ycal \times \Zcal,\qquad
            (\omega_2,\omega_1) \mapsto ( Y(\omega_1)(\omega_2), Z(\omega_1)).\]
           Together with the isomorphism $\Theta:\; \Omega_4 \cong \Omega_2 \times \Omega_1$ we get the quasi-measurable map:
           \[( Y\times Z) \circ \Theta:\; \Omega_4 \xrightarrow{\Theta} \Omega_2 \times \Omega_1 \xrightarrow{ Y \times  Z} \Ycal \times \Zcal,
            \qquad \omega_4 \mapsto ( Y(\Theta_2(\omega_4))(\Theta_1(\omega_4)), Z(\Theta_2(\omega_4)).\]
            Similarly to before by composing the latter map with $\mu \in\PrPF(\Xcal)^{\Ycal \times \Zcal}$  we have that:
            \[ \mu  \circ ( Y \times  Z) \circ \Theta  \in \PrPF(\Xcal)^{\Omega_4}. \]
            So by definition of $\PrPF(\Xcal)^{\Omega_4}$ there exist
            $\varepsilon \in \lp\Xcal^{\Omega_3}\rp^{\Omega_4}$ and $\lambda \in \PrQm(\Omega_3)^{\Omega_4}$
            such that for all $\omega_4 \in \Omega_4$:
            \[ \mu \circ ( Y \times  Z) \circ \Theta(\omega_4) = \varepsilon(\omega_4)_*\lambda(\omega_4).  \]
            For every $(\omega_2,\omega_1) \in \Omega_2\times \Omega_1$ we take
            $\omega_4 := \Theta^{-1}(\omega_2,\omega_1)$ in the above equation and get:
       \[ \mu \circ ( Y \times  Z)(\omega_2,\omega_1) = \varepsilon(\Theta^{-1}(\omega_2,\omega_1))_*\lambda(\Theta^{-1}(\omega_2,\omega_1)).  \]
        We then define  $\kappa \in \PrQm(\Omega_3)^{\Omega_2 \times \Omega_1}$ via the composition of quasi-measurable maps:
    \[\kappa := \lambda \circ \Theta^{-1} :\;\Omega_2 \times \Omega_1 \to \Omega_4 \to \PrQm(\Omega_3),\]
    and $ X \in \lp\Xcal^{\Omega_3} \rp^{\Omega_2 \times \Omega_1}$ via the composition of quasi-measurable maps:
    \[   X:=\varepsilon \circ \Theta^{-1}:\; \Omega_2 \times \Omega_1 \to \Omega_4 \to \Xcal^{\Omega_3}.\]
    We then get for all $(\omega_2,\omega_1)\in \Omega_2\times\Omega_1$:
    \[\mu \circ ( Y\times Z)(\omega_2,\omega_1) = X(\omega_2,\omega_1)_*\kappa(\omega_2,\omega_1).\]
    With this so defined $\kappa \in \PrQm(\Omega_3)^{\Omega_2 \times \Omega_1}$ and $\rho \in \PrQm(\Omega_2)^{\Omega_1}$
    we get that $\kappa \otimes \rho \in \PrQm(\Omega_3 \times \Omega_2)^{\Omega_1}$.\\
  We then again use the isomorphism $\Theta:\; \Omega_5 \cong \Omega_3 \times \Omega_2$ and define:
     \begin{align*}
         \pi&:=\Theta_*^{-1}(\kappa\otimes\rho) \in \PrQm(\Omega_5)^{\Omega_1}.
     \end{align*}
    Furthermore, we define the quasi-measurable map
    $( X \times  Y) \in \lp(\Xcal \times \Ycal)^{\Omega_3 \times \Omega_2} \rp^{\Omega_1}$ via:
    \[  ( X \times  Y)(\omega_1)(\omega_3,\omega_2) :=  ( X(\omega_2,\omega_1)(\omega_3), Y(\omega_1)(\omega_2)) \in \Xcal \times \Ycal,
    \]
    and the quasi-measurable map $\psi  \in \lp (\Xcal \times \Ycal)^{\Omega_5}\rp^{\Omega_1} $ via:
    \[\psi(\omega_1)(\omega_5) := ( X \times  Y)(\omega_1) \circ \Theta(\omega_5) = \lp  X(\Theta_2(\omega_5),\omega_1)(\Theta_3(\omega_5)), Y(\omega_1)(\Theta_2(\omega_5))     \rp \in \Xcal \times \Ycal.\]
    With all these settings and $D \in \Bcal_{\Xcal \times \Ycal}$ and $\omega_1 \in \Omega_1$ we get the calculations:
    \begin{align*}
        (\mu \otimes \nu)( Z(\omega_1))(D)
        &= \int \mu(y, Z(\omega_1))(D^y)\,\nu( Z(\omega_1))(dy) \\
        &= \int \mu(y, Z(\omega_1))(D^y)\, ( Y(\omega_1)_*\rho(\omega_1))(dy) \\
        &= \int \mu( Y(\omega_1)(\omega_2), Z(\omega_1))(D_{ Y(\omega_1)(\omega_2)})\, \rho(\omega_1)(d\omega_2) \\
        &= \int\mu\lp( Y \times  Z)(\omega_2,\omega_1)\rp(D_{ Y(\omega_1)(\omega_2)})\, \rho(\omega_1)(d\omega_2) \\
        &= \int  X(\omega_2,\omega_1)_*\kappa(\omega_2,\omega_1)(D_{ Y(\omega_1)(\omega_2)})\, \rho(\omega_1)(d\omega_2) \\
        &= \int \kappa(\omega_2,\omega_1)( X(\omega_2,\omega_1)^{-1}(D_{ Y(\omega_1)(\omega_2)}))\, \rho(\omega_1)(d\omega_2) \\
        &= \int \kappa(\omega_2,\omega_1)(( X \times  Y)(\omega_1)^{-1}(D)_{\omega_2})\, \rho(\omega_1)(d\omega_2) \\
        &= (\kappa \otimes \rho)(\omega_1) (( X \times  Y)(\omega_1)^{-1}(D))\\
        &= ( X \times  Y)(\omega_1)_* (\kappa \otimes \rho)(\omega_1) (D)\\
        &= ( X \times  Y)(\omega_1)_* \circ \Theta_* \circ \Theta_*^{-1} (\kappa \otimes \rho)(\omega_1) (D) \\
        &= \psi(\omega_1)_* \pi(\omega_1)(D).
    \end{align*}
    So we get that:
    \[ (\mu\otimes\nu) \circ  Z = \psi_* \pi \in \PrPF(\Xcal \times \Ycal)^{\Omega_1}.\]
    This shows that $\otimes$ is a well-defined map.\\
    To show that $\otimes$ is quasi-measurable we need to use the exact same arguments as above but where we replace $\Zcal$ with
    $\Omega \times \Zcal$ and $ Z \in \Zcal^\Omega$ with $\Phi \times  Z \in \Omega^\Omega \times \Zcal^\Omega$
    everywhere.
    This then shows the claim for $\PrPF(\Xcal \times \Ycal)^\Zcal$.\\
    To show the claim for $\PrpF(\Xcal \times \Ycal)^\Zcal$ the same proof applies by noticing that if $ Y$ is not dependent on $\omega_1$ and $\varepsilon$ not dependent on $\omega_4$ then $ X$ does not depend on $(\omega_2,\omega_1)$
    and thus $ X \times  Y$ and $\psi$ will not depend on $\omega_1$. It then follows:
    \[ (\mu\otimes\nu) \circ  Z = \psi_* \pi \in \PrpF(\Xcal \times \Ycal)^{\Omega_1}.\]
    Similarly the claim for $\PrPf(\Xcal \times \Ycal)^\Zcal$ follows by the same arguments --- with $\Theta_*^{-1}(\kappa\otimes\rho)$ rewritten by \Cref{lem:act3-absorb} as a push-forward of a member of the act-3 pool, see \footref{fn:absorb}, which is what the main-text proof of \Cref{qms:thm:product-PrPf} does directly --- by noticing that if
    $\rho$ is not dependent on $\omega_1$ and $\lambda$ not dependent on $\omega_4$ then $\kappa$ is not dependent on $(\omega_2,\omega_1)$ and thus $\kappa \otimes \rho$ and $\pi$ will not depend on $\omega_1$. It then follows that:
        \[ (\mu\otimes\nu) \circ  Z = \psi_* \pi \in \PrPf(\Xcal \times \Ycal)^{\Omega_1}.\]
        Finally, the claim for $\Prpf$ follows from the the case of $\PrPf$, see Lemma \ref{lem:K-P-S-R-eq}.\\

The case for $\PrpF(\Xcal \times \Ycal)^\Zcal$ can be visualized via the commutative diagrams:
        \[\xymatrix{
       \Omega \ar_-{\rho}[d] \ar^-{\Phi \times  Z}[r] & \Omega \times \Zcal \ar^{\nu}[d]   \\
        \PrQm(\Omega) \ar_-{ Y_*}[r] & \PrpF(\Ycal)
    }\qquad
    \xymatrix{
            \Omega \ar^-{\Theta}_-{\sim}[r] \ar_-{\lambda}[dr]&   \Omega \times  \Omega \ar_-{\kappa}[d] \ar^-{ Y\times(\Phi \times Z)}[rr]
                                            && \Ycal \times \Omega \times \Zcal \ar^{\mu}[d]   \\
                                            &  \PrQm(\Omega) \ar_-{ X_*}[rr] && \PrpF(\Xcal)
    }\]
  \[ \xymatrix{
           &\Omega \ar_-{\kappa\otimes\rho}[d] \ar^-{\Phi \times  Z}[r] \ar_-{\pi}[dl] & \Omega \times\Zcal \ar^{\mu\otimes \nu}[d]  \\
\PrQm(\Omega) \ar_-{\Theta_*}^-\sim[r] &
         \PrQm(\Omega \times \Omega) \ar_-{( X \times  Y)_*}[r] & \PrpF(\Xcal \times \Ycal).
            }\]
    \end{proof}
\end{Thm}

\begin{Thm}
    \label{thm:K-P-R-strong-monad}
    Assume that there exists an isomorphism of quasi-measurable spaces:
                 \[ \Theta:\; \Omega \cong \Omega \times \Omega.  \]
    Then the triples $(\PrPF,\delta,\Mbb)$, $(\PrpF,\delta,\Mbb)$, $(\PrPf,\delta,\Mbb)$, $(\Prpf,\delta,\Mbb)$ each define a \emph{strong monad} on the cartesian monoidal quasitopos $(\QMS,\times,\one)$. They are all \emph{affine} by \Cref{qms:lem:affine}, and all \emph{commutative} as soon as the sample space has the Fubini property, by \Cref{qms:thm:commutative}.
    \begin{proof}
        We define the strength of the monad as follows:
            \[ \tau_{\Xcal,\Ycal}:\; \Xcal \times \PrpF(\Ycal) \to \PrpF(\Xcal \times \Ycal),\qquad
         (x,\mu) \mapsto  \delta_x \otimes \mu,\]
    which is a well-defined quasi-measurable map by Theorem \ref{thm:K-P-Q-markov-prod}.
    The rest follows the same steps as in Theorem \ref{thm:Q-strong-monad}. For this note that all involved maps are already shown to be well-defined and quasi-measurable. One then only needs to check those coherence equations on elements. This follows exactly the same as in Theorem \ref{thm:Q-strong-monad}.\\
    Also note that the \emph{bind/extension operation}:
    \[ (\_)_\circ:\; \PrpF(\Ycal)^\Xcal \to \PrpF(\Ycal)^{\PrpF(\Xcal)},\quad \kappa \mapsto \kappa_\circ,\quad \kappa_\circ(\mu)(B):= (\kappa \circ \mu)(B) := \int \kappa(x)(B)\,\mu(dx),\]
    is a well-defined quasi-measurable map (as the adjunction of the push-forward of $\otimes$,
    see Theorem \ref{thm:K-P-Q-markov-prod}), which is functorial in $\Xcal$ and $\Ycal$.\\
    The same arguments hold for the other monads.
    \end{proof}
\end{Thm}

\begin{Lem}
    \label{lem:push-forward-R-K}
    Let $(\Xcal,\Xcal^\Omega)$, $(\Ycal,\Ycal^\Omega)$ be quasi-measurable spaces.
    Then the push-forward map:
    \[ \pf:\; \Ycal^\Xcal \times \PrPF(\Xcal) \to \PrPF(\Ycal), \qquad (f,\mu) \mapsto f_*\mu, \]
    is a well-defined and quasi-measurable map. The same holds true for $\PrPf$ in place of $\PrPF$.
    \begin{proof}
        We already know from Lemma \ref{lem:K-P-R-Q-G-qm-incl} that the map is well-defined.
        Consider the following commutative diagram of quasi-measurable maps, see \ref{lem:push-forward-qm}:
            \[\xymatrix{
            \Ycal^\Xcal \times \Xcal^\Omega \times \PrQm(\Omega) \ar@{->>}_-{\id\times \pf}[d] \ar^-{(\_\circ\_)\times\id}[rrr] &&& \Ycal^\Omega \times \PrQm(\Omega) \ar^-{\pf}[d]\\
            \Ycal^\Xcal \times \PrPF(\Xcal) \ar_-{\pf}[rrr] &&& \PrQm(\Ycal).
    }\]
    This then induces the following commutative diagram of quasi-measurable maps:
        \[\xymatrix{
            \lp\Ycal^\Xcal\rp^{\Omega_1} \times \lp\Xcal^{\Omega_2}\rp^{\Omega_1} \times \PrQm(\Omega_2)^{\Omega_1} \ar@{->>}_-{\id \times \pf}[d] \ar^-{(\_\circ\_)\times\id}[rrr] &&& \lp\Ycal^{\Omega_2}\rp^{\Omega_1} \times \PrQm(\Omega_2)^{\Omega_1} \ar^-{\pf}[d]\\
            \lp\Ycal^\Xcal\rp^{\Omega_1} \times \PrPF(\Xcal)^{\Omega_1} \ar_-{\pf}[rrr] &&& \PrQm(\Ycal)^{\Omega_1}.
    }\]
    The vertical left map is surjective by definition of $\PrPF(\Xcal)^{\Omega_1}$.
    This shows that:
    \[\pf\lp\lp \Ycal^\Xcal\rp^{\Omega_1} \times \PrPF(\Xcal)^{\Omega_1}\rp \ins \pf\lp \lp\Ycal^{\Omega_2}\rp^{\Omega_1} \times \PrQm(\Omega_2)^{\Omega_1}\rp =:\PrPF(\Ycal)^{\Omega_1}. \]
    This shows that the bottom push-forward map lands in $\PrPF(\Ycal)^{\Omega_1}$, showing our claim.\\
    The same arguments hold for $\PrPf$ instead of $\PrPF$ via the diagram:
        \[\xymatrix{
            \lp\Ycal^\Xcal\rp^{\Omega_1} \times \lp\Xcal^{\Omega_2}\rp^{\Omega_1} \times \PrQm(\Omega_2)^{\one} \ar@{->>}_-{\id \times \pf}[d] \ar^-{(\_\circ\_)\times\id}[rrr] &&& \lp\Ycal^{\Omega_2}\rp^{\Omega_1} \times \PrQm(\Omega_2)^{\one} \ar^-{\pf}[d]\\
            \lp\Ycal^\Xcal\rp^{\Omega_1} \times \PrPf(\Xcal)^{\Omega_1} \ar_-{\pf}[rrr] &&& \PrQm(\Ycal)^{\Omega_1}.
    }\]
    This shows the claims.
    \end{proof}
\end{Lem}

\begin{Rem}
    \label{rem:push-forward-R-K}
    Note that the proof of Lemma \ref{lem:push-forward-R-K} would not work for $\PrpF$ (or $\Prpf$) instead of $\PrPF$
    because in the following diagram the top composition map would not be well-defined:
        \[\xymatrix{
            \lp\Ycal^\Xcal\rp^{\Omega_1} \times \lp\Xcal^{\Omega_2}\rp^{\one} \times \PrQm(\Omega_2)^{\Omega_1} \ar@{->>}_-{\id \times \pf}[d] \ar^-{(\_\circ\_)\times\id}[rrr] &&& \lp\Ycal^{\Omega_2}\rp^{\one} \times \PrQm(\Omega_2)^{\Omega_1} \ar^-{\pf}[d]\\
            \lp\Ycal^\Xcal\rp^{\Omega_1} \times \PrpF(\Xcal)^{\Omega_1} \ar_-{\pf}[rrr] &&& \PrQm(\Ycal)^{\Omega_1},
    }\]
 since for $f \in \lp\Ycal^\Xcal\rp^{\Omega_1}$ and $ X \in \lp\Xcal^{\Omega_2}\rp^{\one}$ one would have that:
 \[ (f \circ  X)(\omega_1)(\omega_2)= f(\omega_1)( X(\omega_2))  \]
 depended on $\omega_1 \in \Omega_1$, which one needs to avoid for the image to land in:
 \[\pf\lp \lp\Ycal^{\Omega_2}\rp^{\one} \times \PrQm(\Omega_2)^{\Omega_1}\rp =:\PrpF(\Ycal)^{\Omega_1}.\]
\end{Rem}

\begin{Thm}
    \label{thm:qus-S=P-in-R}
    Let $(\Xcal,\Xcal^\Omega)$ be a quasi-universal space. Then we have the equality of quasi-universal spaces:
    \[\Prpf(\Xcal,\Xcal^\Omega) = \PrPf(\Xcal,\Xcal^\Omega)  = \PrpF(\Xcal,\Xcal^\Omega) = \PrPF(\Xcal,\Xcal^\Omega).\]
    If, furthermore, $\nu=U[0,1]$ is the uniform distribution on $\Omega\cong[0,1]$ then we can represent:
    \[ \PrPf(\Xcal,\Xcal^\Omega) =\lC  X_*\nu\st X \in \Xcal^\Omega \rC, \qquad
    \PrPf(\Xcal,\Xcal^\Omega)^\Omega =\lC  X_*\phi_*\nu \st  X \in \Xcal^\Omega,
    \phi \in \lp\Omega^\Omega \rp^\Omega\rC\footnotemark.  \]
    In particular, the following map is a well-defined (surjective) quotient map of quasi-universal spaces:
    \footnotetext{We could here also directly absorb $\phi$ into $ X$ and write $ X \in \lp\Xcal^\Omega \rp^\Omega$. We could also require $ X \in \Xcal^\Omega$ and restrict $\phi$ further to be a conditional quantile function, see proof.}
      \[ \Xcal^\Omega \srj \PrPf(\Xcal,\Xcal^\Omega) , \qquad  X \mapsto  X_*\nu.\]
    \begin{proof}
        First, it is important to note that by Lemma \ref{lem:univ-qus-G-Q} we have that
        $\nu \in \Giry(\Omega)=\PrQm(\Omega)$.\\
        Next we apply Lemma \ref{lem:K-P-S-R-eq}:
        By Lemma \ref{lem:qus-sample-prod-iso} we have an isomorphism of quasi-universal spaces:
        \[ \Omega \cong \Omega \times \Omega.\]
        By Lemma \ref{lem:K-P-S-R-eq} we get the inclusions and equalities:
        \[ \Prpf(\Xcal)^\Omega=\PrpF(\Xcal)^\Omega \ins \PrPf(\Xcal)^\Omega = \PrPF(\Xcal)^\Omega. \]
        Next we will argue that the following push-forward map is surjective:
        \[ \pf:\; \lp\Omega^\Omega\rp^\Omega \times \lC \nu \rC \inj \lp\Omega^\Omega\rp^\Omega \times \PrQm(\Omega) \to \PrQm(\Omega)^\Omega.\]
        Let $\kappa \in \PrQm(\Omega)^\Omega$ and $\phi \in \lp\Omega^\Omega\rp^\Omega$
        its conditional quantile function (using the fixed isomorphism $\Omega \cong [0,1]$).
        Then it is a well known fact that $\kappa = \phi_*\nu$, see e.g.\ \cite{For21} Appendix G or \cite{Cen82,Kal21}. So the map from above is surjective.\\
        By Lemma \ref{lem:K-P-S-R-eq} we then get the inclusions and equalities:
      \[ \Prpf(\Xcal)^\Omega=\PrPf(\Xcal)^\Omega \ins \PrpF(\Xcal)^\Omega = \PrPF(\Xcal)^\Omega. \]
        This already shows the equalities of quasi-universal spaces:
      \[\Prpf(\Xcal) = \PrPf(\Xcal)  = \PrpF(\Xcal) = \PrPF(\Xcal). \]
    To show that the mentioned map is a quotient map just note that:
    \[  X_*\kappa =  X_*\phi_*\nu = ( X \circ \phi)_*\nu,\]
    for $ X \in \Xcal^\Omega$, $\kappa \in \PrQm(\Omega)^\Omega$, $\phi \in \lp\Omega^\Omega\rp^\Omega$ and $ X \circ \phi \in \lp\Xcal^\Omega\rp^\Omega$ given by:
    \[ ( X \circ \phi)(\omega_1)(\omega_2) =  X\lp\phi(\omega_1)(\omega_2)\rp.\]
    This shows all the claims.
    \end{proof}
\end{Thm}

\begin{Thm}
    \label{thm:univ-qus-GPQKR}
        Let $(\Xcal,\Xcal^\Omega)$ be a standard quasi-universal space (e.g.\ $(\Xcal,\Xcal^\Omega)=(\Omega,\Omega^\Omega)$)
    then we have the equality of quasi-universal spaces:
    \[  \Giry(\Xcal,\Xcal^\Omega)=\PrQm(\Xcal,\Xcal^\Omega)=\PrPf(\Xcal,\Xcal^\Omega)=\PrPF(\Xcal,\Xcal^\Omega)= \PrpF(\Xcal,\Xcal^\Omega)= \Prpf(\Xcal,\Xcal^\Omega).
 \]
 Furthermore, $\Giry(\Xcal,\Xcal^\Omega)$ is also a standard quasi-universal space.
 \begin{proof}
 By Lemma \ref{lem:univ-qus-G-Q} we have the equality of quasi-universal spaces:
        \[\Giry(\Xcal,\Xcal^\Omega)=\PrQm(\Xcal,\Xcal^\Omega).\]
       Then by Lemma \ref{lem:univ-qus-G-P} and Lemma \ref{lem:S-P-R-K-incl} we have the equality and the inclusion of sets:
         \[ \PrPf(\Xcal,\Xcal^\Omega)=\Giry(\Xcal,\Xcal^\Omega), \qquad  \PrPf(\Xcal,\Xcal^\Omega)^\Omega \ins\Giry(\Xcal,\Xcal^\Omega)^\Omega. \]
        By Theorem \ref{thm:qus-S=P-in-R} we already have the equalities of quasi-universal spaces:
    \[\Prpf(\Xcal,\Xcal^\Omega) = \PrPf(\Xcal,\Xcal^\Omega)  =
    \PrpF(\Xcal,\Xcal^\Omega) = \PrPF(\Xcal,\Xcal^\Omega).\]
    To show the reverse inclusion  $\Giry(\Xcal,\Xcal^\Omega)^\Omega \ins \PrPf(\Xcal,\Xcal^\Omega)^\Omega$
    let $\kappa \in \Giry(\Xcal,\Xcal^\Omega)^\Omega$.
    We can then consider, similar to Lemma \ref{lem:univ-meas-G}, an inclusion $i:\; \Xcal \inj \Omega$ with $\Bcal_\Xcal=\Bcal_{\Omega|\Xcal}$ and:
    \[  X:\;\Omega \to \Xcal,\qquad  X|_{\Xcal}:=\id_\Xcal,\qquad  X|_{\Omega\sm\Xcal}:= x_0 \in \Xcal.\]
    Since $\Xcal,\Omega\sm \Xcal \in \Bcal_\Omega$ we have that $ X \in \Fcal(\Bcal_\Xcal)=\Xcal^\Omega$ and $ X \circ i = \id_\Xcal$.
Then $i_*\kappa \in \Giry(\Omega,\Omega^\Omega)^\Omega=\PrQm(\Omega,\Omega^\Omega)^\Omega$ and:
\[ \kappa= X_*(i_*\kappa) \in \PrPf(\Xcal,\Xcal^\Omega)^\Omega.\]
This then shows the first claim.\\
For the claim that $\Giry(\Xcal,\Xcal^\Omega)$ is also a standard quasi-universal space
first note that $\Giry(\Xcal,\Xcal^\Omega)^\Omega=\Fcal\lp\Bcal_{\Giry(\Xcal,\Bcal_\Xcal)} \rp$, where the $\sigma$-algebra is the smallest $\sigma$-algebra that makes all evaluation maps measurable.
It was  shown in \cite{For21} Cor.\ B.49 that $(\Giry(\Xcal,\Bcal_\Xcal),\Bcal_{\Giry(\Xcal,\Bcal_\Xcal)})$ is a universal measurable spaces if $(\Xcal,\Bcal_\Xcal)$ is one.
Then Definition/Lemma \ref{def:univ-qus} shows that:
\[\lp\Giry(\Xcal,\Xcal^\Omega),\Giry(\Xcal,\Xcal^\Omega)^\Omega\rp= \lp\Giry(\Xcal,\Bcal_\Xcal),\Fcal\lp\Bcal_{\Giry(\Xcal,\Bcal_\Xcal)} \rp\rp\]
is a standard quasi-universal space.
    \end{proof}
\end{Thm}
\subsection{A Variant: the Monad of Perfect Probability Measures}
To avoid requiring that all probability measures of $\PrQm(\Omega)$ are perfect one could restrict to
the subset $\Vcal(\Omega) \ins \PrQm(\Omega)$ of all \emph{perfect} probability measures $\nu$
on $(\Omega,\Bcal_\Omega)$ that also lie in $\PrQm(\Omega)$:
\[ \Vcal(\Omega) := \lC \nu:\;\Bcal_\Omega \to [0,1] \text{ perfect prob.\ measure} \st
\forall D \in \Bcal_{\Omega \times \Omega}.\, \lp\omega \mapsto \nu(D_\omega) \rp \in [0,1]^\Omega\rC.\]
We can then define, similar to $\PrPf$:
\begin{align*}
    \Vcal(\Xcal) &:= \lC X_*\nu \st X \in \Xcal^\Omega, \nu \in \Vcal(\Omega) \rC,\\
    \Vcal(\Xcal)^\Omega &:= \lC X_*\nu \st X \in \lp\Xcal^\Omega \rp^\Omega, \nu \in \Vcal(\Omega) \rC.
\end{align*}
Since push-forwards of perfect measures are perfect again, see \cite{Fremlin} 451E, there is no clash in notations here.

One would similar to Lemma \ref{lem:push-forward-R-K} also have a well-defined quasi-measurable push-forward map:
    \[ \pf:\; \Ycal^\Xcal \times \Vcal(\Xcal) \to \Vcal(\Ycal), \qquad (f,\mu) \mapsto f_*\mu.\]

    Furthermore, under $\Omega \times \Omega \cong \Omega$ (and if $\QMS$-products of perfect measures on $\Omega$ are perfect again) we would also have the quasi-measurable map, see Theorem \ref{thm:K-P-Q-markov-prod}:
    \[ \otimes:\; \Vcal(\Xcal)^{\Ycal \times \Zcal} \times \Vcal(\Ycal)^\Zcal \to \Vcal(\Xcal \times \Ycal)^\Zcal,\quad
    (\mu \otimes \nu)(z)(D)  := \int  \mu(y,z)(D^y)\, \nu(z)(dy). \]
Under the same condition $(\Vcal,\delta,\Mbb)$ would be a strong probability monad by the same proofs
as used for $\PrPf$, see Theorem \ref{thm:K-P-R-strong-monad}.

If one, in addition, wanted to have a similar description of the induced $\sigma$-algebras as in Lemma \ref{lem:qus-universally-complete}  and Theorem \ref{thm:countably-separated}
one then would need to require that $\Bcal_\Omega$ is complete w.r.t.\
$\Vcal(\Omega)$: %
\[ \Bcal_\Omega = (\Bcal_\Omega)_{\Vcal(\Omega)}. \]
Finally, the same arguments would hold if we replaced the word ``perfect'' with ``countably compact'' everywhere, see \cite{Fremlin} 451J, 451K, 452R, 454A.

If one wants to make use of the disintegration theorems, see Theorem \ref{cor:markov-kernel-disint-meas-II},
one would need to require that $(\Omega,\Bcal_\Omega)$ is universally countably generated and that all probability measures are perfect. If one also wants $\Bcal_\Omega$ to separate the points of $\Omega$ one
necessarily arrives at that $(\Omega,\Bcal_\Omega)$ is a universal measurable space, see Lemma \ref{deflem:univ-meas}.
If one also wants the convenient descriptions of the induced $\sigma$-algebras, see Lemma \ref{lem:qus-universally-complete} and Theorem \ref{thm:countably-separated}, one needs $(\Omega,\Bcal_\Omega)$ to be universally complete, i.e.\ $\Bcal_\Omega =(\Bcal_\Omega)_\Giry$.
So we arrive at requiring that $(\Omega,\Bcal_\Omega)$ is a universally closed universal measurable space,
under which certainly the universal completions of Polish spaces are the best behaved ones.
The isomorphism $\Omega \times \Omega \cong \Omega$ shows that $\Omega$ must be infinite.
If one wants to allow for non-discrete probability measures one even arrives at uncountable Polish spaces $(\Omega,\Bcal_\Omega)$ endowed with their $\sigma$-algebras of all universally measurable subsets,
and thus at the \emph{category of quasi-universal spaces} $\QUS$.

\subsection{Affinity and Commutativity of the Probability Monads}
\label{qms:sec:app-affine-comm}

\Cref{qms:rem:monad-properties} fixed the four properties this paper tracks for each of its
probability monads. Strength is proved for each monad where it is stated. This subsection supplies
the two remaining ones --- affinity, which is uniform and elementary, and commutativity, which for
every monad on $\QMS$ rests on the Fubini property.

\begin{Lem}[Affinity, see \Cref{qms:lem:affine:prf}]
    \label{qms:lem:affine}
    Let $\one$ denote the terminal object, i.e.\ the one-point quasi-measurable space
    $(\lC * \rC, \lC * \rC^\Omega)$, with induced $\sigma$-algebra
    $\Bcal_\one = \lC \emptyset, \lC * \rC\rC$. Then:
    \begin{align}
        \Giry(\one) = \PrQm(\one) = \PrPF(\one) = \PrpF(\one) = \PrPf(\one) = \Prpf(\one)
        = \lC \delta_* \rC \cong \one .
    \end{align}
    In particular each of the monads $\Giry$, $\PrQm$, $\PrPF$, $\PrpF$, $\PrPf$, $\Prpf$ is
    \emph{affine}.
\end{Lem}

\begin{Prf}[Affinity]{qms:lem:affine}
    A probability measure on $(\lC * \rC, \Bcal_\one)$ is determined by its values on
    $\emptyset$ and $\lC * \rC$, which are forced to be $0$ and $1$. So
    $\Giry(\one)$ is the one-element set $\lC \delta_* \rC$.
    By \Cref{lem:K-P-R-Q-G-qm-incl} every one of $\PrQm(\one)$, $\PrPF(\one)$, $\PrpF(\one)$,
    $\PrPf(\one)$, $\Prpf(\one)$ is a subset of $\Giry(\one)$, and each contains
    $\delta_* = \delta_\one(*)$ because $\delta$ is the unit of the corresponding monad. Hence all
    six sets coincide with $\lC \delta_* \rC$.

    For the quasi-measurable structure: the only map $\Omega \to \lC \delta_* \rC$ is the constant
    one, which is admissible in each case, so all six spaces are the one-point quasi-measurable
    space, which is terminal in $\QMS$ by \Cref{qms:thm:co-complete}. A monad $T$ with
    $T(\one) \cong \one$ is affine by definition.
\end{Prf}

\begin{Rem}[Why affinity is not automatic]
    \label{qms:rem:affine-not-automatic}
    The argument uses that the elements of $\Giry(\Xcal)$ are \emph{probability} measures. Had we
    worked with sub-probability measures, or with all finite measures, $\Giry(\one)$ would be
    $[0,1]$, respectively $[0,\infty]$, and the monad would fail to be affine. Affinity is what
    makes the marginalisation map $\PrPf(\Xcal \times \Ycal) \to \PrPf(\Xcal)$, induced by the
    projection, agree with integrating out the second coordinate.
\end{Rem}

\begin{Thm}[Commutativity, see \Cref{qms:thm:commutative:prf}]
    \label{qms:thm:commutative}
    Let $(\Omega,\Omega^\Omega,\Bcal_\Omega,\PrPf)$ be a sample space satisfying
    \Cref{qms:def:sample-space-act-3}, the Fubini property \Cref{qms:asm:fubini}, and the product
    isomorphism property \Cref{qms:asm:prod-iso}. Then the four push-forward monads
    \begin{align}
        (\PrPF,\delta,\Mbb), \quad (\PrpF,\delta,\Mbb), \quad
        (\PrPf,\delta,\Mbb), \quad (\Prpf,\delta,\Mbb)
    \end{align}
    on $(\QMS,\times,\one)$ are all \emph{commutative}, equivalently \emph{symmetric monoidal}, in
    the sense of \Cref{qms:def:commutative-monad}.
    If, in addition, $\PrQm(\Xcal)=\PrPf(\Xcal)$ for all spaces involved --- which holds for
    standard quasi-universal spaces by \Cref{thm:univ-qus-GPQKR} --- then $(\PrQm,\delta,\Mbb)$ is
    commutative as well.
\end{Thm}

\begin{Rem}[Why $\PrQm$ needs the extra hypothesis]
    \label{qms:rem:comm-PrQm}
    The Fubini property \Cref{qms:asm:fubini} quantifies over $P_1,P_2 \in \PrPf$, and
    \Cref{qms:thm:fubini} is correspondingly stated for push-forward measures
    $P(X) \in \PrPf(\Xcal)$. For $\PrQm$ the two measures involved need not be push-forwards of
    members of $\PrPf$, and classical Tonelli does not close the gap either, because the integrand
    is only $\Bcal_{\Xcal \times \Ycal}$-measurable and $\Bcal_{\Xcal \times \Ycal}$ can be
    strictly larger than $\Bcal_\Xcal \otimes \Bcal_\Ycal$. The hypothesis that removes the difference is
    $\PrQm(\Xcal)=\PrPf(\Xcal)$ for the spaces involved, not the weaker $\PrQm(\Omega)=\PrPf$:
    equality of the \emph{pools} on $\Omega$ does not make every member of $\PrQm(\Xcal)$ a
    push-forward. For a witness take $\Xcal=\lC 0,1\rC$ carrying only its two constant random
    variables, inside the sample space of \Cref{qms:def:sample-space-act-4}, so that
    $\PrQm(\Omega)=\PrPf$ holds; then $\Bcal_\Xcal$ is the power set and
    $\PrPf(\Xcal)=\lC \delta_0,\delta_1\rC$, while
    $\tfrac12\lp\delta_0+\delta_1\rp \in \PrQm(\Xcal)$.
    On \emph{standard} quasi-universal spaces the required equality does hold, by
    \Cref{thm:univ-qus-GPQKR}, so nothing is lost there.
\end{Rem}

\begin{Prf}[Commutativity]{qms:thm:commutative}
    Write $T$ for any one of the four push-forward monads, or for $\PrQm$ under the additional
    hypothesis that $\PrQm(\Xcal)=\PrPf(\Xcal)$ for the spaces involved, under which every measure
    the argument meets is a push-forward of a member of $\PrPf$ and the argument below applies
    verbatim. By \Cref{thm:Q-strong-monad} and
    \Cref{thm:K-P-R-strong-monad} each $T$ is strong, with left strength
    $\tau_{\Xcal,\Ycal}(x,\mu) = \delta_x \otimes \mu$ and right strength
    $\rho_{\Xcal,\Ycal}(\mu,y) = \mu \otimes \delta_y$.

    By \Cref{prp:Q-product-kernels} for $T=\PrQm$, and by \Cref{thm:K-P-Q-markov-prod} for the four
    push-forward monads --- the latter using \Cref{qms:asm:prod-iso} --- the product of Markov
    kernels
    \begin{align}
        \otimes:\; T(\Xcal)^{\Ycal \times \Zcal} \times T(\Ycal)^\Zcal \to
        T(\Xcal \times \Ycal)^\Zcal
    \end{align}
    is well defined and quasi-measurable. Specialising to $\Zcal=\one$ and to a $\mu$ that does not
    depend on $\Ycal$ gives a natural transformation
    $\otimes_{\Xcal,\Ycal}:\, T(\Xcal) \times T(\Ycal) \to T(\Xcal \times \Ycal)$ with
    \begin{align}
        (\mu \otimes \nu)(D) = \int \mu(D^y)\, \nu(dy), \qquad D \in \Bcal_{\Xcal \times \Ycal}.
    \end{align}

    The two composites in \Cref{qms:def:commutative-monad} evaluate, on $(\mu,\nu)$ and
    $D \in \Bcal_{\Xcal\times\Ycal}$, to the two iterated integrals
    \begin{align}
        \begin{split}
            \int \mu(D^y)\, \nu(dy) &= \iint \I_D(x,y)\, \mu(dx)\, \nu(dy),\\
            \int \nu(D_x)\, \mu(dx) &= \iint \I_D(x,y)\, \nu(dy)\, \mu(dx),
        \end{split}
    \end{align}
    which agree by the Fubini theorem \Cref{qms:thm:fubini}, applied to the non-negative function
    $h = \I_D$ on $\Xcal \times \Ycal$. Note that this is exactly the point at which
    \Cref{qms:asm:fubini} is needed: $D$ ranges over the induced $\sigma$-algebra
    $\Bcal_{\Xcal \times \Ycal}$, which can be strictly larger than
    $\Bcal_\Xcal \otimes \Bcal_\Ycal$, so classical Tonelli does not by itself suffice.

    Associativity and unitality of $\otimes$ follow from the same theorem applied to triple
    products, and $\delta:\, \one \to T(\one)$ is an isomorphism by \Cref{qms:lem:affine}. So by
    \Cref{qms:rem:otimes-commutative} and \cite{Koc72} Thm.\ 2.3 the induced strengths make $T$ a
    commutative monad, and by \cite{Koc70} commutative and symmetric monoidal agree for a strong
    monad on a symmetric monoidal category.
\end{Prf}

\begin{Rem}[Commutativity is exactly Fubini]
    \label{qms:rem:comm-is-fubini}
    The proof above says something worth stating on its own: for a strong probability monad on
    $\QMS$, \emph{commutativity and the Fubini property are the same condition}, read once
    categorically and once in integrals. This is why \Cref{qms:asm:fubini} was isolated as an
    assumption on the sample space in \Cref{sec:pf-monad} rather than left implicit: it is not a
    technical convenience but the precise hypothesis under which
    $\PrPf(\Xcal) \times \PrPf(\Ycal) \to \PrPf(\Xcal \times \Ycal)$ is well defined independently
    of the order of the factors, and hence under which the Kleisli category of $\PrPf$ can be a
    Markov category, see \Cref{thm:Markov-category}.
\end{Rem}

\section{Appendix - Countably Separated Spaces}
\begin{Prp}[See \Cref{prp:perf-count-sep:prf}]
    \label{prp:perf-count-sep}
    Let $(\Omega,\Bcal_\Omega)$ be a measurable space and
     $\PrPf$ be a set of \emph{perfect} probability measures
    on $(\Omega,\Bcal_\Omega)$.
    Let $(\Xcal,\Bcal_\Xcal)$ be a \emph{countably separated} measurable space and $\Ecal \ins \Bcal_\Xcal$
    be any countable subset that separates the points of $\Xcal$.
    Let $\Fcal \ins \Meas\lp (\Omega,\Bcal_\Omega), (\Xcal,\Bcal_\Xcal)\rp$ be any subsets of measurable maps.
    We abbreviate the set of push-forward measures:
    \[\PrPf(\Fcal):= \lC  X_*\mu:\, \Bcal_\Xcal \to [0,1]\st X \in \Fcal, \mu \in \PrPf \rC.\]
    Then we have the equalities:
    \[ \lp \Ecal \rp_{\PrPf(\Fcal)}=\lp \Bcal_\Xcal \rp_{\PrPf(\Fcal)} = \Bcal(\Fcal)_{\PrPf(\Fcal)},\]
    where the index $\PrPf(\Fcal)$ refers to the intersections of all (Lebesgue) completions w.r.t.\
    $\nu \in \PrPf(\Fcal)$.\\
    If, furthermore, we have that $\Bcal_\Omega$ is complete w.r.t.\ $\PrPf$, i.e. that:
    \[\Bcal_\Omega = \lp \Bcal_\Omega \rp_\PrPf:= \bigcap_{\mu \in \PrPf} \lp\Bcal_\Omega \rp_\mu,\]
        then we also get the equality:
    \[ \Bcal(\Fcal) = \lp \Ecal \rp_{\PrPf(\Fcal)}.\]

\end{Prp}
\label{qms:sec:app-count-sep}

\begin{Lem}
\label{qms:lem:push-pull-s-alg}
Let $j:\, \Xcal \to \Ycal$ be an injective map and $\Bcal$ a $\sigma$-algebra on $\Xcal$.
Then we have:
        \begin{align} j^*j_*\Bcal = \Bcal.\end{align}
  \begin{proof}
        First, let $C \in j^*j_*\Bcal$. Then there exists $B \in j_*\Bcal$ such that $C=j^{-1}(B)$.
        By definition of $j_*\Bcal$ we have that:
        \begin{align}C=j^{-1}(B) \in \Bcal.\end{align}
        This shows:
        \begin{align}j^*j_*\Bcal \ins \Bcal.\end{align}
        Now let $A \ins \Bcal$ then by the injectivity of $j$:
        \begin{align} j^{-1}(j(A)) = A \in \Bcal.\end{align}
        So by definition of $j_*\Bcal$ we get that $j(A) \in j_*\Bcal$ and thus:
        \begin{align} A = j^{-1}(j(A)) \in j^*j_*\Bcal.\end{align}
         This shows:
          \begin{align}\Bcal \ins j^*j_*\Bcal.\end{align}
         and thus the claim.
  \end{proof}
\end{Lem}

\begin{Lem}
    \label{qms:lem:perf-count-sep}
    Let $(\Omega,\Bcal_\Omega,P)$ be a perfect probability space and $(\Xcal,\Bcal_\Xcal)$
    be a countably separated measurable space and $X:\; \Omega \to \Xcal$ a measurable map.
    Then we have:
    \begin{align} X_*(\Bcal_\Omega)_P = (\Bcal_\Xcal)_{X_*P},  \end{align}
    where the index refers to the (Lebesgue) completion w.r.t.\ the corresponding probability measure.
    \begin{proof}
        By the measurability of $X$ we always have the inclusion:
        \begin{align} X_*(\Bcal_\Omega)_P \sni (\Bcal_\Xcal)_{X_*P}.\end{align}
        Since $(\Xcal,\Bcal_\Xcal)$ is countably separated by \cite{Bog07} Thm.\ 6.5.7
        there exists an injective $\Bcal_\Xcal$-$\Bcal_\R$-measurable map:
        \begin{align}j:\; \Xcal \to \R.\end{align}
        So we have the chain of inclusions:
        \begin{align} j_*X_*(\Bcal_\Omega)_P \sni j_*(\Bcal_\Xcal)_{X_*P} \sni (\Bcal_\R)_{j_*X_*P}. \end{align}
        Since $P$ is perfect we have the equality:
        \begin{align} j_*X_*(\Bcal_\Omega)_P =(\Bcal_\R)_{j_*X_*P},\end{align}
        which then implies the equality:
        \begin{align} j_*X_*(\Bcal_\Omega)_P = j_*(\Bcal_\Xcal)_{X_*P}. \end{align}
        Since $j$ is injective we get by \Cref{qms:lem:push-pull-s-alg}:
        \begin{align}X_*(\Bcal_\Omega)_P = j^*j_*X_*(\Bcal_\Omega)_P = j^*j_*(\Bcal_\Xcal)_{X_*P}=(\Bcal_\Xcal)_{X_*P}. \end{align}
        This shows the claim.
    \end{proof}
\end{Lem}

\begin{Lem}
    \label{qms:lem:push-completion}
    Let $(\Omega,\Bcal_\Omega)$ be a measurable space and
     $\Pcal$ be a set of  probability measures
    on $(\Omega,\Bcal_\Omega)$ such that:
    \begin{align} \Bcal_\Omega \sni \lp \Bcal_\Omega \rp_\Pcal:= \bigcap_{P \in \Pcal} \lp\Bcal_\Omega \rp_P.  \end{align}
    Let $(\Xcal,\Bcal_\Xcal)$ be a measurable space
    and $\Hcal \ins \Meas\lp (\Omega,\Bcal_\Omega), (\Xcal,\Bcal_\Xcal)\rp$ be any subset of measurable maps.
    Let $\Qcal$ be any set of probability measures on
    $(\Xcal,\Bcal_\Xcal)$ that contains the push-forward probability measures:
    \begin{align}\Qcal \sni \Hcal_*\Pcal:= \lC  X_*P:\, \Bcal_\Xcal \to [0,1]\st X \in \Hcal, P \in \Pcal \rC.\end{align}
    Then we have the inclusions:
    \begin{align}
      \Bcal_\Xcal \ins  (\Bcal_\Xcal)_\Qcal \ins \Bcal(\Hcal) := \lC A \ins \Xcal\st \forall  X \in \Hcal.\, X^{-1}(A) \in \Bcal_\Omega  \rC,
    \end{align}
    and thus:
    \begin{align}
        \Hcal &\ins \Meas\lp (\Omega,\Bcal_\Omega), (\Xcal,\Bcal(\Hcal))\rp  \\
              &\ins \Meas\lp (\Omega,\Bcal_\Omega), (\Xcal,(\Bcal_\Xcal)_\Qcal)\rp \\ &
              \ins \Meas\lp (\Omega,\Bcal_\Omega), (\Xcal,\Bcal_\Xcal)\rp.
    \end{align}
    \begin{proof}
        We have the following chain of inclusion:
        \begin{align*}
             \Bcal(\Fcal)
            &= \lC A \ins \Xcal\st \forall  X \in \Fcal.\, X^{-1}(A) \in \Bcal_\Omega  \rC \\
            &= \bigcap_{ X \in \Fcal}  X_*\Bcal_\Omega \\
            &\sni \bigcap_{ X \in \Fcal}  X_*\bigcap_{P \in \PrPf} \lp\Bcal_\Omega \rp_P \\
            &= \bigcap_{ X \in \Fcal} \bigcap_{P \in \PrPf}  X_*\lp\Bcal_\Omega \rp_P \\
            &\sni \bigcap_{ X \in \Fcal} \bigcap_{P \in \PrPf} \lp\Bcal_\Xcal \rp_{ X_*P} \\
            &= \bigcap_{Q \in \Fcal_*\PrPf} \lp\Bcal_\Xcal \rp_{Q}\\
            &\sni \bigcap_{Q \in \PrQm} \lp\Bcal_\Xcal \rp_{Q}\\
            &= \lp \Bcal_\Xcal \rp_{\PrQm} \\
            & \sni \Bcal_\Xcal.
    \end{align*}
    This implies all claims.
    \end{proof}
\end{Lem}

\section{Appendix - Construction of Well-Behaved Sample Spaces}
\label{qms:sec:app-sample-spaces}

\begin{Prf}[Sample spaces from measurable spaces]{qms:thm:spl-spc-gen}
    It is clear that \Cref{qms:def:sample-space-act-1} and \Cref{qms:def:sample-space-act-2} are satisfied.
    Now, applying the functor $\Fcal$ to $\Theta$ shows that we get an isomorphism of quasi-measurable spaces:
    \begin{align}
        \Theta:\, (\Omega,\Fcal(\Bcal_\Omega)) \cong (\Omega \times \Omega, \Fcal(\Bcal_\Omega)\times\Fcal(\Bcal_\Omega)).
    \end{align}
    Since, by assumption, we have $\Omega^\Omega=\Fcal(\Bcal_\Omega)$ we get the isomorphism of quasi-measurable spaces:
    \begin{align}
        \Theta:\, (\Omega,\Omega^\Omega) \cong (\Omega \times \Omega, \Omega^\Omega\times\Omega^\Omega),
    \end{align}
    which already implies \Cref{qms:asm:prod-iso}.
   Further applying the functor $\Bcal$ gives us the isomorphism of measurable spaces:
    \begin{align}
        \Theta:\, (\Omega,\Bcal_\Omega) \cong (\Omega \times \Omega, \Bcal_{\Omega\times\Omega}).
    \end{align}
    Comparing this with the first isomorphism of measurable spaces shows that:
    \begin{align}
        \Bcal_{\Omega\times\Omega} &= \Bcal_\Omega \otimes \Bcal_\Omega.
    \end{align}
    Furthermore, \Cref{qms:def:sample-space-act-3} \Cref{qms:eq:sample-space-act-3-a} and \Cref{qms:eq:sample-space-act-3-b} are also clear. For \Cref{qms:eq:sample-space-act-3-c} we put:
    $P:=\Theta_*^{-1}(P_1 \otimes P_2) \in \PrPf$ and $\Phi:=\Theta \in (\Omega \times \Omega)^\Omega$,
    which shows:
    \begin{align}
        P_1 \otimes P_2 = P(\Phi) \in \PrPf(\Omega\times\Omega).
    \end{align}
    So all points of \Cref{qms:def:sample-space-act-3} are shown.
    Finally, by \Cref{qms:prp:fubini-crit} we get \Cref{qms:asm:fubini}.
\end{Prf}

\begin{Prf}[Sample spaces from completing measurable spaces]{qms:thm:spl-spc-cmpl}
        It is clear that \Cref{qms:def:sample-space-act-1} and \Cref{qms:def:sample-space-act-2} are satisfied.

        First note that we trivially have the inclusion:
        \begin{align}
            \Fcal(\Bcal_\Omega) =  \Omega^\Omega &=  \Meas\lp (\Omega,(\Bcal_\otimes)_\PrPf),(\Omega,(\Bcal_\otimes)_\PrPf) \rp \\ &\ins \Meas\lp (\Omega,(\Bcal_\otimes)_\PrPf),(\Omega,\Bcal_\otimes) \rp \\
                                                 &= \Fcal(\Bcal_\otimes).
        \end{align}
        For showing the reverse inclusion let $\Phi \in \Meas\lp (\Omega,(\Bcal_\otimes)_\PrPf),(\Omega,\Bcal_\otimes) \rp$. By \Cref{qms:eq:PrPf-compl-cond} we have that $\Phi_*\PrPf \ins \PrPf$.
        By \Cref{qms:lem:map-completetions} this implies that $\Phi$ is also $\lp (\Bcal_\otimes)_\PrPf\rp_\PrPf$-$(\Bcal_\otimes)_\PrPf$-measurable. Since $\lp (\Bcal_\otimes)_\PrPf\rp_\PrPf = (\Bcal_\otimes)_\PrPf$ this shows that $\Phi \in \Meas\lp (\Omega,(\Bcal_\otimes)_\PrPf),(\Omega,(\Bcal_\otimes)_\PrPf) \rp$, which we wanted to show.
        So we have the equality:
        \begin{align}
            \Fcal(\Bcal_\Omega) =  \Omega^\Omega &= \Fcal(\Bcal_\otimes).
        \end{align}
            Now, applying the functor $\Fcal$ to $\Theta$ shows that we get an isomorphism of quasi-measurable spaces:
    \begin{align}
        \Theta:\, (\Omega,\Fcal(\Bcal_\otimes)) \cong (\Omega \times \Omega, \Fcal(\Bcal_\otimes)\times\Fcal(\Bcal_\otimes)).
    \end{align}
    Since we have shown that $\Omega^\Omega=\Fcal(\Bcal_\otimes)$ we get the isomorphism of quasi-measurable spaces:
    \begin{align}
        \Theta:\, (\Omega,\Omega^\Omega) \cong (\Omega \times \Omega, \Omega^\Omega\times\Omega^\Omega),
    \end{align}
    which already implies \Cref{qms:asm:prod-iso} and that $\Theta \in (\Omega \times \Omega)^\Omega$.
   Further applying the functor $\Bcal$ gives us the isomorphism of measurable spaces:
    \begin{align}
        \Theta:\, (\Omega,\Bcal_\Omega) \cong (\Omega \times \Omega, \Bcal_{\Omega\times\Omega}).
    \end{align}
   \Cref{qms:lem:map-completetions} together with the original isomorphism $\Theta$ of measurable spaces implies the isomorphism of measurable spaces:
    \begin{align}
        \Theta:\, (\Omega,\Bcal_\Omega) = (\Omega,(\Bcal_\otimes)_{\PrPf}) \cong (\Omega \times \Omega, (\Bcal_\otimes\otimes\Bcal_\otimes)_{\Theta_*\PrPf}).
    \end{align}
    Both isomorphisms together show the equality:
    \begin{align}
            \Bcal_{\Omega\times\Omega} = (\Bcal_\otimes\otimes\Bcal_\otimes)_{\Theta_*\PrPf}.
    \end{align}
    Since $\Theta \in (\Omega \times \Omega)^\Omega$ we have the inclusion: $\Theta_*\PrPf \ins \PrPf(\Omega \times \Omega)$. To also show the reverse inclusion let $P(\Psi) \in \PrPf(\Omega \times \Omega)$ with $P \in \PrPf$ and $\Psi \in (\Omega \times \Omega)^\Omega$. Then we get:
    \begin{align}
       P(\Psi) = \Theta_*P(\Theta^{-1} \circ \Psi),
    \end{align}
    with $\Theta^{-1} \circ \Psi \in \Omega^\Omega=\Fcal(\Bcal_\otimes)$. With \Cref{qms:eq:PrPf-compl-cond} we see that $P(\Theta^{-1} \circ \Psi) \in \PrPf$, which then implies: $P(\Psi) \in \Theta_*\PrPf$.
    This shows the equality of sets of probability measures living on $\Bcal_{\Omega\times\Omega} = (\Bcal_\otimes\otimes\Bcal_\otimes)_{\Theta_*\PrPf}$:
    \begin{align}
        \Theta_*\PrPf = \PrPf(\Omega \times \Omega).
    \end{align}
    So we get the chain of inclusions:
    \begin{align}
        \Bcal_\otimes\otimes\Bcal_\otimes \ins \Bcal_\Omega \otimes \Bcal_\Omega &\ins \Bcal_{\Omega \times \Omega} = (\Bcal_\otimes\otimes\Bcal_\otimes)_{\Theta_*\PrPf} = (\Bcal_\otimes\otimes\Bcal_\otimes)_{\PrPf(\Omega\times\Omega)}.
\end{align}
This shows that we get the  equality:
    \begin{align}
 \Bcal_{\Omega \times \Omega} = (\Bcal_\Omega\otimes\Bcal_\Omega)_{\PrPf(\Omega\times\Omega)}.
    \end{align}
    We are left to show the points from \Cref{qms:def:sample-space-act-3}.
    Since $\Bcal_\Omega=(\Bcal_\otimes)_\PrPf \ins (\Bcal_\otimes)_P$ for every $P \in \PrPf$ we get
    \Cref{qms:eq:sample-space-act-3-a}.
    To show \Cref{qms:eq:sample-space-act-3-b} we want to apply \Cref{qms:eg:markov-kernel-completions}. For this we need to show that $P_1 \otimes P_2 \in \Theta_*\PrPf$ for $P_1,P_2 \in \PrPf$,
    which is equivalent to
    $P:=\Theta_*^{-1}(P_1\otimes P_2) \in \PrPf$, which holds by \Cref{qms:eq:PrPf-prod-cond}.
    So \Cref{qms:eq:sample-space-act-3-b} holds. We further get that $P_1\otimes P_2 = P(\Theta)$ as probability measures on $\Bcal_\otimes\otimes\Bcal_\otimes$. But then they also must agree on the completion $(\Bcal_\otimes\otimes\Bcal_\otimes)_{\Theta_*\PrPf}=\Bcal_{\Omega\times\Omega}$. This shows that \Cref{qms:eq:sample-space-act-3-c} also holds.
    So all points from \Cref{qms:def:sample-space-act-3} are shown.
    Finally, by \Cref{qms:prp:fubini-crit} we get \Cref{qms:asm:fubini}.
\end{Prf}

\begin{Prf}[Countably separated product spaces as sample space for perfect probability distributions]{qms:cor:sample-space-count-sep-perf}
This follows directly from \Cref{qms:thm:spl-spc-cmpl} with the following abbreviations:
    \begin{align}
       \Bcal_0 &:= \sigma(\Ecal_0), & \Bcal_\otimes &:= \bigotimes_{n \in \N} \Bcal_0.
    \end{align}
Note that we have:
    \begin{align}
       \Ecal \ins \Bcal_\otimes \ins \Bcal_\Omega.
    \end{align}
and that $\Ecal$ is countable and separates the points of $\Omega$.

With the above notations we acknowledge that the zip-locker map $\Theta:\, \Omega \to \Omega \times \Omega$ is a bijection, where $\Theta$ is $\Bcal_\otimes$-$(\Bcal_\otimes \otimes \Bcal_\otimes)$-measurable and $\Theta^{-1}$ is $(\Bcal_\otimes \otimes \Bcal_\otimes)$-$\Bcal_\otimes$-measurable, which can be checked on finite cylinder sets.
So $\Theta$ constitutes an isomorphism of measurable spaces:
\begin{align}
        \Theta:\, (\Omega,\Bcal_\otimes) \cong (\Omega \times \Omega, \Bcal_\otimes \otimes \Bcal_\otimes).
\end{align}
The rest follows directly from \Cref{qms:thm:spl-spc-cmpl}.
\end{Prf}

\section{Appendix - \Patchable\ and \Sturdy\ Quasi-Measurable Spaces}

\label{sec:patchable-qms}

The point of this appendix, stated once as a proposition, is that \patchability\ is exactly the
property that was dropped from \cite{Heu17} Def.\ 7, and that putting it back recovers quasi-Borel
spaces on the nose.

\begin{Prp}[Recovery of quasi-Borel spaces]
    \label{prp:pqms-is-qbs}
    Let the sample space be $(\Omega,\Omega^\Omega)=(\R,\Fcal(\Bcal_\R))$, i.e.\ the real numbers
    with the Borel $\sigma$-algebra and all Borel measurable maps as admissible random variables,
    see \Cref{qms:eg:meas-sample-space}. Then the category $\PQMS$ of \patchable\ quasi-measurable
    spaces is isomorphic to the category $\QBS$ of quasi-Borel spaces of \cite{Heu17} Def.\ 7, and
    the isomorphism is the identity on underlying sets and on admissible random variables.
    \begin{proof}
        A quasi-Borel space in the sense of \cite{Heu17} Def.\ 7 is a set $\Xcal$ together with a
        set $\Xcal^\R \ins \Sets(\R,\Xcal)$ satisfying three conditions: it contains the constant
        maps; it is closed under precomposition with Borel measurable $\Phi:\,\R \to \R$; and it is
        closed under gluing along countable Borel partitions of $\R$. The first two are
        \Cref{qms:def:quasi-measurable-spaces} for this sample space, and the third is
        \Cref{def:patchable}. Morphisms agree on both sides, being the maps that preserve
        admissible random variables. So the two definitions coincide term by term.
    \end{proof}
\end{Prp}

\begin{Rem}[Which monad is the one of \cite{Heu17}]
    \label{rem:qbs-monad}
    Under \Cref{prp:pqms-is-qbs} the probability monad of \cite{Heu17} corresponds to $\PrpF$: its
    elements are the push-forwards $\alpha_*\mu$ of a Borel probability measure $\mu$ on $\R$ along
    an admissible $\alpha \in \Xcal^\R$, and its admissible random variables hold $\alpha$ fixed
    while letting $\mu$ vary --- which is exactly the row of $\PrpF$ in \Cref{qms:tab:monads}. This
    is the sense in which the present paper generalises \cite{Heu17}: the sample space is no longer
    fixed to $\R$, \patchability\ is not imposed, and the monad in the foreground is $\PrPf$ rather
    than $\PrpF$, the two agreeing on $\QUS$ by \Cref{thm:qus-S=P-in-R}.
\end{Rem}

In this section we will study the 3rd property in the definition of quasi-Borel spaces from Def.\ 7 in \cite{Heu17} transferred to quasi-measurable spaces. We will call quasi-measurable spaces that satisfy that property \emph{\patchable}. Note that until now and throughout the whole Sections \ref{sec:cat-qms} and \ref{qms:sec:app-other-monads} the requirement of \patchability\ was never needed.
The main reason one would anyways impose \patchability\ onto quasi-measurable spaces is stated in Lemma \ref{lem:patch-coprod}, namely to preserve countable coproducts when going from measurable spaces to quasi-measurable spaces via the functor $\Fcal$.
We will see in Theorem \ref{thm:patchable-reflexive} that the category of \patchable\ quasi-measurable spaces $\PQMS$ is a reflective subcategory of the category of quasi-measurable spaces $\QMS$. Furthermore,
it will turn out that $\PQMS$ is an exponential ideal of $\QMS$ and thus cartesian closed on its own.

\begin{Rem}
For the purpose of this section we will in the following assume
that the sample space $(\Omega,\Omega^\Omega,\Bcal_\Omega)$ satisfies:
\[ \Omega^\Omega=\Fcal(\Bcal_\Omega) = \Meas\lp(\Omega,\Bcal_\Omega),(\Omega,\Bcal_\Omega)\rp.  \]
\end{Rem}

\subsection{\Patchable\ Quasi-Measurable Spaces}

Here we shortly define the notion of \emph{\patchable} quasi-measurable spaces. It resembles the 3rd property of quasi-Borel spaces from Def.\ 7 in \cite{Heu17}.

\begin{Def}[\Patchable\ quasi-measurable spaces]
    \label{def:patchable}
    We call a quasi-measurable space $(\Xcal,\Xcal^\Omega)$ \emph{\patchable} if
  for every countable disjoint decomposition: $\Omega = \bigdcup_{i \in I} C_i$ with $C_i \in \Bcal_\Omega$,  $i \in I \ins \N$,
        and for  $X_i \in \Xcal^\Omega$, $i \in I$, we have that:
        \[ X:= \bigdcup_{i \in I} X_i|_{C_i} \in \Xcal^\Omega,\]
        where, more precisely, $X$ is defined as:
        \[X(\omega) := X_i(\omega) \quad \text{ for } \quad \omega \in C_i.\]
\end{Def}

\begin{Eg}
    \label{eg:patchable}
    \begin{enumerate}
        \item If $\Xcal$ is a set and $\Xcal^\Omega:=\Sets(\Omega,\Xcal)$ then
    $(\Xcal,\Xcal^\Omega)$ is a \patchable\ quasi-measurable space.
        \item If $(\Xcal,\Xcal^\Omega)$ is a quasi-measurable space with $\Xcal^\Omega=\Fcal(\Bcal_\Xcal)$
            for any $\sigma$-algebra $\Bcal_\Xcal$ then $(\Xcal,\Xcal^\Omega)$ is patchable.
        \item[] Indeed,
            if $X:= \bigdcup_{i \in I} X_i|_{C_i} \in \Xcal^\Omega$ with $X_i \in \Xcal^\Omega$
            and $D \in \Bcal_\Xcal$ and $C_i \in \Bcal_\Omega$, $i\in I \ins \N$, then:
            \[ X^{-1}(D) = \bigcup_{i \in I} \lp C_i \cap X_i^{-1}(D) \rp \in \Bcal_\Omega.\]
           So $X$ is measurable and thus $X \in \Fcal(\Bcal_\Xcal)=\Xcal^\Omega$.
    \end{enumerate}
\end{Eg}

\subsection{Heredity of \Patchability}

In this subsection we study under which conditions the property of being \patchable\ is inherited by categorical constructions, like products or quotients, etc.

\begin{Lem}
    \label{lem:pull-back-pqms}
    Let $(\Ycal,\Ycal^\Omega)$ be a \patchable\ quasi-measurable space.
    \begin{enumerate}
        \item Let $g:\; \Xcal \to \Ycal$ be any map.
            Then $(\Xcal,g^*\Ycal^\Omega)$ is also \patchable\ quasi-measurable space.
        \item Let $h:\;\Ycal \srj \Zcal$ be any surjective map.
            Then $(\Zcal,h \circ \Ycal^\Omega)$ is also a \patchable\ quasi-measurable space.
    \end{enumerate}
    \begin{proof}
        1.) We clearly have $\Xcal^\one \ins g^*\Ycal^\Omega$.
        If $\Phi \in \Omega^\Omega$ then:
        \[X \in g^*\Ycal^\Omega \quad\implies\quad g \circ X \in \Ycal^\Omega\quad \implies \quad  g \circ X \circ \Phi \in \Ycal^\Omega
        \quad \implies \quad X \circ \Phi \in g^*\Ycal^\Omega.\]
        Let $X_i \in g^*\Ycal^\Omega$ then $g \circ X_i \in \Ycal^\Omega$, which implies
        for every countable measurable disjoint union decomposition $\Omega = \bigdcup_{i\in I} C_i$
        that:
        \[g \circ \bigdcup_{i \in I} X_i|_{C_i}= \bigdcup_{i \in I} (g \circ X_i)|_{C_i} \in \Ycal^\Omega. \]
        Thus: $\bigdcup_{i \in I} X_i|_{C_i} \in g^*\Ycal^\Omega$. This shows the first claim.\\
        2.) Since $h$ is surjective we have $\Zcal^\one \ins h \circ \Ycal^\Omega$.
            We also have:
            \[ (h \circ \Ycal^\Omega) \circ \Omega^\Omega = h \circ (\Ycal^\Omega \circ \Omega^\Omega)  \ins h \circ \Ycal^\Omega.\]
            Similarly to the first point we have that:
        \[\bigdcup_{i \in I} (h \circ Y_i)|_{C_i}= h \circ \bigdcup_{i \in I} Y_i|_{C_i} \in h \circ\Ycal^\Omega, \]
        because if $Y_i \in \Ycal^\Omega$, $i \in I$, then also $\bigdcup_{i \in I} Y_i|_{C_i} \in \Ycal^\Omega$.
        This shows the second claim.
    \end{proof}
\end{Lem}

\begin{Lem}[Products of \patchable\ quasi-measurable spaces are \patchable]
    \label{lem:patchable-product}
        If $(\Xcal_i,\Xcal_i^\Omega)$ are \patchable\ quasi-measurable spaces for $i \in I$ and any index set $I$
            then also their product of quasi-measurable spaces is \patchable:
            \[ \prod_{i\in I}(\Xcal_i,\Xcal_i^\Omega).  \]
            \begin{proof}
                If we consider a measurable disjoint decomposition $\Omega = \bigdcup_{j \in J} C_j$
                with countable $J \ins \N$ and:
                \[ X^j = (X_i^j)_{i \in I} \in \prod_{i\in I} \Xcal_i^\Omega,   \]
                for $j \in J$, then:
                \[ X^j|_{C_j} = (X_i^j|_{C_j} )_{i \in I}, \]
                and thus:
                \[ \bigdcup_{j \in J} X^j|_{C_j} = \lp \bigdcup_{j \in J} X_i^j|_{C_j} \rp_{i \in I}. \]
                Since by assumption $\bigdcup_{j \in J} X_i^j|_{C_j} \in \Xcal_i^\Omega$, for every $i \in I$ we get:
                \[\bigdcup_{j \in J} X^j|_{C_j} = \lp \bigdcup_{j \in J} X^j_i|_{C_j}\rp_{i\in I} \in \prod_{i \in I} \Xcal_i^\Omega.\]
    \end{proof}
\end{Lem}

\begin{Rem}
    \begin{enumerate}
        \item Also (small) limits and equalizers of \patchable\ quasi-measurable spaces (taken in $\QMS$) are \patchable.
        \item Note that colimits like  coproducts of \patchable\ quasi-measurable spaces might not be \patchable.
    \end{enumerate}
\end{Rem}

\begin{Lem}[Exponential ideal property]
    \label{lem:exponential-ideal}
Let $(\Xcal,\Xcal^\Omega)$ and $(\Ycal,\Ycal^\Omega)$ be two quasi-measurable spaces, where $(\Ycal,\Ycal^\Omega)$ is \patchable. Then also $\lp\Ycal^\Xcal,\lp \Ycal^\Xcal \rp^\Omega\rp$ is \patchable.
\begin{proof} Consider a measurable disjoint union decomposition $\Omega = \bigdcup_{j \in J} C_j$ with countable
            $J \ins \N$ and $Y_j \in \lp \Ycal^\Xcal \rp^\Omega$ for $j \in J$.
            We then consider:
            \[ Y := \bigdcup_{j \in J} Y_j|_{C_j} \in \Sets(\Omega,\Ycal^\Xcal). \]
            Let $\Phi \in \Omega^\Omega$ and $X \in \Xcal^\Omega$. We then have for $\omega \in \Phi^{-1}(C_j)$:
            \[ Y(\Phi(\omega))(X(\omega)) = Y_j(\Phi(\omega))(X(\omega)). \]
            Note that we get the measurable disjoint union decomposition $\Omega = \bigdcup_{j \in J} \Phi^{-1}(C_j)$.
            Furthermore, we know that $Y_j(\Phi)(X) \in \Ycal^\Omega$ for all $j \in J$ by assumption.
            Because $\Ycal$ is \patchable\ we then get that:
            \[ Y(\Phi)(X)= \bigdcup_{j \in J} Y_j(\Phi)(X)|_{\Phi^{-1}(C_j)} \in \Ycal^\Omega. \]
            Since this holds for all $\Phi$ and $X$ the above shows that:
            \[  \bigdcup_{j \in J} Y_j|_{C_j} = Y \in \lp \Ycal^\Xcal \rp^\Omega.\]
    \end{proof}
\end{Lem}

\begin{Eg} Let $(\Xcal,\Xcal^\Omega)$ be a (not necessarily \patchable) quasi-measurable space.
    \begin{enumerate}
        \item The space of indicator functions $(\twob^\Xcal,(\twob^\Xcal)^\Omega)$ is always a \patchable\ quasi-measurable space.
        \item The corresponding $\sigma$-algebra $(\Bcal_\Xcal,(\Bcal_\Xcal)^\Omega)$ is always a \patchable\ quasi-measurable space.
    \end{enumerate}
    \begin{proof}
        $(\twob,\twob^\Omega) = (\twob, \Fcal(\Bcal_\twob))$ is \patchable\ by Example \ref{eg:patchable}
        and thus $(\twob^\Xcal,(\twob^\Xcal)^\Omega)$ as well
        by Lemma \ref{lem:exponential-ideal}. \\
        Since $(\Bcal_\Xcal,(\Bcal_\Xcal)^\Omega) \cong (\twob^\Xcal,(\twob^\Xcal)^\Omega)$
        by Lemma \ref{qms:lem:s-algebra-indicator} also $(\Bcal_\Xcal,(\Bcal_\Xcal)^\Omega)$ is \patchable.
    \end{proof}
\end{Eg}

\subsection{Reflector to \Patchable\ Quasi-Measurable Spaces}

In this subsection we construct the \emph{reflector} $\Pch$ from the category of quasi-measurable spaces $\QMS$ to the category of \emph{\patchable} quasi-measurable spaces $\PQMS$.

\begin{Lem}[Reflector to \patchable\ quasi-measurable spaces]
    \label{lem:reflector-pqms}
    Let $\Xcal$ be a set and $\Hcal \ins \Sets(\Omega,\Xcal)$ any subset of functions, $\Xcal^\one:=\Sets(\one,\Xcal) \circ \Sets(\Omega,\one)$ the set of constant maps to $\Xcal$.
    Then define:
    \[ \Pch(\Hcal) : = \bigcap_{\substack{\Gcal \ins \Sets(\Omega,\Xcal)\\\Xcal^\one \ins \Gcal\\\Gcal \circ \Omega^\Omega \ins \Gcal\\ \Gcal \text{Def.\ } \ref{def:patchable}  }} \Gcal   \qquad \ins \Sets(\Omega,\Xcal), \]
    the intersection of all set of functions $\Gcal \ins \Sets(\Omega,\Xcal)$ that:
    \begin{enumerate}
        \item contains $\Hcal$,
        \item contains the constant maps $\Xcal^\one$,
        \item is closed under right composition with elements  $\Phi \in \Omega^\Omega$,
        \item contains every $X =\bigdcup_{i \in I} X_i|_{C_i}$ if every $X_i \in \Gcal$ and
            $ \Omega = \bigdcup_{i \in I} C_i$ is a countable measurable disjoint union decomposition,
            $i \in I \ins \N$.
    \end{enumerate}
    Then $\Pch(\Hcal)$ satisfies all these points as well and $(\Xcal,\Pch(\Hcal))$ is a well-defined \emph{\patchable} quasi-measurable space.
    \begin{proof}
        First, the intersection is non-empty as $\Gcal=\Sets(\Omega,\Xcal)$ satisfies all points.\\
        We clearly have $\Xcal^\one \cup \Hcal \ins \Pch(\Hcal)$ if $\Xcal^\one \cup \Hcal \ins \Gcal$
        for all $\Gcal$ of consideration. \\
        If $X \in \Pch(\Hcal)$ and $\Phi \in \Omega^\Omega$ then $X \circ \Phi \in \Gcal$ for all those $\Gcal$ and thus $X \circ \Phi \in \Pch(\Hcal)$.\\
        If $X_i \in \Pch(\Hcal)$ for $i \in I \ins \N$ then $X_i \in \Gcal$ for $i \in I$.
        So we get $X :=\bigdcup_{i \in I} X_i|_{C_i} \in \Gcal$ for all $\Gcal$ of consideration, which implies
        $X=\bigdcup_{i \in I} X_i|_{C_i} \in \Pch(\Hcal)$.
    \end{proof}
\end{Lem}

\begin{Thm}
    \label{thm:patchable-reflexive}
    The full subcategory $\PQMS$ of all \patchable\ quasi-measurable spaces inside the quasitopos $\QMS$ of all quasi-measurable spaces is a reflective subcategory, an exponential ideal and cartesian closed in itself.
    In particular, the reflector $\Pch:\; \QMS \to \PQMS$ preserves finite products and arbitrary (small) colimits and coproducts.\\
    Furthermore, $\PQMS$ has all (small) limits and products: the $\QMS$-products/limits of diagrams in $\PQMS$ are already in $\PQMS$. $\PQMS$ has all (small) colimits and coproducts, which are given by applying the reflector $\Pch$ to the corresponding
    $\QMS$-colimits/coproducts:
    \[\coprod_{i \in I}^\PQMS(\Xcal_i,\Xcal_i^\Omega) = \Pch \lp \coprod_{i \in I}^\QMS(\Xcal_i,\Xcal_i^\Omega)\rp
    = \lp \coprod_{i \in I} \Xcal_i, \Pch\lp \coprod_{i \in I} \Xcal_i^\Omega \rp \rp \]
    \begin{proof}
        The reflector is given via Lemma \ref{lem:reflector-pqms}:
        \[\Pch:\; \QMS \to \PQMS,\qquad (\Xcal,\Xcal^\Omega) \mapsto (\Xcal,\Pch(\Xcal^\Omega)), \qquad \Pch(g)=g.\]
        We have to show that $\Pch$ is left-adjoint to the forgetful functor $\PQMS \inj \QMS$, i.e.\ that
        for every quasi-measurable $(\Xcal,\Xcal^\Omega)$ and every \patchable\ quasi-measurable space $(\Ycal,\Ycal^\Omega)$ we have the natural identification:
        \[ \QMS\lp(\Xcal,\Pch(\Xcal^\Omega)),(\Ycal,\Ycal^\Omega)\rp = \QMS\lp(\Xcal,\Xcal^\Omega),(\Ycal,\Ycal^\Omega)\rp.  \]
    Indeed, if $f \in \QMS\lp(\Xcal,\Pch(\Xcal^\Omega)),(\Ycal,\Ycal^\Omega)\rp$ then:
    \[ f \circ \Xcal^\Omega \ins f \circ \Pch(\Xcal^\Omega) \ins \Ycal^\Omega,  \]
    which implies: $f \in \QMS\lp(\Xcal,\Xcal^\Omega),(\Ycal,\Ycal^\Omega)\rp$.\\
    For the reverse inclusion let $g \in \QMS\lp(\Xcal,\Xcal^\Omega),(\Ycal,\Ycal^\Omega)\rp$.
    Then we get:
    \[ g \circ \Xcal^\Omega \ins \Ycal^\Omega,\]
    which is equivalent to:
    \[  \Xcal^\Omega \ins g^*\Ycal^\Omega := \lC X:\,\Omega \to \Xcal\st g \circ X \in \Ycal^\Omega \rC.\]
    Since $g^*\Ycal^\Omega$ satisfies all 4 points from Lemma \ref{lem:reflector-pqms} by Lemma \ref{lem:pull-back-pqms}
    we get:
    \[ \Pch(\Xcal^\Omega) \ins g^*\Ycal^\Omega,\]
    which is equivalent to:
    \[ g \circ \Pch(\Xcal^\Omega)  \ins \Ycal^\Omega. \]
    This implies: $g \in \QMS\lp(\Xcal,\Pch(\Xcal^\Omega)),(\Ycal,\Ycal^\Omega)\rp$.\\
    To show that $\PQMS$ is an exponential ideal of $\QMS$ we need to show that for every $\Xcal \in \QMS$
    and $\Ycal \in \PQMS$ also $\Ycal^\Xcal \in \PQMS$. This was already shown in Lemma \ref{lem:exponential-ideal}.\\
    That $\Pch$ preserves finite products follows from the last statement and \cite{Joh02} Prp.\ A4.3.1.
    \[ \Pch\lp \prod_{n=1}^N (\Xcal_n,\Xcal_n^\Omega) \rp = \prod_{n=1}^N  (\Xcal_n,\Pch(\Xcal_n^\Omega)).   \]
    Since $\Pch$ is a left-adjoint it automatically preserves all colimits like coproducts.\\
    Since $\QMS$ has all small colimits by Theorem \ref{qms:thm:co-complete} also $\PQMS$ has all small colimits
    by applying $\Pch$ afterwards.\\
    Since $\PQMS$ inherits all products, see Lemma \ref{lem:patchable-product}, and equalizers,
    see Lemma \ref{qms:lem:qms-equalizer-up} and check that equalizers of \patchable\ quasi-measurable spaces
    are \patchable, it then also inherits all limits from $\QMS$.
        \end{proof}
\end{Thm}

\begin{Rem}
    \begin{enumerate}
        \item For a quasi-measurable space $(\Xcal,\Xcal^\Omega)$ we always have:
            \[ \Xcal^\Omega \ins \Pch(\Xcal^\Omega) \ins \Fcal\Bcal(\Xcal^\Omega) \ins \Sets(\Omega,\Xcal).\]
        \item $\SQMS \ins \PQMS \ins \QMS$ are reflective full subcategories.
        \item For quasi-universal spaces we even have the full subcategories:
            \[ \StdQUS \ins \SQUS \ins \PQUS \ins \QUS,  \]
            where $\StdQUS$ is the full subcategory of all standard quasi-universal spaces, see Definition \ref{def:univ-qus} and Corollary \ref{cor:umeas-uqus-eq} later on, which might not be reflective inside $\QUS$.
    \end{enumerate}
\end{Rem}

\begin{Thm}[$\PQMS$ is a quasitopos]
    \label{thm:pqms-quasitopos}
    The category $\PQMS$ of \patchable\ quasi-measurable spaces and quasi-measurable maps is a
    \emph{quasitopos}, and in particular locally cartesian closed.
    \begin{proof}
        We check the four ingredients, each of which is already available.

        \emph{Finite limits.} By \Cref{lem:patchable-product} products and equalizers of
        \patchable\ quasi-measurable spaces are again \patchable\ and are the ones of $\QMS$, so
        $\PQMS$ has all (finite) limits and they are computed in $\QMS$.

        \emph{Finite colimits.} By \Cref{thm:patchable-reflexive} colimits exist in $\PQMS$ and are
        obtained by applying the reflector to the corresponding $\QMS$-colimits.

        \emph{Local cartesian closedness.} Let $(\Tcal,\Tcal^\Omega)$ be \patchable\ and let
        $\Ycal,\Zcal$ be \patchable\ quasi-measurable spaces over $\Tcal$. The fibre product
        $\Ycal \ftimes_\Tcal \Zcal$ is a limit, hence \patchable\ by
        \Cref{lem:patchable-product}. The internal hom $\ihom_\Tcal(\Ycal,\Zcal)$ of
        \Cref{qms:sec:int-homs} is built from function spaces $\Zcal^\Ycal$ and fibre products; the
        former are \patchable\ because $\PQMS$ is an exponential ideal,
        \Cref{lem:exponential-ideal}, and the latter by the preceding sentence. So
        $\ihom_\Tcal(\Ycal,\Zcal)$ lies in $\PQMS_\Tcal$, and since $\PQMS_\Tcal$ is a full
        subcategory of $\QMS_\Tcal$ the adjunction of \Cref{qms:thm:lcc} restricts:
        \begin{align}
            \PQMS_\Tcal\lp \Xcal \ftimes_\Tcal \Ycal, \Zcal \rp
            &\cong \PQMS_\Tcal\lp \Xcal, \ihom_\Tcal\lp\Ycal,\Zcal\rp\rp.
        \end{align}
        Every slice of $\PQMS$ is therefore cartesian closed.

        \emph{Classifier for strong monomorphisms.} The indiscrete two-point space $\twop$ is
        \patchable, for the trivial reason that $\twop^\Omega=\Sets(\Omega,\twop)$ consists of
        \emph{all} maps $\Omega \to \twop$ and is thus closed under the gluings of
        \Cref{def:patchable}. Embeddings of \patchable\ spaces are again \patchable, being
        equalizers, and the classifying square $\Wcal \cong \Xcal \ftimes_\twop \one$ of
        \Cref{qms:thm:strong-mono-classifier} is a fibre product, hence computed in $\QMS$ and
        again \patchable. So $\twop$ classifies the strong monomorphisms of $\PQMS$ as well.
    \end{proof}
\end{Thm}

\begin{Rem}[Two comments on the preceding proof]
    \label{rem:pqms-quasitopos}
    \leavevmode
    \begin{enumerate}
        \item The proof does \emph{not} go through the general principle that a reflective
            subcategory of a quasitopos is again one; that principle requires the reflector to
            preserve finite limits, and \Cref{thm:patchable-reflexive} gives only that $\Pch$
            preserves finite \emph{products}. Indeed, for a subspace $\Wcal \ins \Xcal$ the space
            $\Pch(\Wcal)$ carries the closure of $\Wcal^\Omega$ under the gluings of
            \Cref{def:patchable}, whereas the equalizer computed in $\PQMS$ carries
            $\Pch(\Xcal^\Omega) \cap \Sets(\Omega,\Wcal)$, which may be strictly larger: a glued
            $X=\bigdcup_n X_n|_{C_n}$ may take values in $\Wcal$ without any $X_n$ lying in
            $\Wcal^\Omega$, only its restriction to $C_n$ doing so. The argument above avoids this
            by never asking $\Pch$ to preserve limits; it only uses that $\PQMS$ is \emph{closed
            under} them, which is \Cref{lem:patchable-product}.
        \item What does change on passing from $\QMS$ to $\PQMS$ is the colimits, not the limits.
            \Cref{def:patchable} is a gluing condition, and gluing conditions are preserved by
            limits but not by colimits; the countable coproduct has to be reflected, see
            \Cref{lem:patch-coprod}, and this is exactly what is behind the coproducts of
            \Cref{thm:Markov-category} being those of $\PQUS$ rather than those of $\QUS$.
        \item There is a second route to \Cref{thm:pqms-quasitopos}, for readers who want it, and
            it takes three lines. Enlarge the site $\Ccal_\Omega$ of
            \Cref{qms:rem:concrete-presheaf} to the category $\Ccal_\Omega^\flat$ whose objects are
            $\one$ and the events $D \in \Bcal_\Omega$, with $\Ccal_\Omega^\flat(\one,D)=D$ and
            $\Ccal_\Omega^\flat(D,E) = \lC \Phi|_D \st \Phi \in \Omega^\Omega,\, \Phi(D) \ins E \rC$,
            and equip it with the coverage whose covering families of $D$ are the countable
            measurable partitions $D = \bigdcup_{n \in \N} D_n$. For the archetypical sample
            space $\Omega^\Omega=\Fcal(\Bcal_\Omega)$ of \Cref{qms:eg:meas-sample-space} this is
            again a concrete site in the sense of \cite{BH11} Def.\ 4.9: the covering families
            are jointly surjective on points, $\Ccal_\Omega^\flat(\one,-)$ is faithful, and the
            coverage is subcanonical, because gluing measurable maps along a countable
            measurable partition gives a measurable map, \Cref{eg:patchable}.
            A concrete presheaf on $\Ccal_\Omega^\flat$ that is a \emph{sheaf} for this coverage is
            exactly a quasi-measurable space whose admissible random variables are closed under the
            gluings of \Cref{def:patchable}, so
            \begin{align}
                \PQMS &\simeq \mathrm{Conc}\lp\Ccal_\Omega^\flat,\text{countable measurable partitions}\rp,
            \end{align}
            and \cite{BH11} Thm.\ 5.25 gives \Cref{thm:pqms-quasitopos} again. For
            $\Omega=(\R,\Fcal(\Bcal_\R))$ this is the presentation of $\QBS$ as the concrete
            sheaves on standard Borel spaces with countable measurable covers, see \cite{VO18}
            Section 11 and \Cref{prp:pqms-is-qbs}.
    \end{enumerate}
\end{Rem}

\begin{Lem}[Countable coproducts]
    \label{lem:patch-coprod}
    Let $I$ be a countable set and $(\Xcal_i,\Bcal_i)$ measurable spaces for $i \in I$.
    Then we have:
    \[ \Pch \lp\coprod_{i \in I} (\Xcal_i,\Fcal(\Bcal_i)) \rp = \Fcal\lp \coprod_{i \in I}(\Xcal_i,\Bcal_i) \rp,  \]
    where on the lhs the reflector $\Pch$ is applied to the $\QMS$-coproduct, which together is the $\PQMS$-coproduct, of the individual spaces after applying $\Fcal$, and on the rhs we have $\Fcal$ applied to the coproduct in $\Meas$. \\
    In other words, $\Fcal$, when seen as a functor from $\Meas$ to $\PQMS$, preserves countable coproducts.
    \begin{proof}
    Remember that the coproduct $\sigma$-algebra is:
    \[\Bcal_{\coprod_{i \in I}(\Xcal_i,\Bcal_i)}:= \lC D \ins \coprod_{i \in I} \Xcal_i \st \forall i \in I.\,
    \incl_i^{-1}(D) \in \Bcal_i\rC. \]
    Since the inclusion are measurable:
    \[ \incl_k:\; (\Xcal_k,\Bcal_k) \to  \coprod_{i \in I}(\Xcal_i,\Bcal_i), \]
    we get by applying $\Fcal$ the quasi-measurable maps:
    \[ \incl_k:\; (\Xcal_k,\Fcal(\Bcal_k)) \to \Fcal\lp \coprod_{i \in I}(\Xcal_i,\Bcal_i)\rp.\]
    By the universal property of the coproduct we then have a quasi-measurable (identity) map:
    \[ \id:\; \coprod_{i\in I} (\Xcal_i,\Fcal(\Bcal_i)) \to \Fcal\lp \coprod_{i \in I}(\Xcal_i,\Bcal_i)\rp.\]
    Since the rhs is a \patchable\ quasi-measurable space we then also get the well-defined quasi-measurable map:
    \[ \id:\; \Pch\lp\coprod_{i\in I} (\Xcal_i,\Fcal(\Bcal_i))\rp \to \Fcal\lp \coprod_{i \in I}(\Xcal_i,\Bcal_i)\rp.\]
    This shows the inclusion:
    \[ \Lcal\lp\coprod_{i \in I} \Fcal(\Bcal_i)\rp \ins \Fcal\lp \Bcal_{\coprod_{i \in I}(\Xcal_i,\Bcal_i)}\rp. \]
To show the reverse inclusion consider:
\[ X \in \Fcal\lp \Bcal_{\coprod_{i \in I}(\Xcal_i,\Bcal_i)}\rp= \Meas\lp(\Omega,\Bcal_\Omega),\lp
\coprod_{i \in I} \Xcal_i,\Bcal_{\coprod_{i \in I}\Xcal_i} \rp\rp.  \]
Since $\Xcal_k \in \Bcal_{\coprod_{i \in I}(\Xcal_i,\Bcal_i)}$ we get that:
\[ C_k := X^{-1}\lp \Xcal_k\rp \in \Bcal_\Omega.\]
We can then define:
\[ X_k:\; \Omega \to \Xcal_k,\qquad X_k|_{C_k}:=X|_{C_k},\qquad X_k|_{\Omega\sm C_k}:= x_k \in \Xcal_k.\]
Since $X$ is measurable and $C_k \in \Bcal_\Omega$ then also $X_k$ is measurable.
So $X_k \in \Fcal(\Bcal_k)$ for $k \in I$.
By definition of the coproduct in $\QMS$ we get for every $k \in I$ that:
\[   \incl_k \circ X_k \in \coprod_{i \in I} \Fcal(\Bcal_i).\]
Furthermore, we get by definition of $\Pch$ that:
\[ X = \bigdcup_{i \in I} (\incl_i \circ X_i|_{C_i}) \in \Lcal\lp\coprod_{i \in I} \Fcal(\Bcal_i)\rp.\]
This shows that:
\[ \Fcal\lp \Bcal_{\coprod_{i \in I}(\Xcal_i,\Bcal_i)}\rp \ins \Lcal\lp\coprod_{i \in I} \Fcal(\Bcal_i)\rp, \]
and thus the claim.
    \end{proof}
\end{Lem}

\begin{Rem}
    The strongest argument to work inside the category of \patchable\ quasi-measurable spaces $\PQMS$
    rather than inside the category of all quasi-measurable spaces $\QMS$
    is in fact Lemma \ref{lem:patch-coprod}, which asserts that when going from $\Meas$ to $\PQMS$ via the functor $\Fcal$ one preserves countable coproducts (as well as all products):
    \[ \coprod^\PQMS_{n \in \N} (\Xcal_n,\Fcal(\Bcal_n)) = \Pch \lp\coprod^\QMS_{n \in \N} (\Xcal_n,\Fcal(\Bcal_n)) \rp = \Fcal\lp \coprod^\Meas_{n \in \N}(\Xcal_n,\Bcal_n) \rp.  \]
    Furthermore, one does not lose anything over $\QMS$ as one also has all (small) limits and colimits in $\PQMS$, which is also cartesian closed.
    On the downside one needs to require and check the condition from Definition \ref{def:patchable}
    through all constructions and examples,
    although that condition does not seem to play any role besides Lemma \ref{lem:patch-coprod}, which
    also states how to recover from the disalignment by applying the reflector $\Pch$.
\end{Rem}

\subsection{\Patchable\ Probability Spaces}

In this subsection we study under which conditions our probability monads become \patchable.

\begin{Lem}[Spaces of probability measures $\Giry$, $\PrQm$, $\PrPF$]
    \label{lem:patchable-G-Q-K}
    Let $(\Xcal,\Xcal^\Omega)$ be a (not necessarily \patchable) quasi-measurable space.
    Then the following points hold:
    \begin{enumerate}
        \item $(\Giry(\Xcal),\Giry(\Xcal)^\Omega)$ is \patchable.
        \item $(\PrQm(\Xcal),\PrQm(\Xcal)^\Omega)$ is \patchable.
    \end{enumerate}
    If $(\Xcal,\Xcal^\Omega)$ is \patchable\ then also the following holds:
    \begin{enumerate}[resume]
        \item $(\PrPF(\Xcal),\PrPF(\Xcal)^\Omega)$ is \patchable.
    \end{enumerate}
    \begin{proof}
        Since $\Giry(\Xcal)^\Omega =\Fcal(\Bcal_{\Giry(\Xcal)})$, where:
        \[ \Bcal_{\Giry(\Xcal)} := \sigma\lp\lC \ev_A^{-1}((t,1])\st A \in \Bcal_\Xcal, t \in \R  \rC\rp, \]
        it is clear by Example \ref{eg:patchable} that $(\Giry(\Xcal),\Giry(\Xcal)^\Omega)$ is \patchable.\\
        For $\PrQm(\Xcal)$ first note that $([0,1],[0,1]^\Omega)=([0,1],\Fcal(\Bcal_{[0,1]}))$ is \patchable\ by
        Example \ref{eg:patchable}, and thus by Lemma \ref{lem:exponential-ideal} also $[0,1]^{\Bcal_\Xcal}$.
        Since $\PrQm(\Xcal)$ carries the subspace quasi-measurable space structure of $[0,1]^{\Bcal_\Xcal}$ by
        Lemma \ref{lem:pull-back-pqms} also $(\PrQm(\Xcal),\PrQm(\Xcal)^\Omega)$ is \patchable.\\
        Now assume that $(\Xcal,\Xcal^\Omega)$ is \patchable.\\
        $\PrPF(\Xcal)$ carries the quotient quasi-measurable space structure of $\Xcal^\Omega \times \PrQm(\Omega)$ along
        the push-forward map $\pf$. Since $\Xcal$ is \patchable\ so is $\Xcal^\Omega$ by Lemma \ref{lem:exponential-ideal}.
        By the point above $\PrQm(\Omega)$ is \patchable\ as well. By Lemma \ref{lem:patchable-product} then also their product
        $\Xcal^\Omega \times \PrQm(\Omega)$ is \patchable. As a quotient of the latter by Lemma \ref{lem:pull-back-pqms} also the
        space $(\PrPF(\Xcal),\PrPF(\Xcal)^\Omega)$ becomes \patchable.\\
    \end{proof}
\end{Lem}

\begin{Lem}[The space of probability measures $\PrpF$]
    \label{lem:patchable-P}
    Let $(\Xcal,\Xcal^\Omega)$ be a \patchable\ quasi-measurable space and assume that there exists an isomorphism
    of measurable spaces:
    \[ \coprod_{n \in \N} (\Omega,\Bcal_\Omega) \cong (\Omega,\Bcal_\Omega).\]
    Further assume that $\Omega^\Omega = \Fcal(\Bcal_\Omega)$. Then the space of push-forward probability measures $(\PrpF(\Xcal),\PrpF(\Xcal)^\Omega)$
    is a \patchable\ quasi-measurable space.
    \begin{proof}
        Consider a measurable disjoint union decomposition $\Omega = \bigdcup_{j \in J} C_j$ with countable $J \ins \N$ and
                  $\nu_j=X_{j,*}\kappa_j \in \PrpF(\Xcal)^\Omega$ with $X_j \in \Xcal^\Omega$ and $\kappa_j \in \PrQm(\Omega)^\Omega$.
                  By the made assumption we also have a countable measurable disjoint union decomposition:
                  \[\Omega = \bigdcup_{n \in \N} D_n,\]
                  with measurable isomorphisms: $\Omega \cong D_n$ for $n \in \N$.
                  Consider those measurable isomorphisms composed with the measurable inclusion map:
                  \[ d_n :\; \Omega \cong D_n \inj \Omega.\]
                  Then we get that $d_n \in \Fcal(\Bcal_\Omega)=\Omega^\Omega$.\\
                  Define the map:
                  \[ \kappa := \bigdcup_{j \in J} d_{j,*}\kappa_j|_{C_j}:\; \Omega =
                  \bigdcup_{j \in J} C_j \to \PrQm\lp\bigdcup_{n \in \N} D_n \rp = \PrQm(\Omega).\]
                  Note that each $d_{j,*}\kappa_j$ is supported on $D_j$. \\
                  To check that $\kappa \in \PrQm(\Omega)^\Omega$ let $\Phi \in \Omega^\Omega$ and $E \in \Bcal_{\Omega \times \Omega}$.
                  Then the map:
                  \[ \Phi^{-1}(C_i) \ni \omega \mapsto \kappa(\Phi(\omega))(E_\omega) = d_{i,*}\kappa_i(\Phi(\omega))(E_\omega)
                    = \kappa_i(\Phi(\omega))((\id_\Omega,d_i)^{-1}(E)_\omega), \]
                    is seen to be measurable by the assumption that $\kappa_i \in \PrQm(\Omega)^\Omega$, thus $\kappa \in \PrQm(\Omega)^\Omega$.\\
                  We now pick a fixed point $\omega_0 \in \Omega$ and $x_0 \in \Xcal$ and define for $i\in J$:
                  \[ d_i^-: \Omega =\bigdcup_{n \in \N} D_n \to \Omega,\qquad d_i^-|_{D_i}:=d_i^{-1},\qquad d_i^-|_{D_j}:=\omega_0, \; j\neq i.   \]
                  For $n \in \N \sm J$ we define the constant maps: $d_n^-:=\omega_0$ and $X_n:=x_0$.
                  Then all $d_n^-$ are measurable and thus $d_n^- \in \Fcal(\Bcal_\Omega)=\Omega^\Omega$. Since all $X_n \in \Xcal^\Omega$ we get that
                  $X_n \circ d_n^- \in \Xcal^\Omega$ for all $n \in \N$.
                  With these settings we can then put:
                  \[ X :=  \bigdcup_{n \in \N} (X_n \circ d_n^-)|_{D_n} :\;\Omega = \bigdcup_{n \in \N} D_n \to \Xcal.\]
                  It is then $X \in \Xcal^\Omega$ since $\Xcal^\Omega$ is \patchable.
                  Furthermore, for $\omega \in C_i$, $i \in J$, we get:
                  \[ X_*\kappa(\omega) = X_{i,*} \, d_{i,*}^- \, d_{i,*}\,\kappa_i(\omega) = X_{i,*}\kappa_i(\omega) = \nu_i(\omega).\]
                  So we get that:
                  \[ \bigdcup_{j \in J} \nu_j|_{C_j} = X_*\kappa \in \PrpF(\Xcal)^\Omega, \]
                  with $X \in \Xcal^\Omega$ and $\kappa \in \PrQm(\Omega)^\Omega$,
                  which shows the claim.
    \end{proof}
\end{Lem}

\begin{Rem}
    It has not yet been investigated under which conditions $(\PrPf(\Xcal),\PrPf(\Xcal)^\Omega)$
    or $(\Prpf(\Xcal),\Prpf(\Xcal)^\Omega)$ become \patchable.
    We will see in Theorem \ref{thm:qus-S=P-in-R} however that for the category of quasi-universal spaces $\QUS$ we have:
    \[\Prpf = \PrpF = \PrPf=\PrPF,\]
    which then will be covered by the results above.
\end{Rem}

\section{Appendix - Quasi-Universal Spaces}
\subsection{Disintegration of Markov Kernels}
\label{qms:sec:app-disint-meas}

The disintegration results for quasi-universal spaces reduce to the following two
results about measurable spaces, which we quote. Together they are the only
external input to \Cref{thm:disint-X}. Each lists three alternative hypotheses; none of them
refers to the Markov kernel that is to be disintegrated, so each makes the triple of spaces a
\emph{disintegration triple} in the following sense.

\begin{Def}[Disintegration triple]
    \label{def:disintegration-triple}
    A triple $(\Xcal,\Bcal_\Xcal)$, $(\Ycal,\Bcal_\Ycal)$, $(\Zcal,\Bcal_\Zcal)$ of measurable
    spaces is called a \emph{disintegration triple} if \emph{every} Markov kernel
    \[ K(X,Y|Z):\; (\Zcal,\Bcal_\Zcal) \to \Giry(\Xcal\times\Ycal, \Bcal_\Xcal\otimes\Bcal_\Ycal) \]
    admits a \emph{conditional Markov kernel}
    \[ K(X|Y,Z):\; (\Ycal\times\Zcal, \Bcal_\Ycal\otimes\Bcal_\Zcal)\to \Giry(\Xcal, \Bcal_\Xcal), \]
    i.e.\ one with
    \[ K(X,Y|Z) = K(X|Y,Z) \otimes K(Y|Z), \]
    where $K(Y|Z)$ denotes the marginal Markov kernel.
    It is called a \emph{universal disintegration triple} if the same holds with
    $\Bcal_\Ycal \otimes \Bcal_\Zcal$ replaced by its universal completion
    $(\Bcal_\Ycal \otimes \Bcal_\Zcal)_\Giry$ as the domain of $K(X|Y,Z)$.
    Note that neither notion refers to a particular Markov kernel; both are properties of the
    triple of spaces alone.
\end{Def}

\begin{Thm}[See \cite{For21} Appendix C]
    \label{thm:markov-kernel-disint-meas-II}
    Let $(\Xcal,\Bcal_\Xcal)$, $(\Ycal,\Bcal_\Ycal)$ and $(\Zcal,\Bcal_\Zcal)$ be measurable spaces
    satisfying at least one of the following:
    \begin{enumerate}
        \item $(\Xcal,\Bcal_\Xcal)$ is standard Borel, $(\Ycal,\Bcal_\Ycal)$ is countably generated
            and $(\Zcal,\Bcal_\Zcal)$ is arbitrary; or:
        \item $(\Xcal,\Bcal_\Xcal)$ is standard Borel, $(\Ycal,\Bcal_\Ycal)$ is arbitrary and
            $(\Zcal,\Bcal_\Zcal)$ is discrete (i.e.\ $\Zcal$ countable, $\Bcal_\Zcal = \two^\Zcal$); or:
        \item $(\Xcal,\Bcal_\Xcal)$ is arbitrary, $(\Ycal,\Bcal_\Ycal)$ is discrete (i.e.\ $\Ycal$
            countable, $\Bcal_\Ycal = \two^\Ycal$) and $(\Zcal,\Bcal_\Zcal)$ is arbitrary.
    \end{enumerate}
    Then $(\Xcal,\Ycal,\Zcal)$ is a \emph{disintegration triple} in the sense of
    \Cref{def:disintegration-triple}.
\end{Thm}

\begin{Rem}[The three cases, and what else is available]
    \label{rem:disint-meas-cases}
    None of the three conditions refers to the kernel $K(X,Y|Z)$; each is a condition on the triple
    of spaces alone, so each of them makes $(\Xcal,\Ycal,\Zcal)$ a \emph{disintegration triple}.
    They are genuinely different. Cases 1 and 2 are the two familiar ones and both require
    $(\Xcal,\Bcal_\Xcal)$ standard Borel. Case 3 requires nothing of $\Xcal$ at all: when $\Ycal$ is
    discrete the conditional kernel is the elementary quotient
    \begin{align}
        K(X \in A|Y=y,Z=z) &= \frac{K(X \in A,Y=y|Z=z)}{K(Y=y|Z=z)}
        \label{eq:disint-discrete-quotient}
    \end{align}
    wherever $K(Y=y|Z=z)>0$, and an arbitrary fixed probability measure on $(\Xcal,\Bcal_\Xcal)$
    elsewhere; measurability in $z$ is inherited from $K(X,Y|Z)$ and measurability in $y$ is
    automatic, every subset of $\Ycal$ being measurable.
    Beyond these, \cite{For21} also gives criteria that are conditions on the given kernel rather
    than on the spaces --- a joint density w.r.t.\ a product of $\sigma$-finite measures, and two
    determinism conditions --- which we do not need here.
\end{Rem}

The version we actually apply extends this from standard Borel spaces to their universal
completions, and correspondingly replaces countable generation by universal countable generation.

\begin{Cor}[See \cite{For21} Cor.\ C.8]
    \label{cor:markov-kernel-disint-meas-II}
    Let $(\Xcal,\Bcal_\Xcal)$, $(\Ycal,\Bcal_\Ycal)$ and $(\Zcal,\Bcal_\Zcal)$ be measurable spaces
    satisfying at least one of the following:
    \begin{enumerate}
        \item $(\Xcal,\Bcal_\Xcal)$ is perfect (i.e.\ every probability measure on
            $(\Xcal,\Bcal_\Xcal)$ is perfect) and universally countably generated (i.e.\ its
            universal completion is the universal completion of a countably generated measurable
            space), e.g.\ a universal measurable space; $(\Ycal,\Bcal_\Ycal)$ is universally
            countably generated; and $(\Zcal,\Bcal_\Zcal)$ is arbitrary; or:
        \item $(\Xcal,\Bcal_\Xcal)$ is perfect and universally countably generated,
            $(\Ycal,\Bcal_\Ycal)$ is arbitrary and $(\Zcal,\Bcal_\Zcal)$ is discrete; or:
        \item $(\Xcal,\Bcal_\Xcal)$ is arbitrary, $(\Ycal,\Bcal_\Ycal)$ is discrete and
            $(\Zcal,\Bcal_\Zcal)$ is arbitrary.
    \end{enumerate}
    Then $(\Xcal,\Ycal,\Zcal)$ is a \emph{universal disintegration triple} in the sense of
    \Cref{def:disintegration-triple}. In case 3 it is even a disintegration triple, no completion
    being needed.
\end{Cor}

\begin{Rem}[Why case 3 needs no completion]
    \label{rem:disint-cor-case3}
    Cases 1 and 2 are the universal-completion versions of cases 1 and 2 of
    \Cref{thm:markov-kernel-disint-meas-II}: one passes to a standard Borel model on a set of full
    measure, which is what perfectness and universal countable generation are for, and the
    conditional kernel is then only measurable for the universal completion
    $(\Bcal_\Ycal \otimes \Bcal_\Zcal)_\Giry$.
    Case 3 is stated verbatim as in \Cref{thm:markov-kernel-disint-meas-II}, because there is
    nothing to complete: $\Bcal_\Ycal = \two^\Ycal$ is its own universal completion, the
    construction \Cref{eq:disint-discrete-quotient} is explicit, and it produces a kernel that is
    measurable for $\Bcal_\Ycal \otimes \Bcal_\Zcal$ on the nose. In particular case 3 needs
    neither perfectness nor universal countable generation of $(\Xcal,\Bcal_\Xcal)$.
\end{Rem}

\label{qms:sec:app-qus}
\begin{Prf}{qms:lem:map-completetions}
First remember that for every $\mu \in \Giry(\Xcal,\Bcal_\Xcal)$ we have that $f$ is
        $(\Bcal_\Xcal)_\mu$-$(\Bcal_\Ycal)_{f_*\mu}$-measurable.
        Then for $\mu \in \Pcal$ we have by assumption that $f_*\mu \in \Qcal$. So for
        $B \in (\Bcal_\Ycal)_\Qcal \ins (\Bcal_\Ycal)_{f_*\mu}$ we get that
        $f^{-1}(B) \in (\Bcal_\Xcal)_\mu$ for all $\mu \in \Pcal$, thus
        $f^{-1}(B) \in (\Bcal_\Xcal)_\Pcal$.
\end{Prf}

\begin{Prf}{qms:thm:markov-kernel-completions}
        We use Lemma \ref{qms:lem:map-completetions} and Proposition \ref{prp:prob-prob-meas} with the choices:
        \[ \Bcal_\Pcal := \Bcal_{\Giry(\Xcal,\Bcal_\Xcal)|\Pcal}, \qquad \Mcal := \Mbb_\iota^{-1}(\Pcal) \ins \Giry(\Pcal,\Bcal_\Pcal).\]
            Since $\kappa(\Zcal) \ins \Pcal$ the map $\kappa$ factorizes into the following measurable maps:
            \[ \kappa:\; (\Zcal,\Bcal_\Zcal) \xrightarrow{\kappa} (\Pcal,\Bcal_\Pcal) \xhookrightarrow{\iota}
    \lp \Giry(\Xcal,\Bcal_\Xcal), \Bcal_{\Giry(\Xcal,\Bcal_\Xcal)}\rp.\]
This induces the following measurable factorization of $\kappa_\circ$:
\[ \kappa_\circ:\; \Giry(\Zcal,\Bcal_\Zcal) \xrightarrow{\kappa_*} \Giry(\Pcal,\Bcal_\Pcal) \xrightarrow{\iota_*}
\Giry\lp \Giry(\Xcal,\Bcal_\Xcal), \Bcal_{\Giry(\Xcal,\Bcal_\Xcal)}\rp \xrightarrow{\Mbb} \Giry(\Xcal,\Bcal_\Xcal),\]
showing that: $\Qcal:= (\kappa_\circ)^{-1}(\Pcal) =(\kappa_*)^{-1}(\Mcal)$.\\
Then we get a commutative diagram of measurable maps:
    \[\xymatrix{
            (\Zcal,(\Bcal_\Zcal)_\Qcal) \ar@{^(->}_-{\id_\Zcal}[d] \ar^-{\kappa}[rr] && (\Pcal,(\Bcal_\Pcal)_\Mcal) \ar@{^(->}^-{\id_\Pcal}[d]
            \ar@{^(->}^-{\incl}[rr]&& \Giry(\Xcal,(\Bcal_\Xcal)_\Pcal) \ar^-{\res}[d]\\
            (\Zcal,\Bcal_\Zcal) \ar_-{\kappa}[rr] &&(\Pcal,\Bcal_\Pcal) \ar@{^(->}_-{\iota}[rr]&& \Giry(\Xcal,\Bcal_\Xcal).
    }\]
    Indeed, the upper left map is well-defined and measurable by Lemma \ref{qms:lem:map-completetions} using $\kappa_*(\Qcal) \ins \Mcal$ and
    and the upper right map by Proposition \ref{prp:prob-prob-meas} using $\Mcal \ins \Mbb_\iota^{-1}(\Pcal)$.
    Since the composition of the upper maps is our induced Markov kernel, the claim is shown.

\end{Prf}

\subsection{The Sigma-Algebras of Quasi-Universal Spaces}
\begin{Lem}[The $\sigma$-algebras of quasi-universal spaces]
    \label{lem:qus-universally-complete}
    Let $(\Xcal,\Xcal^\Omega)$ be a quasi-universal space.
    Then $\Bcal(\Xcal^\Omega)$ is universally complete and complete w.r.t.\ the set of push-forward probability measures:
    \[\PrPf(\Xcal,\Xcal^\Omega):= \lC  X_*\mu :\; \Bcal(\Xcal^\Omega) \to [0,1]\st  X \in \Xcal^\Omega, \mu \in
    \Giry(\Omega,\Bcal_\Omega)\footnotemark\rC,
\]
    i.e.\ we have:\footnotetext{\label{fn:G-Q-1}We will see in Lemma \ref{lem:univ-qus-G-Q} that $\Giry(\Omega,\Bcal_\Omega) =\PrQm(\Omega,\Omega^\Omega)$ as sets and quasi-universal spaces and the ambiguity in the notation for
$\PrPf(\Xcal,\Xcal^\Omega)$ will disappear.}
    \[\Bcal(\Xcal^\Omega) =  \Bcal(\Xcal^\Omega)_\Giry = \Bcal(\Xcal^\Omega)_{\PrPf(\Xcal,\Xcal^\Omega)}
    = \bigcap_{\substack{\mu \in \Giry(\Omega,\Bcal_\Omega)\\ X \in \Xcal^\Omega}} X_*(\Bcal_\Omega)_\mu.\]
    Note that each $ X_*(\Bcal_\Omega)_\mu$ is a complete $\sigma$-algebra on $\Xcal$ w.r.t.\
    the push-forward probability measure $ X_*\mu$.
    \begin{proof}
        Since the set of all probability measures contains all push-forward probability measures:
        \[\Giry(\Xcal,\Bcal(\Xcal^\Omega)) =: \Giry(\Xcal,\Xcal^\Omega)  \sni \PrPf(\Xcal,\Xcal^\Omega),\]
        we have the chain of  inclusions:
        \begin{align*}
            \Bcal(\Xcal^\Omega)_{\PrPf(\Xcal,\Xcal^\Omega)} \sni \Bcal(\Xcal^\Omega)_\Giry
            &\sni     \Bcal(\Xcal^\Omega) \\
          &= \bigcap_{ X \in \Xcal^\Omega}  X_*\Bcal_\Omega \\
          &= \bigcap_{ X \in \Xcal^\Omega}  X_*(\Bcal_\Omega)_\Giry \\
        &= \bigcap_{ X \in \Xcal^\Omega}  X_* \bigcap_{\mu \in \Giry(\Omega,\Bcal_\Omega)}(\Bcal_\Omega)_\mu \\
            &= \bigcap_{ X \in \Xcal^\Omega} \bigcap_{\mu \in \Giry(\Omega,\Bcal_\Omega)} X_*(\Bcal_\Omega)_\mu \\
            &\sni \bigcap_{ X \in \Xcal^\Omega} \bigcap_{\mu \in \Giry(\Omega,\Bcal_\Omega)} (\Bcal(\Xcal^\Omega) )_{ X_*\mu} \\
            &= \bigcap_{\nu \in \PrPf(\Xcal,\Xcal^\Omega)} (\Bcal(\Xcal^\Omega) )_{\nu}\\
            &= \Bcal(\Xcal^\Omega)_{\PrPf(\Xcal,\Xcal^\Omega)}.
        \end{align*}
        This shows equality and thus the claim.
    \end{proof}
\end{Lem}
With these in hand we can prove the disintegration theorem for quasi-universal spaces.

\begin{Prf}{thm:disint-X}
        Let $K(X,Y|Z) \in \PrPf(\Xcal \times \Ycal)^\Zcal$. Then the marginal $K(Y|Z) \in \PrPf(\Ycal)^\Zcal$.\\
        1.) Assume that $(\Xcal,\Xcal^\Omega)$ is a standard quasi-universal space and
        $(\Ycal,\Ycal^\Omega)$ countably separated.
        Then by Theorem \ref{thm:univ-qus-GPQKR} we have the identification of quasi-universal spaces:
        \[\PrPf(\Xcal)=\PrPF(\Xcal) = \PrQm(\Xcal)= \Giry(\Xcal).\]
        In particular, $\PrPf(\Xcal)^\Omega=\Giry(\Xcal)^\Omega$ consist of all measurable maps
        $\Omega \to \Giry(\Xcal,\Bcal_\Xcal)$.\\
        Also by Lemma \ref{lem:qus-universally-complete} and Theorem \ref{thm:countably-separated} there exists a
    countably generated $\sigma$-algebras $\Ecal_\Xcal \ins \Bcal_\Xcal$ and $\Ecal_\Ycal \ins \Bcal_\Ycal$ such that:
    \begin{enumerate}
        \item[a.)] $\Ecal_\Xcal$ separates the points of $\Xcal$ such that
            $(\Ecal_\Xcal)_\Giry = \Bcal_\Xcal$, and:
        \item[b.)] $\Ecal_\Ycal$ separates the points of $\Ycal$ such that
    $(\Ecal_\Ycal)_{\PrPf(\Ycal,\Ycal^\Omega)} = \Bcal_\Ycal$.
    \end{enumerate}
    Then consider the composition with the measurable restriction maps:
    \[K(X,Y|Z):\; \Zcal \to \PrPf(\Xcal \times \Ycal)\to \Giry(\Xcal \times \Ycal, \Bcal_{\Xcal \times \Ycal}) \to
    \Giry(\Xcal \times \Ycal, \Bcal_\Xcal \otimes \Ecal_\Ycal).\]
    By Theorem \ref{cor:markov-kernel-disint-meas-II} and the fact that $(\Xcal,\Bcal_\Xcal)$ is a universal measurable space  there exists a measurable:
    \[ K(X|Y,Z):\; (\Ycal \times \Zcal, (\Ecal_\Ycal \otimes \Bcal_\Zcal)_\Giry) \to \Giry(\Xcal,\Bcal_\Xcal), \]
    such that for all $z \in \Zcal$:
    \[  K(X,Y|Z=z) = K(X|Y,Z=z) \otimes K(Y|Z=z). \]
    as measures on $\Bcal_\Xcal \otimes \Ecal_\Ycal$.
    By Lemma \ref{lem:qus-universally-complete} we get that:
    \[ (\Ecal_\Ycal \otimes \Bcal_\Zcal)_\Giry \ins (\Bcal_{\Ycal \times \Zcal})_\Giry=\Bcal_{\Ycal \times \Zcal}.\]
    Let $ Y \in \Ycal^\Omega$ and $ Z \in \Zcal^\Omega$ then the composition:
    \[ (\Omega,\Bcal_\Omega) \xrightarrow{ Y \times  Z} (\Ycal \times \Zcal,\Bcal_{\Ycal \times \Zcal})
    \xrightarrow{\id_{\Ycal \times\Zcal}} (\Ycal \times \Zcal,(\Ecal_\Ycal \otimes \Bcal_\Zcal)_\Giry) \xrightarrow{K(X|Y,Z)}\Giry(\Xcal,\Bcal_\Xcal)=\PrPf(\Xcal)  \]
    is measurable and thus an element in $\PrPf(\Xcal)^\Omega$ by the remarks above.
    Since this holds for all $ Y \in \Ycal^\Omega$ and $ Z \in \Zcal^\Omega$ we get that:
    \[ K(X|Y,Z) \in \PrPf(\Xcal)^{\Ycal \times \Zcal}.  \]
    So for every $z \in \Zcal$ we have that: $K(X,Y|Z=z) \in \PrPf(\Xcal \times \Ycal)$ and
    $K(X|Y,Z=z)\otimes K(Y|Z=z) \in \PrPf(\Xcal \times \Ycal)$ and both agree on $\Bcal_\Xcal \otimes \Ecal_\Ycal$.
    Since $\Bcal_\Xcal$ and $\Ecal_\Ycal$ are countably separated also the product $\Bcal_\Xcal \otimes \Ecal_\Ycal$ is.
    Then again by Theorem \ref{thm:countably-separated} we get that:
    \[ \Bcal_{\Xcal \times \Ycal} = (\Bcal_\Xcal \otimes \Ecal_\Ycal)_{\PrPf(\Xcal \times \Ycal)}.\]
    This shows that the equality:
    \[K(X,Y|Z=z) = K(X|Y,Z=z) \otimes K(Y|Z=z) \]
    uniquely extends from $\Bcal_\Xcal \otimes \Ecal_\Ycal$ to $(\Bcal_\Xcal \otimes \Ecal_\Ycal)_{\PrPf(\Xcal \times \Ycal)} =\Bcal_{\Xcal \times \Ycal}$.\footnote{Even though $\PrQm(\Xcal \times \Ycal)$ is defined on the same $\sigma$-algebra $\Bcal_{\Xcal \times \Ycal}$ its probability measures might not come as completions of $\Bcal_\Xcal\otimes\Ecal_\Ycal$ unless they are from $\PrPf(\Xcal \times \Ycal)$. This is where the arguments for $\PrQm$  would break.\label{fn:Q-compl}}
    This shows the surjectivity of the map under the first assumptions. \\
    Since the same arguments hold also when $\Zcal$ is replaced by $\Zcal \times \Omega$ the map is even a quotient map, i.e.\ surjective on the level of quasi-measurable functions.\\
    2.) Now assume that $(\Zcal,\Zcal^\Omega)$ is a standard quasi-universal space and
    $(\Ycal,\Ycal^\Omega)$ countably separated. Then w.l.o.g.\ we can assume $\Zcal \in \Bcal_\Omega$ and
    $\Bcal_\Zcal=\Bcal_{\Omega|\Zcal}$.
    Let $\iota:\; \Zcal \inj \Omega$ be the (quasi-)measurable inclusion map and define the measurable projection map:
    \[ Z:\; \Omega \srj \Zcal,\qquad  Z_{|\Zcal}:=\id_\Zcal,\qquad  Z_{|\Omega\sm\Zcal}:= z_0 \in \Zcal.\]
    Then $ Z \in \Fcal(\Bcal_\Zcal)=\Zcal^\Omega$, since $\Zcal$ is a standard quasi-universal space, and: $ Z \circ \iota = \id_\Zcal$.\\
    We then consider the composition of quasi-measurable maps:
    \[ K(X,Y|W):\; \Omega \xrightarrow{ Z} \Zcal \xrightarrow{K(X,Y|Z)} \PrPf(\Xcal \times \Ycal),\]
    which shows that $K(X,Y|W) \in \PrPf(\Xcal \times \Ycal)^\Omega$.
    By definition there exist $ X \in \Xcal^\Omega$ and $ Y \in \Ycal^\Omega$
    and $K(U|W) \in \Giry(\Omega)^\Omega$
    such that:
    \[ K(X,Y|Z= Z)=K(X,Y|W) \stackrel{!}{=} K( X(U), Y(U)|W) := ( X\times Y)_*K(U|W).\]
    Then consider:
    \[ K(U,Y|Z):= K(U, Y(U)|W=\iota(Z)) := (\id_\Omega\times Y\iota)_*K(U|W=\iota(Z))
    \in \PrPf(\Omega \times \Ycal)^\Zcal.\]
    Note that we have:
    \[ ( X\iota\times \id_\Ycal)_*K(U,Y|Z) = K(X,Y|Z).  \]
    Then the $Y$-marginal of $K(U,Y|Z)$ agrees with the $Y$-marginal of $K(X,Y|Z)$, because:
    \[ K(Y|W=\iota(z))=K(Y|Z= Z(\iota(z)))=K(Y|Z=z).\]
    Since $\Omega$ is a standard quasi-universal space by the first point we now have:
    \[ K(U,Y|Z) = K(U|Y,Z) \otimes K(Y|Z) \in \PrPf(\Omega \times \Ycal)^\Omega,  \]
    with $K(U|Y,Z) \in \PrPf(\Omega)^{\Ycal \times \Zcal}$ and $K(Y|Z) \in \PrPf(\Ycal)^\Zcal$.
    We now put:
    \[ K(X|Y,Z) := ( X\iota)_*K(U|Y,Z) \in \PrPf(\Xcal)^{\Ycal \times \Zcal}.\]
    With this we get
    \begin{align*}
        K(X,Y|Z) &= ( X\iota \times \id_\Ycal)_*K(U,Y|Z) = ( X\iota)_*K(U|Y,Z) \otimes K(Y|Z)\\
                 &= K(X|Y,Z) \otimes K(Y|Z).
    \end{align*}
    This shows the surjectivity for the second case.\\
    To see that these maps are quotient maps again note that we can use the same proof by replacing $\Zcal$ with $\Zcal \times \Omega$. For this note that $\Zcal \times \Omega$ is again a standard quasi-universal space if $\Zcal$ is, see Lemma \ref{lem:uqus-count-prod}.

    3.) Assume finally that $(\Ycal,\Ycal^\Omega)$ is discrete, with $(\Xcal,\Xcal^\Omega)$ and
    $(\Zcal,\Zcal^\Omega)$ arbitrary quasi-universal spaces. This is the only case in which no
    space is required to be standard, and it is also the cleanest, because by
    \Cref{lem:prod-discrete} we have the equality of $\sigma$-algebras
    \[ \Bcal_{\Xcal \times \Ycal} = \Bcal_\Xcal \otimes \Bcal_\Ycal, \]
    so the extension step that the first three cases need does not arise.

    Composing $K(X,Y|Z)$ with the restriction map as before gives a measurable Markov kernel
    \[ K(X,Y|Z):\; (\Zcal,\Bcal_\Zcal) \to \Giry(\Xcal \times \Ycal, \Bcal_\Xcal \otimes \Bcal_\Ycal), \]
    and \Cref{cor:markov-kernel-disint-meas-II} in its third case --- which requires nothing of
    $(\Xcal,\Bcal_\Xcal)$ and nothing of $(\Zcal,\Bcal_\Zcal)$, and no completion, see
    \Cref{rem:disint-cor-case3} --- yields a measurable
    \[ K(X|Y,Z):\; (\Ycal \times \Zcal, \Bcal_\Ycal \otimes \Bcal_\Zcal) \to \Giry(\Xcal,\Bcal_\Xcal), \]
    with $K(X,Y|Z=z)=K(X|Y,Z=z) \otimes K(Y|Z=z)$ on $\Bcal_\Xcal \otimes \Bcal_\Ycal$, hence on all
    of $\Bcal_{\Xcal \times \Ycal}$, for every $z \in \Zcal$. Explicitly, by
    \Cref{eq:disint-discrete-quotient},
    \begin{align}
        K(X \in A|Y=y,Z=z) &= \frac{K(X \in A,Y=y|Z=z)}{K(Y=y|Z=z)},
        \quad \text{whenever } K(Y=y|Z=z)>0.
        \label{eq:disint-qus-quotient}
    \end{align}

    It remains to see that $K(X|Y,Z)$ lands in $\PrPf(\Xcal)$ and is quasi-measurable, i.e.\ that
    $K(X|Y,Z) \in \PrPf(\Xcal)^{\Ycal \times \Zcal}$. Here we cannot argue as in the first case,
    since $(\Xcal,\Xcal^\Omega)$ is not assumed standard and $\PrPf(\Xcal) \subsetneq \Giry(\Xcal)$
    in general. Instead we use \Cref{eq:disint-qus-quotient} directly.
    Let $(Y,Z) \in (\Ycal \times \Zcal)^\Omega$. Since we are in $\QUS$ we have
    $\PrPf = \PrPF$ by \Cref{thm:qus-S=P-in-R}, so from
    $K(X,Y|Z=Z(\cdot)) \in \PrPf(\Xcal \times \Ycal)^\Omega$ we obtain
    $V=(V_\Xcal,V_\Ycal) \in \lp(\Xcal \times \Ycal)^\Omega\rp^\Omega$ and
    $\mu \in \PrQm(\Omega)^\Omega$ with
    \[ K(X,Y|Z=Z(\omega)) = V(\omega)_*\mu(\omega) \qquad \text{for all }\omega \in \Omega. \]
    Put $B(\omega) := V_\Ycal(\omega)^{-1}\lp \lC Y(\omega) \rC \rp \ins \Omega$ and
    \[ \tilde\mu(\omega) := \begin{cases}
        \mu(\omega)\lp\, \cdot \mid B(\omega)\rp, & \text{if } \mu(\omega)(B(\omega))>0,\\
        \mu_0, & \text{otherwise,}
    \end{cases} \]
    for a fixed $\mu_0 \in \PrQm(\Omega)$. By \Cref{eq:disint-qus-quotient} we then have
    $K(X|Y=Y(\omega),Z=Z(\omega)) = V_\Xcal(\omega)_*\tilde\mu(\omega)$, with
    $V_\Xcal \in (\Xcal^\Omega)^\Omega$, so it suffices to show $\tilde\mu \in \PrQm(\Omega)^\Omega$.
    Each $\tilde\mu(\omega)$ lies in $\PrQm(\Omega)=\Giry(\Omega,\Bcal_\Omega)$ by
    \Cref{lem:univ-qus-G-Q}, $(\Omega,\Omega^\Omega)$ being standard. For the measurability, note
    first that
    \[ \tilde B := \lC (\tilde\omega,\omega) \in \Omega \times \Omega \st \tilde\omega \in B(\omega) \rC
       = \lp \tilde\omega,\omega \rp \mapsto \lp V_\Ycal(\omega)(\tilde\omega), Y(\omega) \rp
       \text{ pulled back along } \Delta_\Ycal \]
    lies in $\Bcal_{\Omega \times \Omega}$: the displayed map is quasi-measurable
    $\Omega \times \Omega \to \Ycal \times \Ycal$ and $\Delta_\Ycal \in \Bcal_{\Ycal \times \Ycal}$
    by \Cref{lem:prod-discrete}. Hence for every $D \in \Bcal_{\Omega \times \Omega}$ also
    $D \cap \tilde B \in \Bcal_{\Omega \times \Omega}$, and
    $\mu \in \PrQm(\Omega)^\Omega$ makes both
    \[ \omega \mapsto \mu(\omega)\lp (D \cap \tilde B)^\omega \rp
       \qquad\text{and}\qquad
       \omega \mapsto \mu(\omega)\lp \tilde B^\omega \rp = \mu(\omega)(B(\omega)) \]
    $\Bcal_\Omega$-measurable, see \Cref{qms:def:PrQm-maps}. The set
    $\lC \omega \st \mu(\omega)(B(\omega))>0 \rC$ is therefore in $\Bcal_\Omega$, and on it
    $\omega \mapsto \tilde\mu(\omega)(D^\omega)$ is the quotient of the two, while off it the map is
    constant. So $\tilde\mu \in \PrQm(\Omega)^\Omega$ and thus
    $K(X|Y,Z) \circ (Y,Z) \in \PrPf(\Xcal)^\Omega$. As $(Y,Z)$ was arbitrary,
    $K(X|Y,Z) \in \PrPf(\Xcal)^{\Ycal \times \Zcal}$, which shows surjectivity.\\
    Since $(\Zcal,\Zcal^\Omega)$ carries no hypothesis in this case, the same argument applies with
    $\Zcal$ replaced by $\Zcal \times \Omega$, so the map is again a quotient map.

\end{Prf}

\begin{Prf}{lem:prod-discrete}
    $\supseteq$: Let $A_y \in \Bcal_\Xcal$ for $y \in \Ycal$ and put
    $D:=\bigdcup_{y \in \Ycal} A_y \times \lC y\rC$. Let $(X,Y) \in (\Xcal \times \Ycal)^\Omega$, so
    $X \in \Xcal^\Omega$ and $Y \in \Ycal^\Omega = \Fcal(\two^\Ycal)$. Then
    $Y^{-1}(\lC y\rC) \in \Bcal_\Omega$, because $Y$ is $\Bcal_\Omega$-$\two^\Ycal$-measurable, and
    $X^{-1}(A_y) \in \Bcal_\Omega$ by the definition of $\Bcal_\Xcal=\Bcal(\Xcal^\Omega)$. Since
    $\Ycal$ is countable,
    \[ (X,Y)^{-1}(D) = \bigdcup_{y \in \Ycal} \lp X^{-1}(A_y) \cap Y^{-1}(\lC y\rC) \rp \in \Bcal_\Omega, \]
    as a countable disjoint union. As $(X,Y)$ was arbitrary, $D \in \Bcal_{\Xcal \times \Ycal}$.\\
    $\subseteq$: Let $D \in \Bcal_{\Xcal \times \Ycal}$ and $y \in \Ycal$, and put
    $D^y := \lC x \in \Xcal \st (x,y) \in D \rC$. The constant map $Y_y:\, \Omega \to \Ycal$,
    $\omega \mapsto y$, lies in $\Ycal^\Omega$. So for every $X \in \Xcal^\Omega$ we have
    $(X,Y_y) \in (\Xcal \times \Ycal)^\Omega$ and
    \[ X^{-1}(D^y) = (X,Y_y)^{-1}(D) \in \Bcal_\Omega. \]
    Since $X$ was arbitrary this gives $D^y \in \Bcal_\Xcal$, and
    $D = \bigdcup_{y \in \Ycal} D^y \times \lC y \rC$ has the required form.\\
    For the last claim apply the equality twice to get
    $\Bcal_{\Ycal \times \Ycal}=\two^\Ycal \otimes \two^\Ycal$, which contains
    $\Delta_\Ycal = \bigdcup_{y \in \Ycal} \lC y \rC \times \lC y \rC$.
\end{Prf}

\subsection{Countably Separated Quasi-Universal Spaces}

\subsubsection{Countably Separated Spaces}
\begin{Prf}{qms:qus:rem:auto-4}
        The first claim is clear. \\
        For the second claim note that by \cite{Bog07} Thm.\ 6.5.7 for countably separated spaces we have:
        \[ \Delta_\Xcal \in \Bcal_\Xcal \otimes \Bcal_\Xcal \ins \Bcal_{\Xcal \times \Xcal}. \]
        For the third claim take $\Ecal:= \lC \pr_n^{-1}(A_n)\st n \in \N, A_n \in \Ecal_n \rC$ on the product if $\Ecal_n$ is countable and separates the points of $\Xcal_n$, $n \in \N$.

\end{Prf}
\begin{Lem}
    \label{lem:count-sep-qms-mono}
    Let $(\Ycal,\Ycal^\Omega)$ be a quasi-universal space. Then the following statements are equivalent:
    \begin{enumerate}
        \item $(\Ycal,\Ycal^\Omega)$ is countably separated.
        \item There exists an injective measurable map:
            \[\iota:\; (\Ycal,\Bcal(\Ycal^\Omega)) \inj (\R,\Bcal_\R).\]
        \item There exists an injective measurable map:
            \[\iota:\; (\Ycal,\Bcal(\Ycal^\Omega)) \inj (\Omega,\Bcal_\Omega).\]
        \item There exists an injective quasi-measurable map:
            \[\iota:\; (\Ycal,\Ycal^\Omega) \inj (\Omega,\Omega^\Omega).\]
        \item There exists a mono-morphism $(\Ycal,\Ycal^\Omega) \inj (\Omega,\Omega^\Omega)$ in $\QUS$.
    \end{enumerate}
        \begin{proof}
        4. $\iff$ 5.: These are the same statements.\\
        1. $\iff$ 2. $\iff$ 3.: This follows from \cite{Bog07} Thm.\ 6.5.7 and the fact that
        $\Bcal(\Ycal^\Omega)$ is universally complete by Lemma \ref{lem:qus-universally-complete}.\\
        3. $\implies$ 4.: Apply functor $\Fcal$ and use $\Ycal^\Omega \ins \Fcal\Bcal(\Ycal^\Omega)$.\\
        4. $\implies$ 3.: Apply functor $\Bcal$.
        \end{proof}
\end{Lem}\begin{Lem}
    \label{lem:sturdily-countably-separated-spaces}
    Let $(\Ycal,\Ycal^\Omega)$ be a quasi-universal space. Then the following statements are equivalent:
    \begin{enumerate}
        \item There exists a countably generated $\sigma$-algebra $\Ecal_\Ycal$
    that separates the points of $\Ycal$ such that $\Fcal(\Ecal_\Ycal)=\Ycal^\Omega$.
        \item There exists a universally countably generated and separated $\sigma$-algebra $\Ecal_\Ycal$  on $\Ycal$ such that $\Fcal(\Ecal_\Ycal)=\Ycal^\Omega$.
        \item There exists an injective quasi-measurable map:
            \[\iota:\; (\Ycal,\Ycal^\Omega) \inj (\Omega,\Omega^\Omega), \]
            such that: $\Ycal^\Omega = i^*\Omega^\Omega$.
        \item There exists an embedding $(\Ycal,\Ycal^\Omega) \inj (\Omega,\Omega^\Omega)$ in $\QUS$.
    \end{enumerate}
    \begin{proof}
        3. $\iff$ 4.: The first is the definition of the latter.\\
        Now let $\Bcal_\R \ins \Bcal_\Omega$ be the Borel $\sigma$-algebra of $\R$.\\
        2. $\iff$ 3.: Since $\Bcal_\Omega$ is universally complete we have:
        $\Fcal(\Ecal_\Ycal)=\Fcal((\Ecal_\Ycal)_\Giry)$.\\
        1. $\implies$ 3.: $\Ecal_\Ycal$ induces an injective measurable map:
        \[\iota:\; (\Ycal,\Ecal_\Ycal) \inj (\R,\Bcal_\R),\]
        such that $\Ecal_\Ycal=\iota^*\Bcal_\R$, see \cite{Bog07} Thm.\ 6.5.8.
        Applying the functor $\Fcal$ gives us an injective quasi-measurable map:
        \[\iota:\; (\Ycal,\Ycal^\Omega)=(\Ycal,\Fcal(\Ecal_\Ycal)) \inj (\R,\Fcal(\Bcal_\R)) = (\Omega,\Omega^\Omega).\]
        So we have the inclusion $\Ycal^\Omega \ins \iota^*\Omega^\Omega$.
        For the reverse inclusion let: $ Y \in \iota^*\Omega^\Omega$.
        Then $ Y \in \Sets(\Omega,\Ycal)$ with:
        \[\iota \circ  Y \in \Omega^\Omega=\Fcal(\Bcal_\R) = \Meas\lp(\Omega,\Bcal_\Omega),(\R,\Bcal_\R)\rp.\]
        But this means that $ Y$ is $\Bcal_\Omega$-$\iota^*\Bcal_\R$-measurable.
        Since $\Ecal_\Ycal=\iota^*\Bcal_\R$ we have that $ Y \in \Fcal(\Ecal_\Ycal)$.
        So we get $\iota^*\Omega^\Omega \ins \Fcal(\Ecal_\Ycal)=\Ycal^\Omega$.\\
        3. $\implies$ 1.: We assume that we have an injective quasi-measurable map:
            \[\iota:\; (\Ycal,\Ycal^\Omega) \inj (\Omega,\Omega^\Omega), \]
            such that: $\Ycal^\Omega = i^*\Omega^\Omega$. Applying the functor $\Bcal$ to it gives an
            injective measurable map:
            \[\iota:\; (\Ycal,\Bcal(\Ycal^\Omega)) \inj (\Omega,\Bcal_\Omega). \]
        We put:
        \[\Ecal_\Ycal:=\iota^*\Bcal_\R \ins \iota^*\Bcal_\Omega \ins \Bcal(\Ycal^\Omega).\]
        Then clearly $\Ecal_\Ycal$ separates the points of $\Ycal$.
        Furthermore, we have:
        \[ \Ycal^\Omega \ins \Fcal\Bcal(\Ycal^\Omega) \ins \Fcal(\Ecal_\Ycal) = \Fcal((\Ecal_\Ycal)_\Giry).\]
        For the reverse inclusion let $ Y \in \Fcal(\Ecal_\Ycal) = \Fcal((\Ecal_\Ycal)_\Giry)$.
        Since $\iota$ is $\Ecal_\Ycal$-$\Bcal_\R$-measurable it is also
         $(\Ecal_\Ycal)_\Giry$-$\Bcal_\Omega$-measurable.
         So the composition $\iota \circ  Y$ is $\Bcal_\Omega$-$\Bcal_\Omega$-measurable.
         So $\iota \circ  Y \in \Omega^\Omega$ and thus $ Y \in \iota^*\Omega^\Omega$.
         This shows: $\Fcal(\Ecal_\Ycal) \ins \iota^*\Omega^\Omega = \Ycal^\Omega$.
    \end{proof}
\end{Lem}
\subsubsection{The Sigma-Algebra of Countably Separated Quasi-Universal Spaces}
\begin{Prf}{prp:perf-count-sep}
        We clearly have the chain of inclusions and equalities:
        \begin{align*}
          \Ecal \ins \Bcal_\Xcal \ins   \Bcal(\Fcal)
            &:= \lC A \ins \Xcal\st \forall  X \in \Fcal.\, X^{-1}(A) \in \Bcal_\Omega  \rC \\
            &= \bigcap_{ X \in \Fcal}  X_*\Bcal_\Omega \\
            &\ins \bigcap_{ X \in \Fcal}  X_*(\Bcal_\Omega)_\PrPf \\
            &= \bigcap_{ X \in \Fcal}  X_*\bigcap_{\mu \in \PrPf} \lp\Bcal_\Omega \rp_\mu \\
            &= \bigcap_{ X \in \Fcal} \bigcap_{\mu \in \PrPf}  X_*\lp\Bcal_\Omega \rp_\mu \\
            &= \bigcap_{ X \in \Fcal} \bigcap_{\mu \in \PrPf} \lp\Ecal \rp_{ X_*\mu} \\
            &= \bigcap_{\nu \in \PrPf(\Fcal)} \lp\Ecal \rp_{\nu}\\
            &= \lp \Ecal \rp_{\PrPf(\Fcal)},
    \end{align*}
    where the equality: $ X_*\lp\Bcal_\Omega \rp_\mu = \lp\Ecal \rp_{ X_*\mu}$
    comes from Lemma \ref{qms:lem:perf-count-sep} and the fact that $\Ecal$ is countably separated
    and that every $ X \in \Fcal$ is a $\Bcal_\Omega$-$\sigma(\Ecal)$-measurable map and that
    every $\mu \in \PrPf$ is perfect.
    Together we get the inclusions:
    \[\Ecal \ins \Bcal_\Xcal \ins \Bcal(\Fcal) \ins \lp \Ecal \rp_{\PrPf(\Fcal)}, \]
    where the latter is an equality if $\Bcal_\Omega = \lp \Bcal_\Omega \rp_\PrPf$.
    Taking the completions w.r.t.\ $\PrPf(\Fcal)$ we get:
    \[ \lp \Ecal \rp_{\PrPf(\Fcal)} \ins \lp \Bcal_\Xcal \rp_{\PrPf(\Fcal)} \ins \Bcal(\Fcal)_{\PrPf(\Fcal)} \ins \lp \Ecal \rp_{\PrPf(\Fcal)},  \]
    which imply the equalities:
    \[ \lp \Ecal \rp_{\PrPf(\Fcal)} = \lp \Bcal_\Xcal \rp_{\PrPf(\Fcal)} = \Bcal(\Fcal)_{\PrPf(\Fcal)}. \]
    This shows the claim.

\end{Prf}
\begin{Prf}{thm:countably-separated}
        This follows directly from Proposition \ref{prp:perf-count-sep} together with the fact that every probability measure on $(\Omega,\Bcal_\Omega)$ is perfect and $(\Omega,\Bcal_\Omega)$ is universally complete.
\end{Prf}
\begin{Prf}{qms:qus:cor:auto-13}
        By Theorem \ref{thm:countably-separated} we have that $\Bcal(\Xcal^\Omega) = \lp \Ecal \rp_{\PrPf(\Xcal,\Xcal^\Omega)}$ for a countable algebra $\Ecal \ins \Bcal(\Xcal^\Omega)$.
        Then we have that:
        \[ \Ecal_\PrPf:= \lC \ev_A^{-1}\lp(t,1]\rp\st A \in \Ecal, t \in \Q\rC \ins \Bcal\lp\PrPf(\Xcal,\Xcal^\Omega)^\Omega \rp\]
        is countable and separates the points of $\PrPf(\Xcal,\Xcal^\Omega)$.
        Indeed, two measures $\mu_1$, $\mu_2$ agree on $\sigma(\Ecal)$
        if they agree on the generating $\pi$-system $\Ecal$. Since
        $\mu_1,\mu_2 \in \PrPf(\Xcal,\Xcal^\Omega)$ they uniquely extend to $\Bcal(\Xcal^\Omega) = \lp \Ecal \rp_{\PrPf(\Xcal,\Xcal^\Omega)}$ and thus are then equal there as well.\\
        Note that the above set $\Ecal_\PrPf$ really lies inside the $\sigma$-algebra
        $\Bcal\lp\PrPf(\Xcal,\Xcal^\Omega)^\Omega \rp$
        as the evaluation maps $\ev_A$ are measurable.
\end{Prf}

\subsection{Universal Measurable and Standard Quasi-Universal Spaces}

\subsubsection{Universal Measurable Spaces}
\begin{Prf}{deflem:univ-meas}
    2.\ $\iff$ 3.: This is just an explicit reformulation. Note that if $(\Ecal)_\Giry = (\Bcal_\Xcal)_\Giry$ then we have the equality of sets:
    \[ \Giry(\Xcal,\Ecal)= \Giry(\Xcal,(\Ecal)_\Giry) = \Giry(\Xcal,(\Bcal_\Xcal)_\Giry) = \Giry(\Xcal,\Bcal_\Xcal).\]
    Also note that a probability measure is perfect if and only if its completion is perfect.\\
    1.\ $\implies$ 4.: Let $i$ be the composition of the measurable isomorphism and the measurable inclusion map:
    \[i:\; (\Xcal,(\Bcal_\Xcal)_{\Giry})\cong (\Zcal,\Bcal_{\Omega|\Zcal}) \inj (\Omega,\Bcal_\Omega). \]
    Define the map $r$ via:
    \[ r:\; (\Omega,\Bcal_\Omega) \to (\Xcal,(\Bcal_\Xcal)_{\Giry}),
       \qquad r|_{\Zcal}:=i^{-1},\qquad r|_{\Omega \sm \Zcal}:=x_0 \in \Xcal. \]
    Since $\Zcal \in \Bcal_\Omega$ the map $r$ is measurable. It is clear that $r \circ i = \id_\Xcal$.
    This shows that $(\Xcal,(\Bcal_\Xcal)_{\Giry})$ is a retract of $(\Omega,\Bcal_\Omega)$.\\
    4.\ $\implies$ 3.: Put $\Ecal:=i^*\Bcal_\R \ins (\Bcal_\Xcal)_\Giry$.
    Then $\Ecal$ is countably generated and separates the points of $\Xcal$ since $i$ is injective.
    So $(\Bcal_\Xcal)_\Giry$ is universally countably separated. Let $\mu \in \Giry(\Xcal,\Ecal)$.
    Furthermore, $i$ is then $\Ecal$-$\Bcal_\R$-measurable and thus $(\Ecal)_\Giry$-$\Bcal_\Omega$-measurable.
    Note that $r$ is $\Bcal_\Omega$-$(\Bcal_\Xcal)_\Giry$-measurable.
    Together we get for $A \in (\Bcal_\Xcal)_\Giry$ that $r^{-1}(A) \in \Bcal_\Omega$ and thus:
    \[ A = \id_\Xcal^{-1}(A) = i^{-1}(r^{-1}(A)) \in (\Ecal)_\Giry.\]
    This shows the inclusions:
    \[ (\Bcal_\Xcal)_\Giry \ins (\Ecal)_\Giry \ins (\Bcal_\Xcal)_\Giry,  \]
    and thus equality. This shows that $(\Xcal,\Bcal_\Xcal)$ is universally countably generated.
    Let $\mu \in \Giry(\Xcal,(\Bcal_\Xcal)_\Giry)$. Then $i_*\mu \in \Giry(\Omega,\Bcal_\Omega)$. Since all probability measures
    on $(\Omega,\Bcal_\Omega)$ are perfect so is $i_*\mu$. Since also push-forwards of perfect probability measures are
    perfect also:
    \[ \mu = \id_{\Xcal,*}\mu = r_*(i_*\mu), \]
    is perfect. This shows that all probability measures on $(\Xcal,\Bcal_\Xcal)$ are perfect.\\
    3.\ $\implies$ 1.: Let $\Ecal \ins (\Bcal_\Xcal)_\Giry$ be a countably generated $\sigma$-algebra that separates the points of
    $\Xcal$ such that $(\Ecal)_\Giry = (\Bcal_\Xcal)_\Giry$ and $(\Xcal,\Ecal)$ is perfect.
    $\Ecal$ then induces an injective $\Ecal$-$\Bcal_\R$-measurable map $i:\; \Xcal \to \Omega$
  such that $\Ecal=i^*\Bcal_\R$ (see \cite{Bog07} Thm.\ 6.5.8.). If we put $\Zcal:=i(\Xcal)$
  then $i$ induces an isomorphism of measurable spaces:
  \[(\Xcal,\Ecal)\cong (\Zcal,\Bcal_{\R|\Zcal}).\]
  Since $(\Xcal,\Ecal)$ and thus $(\Zcal,\Bcal_{\R|\Zcal})$ is perfect we get that $\Zcal:=i(\Xcal) \in \Bcal_\Omega$
  (see \cite{Bog07} Thm.\ 7.5.7 or \cite{Dar71,Saz62} Lem.\ 3). With $\Zcal \in \Bcal_\Omega$ we get that
  $(\Bcal_{\R|\Zcal})_\Giry = \Bcal_{\Omega|\Zcal}$ and thus an measurable isomorphism:
  \[(\Xcal,(\Ecal)_\Giry)\cong (\Zcal,(\Bcal_{\R|\Zcal})_\Giry) = (\Zcal, \Bcal_{\Omega|\Zcal}), \]
  with $\Zcal \in \Bcal_\Omega$.\\
  This shows all claims.
\end{Prf}
\begin{Lem}
    \label{lem:univ-meas-G}
    Let $(\Xcal,\Bcal_\Xcal)$ be a universal measurable space then  we get that:
    \[(\Xcal,\Bcal\Fcal(\Bcal_\Xcal)) = (\Xcal,(\Bcal_\Xcal)_\Giry) \]
        is the universal completion of $(\Xcal,\Bcal_\Xcal)$.
        Furthermore, we have the equality of sets:
   \begin{align*}
        \Giry(\Xcal,\Bcal_\Xcal) &= \Giry(\Xcal,\Xcal^\Omega) = \PrPf(\Xcal,\Xcal^\Omega)\\
        &:=\lC X_*P :\; \Bcal(\Xcal^\Omega) \to [0,1]\st  X \in \Xcal^\Omega,\; P \in
        \Giry(\Omega,\Bcal_\Omega)\footref{fn:G-Q-1}\rC,
    \end{align*}
    where we put: $\Xcal^\Omega:=\Fcal(\Bcal_\Xcal)$.
 \begin{proof}
    By Definition/Lemma \ref{deflem:univ-meas} there are measurable maps:
    \[ \id_\Xcal:\; (\Xcal,(\Bcal_\Xcal)_\Giry) \xrightarrow{i} (\Omega,\Bcal_\Omega) \xrightarrow{r} (\Xcal,(\Bcal_\Xcal)_\Giry), \]
    which show the inclusions and thus equalities:
    \[  (\Bcal_\Xcal)_\Giry \sni i^*\Bcal_\Omega \sni i^*r^*(\Bcal_\Xcal)_\Giry=\id_\Xcal^*(\Bcal_\Xcal)_\Giry =(\Bcal_\Xcal)_\Giry. \]
    Applying the functor $\Bcal\Fcal$ to the above gives the measurable maps:
    \[ \id_\Xcal:\; (\Xcal,\Bcal\Fcal(\Bcal_\Xcal)) \xrightarrow{i} (\Omega,\Bcal_\Omega) \xrightarrow{r} (\Xcal,\Bcal\Fcal(\Bcal_\Xcal)). \]
            For the latter note that since $\Bcal_\Omega$ is universally complete we have $\Fcal(\Bcal_\Xcal)=\Fcal((\Bcal_\Xcal)_\Giry)$, and also:  $\Bcal\Fcal(\Bcal_\Omega)=\Bcal\Fcal\Bcal(\Omega^\Omega)=\Bcal(\Omega^\Omega)=\Bcal_\Omega$.\\
    So we get the inclusions and thus equalities:
    \[  \Bcal\Fcal(\Bcal_\Xcal)\sni i^*\Bcal_\Omega \sni i^*r^*\Bcal\Fcal(\Bcal_\Xcal)=\id_\Xcal^*\Bcal\Fcal(\Bcal_\Xcal) =\Bcal\Fcal(\Bcal_\Xcal). \]
    Together we then get:
    \[  (\Bcal_\Xcal)_\Giry = i^*\Bcal_\Omega =\Bcal\Fcal(\Bcal_\Xcal).\]
    This shows the first claim.\\
    Now let $\Xcal^\Omega:=\Fcal(\Bcal_\Xcal)$ and $\mu \in \Giry(\Xcal,\Xcal^\Omega)$.
    Then we get:
     \[ \mu = r_*(i_*\mu) \in \PrPf(\Xcal,\Xcal^\Omega) \ins \Giry(\Xcal,\Xcal^\Omega),\]
     since $r \in \Fcal(\Bcal_\Xcal)=\Xcal^\Omega$ and $i_*\mu \in \Giry(\Omega,\Bcal_\Omega)$.
     This shows:
     \[\Giry(\Xcal,\Xcal^\Omega)=\PrPf(\Xcal,\Xcal^\Omega). \]
    This shows the second claim.
    \end{proof}
\end{Lem}
\subsubsection{Standard Quasi-Universal Spaces}
\begin{Prf}{def:univ-qus}
        3.\ $\implies$ 1.: Consider the composition of quasi-measurable maps:
        \[ \id_\Xcal:\;  (\Xcal,\Xcal^\Omega) \xrightarrow{i}  (\Omega,\Omega^\Omega) \xrightarrow{r} (\Xcal,\Xcal^\Omega).\]
        This gives us the inclusions and thus equalities:
        \[ \Xcal^\Omega = \id_\Xcal \circ \Xcal^\Omega = r \circ i \circ \Xcal^\Omega \ins r \circ \Omega^\Omega \ins \Xcal^\Omega.  \]
        Applying the functor $\Bcal$ to the above gives the composition of measurable maps:
        \[ \id_\Xcal:\;  (\Xcal,\Bcal(\Xcal^\Omega)) \xrightarrow{i}  (\Omega,\Bcal_\Omega) \xrightarrow{r} (\Xcal,\Bcal(\Xcal^\Omega)).\]
        So $(\Xcal,\Bcal(\Xcal^\Omega))$ is a retract of $(\Omega,\Bcal_\Omega)$ in $\Meas$ and universally complete by Lemma \ref{lem:qus-universally-complete}, showing that $(\Xcal,\Bcal(\Xcal^\Omega))$ is a (universally complete) universal measurable space.
        Again, applying the functor $\Fcal$ to the above gives the composition of quasi-measurable maps:
        \[ \id_\Xcal:\;  (\Xcal,\Fcal\Bcal(\Xcal^\Omega)) \xrightarrow{i}  (\Omega,\Omega^\Omega) \xrightarrow{r} (\Xcal,\Fcal\Bcal(\Xcal^\Omega)).\]
        This gives us the inclusions and thus equalities:
        \[ \Fcal\Bcal(\Xcal^\Omega) = \id_\Xcal \circ \Fcal\Bcal(\Xcal^\Omega) = r \circ i \circ \Fcal\Bcal(\Xcal^\Omega) \ins r \circ \Omega^\Omega \ins \Fcal\Bcal(\Xcal^\Omega).  \]
        Together we get:
        \[ \Xcal^\Omega = r \circ \Omega^\Omega = \Fcal\Bcal(\Xcal^\Omega).\]
        1.\ $\implies$ 2.: Pick $\Bcal_\Xcal:=\Bcal(\Xcal^\Omega)$.\\
        2.\ $\implies$ 1.: If $(\Xcal,\Bcal_\Xcal)$ is a universal measurable space with $\Xcal^\Omega=\Fcal(\Bcal_\Xcal)$ then
        by Lemma \ref{lem:univ-meas-G} we have that $(\Xcal,\Bcal(\Xcal^\Omega))=(\Xcal,\Bcal\Fcal(\Bcal_\Xcal))$ is a
        (universally complete) universal measurable space and $(\Bcal_\Xcal)_\Giry=\Bcal\Fcal(\Bcal_\Xcal)$.
        Applying $\Fcal$ to the latter we get:
        \[ \Xcal^\Omega = \Fcal(\Bcal_\Xcal) = \Fcal((\Bcal_\Xcal)_\Giry) = \Fcal\Bcal\Fcal(\Bcal_\Xcal) = \Fcal\Bcal(\Xcal^\Omega).\]
        1. $\implies$ 3.: Since $(\Xcal,\Bcal(\Xcal^\Omega))$ is a (universally complete) universal measurable space we by
        Definition/Lemma \ref{deflem:univ-meas} have a composition of measurable spaces:
       \[ \id_\Xcal:\;  (\Xcal,\Bcal(\Xcal^\Omega)) \xrightarrow{i}  (\Omega,\Bcal_\Omega) \xrightarrow{r} (\Xcal,\Bcal(\Xcal^\Omega)).\]
       Applying the functor $\Fcal$ to it and exploiting $\Xcal^\Omega=\Fcal\Bcal(\Xcal^\Omega)$ gives a composition of quasi-measurable maps:
       \[ \id_\Xcal:\;  (\Xcal,\Xcal^\Omega)=(\Xcal,\Fcal\Bcal(\Xcal^\Omega)) \xrightarrow{i}  (\Omega,\Omega^\Omega) \xrightarrow{r} (\Xcal,\Fcal\Bcal(\Xcal^\Omega))= (\Xcal,\Xcal^\Omega).\]
        This shows that $(\Xcal,\Xcal^\Omega)$ is a retract of $(\Omega,\Omega^\Omega)$ in $\QUS$.

\end{Prf}
\begin{Lem}
    \label{lem:univ-qus-G-P}
    Let $(\Xcal,\Xcal^\Omega)$ be a standard quasi-universal space.
    Then we have the equality of sets:
   \[ \Giry(\Xcal,\Xcal^\Omega) =
    \PrPf(\Xcal,\Xcal^\Omega):=\lC X_*\nu :\; \Bcal(\Xcal^\Omega) \to [0,1]\st  X \in \Xcal^\Omega, \nu \in
    \Giry(\Omega,\Bcal_\Omega)\footref{fn:G-Q-1}\rC.\]
    \begin{proof}
        This follows from Definition/Lemma \ref{def:univ-qus} and Lemma \ref{lem:univ-meas-G}.
    \end{proof}
\end{Lem}\begin{Lem}
    \label{lem:qus-sample-prod-iso}
    Let $J$ be a countable set (e.g.\ $J=\N$).
    Then there is an isomorphism of quasi-universal spaces (in $\QUS$):
    \[ (\Omega,\Omega^\Omega) \cong \prod_{j \in J} (\Omega,\Omega^\Omega). \]
    \begin{proof}
        Since $J$ is countable by \cite{Fremlin} 424C there exists an isomorphism of measurable spaces (in $\Meas$):
        \[ (\R,\Bcal_\R) \cong \prod_{j \in J} (\R,\Bcal_\R) = \lp\prod_{j \in J} \R, \bigotimes_{j \in J} \Bcal_\R \rp. \]
        We now apply the functor $\Fcal$ and exploit that $\Fcal$ as a right-adjoint preserves products, see Lemma \ref{qms:lem:F-prod}:
        \[ (\R,\Fcal(\Bcal_\R)) \cong \lp\prod_{j \in J} \R, \Fcal\lp\bigotimes_{j \in J} \Bcal_\R\rp \rp
            = \lp\prod_{j \in J} \R, \prod_{j \in J} \Fcal(\Bcal_\R) \rp.
\]
    Since $(\Omega,\Bcal_\Omega)=(\R,(\Bcal_\R)_\Giry)$ is universally complete we have:
    \[ \Fcal(\Bcal_\R)= \Fcal((\Bcal_\R)_\Giry)=\Fcal(\Bcal_\Omega)=\Omega^\Omega.\]
    This gives us the desired isomorphism of quasi-universal spaces:
    \[(\Omega,\Omega^\Omega)= (\R,\Fcal(\Bcal_\R)) \cong  \lp\prod_{j \in J} \R, \prod_{j \in J} \Fcal(\Bcal_\R) \rp= \lp\prod_{j \in J} \Omega, \prod_{j \in J} \Omega^\Omega \rp = \prod_{j \in J} (\Omega,\Omega^\Omega). \]
    This shows the claim.
    \end{proof}
\end{Lem}\begin{Lem}[Countable products of standard quasi-universal spaces]
    \label{lem:uqus-count-prod}
    Let $(\Xcal_j,\Xcal_j^\Omega)$ be standard quasi-universal spaces for $j \in J$ for a countable index set $J$.
    Then the product of quasi-universal spaces (in $\QUS$):
    \[ \prod_{j \in J}(\Xcal_j,\Xcal_j^\Omega) = \lp \prod_{j \in J} \Xcal_j, \prod_{j \in J} \Xcal_j^\Omega \rp   \]
    is again a standard quasi-universal space.
    \begin{proof}
        By Definition/Lemma \ref{def:univ-qus} we have compositions of quasi-measurable maps for $j \in J$:
        \[ \id_{\Xcal_j}:\; (\Xcal_j,\Xcal_j^\Omega) \xrightarrow{i_j} (\Omega,\Omega^\Omega) \xrightarrow{r_j}(\Xcal_j,\Xcal_j^\Omega). \]
        Taking the product of all those maps for $j \in J$ gives the composition of quasi-measurable maps:
        \[ \id_{\prod_{j \in J}\Xcal_j}:\; \prod_{j \in J}(\Xcal_j,\Xcal_j^\Omega) \xrightarrow{(i_j)_{j \in J}} \prod_{j \in J}(\Omega,\Omega^\Omega) \xrightarrow{(r_j)_{j \in J}}\prod_{j \in J}(\Xcal_j,\Xcal_j^\Omega). \]
        Since $J$ is countable by Lemma \ref{lem:qus-sample-prod-iso} we now have an isomorphism of quasi-universal spaces:
        \[ \prod_{j \in J} (\Omega,\Omega^\Omega) \cong (\Omega,\Omega^\Omega),\]
        showing together with the above that $\prod_{j \in J}(\Xcal_j,\Xcal_j^\Omega)$ is a retract of
    $(\Omega,\Omega^\Omega)$ in $\QUS$.
    This shows the claim that $\prod_{j \in J}(\Xcal_j,\Xcal_j^\Omega)$ is also a standard quasi-universal space.
    \end{proof}
\end{Lem}\begin{Prf}{cor:umeas-uqus-eq}
        This immediately follows from Theorem \ref{thm:qms-meas-adjunction} and Lemmas \ref{deflem:univ-meas}, \ref{lem:univ-meas-G}, and \ref{def:univ-qus}.
        Note that $\Bcal\Fcal$ and $\Fcal\Bcal$ do not change the maps. Furthermore, $\Bcal\Fcal=\id_\UMeas$, since by Lemma
        \ref{lem:univ-meas-G}:
        \[ \Bcal\Fcal(\Bcal_\Xcal) = (\Bcal_\Xcal)_\Giry=\Bcal_\Xcal,  \]
        if $(\Xcal,\Bcal_\Xcal)$ is a universally complete universal measurable space.
        We also have that $\Fcal\Bcal=\id_\StdQUS$, since:
        \[ \Fcal\Bcal(\Xcal^\Omega)=\Xcal^\Omega,  \]
        by Definition/Lemma \ref{def:univ-qus} of standard quasi-universal spaces.

\end{Prf}

\subsection{Spaces of Probability Measures for Quasi-Universal Spaces}
\begin{Lem}
    \label{lem:univ-qus-G-Q}
    Let $(\Xcal,\Xcal^\Omega)$ be a standard quasi-universal space (e.g.\ $(\Xcal,\Xcal^\Omega)=(\Omega,\Omega^\Omega)$)
    then we have the equality of quasi-universal spaces:
    \[ \Giry(\Xcal,\Xcal^\Omega) = \PrQm(\Xcal,\Xcal^\Omega).\]
    \begin{proof}
The inclusion $\PrQm(\Xcal,\Xcal^\Omega) \ins \Giry(\Xcal,\Xcal^\Omega)$ is quasi-measurable by
    Remark \ref{rem:Q-eval}.
    Now let $\mu \in \Giry(\Xcal,\Xcal^\Omega)$ and $\kappa \in \Giry(\Xcal,\Xcal^\Omega)^\Omega$.
    For $D \in \Bcal_{\Omega \times \Xcal}$ and $\Phi \in \Omega^\Omega$
    we need to show that the maps:
    \[ \Omega \to [0,1],\qquad \omega \mapsto \mu(D_\omega);\qquad  \omega \mapsto \kappa(\Phi(\omega))(D_\omega),\]
    are $\Bcal_\Omega$-$\Bcal_{[0,1]}$-measurable. Since the former can be viewed as a special case of the latter we only focus on the $\kappa$ case. The latter map can then be factorized as:
  \[h:\; \Omega \xrightarrow{\Phi \times \kappa} \Omega \times \Giry(\Xcal,\Bcal_\Xcal) \xrightarrow{\delta \times (\id_\Xcal)_*} \Giry(\Omega\times \Xcal,\Bcal_\Omega\otimes\Bcal_\Xcal) \xrightarrow{\ev_D} [0,1],\]
    \[\omega \mapsto  (\delta_\omega\otimes\kappa(\Phi(\omega)))(D) =\kappa(\Phi(\omega))(D_\omega).\]
    The left map is measurable. The middle map is measurable as the strength of the
    Giry monad, see \cite{Gir82} or also \cite{For21} Lem.\ B.39.
    The type of measurability of the evaluation map $\ev_D$ on the right depends on the
    kind of the subset $D$,
    see Proposition \ref{prp:prob-prob-meas}. We need to analyse this further.\\
    Since countable products of universal measurable spaces are universal measurable spaces
    (see Remark \ref{rem:umeas} and Lemma \ref{lem:uqus-count-prod})
    we get from Lemma \ref{lem:univ-meas-G}:
    \[D \in \Bcal_{\Omega \times \Xcal} = \Bcal(\Omega^\Omega \times \Xcal^\Omega) = \Bcal(\Fcal(\Bcal_\Omega) \times \Fcal(\Bcal_\Xcal))
        = \Bcal\Fcal(\Bcal_\Omega \otimes \Bcal_\Xcal)
    = \lp\Bcal_\Omega \otimes \Bcal_\Xcal\rp_\Giry. \]
    Proposition \ref{prp:prob-prob-meas} shows that $\ev_D$ and thus $h$ is universally measurable,
    i.e.\  $(\Bcal_\Omega)_\Giry$-$\Bcal_{[0,1]}$-measurable.
    Since $(\Bcal_\Omega)_\Giry=\Bcal_\Omega$ it is even
    $\Bcal_\Omega$-$\Bcal_{[0,1]}$-measurable. This shows that
    $\kappa \in \PrQm(\Xcal,\Xcal^\Omega)^\Omega$. Similarly, $\mu \in \PrQm(\Xcal,\Xcal^\Omega)$.
    This shows the wanted equality of quasi-measurable spaces:
    \[ \Giry(\Xcal,\Xcal^\Omega) = \PrQm(\Xcal,\Xcal^\Omega),\]
    and thus the claim.
    \end{proof}
\end{Lem}\begin{Prf}{cor:random-fct}
        It is clear that all maps are well-defined and quasi-measurable, the diagram commutes and that the top horizontal map is a quotient map,  see Theorem \ref{thm:qus-S=P-in-R}. It is left to show that the bottom horizontal map is a quotient map, i.e.\ that it is surjective for both $\Zcal$ and
        $\Omega \times \Zcal$ in the exponent.\\
        1.) Assume that   $\Xcal$ is a standard quasi-universal space.
        By Theorem \ref{thm:univ-qus-GPQKR} we then have that:
        $(\PrPf(\Xcal),\PrPf(\Xcal)^\Omega)= (\Giry(\Xcal),\Fcal(\Bcal_{\Giry(\Xcal)}))$.
        By Definition/Lemma \ref{def:univ-qus} We also have quasi-measurable maps $\Xcal \xrightarrow{\iota} \Omega \xrightarrow{r} \Xcal$ with $r \circ \iota =\id_\Xcal$.
        Now let $\rho$ be given by the adjunction, Theorem \ref{thm:qms-meas-adjunction}:
        \[ \rho \in \PrPf(\Xcal)^\Zcal=\QUS(\Zcal,\PrPf(\Xcal))=\Meas(\Zcal,\Giry(\Xcal)). \]
        Let $ Y:\;\Zcal \times \Omega \to \Omega$ be the conditional quantile function
        of the Markov kernel $\iota_*\rho$, which is measurable, thus $ Y \in \lp\Omega^\Omega\rp^\Zcal$. We have that $ Y_*\nu = \iota_*\rho$, see e.g.\ \cite{Cen82,Kal21} or \cite{For21} Appendix G.
        Then $(r \circ  Y)_*\nu = r_*\iota_*\rho= \rho$ and $(r \circ  Y) \in \lp\Xcal^\Omega\rp^\Zcal$.
        This shows the surjectivity of $\lp\Xcal^\Omega\rp^\Zcal  \to \PrPf(\Xcal)^\Zcal$.
        By replacing $\Zcal$ with $\Omega \times \Zcal$ we can use the same arguments to see that this
        is actually a quotient map. \\
        This shows the claim if $(\Xcal,\Xcal^\Omega)$ is a standard quasi-universal space.\\
        2.) Now assume that $(\Zcal,\Zcal^\Omega)$ is a standard quasi-universal space and
         $\rho \in  \PrPf(\Xcal)^\Zcal$.
By Definition/Lemma \ref{def:univ-qus} we get quasi-measurable maps
        $ \Zcal \xrightarrow{\iota} \Omega \xrightarrow{r}  \Zcal$
        with $r \circ \iota =\id_{\Zcal}$.
        Then we have:
        \[\rho \circ r \in \PrPf(\Zcal)^\Omega = \lC  X_*\nu\st  X \in
        \lp\Xcal^\Omega\rp^\Omega\rC. \]
        So $\rho \circ r =  X_*\nu$ and we get by composing with $\iota$:
        \[\rho = \rho \circ r \circ \iota = ( X \circ \iota)_*\nu, \]
        with now $ X \circ \iota \in \lp\Xcal^\Omega\rp^{ \Zcal}$.
        This shows the surjectivity of $\lp\Xcal^\Omega\rp^\Zcal  \to \PrPf(\Xcal)^\Zcal$.
        By replacing $\Zcal$ with $\Omega \times \Zcal$ we can use the same arguments to see that this
        is actually a quotient map.
        For this note that by Lemma \ref{lem:uqus-count-prod} the product $\Omega \times \Zcal$
        is again a standard quasi-universal space.
        This shows the claim.

\end{Prf}

\subsection{Fubini Theorem for Probability Monad $\PrPf$ on Quasi-Universal Spaces}
\begin{Prf}{thm:fubini}
        It was already shown that all maps are quasi-measurable, see Theorem \ref{thm:K-P-Q-markov-prod}
        and Lemma \ref{lem:K-P-R-Q-G-qm-incl}.\footnotetext{One would more consistently use $f^s(y,x):=f(x,y)$ here in this place.\label{fn:fubini}}
        We only need to check the commutativity on elements and fix $u \in \Ucal$ and $z \in \Zcal$.
        So we can w.l.o.g.\ ignore
        $\Ucal$ and $\Zcal$.\\
        Let $\mu= X_*\kappa \in \PrPf(\Xcal)$ with $ X \in \Xcal^{\Omega_2}$,
        $\kappa \in \PrQm(\Omega_2)$ and
        $\nu = Y_*\rho \in \PrPf(\Ycal)$ with $ Y \in \Ycal^{\Omega_1}$, $\rho \in \PrQm(\Omega_1)$ and $D \in \Bcal_{\Ycal \times \Xcal}$, where $\Omega_2$, $\Omega_1$ are copies of $\Omega=\R$,
        with an index for readability. Then we use the superscript $(\_)^s$ to indicate swaps:
        \[ D^s := \lC (x,y) \in \Xcal \times \Ycal\st (y,x) \in D \rC.\]
        We need to show:
        \[ (\nu \otimes \mu)(D) = (\mu \otimes \nu)(D^s).\]
        We first show the analogue for $\kappa$ and $\rho$. For this note that by Lemma \ref{lem:univ-meas-G}:
        \[ \Bcal_{\Omega_1 \times \Omega_2} = \lp \Bcal_{\Omega_1} \otimes \Bcal_{\Omega_2}\rp_\Giry
            = \lp\lC B \times A\st B \in \Bcal_{\Omega_1}, A \in \Bcal_{\Omega_2}  \rC\rp_\Giry.
        \]
        We trivially have the required equality on the elements of the generating $\pi$-system:
        \[ (\rho \otimes \kappa)(B \times A) = \rho(B) \cdot \kappa(A) =
            \kappa(A) \cdot \rho(B) = (\kappa\otimes\rho)(A \times B) =
            (\kappa\otimes\rho)( (B \times A)^s).
        \]
        This equality then extends to its universal completion and we get for all $E \in \Bcal_{\Omega_1 \times \Omega_2}$:
        \[ (\rho \otimes \kappa)(E) = (\kappa\otimes\rho)(E^s).\]
        For the case of $\mu$ and $\nu$ we first note the following:
       \begin{align*}
            X^{-1}(D^s_{ Y(\omega_1)})
        &= \lC \omega_2 \in \Omega_2\st ( X(\omega_2), Y(\omega_1)) \in D^s \rC\\
        &= \lC \omega_2 \in \Omega_2\st (\omega_2,\omega_1) \in ( X \times  Y)^{-1}(D^s) \rC\\
        &= ( X \times  Y)^{-1}(D^s)_{\omega_1}.
       \end{align*}
       Also note that:
       \begin{align*}
          ( X\times Y)^{-1}(D^s)
           &= \lC (\omega_2,\omega_1) \in \Omega_2 \times \Omega_1\st ( Y(\omega_1), X(\omega_2)) \in D \rC \\
           &= ( Y \times  X)^{-1}(D)^s.
       \end{align*}
        Then we get:
        \begin{align*}
            & (\mu \otimes \nu)(D^s)\\
            &= \int \mu(D_y^s)\, \nu(dy)\\
            &= \int ( X_*\kappa)(D^s_y)\, ( Y_*\rho)(dy)\\
            &= \int \kappa\lp X^{-1}(D^s_{ Y(\omega_1)})\rp\, \rho(d\omega_1)\\
            &= \int \kappa\lp( X \times  Y)^{-1}(D^s)_{\omega_1}\rp\, \rho(d\omega_1)\\
            &= (\kappa\otimes\rho)\lp( X \times  Y)^{-1}(D^s)\rp\\
            &= (\kappa\otimes\rho)\lp( Y \times  X)^{-1}(D)^s\rp\\
            &= (\rho\otimes\kappa)\lp( Y \times  X)^{-1}(D)\rp\\
&= \int \rho\lp( Y \times  X)^{-1}(D)_{\omega_2}\rp\, \kappa(d\omega_2)\\
            &= \int \rho\lp Y^{-1}(D_{ X(\omega_2)})\rp\, \kappa(d\omega_2)\\
            &= \int ( Y_*\rho)(D_x)\, ( X_*\kappa)(dx)\\
            &= \int \nu(D_x)\, \mu(dx)\\
            &= (\nu \otimes \mu)(D).
        \end{align*}
        This shows the claim for indicator functions $\I_D(y,x)$.\\
        The claim for $f \in [0,\infty]^{\Xcal \times \Ycal}$ follows from \cite{Kle20} Thm.\ 1.96 giving us the representation:
        \[ f(x,y) = \sum_{n \in \N} c_n\cdot\I_{D^s_{(n)}}(x,y),\]
        with $D_{(n)} \in \Bcal_{\Ycal \times \Xcal}$, $n \in \N$.
        The linearity of the integral implies the claim.

\end{Prf}

\subsection{Essential Uniqueness of Factorizations}

\begin{Prf}{thm:fact-ess-uniq}
     First consider for $D \in \Bcal_{\Xcal \times \Ycal \times \Zcal}$ the set:
       \[ N_D^>:=\lC (y,z) \in \Ycal \times \Zcal\st P_1(X \in D^{y,z}|Y=y,Z=z) > P_2(X \in D^{y,z}|Y=y,Z=z)\rC   \]
       Then define:
     \[ D^>:=D \cap \lp\Xcal \times N_D^>\rp  \in \Bcal_{\Xcal \times \Ycal \times \Zcal}. \]
     With this we then get:
       \begin{align*}
        0 &=  \lp P_1(X|Y,Z) \otimes Q(Y|Z)\rp(D_z^>|z) - \lp P_2(X|Y,Z) \otimes Q(Y|Z)\rp(D_z^>|z)\\
          &= \int\int \I_{N_D^>}(y,z)\cdot \I_D(x,y,z)\, \\
          &\qquad\qquad \lp  P_1(X \in dx|Y=y,Z=z) - P_2(X \in dx|Y=y,Z=z) \rp\, Q(Y \in dy|Z=z)\\
          &= \int \I_{N_D^>}(y,z)\cdot\underbrace{\lp  P_1(X \in D^{y,z}|Y=y,Z=z) - P_2(X \in D^{y,z}|Y=y,Z=z) \rp}_{> 0 \text{ for }(y,z) \in N_D^>}\\
          &\qquad\qquad\qquad\qquad\qquad\qquad\qquad\qquad\qquad\qquad\qquad Q(Y \in dy|Z=z).
       \end{align*}
       It follows from this that $Q(Y \in \lp N_D^>\rp^z|Z=z)=0$ for all $z \in \Zcal$.\\
       We get by a symmetric argument that $Q(Y \in \lp N_D^<\rp^z|Z=z)=0$ and thus $Q(Y \in N_D^z|Z=z)=0$ for all $z \in \Zcal$. This shows the first claim.\\
       If now $(\Xcal,\Xcal^\Omega)$ is countably separated then by Theorem \ref{thm:countably-separated} there exists a countable algebra $\Ecal \ins \Bcal_\Xcal$ that separates the points of $\Xcal$ and such that:
       \[ \Bcal_\Xcal = (\Ecal)_{\PrPf(\Xcal)}.  \]
       Then define the (countable) union:
        \[ N := \bigcup_{E \in \Ecal} N_{E \times \Ycal \times \Zcal} \in \Bcal_{\Ycal \times \Zcal}.\]
       Similar to above we then have that $Q(Y \in N^z|Z=z)=0$ for every $z \in \Zcal$ and
       also that for every $E \in \Ecal$ and $(y,z) \notin N$:
       \[ P_1(X \in E|Y=y,Z=z)= P_2(X \in E|Y=y,Z=z).\]
       Since $\Ecal$ is a $\pi$-system, $\Bcal_\Xcal = (\Ecal)_{\PrPf(\Xcal)}$ and $P_1(X|Y=y,Z=z),P_2(X|Y=y,Z=z) \in \PrPf(\Xcal)$ the above equality uniquely extends to all $A \in \Bcal_\Xcal$ and $(y,z) \notin N$:
          \[ P_1(X \in A|Y=y,Z=z)= P_2(X \in A|Y=y,Z=z).\]
        This shows the second claim.

\end{Prf}

\subsection{Kolmogorov Extension Theorems for Quasi-Universal Spaces}
\begin{Prf}{thm:kol-ext-mar-ker}
    For every fixed $z \in \Zcal$ we individually apply Kolmogorov's extension theorem,
    see \cite{Kol33,Fremlin} 454D-G and we get a unique probability measure:
    $Q_J(X_J|Z=z)$ on $\Bcal_J$ such that for all finite subsets $L \ins J$ we have:
    \[ \pr_{L,*} Q_J(X_J|Z=z) = Q_L(X_L|Z=z).\]
    Furthermore, $Q_J(X_J|Z=z)$ is perfect. \\
    We are left to show that the map $z \mapsto Q_J(X_J|Z=z)$ is measurable.
    For this consider the set:
    \[\Dcal:=\lC D \in \Bcal_J\st (z \mapsto Q_J(X_J \in D|Z=z)) \text{ is $\Bcal_\Zcal$-$\Bcal_{[0,1]}$-measurable} \rC.\]
    Then $\Ecal_J \ins \Dcal$. Indeed, for a finite subset $L \ins J$ and
    $A_L = \times_{l\in L}A_l$ with $A_l \in \Bcal_l$ we have by the consistency relation:
    \[ Q_J(X_J \in \pr_L^{-1}(A_L)|Z)= Q_L(X_L \in A_L|Z),\]
    which by assumption on $Q_L(X_L|Z)$ is measurable.\\
    It is also easily seen that $\Dcal$ is Dynkin system (aka $\lambda$-system),
    as it is closed under complements and disjoint unions and contains $\emptyset$.
    Since $\Ecal_J$ is a closed under intersections (aka a $\pi$-system), we get from
    Dynkin's $\pi$-$\lambda$-theorem, see \cite{Kle20} Theorem 1.19, that:
    \[ \Bcal_J = \sigma(\Ecal_J) \ins \Dcal.\]
    This shows that $Q_J(X_J|Z)$ is measurable and thus the claim.
\end{Prf}
\begin{Prf}{thm:kol-ext-univ}
        By Lemma \ref{def:univ-qus} we have that $\Bcal_n:=\Bcal(\Xcal_n^\Omega)$ makes
        $(\Xcal_n,\Bcal_n)$ a universal measurable spaces with $\Fcal(\Bcal_n)=\Xcal_n^\Omega$.
        Put $\Bcal_\Zcal:=\Bcal(\Zcal^\Omega)$ and note that this is universally complete.
        Since countable products of universal measurable spaces are universal,
        see Lemma \ref{lem:uqus-count-prod}, we have that:
        \[ (\Bcal_{0:n})_\Giry =\lp \bigotimes_{k=0}^n \Bcal_k\rp_\Giry = \Bcal\lp\prod_{k=0}^n \Xcal_k^\Omega \rp = \Bcal(\Xcal_{0:n}^\Omega).   \]
        Then we get the compatible measurable Markov kernels:
        \[Q_n(X_{0:n}|Z):\; (\Zcal,\Bcal_\Zcal) \to \PrPf(\Xcal_{0:n}) =
        \Giry(\Xcal_{0:n},\Bcal_{0:n}),\]
        where the latter identification comes from by Theorem \ref{thm:univ-qus-GPQKR}.
        Then note that every probability measure on a universal measurable space
        $(\Xcal_n,\Bcal_n)$ is perfect.
        Then by Theorem \ref{thm:kol-ext-mar-ker} and Remark \ref{rem:kol-ext}
        we get a unique measurable Markov kernel:
        \[ Q(X_\N|Z):\; (\Zcal,\Bcal_\Zcal) \to \Giry(\Xcal_\N,(\Bcal_\N)_\Giry), \]
        where $\Bcal_\N:=\bigotimes_{n \in \N} \Bcal_n$, such that for all $n \in \N$:
        \[ \pr_{0:n,*}Q(X_\N|Z)=Q_n(X_{0:n}|Z). \]
        We are left to show that $Q(X_\N|Z) \in \PrPf(\Xcal_\N)^\Zcal$.
        For this let $ Z \in \Zcal^\Omega \ins \Fcal(\Bcal_\Zcal)$.
        Then the composition:
        \[Q(X_\N|Z) \circ  Z  \]
        is measurable. We thus get:
        \[ Q(X_\N|Z) \circ  Z  \in \Giry(\Xcal_\N)^\Omega = \PrPf(\Xcal_\N)^\Omega,\]
        since $\Xcal_\N$ is a standard quasi-universal space, again by Theorem \ref{thm:univ-qus-GPQKR}.
        Since this holds for all $ Z \in \Zcal^\Omega$ we get:
        \[Q(X_\N|Z) \in \QUS(\Zcal,\PrPf(\Xcal_\N)) = \PrPf(\Xcal_\N)^\Zcal.\]
        This shows the claim.

\end{Prf}
\begin{Prf}{thm:kol-ext-cs}
        Since $(\Zcal,\Zcal^\Omega)$ is a standard quasi-universal space by Lemma \ref{def:univ-qus}
        there exists $ Z \in \Zcal^\Omega$ and $\iota \in \Omega^\Zcal$ such that
        $  Z \circ \iota = \id_\Zcal$.\\
        We then have that:
        \[ Q_{n}(X_{0:n}|W):=Q_{n}(X_{0:n}|Z) \circ  Z \in \PrPf(\Xcal_{0:n})^\Omega,\]
        which still satify the consistency relations for $n \in \N$.\\
        By the disintegration Theorem \ref{thm:disint-X} the following map is a surjective quotient map:
    \[ \otimes:\; \PrPf(\Xcal_{n+1})^{\Xcal_{0:n} \times \Omega} \times \PrPf(\Xcal_{0:n})^\Omega
        \srj \PrPf(\Xcal_{0:n+1})^{\Omega}.  \]
        So there exist $Q_{n+1}(X_{n+1}|X_{0:n},W) \in \PrPf(\Xcal_{n+1})^{\Xcal_{0:n} \times \Omega}$,
        such that:
        \[ Q_{n+1}(X_{n+1}|X_{0:n},W) \otimes Q_n(X_{0:n}|W) = Q_{n+1}(X_{0:n+1}|W).  \]
        By definition of $\PrPf(\Xcal_{n+1})$ and the isomorphism
        $\Omega^n \cong \Omega$
         we can inductively pick $ X_{n+1} \in \lp\Xcal_{n+1}^{\Omega_{0:n}}\rp^\Omega$ and
        $\kappa_{n+1}(W_{n+1}|W_{0:n},W) \in \PrQm(\Omega_{n+1})^{\Omega_{0:n} \times \Omega}$ such that:
        \begin{align*}  X_{n+1,*}\kappa_{n+1}(W_{n+1}|W_{0:n},W) &= Q_{n+1}(X_{n+1}|W_{0:n},W) \\
            &:= Q_{n+1}(X_{n+1}|X_{0:n},W) \circ ( X_{0:n} \times \id_\Omega)(W_{0:n},W).
        \end{align*}
        By Kolmogorov's extension theorem for standard quasi-universal spaces,
        see Theorem \ref{thm:kol-ext-univ}, there now exists
        a unique $\kappa(W_\N|W) \in \PrQm(\Omega^\N)^\Omega \cong \PrQm(\Omega)^\Omega$ such that for all $n \in N$:
        \[ \pr_{0:n,*}\kappa(W_\N|W) = \bigotimes_{j=n}^0 \kappa_j(W_j|W_{0:j-1},W).  \]
        We then put:
        \[  X := ( X_n)_{n \in \N} \in \lp\Xcal_\N^{\Omega^\N}\rp^\Omega \cong \lp\Xcal^\Omega\rp^\Omega.\]
        Then:
            \[ Q(X_\N|W) : =  X_* \kappa(W_\N|W) \in \PrPf(\Xcal_\N)^\Omega,\]
        and thus:
        \[ Q(X_\N|Z) : = Q(X_\N|W =\iota(Z)):= Q(X_\N|W) \circ \iota \in \PrPf(\Xcal_\N)^\Zcal.\]
        By construction we then have that for all $n \in \N$:
        \[ \pr_{0:n,*}Q(X_\N|Z) = Q_n(X_{0:n}|Z). \]
        We are now left to show the uniqueness.
        Since the spaces $\Xcal_n$ are countably separated there are countably generated
        $\sigma$-algebras $\Ecal_n \ins \Bcal_n:=\Bcal(\Xcal_n^\Omega)$ that separate the points of
        $\Xcal_n$. Then $\Ecal_\N := \bigotimes_{n \in \N} \Ecal_n$
        is countably generated and separates the points of $\Xcal_\N := \prod_{n \in \N} \Xcal_n$.
        Then $Q(X_\N|Z)$, when restricted to $\Ecal_\N$, is the unique extension
        of the $Q_n(X_{0:n}|Z)$, when restricted to $\bigotimes_{j=0}^n\Ecal_j$.
        Since by Theorem \ref{thm:countably-separated} we have:
        \[\Bcal(\Xcal_\N^\Omega) = \lp\Ecal_\N \rp_{\PrPf(\Xcal_\N)},\]
        the Markov kernel $Q(X_\N|Z)$ uniquely extends from $\Ecal_\N$ to $\Bcal(\Xcal_\N^\Omega)$.\\
        This shows the claim.

\end{Prf}

\subsection{Deterministic Markov Kernels between Quasi-Universal Spaces}
\begin{Prf}{def:qms-det-markov-ker}
     1. $\implies$ 2.: Let $D \in \Bcal_{\Ycal \times \Ycal}$ and $z \in \Zcal$. Then:
        \begin{align*}
           & \lp \delta_g(Y_1|Z=z)\otimes\delta_g(Y_2|Z=z)\rp(D) \\
           & =   \int \delta_{g(z)}(D_{y_2})\,\delta_{g(z)}(dy_2) \\
           & =   \int \I_D(g(z),y_2) \,\delta_{g(z)}(dy_2) \\
           & =   \I_D(g(z),g(z)) \\
&= \delta_g((Y_1,Y_2) \in D|Z=z).
        \end{align*}
        2. $\implies$ 3.: We have by assumption for every $B \in \Bcal_\Ycal$:
        \[ K(Y \in B|Z=z) \cdot K(Y \in B|Z=z) =  K(Y \in B,Y \in B|Z=z) = K(Y \in B|Z=z).
 \]
        Solving for $K(Y \in B|Z=z)$ shows:
        \[K(Y \in B|Z=z) \in \lC 0,1\rC.\]
        This shows the claims.

\end{Prf}
\begin{Lem}
    \label{lem:ms-0-1-det-fct-det}
    Let $K(Y|Z):\; (\Zcal,\Bcal_\Zcal) \to \Giry(\Ycal,\Bcal_\Ycal)$ be a measurable and 0-1-deterministic Markov kernel.
    \begin{enumerate}
        \item If $(\Ycal,\Bcal_\Ycal)$ is countably generated then there exists a measurable
            $g:\; (\Zcal,\Bcal_\Zcal) \to (\Ycal,\Bcal_\Ycal)$ such that:
            $ K(Y|Z) = \delta_g(Y|Z)$.
        \item If $(\Ycal,\Bcal_\Ycal)$ is universally countably generated then there exists a universally measurable $g:\; (\Zcal,(\Bcal_\Zcal)_\Giry) \to (\Ycal,(\Bcal_\Ycal)_\Giry)$ such that:
            $ K(Y|Z) = \delta_g(Y|Z)$.
    \end{enumerate}
    Furthermore, if $(\Ycal,\Bcal_\Ycal)$ is separated then the map $g$ from above is unique.
    \begin{proof}
        1.) By assumption we have a countable algebra $\Ecal \ins \Bcal_\Ycal$
        such that $\sigma(\Ecal)=\Bcal_\Ycal$. Then for every $z \in \Zcal$ we pick an arbitrary point
        (using the axiom of choice):
        \[ g(z) \in C_z:=\bigcap_{\substack{B \in \Ecal\\K(Y\in B|Z=z)=1.}} B \in \Bcal_\Ycal.\]
        Note that since $\Ecal$ is countable the intersection $C_z$ lies in $\Bcal_\Ycal$.
        Furthermore, $C_z\neq \emptyset$. Otherwise we would get the contradiction:
        \[ 1 = K(Y \in C_z^\cmpl|Z=z) \le \sum_{\substack{B \in \Ecal\\K(Y\in B|Z=z)=1.}} K(Y \in B^\cmpl|Z=z) = 0.   \]
        To check that $g$ is measurable we only need to check this for $B \in \Ecal$. We then get
        \begin{align*}
        g^{-1}(B) &= \lC z\in \Zcal\st g(z) \in B \rC\\
                  &= \lC z\in \Zcal\st C_z \ins B \rC\\
                  &= \lC z\in \Zcal\st K(Y \in B|Z=z)=1 \rC\\
                  &= K(Y \in B|Z)^{-1}(1) \\
                  & \in \Bcal_\Zcal,
        \end{align*}
        since the Markov kernel is a measurable map for each $B$. This shows the measurability of $g$.
        Similarly, for $B \in \Ecal$ we have:
        \[ \delta_g(Y \in B|Z=z) = \I_B(g(z)) = K(Y \in B|Z=z).  \]
        Since $\Ecal$ is a generating $\pi$-system the equality also holds for all $ B \in \sigma(\Ecal)=\Bcal_\Ycal$.
        This shows the first claim.\\
        2.) When $K(Y|Z)$ is extended to $(\Bcal_\Ycal)_\Giry$ then it is still $(\Bcal_\Zcal)_\Giry$-measurable by Proposition \ref{prp:prob-prob-meas}. For fixed $z \in \Zcal$ and $C \in  (\Bcal_\Zcal)_\Giry$ we can write:
        \[ C = B \triangle N, \]
        with $B \in \Bcal_\Ycal$ and $N$ a $K(Y|Z=z)$-null set. So $K(Y \in C|Z=z) \in \lC0,1 \rC$
        for all $C \in (\Bcal_\Ycal)_\Giry$ and $z \in \Zcal$ as well.\\
        By assumption there now exists a countably generated $\sigma$-algebra
        $\Ecal \ins (\Bcal_\Ycal)_\Giry$ such that $(\Ecal)_\Giry=(\Bcal_\Ycal)_\Giry$.
        By the same arguments as above we get a
        $(\Bcal_\Zcal)_\Giry$-$\Ecal$-measurable map $g$ such that:
        \[ \delta_g(Y|Z) = K(Y|Z),\]
        when restricted to $\Ecal$.
        This equality then uniquely extends to $(\Ecal)_\Giry=(\Bcal_\Ycal)_\Giry$.
        Note that $g$ is also $(\Bcal_\Zcal)_\Giry$-$(\Ecal)_\Giry$-measurable,
        which is the same as being
        $(\Bcal_\Zcal)_\Giry$-$(\Bcal_\Ycal)_\Giry$-measurable.
        This shows the second claim.\\
        For the uniqueness let $g_1,g_2$ be two different (universally) measurable maps
        such that:
        \[ \delta_{g_1}(Y|Z) = \delta_{g_2}(Y|Z).\]
        Then there exists a point $z \in \Zcal$ such that $g_1(z) \neq g_2(z)$.
        Since $(\Ycal,\Bcal_\Ycal)$ is separated there exists $B \in \Bcal_\Ycal$ such that
        $g_1(z) \in B$ and $g_2(z) \notin B$, which  implies:
        \[ \delta_{g_1}(Y \in B|Z=z) = 1 \neq 0 = \delta_{g_2}(Y \in B|Z=z),\]
        in contradiction to the assumption. So $g_1=g_2$.
    \end{proof}
\end{Lem}\begin{Prf}{thm:qms-0-1-det-fct-det}
      Two directions follow from Lemma \ref{def:qms-det-markov-ker}.
        For the remaining claim let $K(Y|Z)$ be 0-1-deterministic.
        We need to be careful with handling the $\sigma$-algebras in the following.
        Let $\Bcal_\Zcal:=\Bcal(\Zcal^\Omega)$ and $\Bcal_\Ycal:=\Bcal(\Ycal^\Omega)$.
        Since $(\Ycal,\Ycal^\Omega)$ is countably separated we have a countably generated $\sigma$-algebra $\Ecal_\Ycal \ins \Bcal_\Ycal$ that separates the points of $\Ycal$ such that $\Bcal_\Ycal=(\Ecal_\Ycal)_{\PrPf(\Ycal)}$, see Theorem \ref{thm:countably-separated}.
        Because of $\Bcal_\Ycal=(\Ecal_\Ycal)_{\PrPf(\Ycal)}$ we know that the following composition of
        (measurable) inclusion and restriction maps is injective:
        \[\res:\; (\PrPf(\Ycal),\Bcal(\PrPf(\Ycal)^\Omega)) \xrightarrow{\incl} (\Giry(\Ycal),\Bcal(\Giry(\Ycal)^\Omega))  \xrightarrow{\res}
            (\Giry(\Ycal,\Ecal_\Ycal),\Bcal_{\Giry(\Ycal,\Ecal_\Ycal)}), \]
            where $\Bcal_{\Giry(\Ycal,\Ecal_\Ycal)}$ is the $\sigma$-algebra generated by the evaluation maps $\ev_B$ for $B \in \Ecal_\Ycal$,
            which is countably generated (since $\Ecal_\Ycal$ is) and separates the points of $\Giry(\Ycal,\Ecal_\Ycal)$.  The $\sigma$-algebra $\Bcal_{\PrPf(\Ycal)} := \res^*\Bcal_{\Giry(\Ycal,\Ecal_\Ycal)} \ins \Bcal(\PrPf(\Ycal)^\Omega)$ is then also countably generated and separates the points of $\PrPf(\Ycal)$, since $\res$ is injective. Again by Theorem \ref{thm:countably-separated} we then have that
            $\Bcal(\PrPf(\Ycal)^\Omega) = \lp\Bcal_{\PrPf(\Ycal)}\rp_{\PrPf(\PrPf(\Ycal))}$.
            Also note that: $(\Bcal_\Zcal)_{\PrPf(\Zcal)} = \Bcal_\Zcal$.\\
        We then have the composition of measurable maps:
        \[(\Zcal,\Bcal_\Zcal) \xrightarrow{K(Y|Z)} (\PrPf(\Ycal),\Bcal_{\PrPf(\Ycal)})
            \xhookrightarrow{\res} (\Giry(\Ycal,\Ecal_\Ycal),\Bcal_{\Giry(\Ycal,\Ecal_\Ycal)}). \]
         Since $K(Y|Z)$ is 0-1-deterministic and $\Ecal_\Ycal$ countably generated by Lemma \ref{lem:ms-0-1-det-fct-det} there exists a (unique)
        $\Bcal_\Zcal$-$\Ecal_\Ycal$-measurable map $g$ such that for $z\in \Zcal$:
        \[ K(Y|Z=z) = \delta_g(Y|Z=z) \quad \in \Giry(\Ycal,\Ecal_\Ycal),\]
        suppressing mentioning $\res$ for brevity. So we get the following commutative diagram of measurable maps:
        \[\xymatrix{
                (\Zcal,\Bcal_\Zcal) \ar_-{g}[d] \ar^-{K(Y|Z)}[rr] && (\PrPf(\Ycal),\Bcal_{\PrPf(\Ycal)}) \ar@{^(->}^-{\res}[d] \\
                (\Ycal,\Ecal_\Ycal) \ar@{^(->}_-{\delta}[rr] && (\Giry(\Ycal,\Ecal_\Ycal), \Bcal_{\Giry(\Ycal,\Ecal_\Ycal)}).
        }\]
        We now want to show that $g$ is even $\Bcal_\Zcal$-$\Bcal_\Ycal$-measurable.
        Since $\Bcal_\Ycal=(\Ecal_\Ycal)_{\PrPf(\Ycal)}$  we want to  complete $\Ecal_\Ycal$ w.r.t.\ $\PrPf(\Ycal)$ and show that $g$ stays measurable even after completion. By Lemma \ref{qms:lem:map-completetions} and $(\Bcal_\Zcal)_{\PrPf(\Zcal)} = \Bcal_\Zcal$ it is enough to show that $g_*(\PrPf(\Zcal)) \ins \PrPf(\Ycal)$. We start by noting that the equality $K(Y|Z) = \delta_g(Y|Z)$ induces the identity:
        \[ \delta_* \circ g_* = \res_* \circ K(Y|X)_*:\; \Giry(\Zcal,\Bcal_\Zcal) \to
        \Giry\lp\Giry(\Ycal,\Ecal_\Ycal),\Bcal_{\Giry(\Ycal,\Ecal_\Ycal)}\rp.\]
        We can apply the monad action:
         \[\Mbb:\;\Giry\lp\Giry(\Ycal,\Ecal_\Ycal),\Bcal_{\Giry(\Ycal,\Ecal_\Ycal)}\rp \to \Giry(\Ycal,\Ecal_\Ycal), \]
         and use the monad rule that $\Mbb \circ \delta_* = \id$ to get the equality (suppressing $\res_*$ for brevity):
         \[  g_* = \Mbb \circ  K(Y|X)_* :\; \Giry(\Zcal,\Bcal_\Zcal) \to \Giry(\Ycal,\Ecal_\Ycal).\]
         Since the map:
         \[ K(Y|Z):\; (\Zcal,\Zcal^\Omega) \to (\PrPf(\Ycal),\PrPf(\Ycal)^\Omega) \]
         is quasi-measurable, we get that $K(Y|Z)_*\lp\PrPf(\Zcal)\rp \ins \PrPf(\PrPf(\Ycal))$. Since $(\PrPf,\delta,\Mbb)$ is a monad we also get that $\Mbb:\; \PrPf(\PrPf(\Ycal)) \to \PrPf(\Ycal)$ is a well-defined map, showing that:
         $\Mbb\lp  \PrPf(\PrPf(\Ycal))\rp \ins \PrPf(\Ycal)$. Note that $\Mbb$ for $\Giry$ and $\Mbb$ for $\PrPf$  commute with $\res$, i.e.\ if $\pi \in \PrPf(\PrPf(\Ycal))$ and $B \in \Ecal_\Ycal$ we get:
         \begin{align*}
             \Mbb_\Giry(\res_*\pi)(B) &= \int_{\Giry(\Ycal,\Ecal_\Ycal)} \ev_B(\mu)\,(\res_*\pi)(d\mu) \\
             &= \int_{\PrPf(\Ycal)} \ev_B(\res(\nu))\,\pi(d\nu) \\
             &= \int_{\PrPf(\Ycal)} \ev_B(\nu)\,\pi(d\nu)\\
             &= \Mbb_\PrPf(\pi)(B).
     \end{align*}
    Without making further distinctions between the $\Mbb$'s we then arrive at:
    \[g_*\lp\PrPf(\Zcal)\rp = \Mbb\circ K(Y|Z)_*\lp\PrPf(\Zcal)\rp \ins \Mbb\lp  \PrPf(\PrPf(\Ycal))\rp \ins \PrPf(\Ycal). \]
    Since $g$ is $\Bcal_\Zcal$-$\Ecal_\Ycal$-measurable it is then, because of $g_*\lp\PrPf(\Zcal)\rp \ins \PrPf(\Ycal)$, by Lemma \ref{qms:lem:map-completetions}, also
    $(\Bcal_\Zcal)_{\PrPf(\Zcal)}$-$(\Ecal_\Ycal)_{\PrPf(\Ycal)}$-measurable.
    Since $(\Bcal_\Zcal)_{\PrPf(\Zcal)}=\Bcal_\Zcal$ and $(\Ecal_\Ycal)_{\PrPf(\Ycal)}=\Bcal_\Ycal$ we
    get that $g$ is $\Bcal_\Zcal$-$\Bcal_\Ycal$-measurable.\\
        Now let $ Z \in \Zcal^\Omega \ins \Fcal(\Bcal_\Zcal)$. Then we
        have $g \circ  Z \in \Fcal(\Bcal_\Ycal)=\Ycal^\Omega$, where the last equality holds by assumption.
        This shows that $g \in \Ycal^\Zcal$ and thus: $\delta_g(Y|Z) \in \PrPf(\Ycal)^\Zcal$.
        Both $K(Y|Z)$ and $\delta_g(Y|Z)$ are in
        $\PrPf(\Ycal)^\Zcal$ and we already have the equality: $K(Y|Z) = \delta_g(Y|Z)$ on $\Ecal_\Ycal$.  Since $(\Ecal_\Ycal)_{\PrPf(\Ycal)}=\Bcal_\Ycal$ the above equality
        uniquely extends to $\Bcal_\Ycal$ and we get:
        \[ K(Y|Z) = \delta_g(Y|Z) \quad \in \PrPf(\Ycal)^\Zcal.\]
        This shows that $K(Y|Z)$ is function-deterministic.
         \[\xymatrix{
                 \PrPf(\Zcal) \ar@{-->}_-{g_*?}[dddd]\ar@{^(->}^-{\incl}[drr] \ar^-{K(Y|Z)_*}[rrrrrr] && && &&
                 \PrPf(\PrPf(\Ycal)) \ar^-{\Mbb}[dddd]\ar@{_(->}_-{\res}[dll]\\
              &&  \Giry(\Zcal,\Bcal_\Zcal) \ar_-{g_*}[d] \ar^-{K(Y|Z)_*}[rr] &&
              \Giry(\PrPf(\Ycal),\Bcal_{\PrPf(\Ycal)}) \ar@{^(->}^-{\res_*}[d] \\
                    &&  \Giry(\Ycal,\Ecal_\Ycal) \ar_-{\id}[rrd] \ar@{^(->}_-{\delta_*}[rr] && \Giry(\Giry(\Ycal,\Ecal_\Ycal), \Bcal_{\Giry(\Ycal,\Ecal_\Ycal)}) \ar^-{\Mbb}[d]\\
              && && \Giry(\Ycal,\Ecal_\Ycal) \\
              \PrPf(\Ycal) \ar_-{\id}[rrrrrr] \ar@{^(->}^-{\res}[rruu]&& && && \PrPf(\Ycal)\ar@{_(->}_-{\res}[llu].
        }\]

\end{Prf}

\subsection{The Markov Category of Markov Kernels on $\StdQUS$}
\begin{Prf}{thm:Markov-category}
        We apply Prop.\ 3.1.\ from \cite{FriG20}.\\
        Strong monad $\PrPf$: See Theorem \ref{thm:K-P-R-strong-monad}. Note also that $\PrPF=\PrPf=\PrpF=\PrQm=\Giry=\Prpf$ for standard quasi-universal spaces by Theorem \ref{thm:univ-qus-GPQKR} and that $\PrPf(\Xcal)$ stays inside $\StdQUS$.\\
        Symmetry/commutativity: See Fubini Theorem \ref{thm:fubini}.\\
        Conditionals: See Disintegration Theorem \ref{thm:disint-X}.\\
        Kolmogorov powers: See Kolmogorov's Extension Theorem \ref{thm:kol-ext-univ}.\\
        Deterministic morphisms/Markov kernels: See Theorem \ref{thm:qms-0-1-det-fct-det}.\\
        A.s.-compatible representability: See Remark \ref{rem:a.s.-det} and essential uniqueness of factorizations Theorem \ref{thm:fact-ess-uniq}.\\
        All other aspects are easy to see or follow from general considerations.

\end{Prf}
\begin{Prf}{cor:de-finetti}
         This directly follows from the synthetic version of de Finetti's Theorem in
         categorical probability theory, see \cite{Fri21},
         together with Theorem \ref{thm:Markov-category}.
         Note that for our case we do not need to require that $(\Zcal,\Zcal^\Omega)$ is
         a \emph{standard} quasi-universal space as all required results work without that assumption.

    We spell out why $(\Zcal,\Zcal^\Omega)$ may be arbitrary, since \Cref{thm:Markov-category}
    itself is stated for $\StdQUS$. The synthetic proof of de Finetti's theorem in
    \cite{Fri21} uses exactly three properties of the ambient Markov category: it has
    \emph{conditionals}, it is \emph{a.s.-compatibly representable}, and it has \emph{countable
    Kolmogorov powers}. In the argument these are applied to the objects being mixed over --- here
    $\Xcal$, the countable power $\Xcal_\N$, and the space of measures $\PrPf(\Xcal)$, all of which
    are standard when $\Xcal$ is --- whereas $\Zcal$ occurs only as the source of the Kleisli
    morphism $Q(X_\N|Z)$, i.e.\ as a parameter carried along. Each of the three steps is available
    in $\QUS$ with $\Zcal$ arbitrary:
    \begin{enumerate}
        \item \emph{Conditionals.} The factorisation of the joint kernel into
            $Q(X|S) \otimes P(S|Z)$ exists because the product of Markov kernels
            \[ \otimes:\; \PrPf(\Xcal)^{\Ycal \times \Zcal} \times \PrPf(\Ycal)^\Zcal \to
               \PrPf(\Xcal \times \Ycal)^\Zcal \]
            is a surjective quotient map by \Cref{thm:disint-X}. That theorem offers three
            alternative hypotheses, and we use the \emph{first} one --- $(\Ycal,\Ycal^\Omega)$
            countably separated and $(\Xcal,\Xcal^\Omega)$ standard --- which leaves
            $(\Zcal,\Zcal^\Omega)$ unconstrained and does give the quotient property, see
            \footref{fn:disint-quotient}.
        \item \emph{Almost-sure uniqueness.} That the factorisation is unique up to a null set,
            which is what makes the resulting $S$ well defined, is \Cref{thm:fact-ess-uniq}, stated
            for arbitrary quasi-universal spaces $\Xcal$, $\Ycal$, $\Zcal$.
        \item \emph{Kolmogorov powers.} Passing from the consistent family of finite-dimensional
            kernels $Q(X_{0:n}|Z)$ to $Q(X_\N|Z)$ is \Cref{thm:kol-ext-univ}, which asks the
            $\Xcal_n$ to be standard and allows $(\Zcal,\Zcal^\Omega)$ to be \emph{any}
            quasi-universal space.
    \end{enumerate}
    With these three substitutions the argument of \cite{Fri21} goes through verbatim, and yields
    the stated equivalence for arbitrary $(\Zcal,\Zcal^\Omega)$. We have given the reduction rather
    than repeating the synthetic argument, which is unchanged.
\end{Prf}

\begin{Prf}[The universal Hilbert cube is a well-behaved sample space]{qms:thm:act4-well-behaved}
    We apply \Cref{qms:thm:spl-spc-cmpl} to $(\Omega,\Bcal_\otimes)$ and $\PrPf=\Giry(\Omega,\Bcal_\otimes)$.

    The ``zip-locker'' map, Hilbert's hotel:
    \begin{align}
        \Theta = \Theta_0 \times \Theta_1:\; \Omega \to \Omega \times \Omega, \qquad
        (e_n)_{n\in\N} \mapsto \lp (e_{2k})_{k \in \N},(e_{2m+1})_{m \in \N}\rp,
    \end{align}
    is a bijection which is $\Bcal_\otimes$-$(\Bcal_\otimes\otimes\Bcal_\otimes)$-bimeasurable, as one checks
    on finite cylinder sets. So it is an isomorphism of measurable spaces:
    \begin{align}
        \Theta:\; (\Omega,\Bcal_\otimes) \cong (\Omega \times \Omega, \Bcal_\otimes \otimes \Bcal_\otimes),
    \end{align}
    which is the hypothesis of \Cref{qms:thm:spl-spc-cmpl}.
    The two conditions \Cref{qms:eq:PrPf-compl-cond} and \Cref{qms:eq:PrPf-prod-cond} on $\PrPf$ hold
    because $\PrPf=\Giry(\Omega,\Bcal_\otimes)$ is \emph{all} Borel probability measures: restrictions of
    completions, push-forwards along $\Bcal_\otimes$-measurable maps and products of Borel probability
    measures are again Borel probability measures.
    \Cref{qms:thm:spl-spc-cmpl} therefore gives all three acts, the Fubini property, the product isomorphism
    property, the $\PrPf(\Omega)$-completeness of $\Bcal_\Omega$, and the two displayed equalities.

    For the remaining points: the Hilbert cube is a Polish space, so $\Bcal_\otimes$ is countably generated
    and separates points --- the cylinder sets $\pr_n^{-1}([0,q])$ with $n \in \N$, $q \in \Q \cap [0,1]$,
    are a countable separating family in $\Bcal_\otimes \ins \Bcal_\Omega$. Hence $(\Omega,\Bcal_\Omega)$ is
    countably separated. Every Borel probability measure on a Polish space is Radon, hence perfect, see
    \cite{Fremlin} 451, so every $P \in \PrPf$ is perfect. The conditions listed in
    \Cref{qms:rem:count-sep-perf-cond} are exactly these three, so they all hold, and
    \Cref{qms:thm:count-sep-perf} and \Cref{qms:cor:count-sep-perf} apply.
\end{Prf}

\section{Appendix - Causal Models}
\label{qms:sec:app-causal}

\begin{Lem}
    \label{lem:indep-det}
    Consider transitional random variables  $\Xk \in \PrPf(\Xcal)^{\Wcal\times\Tcal}$ and $\Yk \in \PrPf(\Ycal)^{\Wcal\times\Tcal}$ and $\Zk \in \PrPf(\Zcal)^{\Wcal\times\Tcal}$. Assume that:
    \[\Delta_\Xcal:=\lC (x,x) \in \Xcal \times \Xcal\st x \in \Xcal \rC \in\Bcal_{\Xcal \times \Xcal}.\footnotemark\]
Then we have the implication:\footnotetext{If $(\Xcal,\Xcal^\Omega)$ is countably separated then we have: $\Delta_\Xcal \in \Bcal_\Xcal \otimes \Bcal_\Xcal \ins \Bcal_{\Xcal \times \Xcal}.$ \label{fn:diagonal}}

    \[ \Xk \ismapof_\Kk \Zk \qquad \implies \qquad \Xk\Indep_\Kk \Yk \given \Zk.\]
\begin{proof}
Assume $\Xk \ismapof_\Kk \Zk$. Then there exists a $\Phi \in \Xcal^\Zcal$ such that:
     \[ \Kk(X,Z|T) = \deltabf_\Phi(X|Z) \otimes \Kk(Z|T).\]
     Since $\Delta_\Xcal \in \Bcal_{\Xcal \times \Xcal}$ we know that:
     \[M:= \lC(x,z) \in \Xcal \times \Zcal\st x \neq \Phi(z)\rC  = (\id_\Xcal \times \Phi)^{-1}(\Delta_\Xcal^\cmpl) \in \Bcal_{\Xcal \times \Zcal}.\]
     Then $M$ is a $\Kk(X,Z|T)$-null set. Indeed:
     \begin{align*}
         \Kk\lp(X,Z) \in M|T\rp
         &= \int \deltabf_\Phi(X \in M_z|Z=z)\, \Kk(Z \in dz|T)\\
         &=\int \I_{M_z}(\Phi(z))\, \Kk(Z \in dz|T)\\
         &= \int 0 \,\Kk(Z \in dz|T)\\
         &=0.
     \end{align*}
     Now consider the Markov kernel $\Kk(X_1,X_2,Z,Y|T) \in \PrPf\lp\Xcal \times \Xcal \times \Zcal \times \Ycal\rp^\Tcal$ with $X_1:=X$ and $X_2:=\Phi(Z)$ and the set:
     \[N:= \lC (x_1,x_2,z,y) \in \Xcal \times \Xcal \times \Zcal \times \Ycal\st x_1 \neq x_2\rC = \Delta_\Xcal^\cmpl \times \Zcal \times \Ycal \in \Bcal_{\Xcal \times \Xcal \times \Zcal \times \Ycal}. \]
     Then we have:
     \[\Kk\lp (X_1,X_2,Z,Y) \in N|T\rp = \Kk\lp(X,Z,Y) \in M \times \Ycal|T\rp =0. \]
     So for every $D \in \Bcal_{\Xcal \times \Xcal \times \Zcal \times \Ycal}$ we get:
     \begin{align*}
         \Kk\lp (X_1,X_2,Z,Y) \in D |T\rp &= \Kk\lp (X_1,X_2,Z,Y) \in D \cap N^\cmpl |T\rp \\
         &= \Kk\lp (X_2,X_1,Z,Y) \in D \cap N^\cmpl |T\rp \\
         &= \Kk\lp (X_2,X_1,Z,Y) \in D |T\rp.
     \end{align*}
     For $D=\Xcal \times E$ with $E \in \Bcal_{\Xcal \times \Zcal \times \Ycal}$ this gives:
     \begin{align*} \lp\deltabf_\Phi(X|Z) \otimes \Kk(Y,Z|T) \rp(E) &= \Kk\lp (X_2,Z,Y) \in E |T\rp\\
     &= \Kk\lp (X_1,Z,Y) \in E |T\rp  \\
     &= \Kk\lp (X,Z,Y) \in E |T\rp.
     \end{align*}
     So this implies that we have:
     \[ \Kk(X,Y,Z|T) = \deltabf_\Phi(X|Z) \otimes \Kk(Y,Z|T),\]
     for every transitional random variable $\Yk$.
     If we pick $\Qk(X|Z):=\deltabf_\Phi(X|Z)$ we showed:
     \[ \Xk\Indep_\Kk \Yk \given \Zk,\]
     which is the claim.
 \end{proof}
\end{Lem}\begin{Prf}{thm:separoid_axioms-tci}
    The constructions follow the same lines as in \cite{For21} Theorem 3.11 and Appendix E.
    Note that for Left Weak Union we need the disintegration $\Qk(X,U|Z)=\Qk(X|U,Z) \otimes \Qk(U|Z)$.
    That is the reason we need to make the extra assumptions on the underlying spaces, see disintegration Theorem \ref{thm:disint-X}. Similarly, for $\Tk$-Restricted Right Redundancy we need the disintegration
    $\Kk(X,Z|T) = \Kk(X|Z,T) \otimes \Kk(Z|T)$.
    The corresponding conditions on the underlying spaces could be weakened if one could find weaker conditions under which disintegrations exist, see again Theorem \ref{thm:disint-X}.
\end{Prf}

\begin{Prf}{lem:absorb-input}
    Write $u:\,\Tcal \to \Ucal$ for the quasi-measurable map with $\Uk(U|W,T)=\deltabf_{u(T)}$.
    Since $\Uk$ does not depend on $W$ and is a Dirac measure, every joint Markov kernel in which
    $\Uk$ occurs factorizes off a Dirac component: for all transitional random variables
    $\Vk_1,\dots,\Vk_n$ we have
    \begin{align}
        \Kk(V_1,\dots,V_n,U|T) &= \Kk(V_1,\dots,V_n|T) \otimes \deltabf_{u(T)}
        \label{eq:absorb-dirac}
    \end{align}
    in $\PrPf(\Vcal_1 \times \cdots \times \Vcal_n \times \Ucal)^\Tcal$, directly from
    \Cref{def:tr-cond-ind} and the definition of the product of Markov kernels,
    \Cref{qms:thm:product-PrPf}.

    1.: Assume $\Xk \otimes \Uk \Indep \Yk \given \Zk \otimes \Tk$, i.e.\ there is a Markov kernel
    $\Qk(X,U|Z,T) \in \PrPf(\Xcal \times \Ucal)^{\Zcal \times \Tcal}$ with
    \begin{align}
        \Kk(X,U,Y,Z|T) &= \Qk(X,U|Z,T) \otimes \Kk(Y,Z|T).
    \end{align}
    Let $\Qk(X|Z,T)$ be its marginal in the $\Xcal$-component, which requires no disintegration.
    Taking the marginal in the $\Ucal$-component on both sides gives
    \begin{align}
        \Kk(X,Y,Z|T) &= \Qk(X|Z,T) \otimes \Kk(Y,Z|T).
        \label{eq:absorb-marginal}
    \end{align}
    Now define
    \begin{align}
        \Qk'(X|U,Z,T) &:= \Qk(X|Z,T),
    \end{align}
    which ignores its $\Ucal$-argument and is therefore quasi-measurable as the composition of
    $\Qk(X|Z,T)$ with the projection $\Ucal \times \Zcal \times \Tcal \to \Zcal \times \Tcal$.
    Using \Cref{eq:absorb-dirac} twice and then \Cref{eq:absorb-marginal} we get
    \begin{align}
        \Kk(X,Y,U,Z|T)
            &= \Kk(X,Y,Z|T) \otimes \deltabf_{u(T)}\\
            &= \Qk(X|Z,T) \otimes \Kk(Y,Z|T) \otimes \deltabf_{u(T)}\\
            &= \Qk'(X|U,Z,T) \otimes \Kk(Y,U,Z|T),
    \end{align}
    where the last step uses that $\Qk'$ does not depend on its $\Ucal$-argument, so that
    re-inserting the Dirac component $\deltabf_{u(T)}$ into the conditioning slot changes nothing.
    This is $\Xk \Indep \Yk \given \Uk \otimes \Zk \otimes \Tk$.
    Note that no disintegration was used: the only operations were a marginalisation and a
    re-indexing, and both are available for arbitrary $\Ucal$.\\

    2.: The first implication is the special case of 1.\ in which the left-hand side does not
    contain $\Uk$ at all; the same two lines apply verbatim.
    For $\Tk$-Restricted Right Redundancy, the statement
    $\Xk \Indep \deltabf_0 \given \Zk \otimes \Uk \otimes \Tk$ asks for a Markov kernel
    $\Qk(X|Z,U,T)$ with $\Kk(X,Z,U|T) = \Qk(X|Z,U,T) \otimes \Kk(Z,U|T)$. By
    \Cref{eq:absorb-dirac} both sides carry the Dirac factor $\deltabf_{u(T)}$, so the equation is
    equivalent to $\Kk(X,Z|T)=\Qk(X|Z,T) \otimes \Kk(Z|T)$ with
    $\Qk(X|Z,U,T):=\Qk(X|Z,T)$, which is exactly the disintegration required for
    $\Xk \Indep \deltabf_0 \given \Zk \otimes \Tk$. So the side condition of
    \Cref{thm:separoid_axioms-tci} b) needs to be checked for $\Zcal$ only.
\end{Prf}

\begin{Prf}{thm:gmp-mI-CBN}
    The graphical induction is the one of \cite{For21} Theorem 6.3 and Appendix J, carried out
    with the \ids-separation of \Cref{def:d-blocked} and the transitional conditional independence
    of \Cref{def:tr-cond-ind} in place of their counterparts in $\Meas$. It uses no property of
    $\Indep$ beyond the (asymmetric) separoid rules of \Cref{thm:separoid_axioms-tci}, so what has
    to be checked here is that every instance of those rules that the induction produces is
    licensed by its side conditions. There are exactly three rules that carry one.

    \emph{Extended Left Redundancy} a) is not needed at all. The only instances the induction
    produces are of the form $X_A \Indep X_B \given X_A$, and these are verified directly by
    taking $\Qk(X_A|X_A):=\deltabf(X_A|X_A)$, which is the unit of the monad $\PrPf$ and is
    quasi-measurable for every quasi-universal space. In particular the hypothesis
    $\Delta_{\Xcal_A} \in \Bcal_{\Xcal_A \times \Xcal_A}$ of that rule is never invoked.

    \emph{$\Tk$-Restricted Right Redundancy} b) and \emph{Left Weak Union} f) are the two rules
    that appeal to the disintegration \Cref{thm:disint-X}, and their side condition is that the
    relevant triple be a disintegration triple, which constrains the
    space in the conditioning slot. In the induction those slots are occupied by the spaces
    $\Xcal_D$ for various $D \ins J \dcup V$, which split as
    \[ \Xcal_D = \Xcal_{D \cap V} \times \Xcal_{D \cap J}. \]
    The second factor consists of input variables. By \Cref{def:cbn} the variable $X_j$, $j \in J$,
    is the projection $\pr_j:\,\Xcal_J \to \Xcal_j$, so $X_{D \cap J}$ is deterministic in the
    input $T = X_J$ in the sense of \Cref{lem:absorb-input}, and $T$ itself is in the conditioning
    slot of every statement of the theorem by the convention recalled there. \Cref{lem:absorb-input}
    then removes that factor from both rules --- by a marginalisation and a re-indexing, with no
    disintegration and hence with no hypothesis on $\Xcal_{D \cap J}$ --- and leaves the side
    conditions to be checked for $\Xcal_{D \cap V}$ alone. Since $\Xcal_v$ is a standard
    quasi-universal space for every $v \in V$ and $D \cap V$ is finite, $\Xcal_{D \cap V}$ is
    standard, see \Cref{cor:umeas-uqus-eq} and \Cref{rem:umeas}, hence countably separated. Both
    side conditions are therefore met, for every $A,B,C \ins J \dcup V$.

    All the remaining rules --- Left and Right Decomposition, $\Tk$-Inverted Right Decomposition,
    Right Weak Union, Left and Right Contraction, Right Cross Contraction and Flipped Left Cross
    Contraction --- carry no side condition, being proved by marginalisation, by re-indexing, or
    by the essential uniqueness of factorizations \Cref{thm:fact-ess-uniq}, which needs none.
\end{Prf}


\begin{thebibliography}{SYW{\etalchar{+}}16}

\bibitem[AKM26]{AKM26}
Danel Ahman, Ohad Kammar, and Rasmus~Ejlers M{\o}gelberg, \emph{{A Convenient
  Fibration for Dependently-Typed Probability Theory}}, {41st Annual ACM/IEEE
  Symposium on Logic in Computer Science (LICS 2026)} (Dagstuhl, Germany)
  (Claudia Faggian and Joost-Pieter Katoen, eds.), Leibniz International
  Proceedings in Informatics (LIPIcs), vol. 380, Schloss Dagstuhl --
  Leibniz-Zentrum f{\"u}r Informatik, 2026, pp.~4:1--4:27.
  \doi{10.4230/LIPIcs.LICS.2026.4}

\bibitem[Aum61]{Aum61}
Robert~J. Aumann, \emph{{Borel Structures for Function Spaces}}, Illinois
  Journal of Mathematics \textbf{5} (1961), no.~4, 614--630.
  \doi{10.1215/ijm/1255631584}

\bibitem[BH11]{BH11}
John~C. Baez and Alexander~E. Hoffnung, \emph{{Convenient categories of smooth
  spaces}}, Transactions of the American Mathematical Society \textbf{363}
  (2011), no.~11, 5789--5825. \doi{10.1090/S0002-9947-2011-05107-X}
  \eprintlink{0807.1704}

\bibitem[Bog07]{Bog07}
Vladimir~I. Bogachev, \emph{{Measure Theory}}, vol. 1+2, Springer, Berlin,
  2007. \doi{10.1007/978-3-540-34514-5}

\bibitem[CD17]{CD17b}
Panayiota Constantinou and A.~Philip Dawid, \emph{{Extended conditional
  independence and applications in causal inference}}, The Annals of Statistics
  \textbf{45} (2017), no.~6, 2618--2653. \doi{10.1214/16-AOS1537}
  \eprintlink{1512.00245}

\bibitem[Chu40]{Chu40}
Alonzo Church, \emph{{A Formulation of the Simple Theory of Types}}, The
  Journal of Symbolic Logic \textbf{5} (1940), no.~2, 56--68.
  \doi{10.2307/2266170}

\bibitem[CJ19]{Cho17}
Kenta Cho and Bart Jacobs, \emph{{Disintegration and Bayesian Inversion via
  String Diagrams}}, Mathematical Structures in Computer Science \textbf{29}
  (2019), no.~7, 938--971. \doi{10.1017/S0960129518000488}

\bibitem[Dar71]{Dar71}
Richard~B. Darst, \emph{{On Universal Measurability and Perfect Probability}},
  The Annals of Mathematical Statistics \textbf{42} (1971), no.~1, 352--354.
  \doi{10.1214/aoms/1177693519}

\bibitem[Daw02]{Daw02b}
A.~Philip Dawid, \emph{{Influence Diagrams for Causal Modelling and
  Inference}}, International Statistical Review \textbf{70} (2002), no.~2,
  161--189, Corrigenda: ibid.\ \textbf{70} (2002), no.~3, 437,
  \doi{10.1111/j.1751-5823.2002.tb00179.x}.
  \doi{10.1111/j.1751-5823.2002.tb00354.x}

\bibitem[DEn06]{DE06}
Eduardo~J. Dubuc and Luis Espa\~nol, \emph{{Quasitopoi over a base category}},
  2006. \eprintlink{math/0612727}

\bibitem[DF37]{deF37}
Bruno De~Finetti, \emph{{La pr{\'e}vision: ses lois logiques, ses sources
  subjectives}}, Annales de l'institut Henri Poincar{\'e} \textbf{7} (1937),
  no.~1, 1--68, Digitised by Numdam.
  \urllink{http://www.numdam.org/item/AIHP_1937__7_1_1_0/}

\bibitem[DF79]{DF79}
Lester~E. Dubins and David~A. Freedman, \emph{{Exchangeable processes need not
  be mixtures of independent, identically distributed random variables}},
  Zeitschrift f\"ur Wahrscheinlichkeitstheorie und Verwandte Gebiete
  \textbf{48} (1979), no.~2, 115--132. \doi{10.1007/BF01886868}

\bibitem[EH99]{EH99}
Mart{\'i}n Escard{\'o} and Reinhold Heckmann, \emph{{On topologies for general
  function spaces}}, Unpublished expository note, dated 19 August 1999, revised
  12 November 1999, 1999.
  \urllink{https://martinescardo.github.io/papers/exponentiablespaces.pdf}

\bibitem[EH02]{EH02}
\bysame, \emph{{Topologies on spaces of continuous functions}}, Topology
  Proceedings \textbf{26} (2002), no.~2, 545--564, Volume for 2001--2002;
  proceedings of the 16th Summer Conference on General Topology and its
  Applications. \urllink{https://martinescardo.github.io/papers/newyork.pdf}

\bibitem[Fad85]{Fad85}
Arnold~M. Faden, \emph{{The Existence of Regular Conditional Probabilities:
  Necessary and Sufficient Conditions}}, The Annals of Probability \textbf{13}
  (1985), no.~1, 288--298. \doi{10.1214/aop/1176993081}

\bibitem[FGP21]{Fri21}
Tobias Fritz, Tom{\'a}{\v{s}} Gonda, and Paolo Perrone, \emph{{De Finetti's
  Theorem in Categorical Probability}}, Journal of Stochastic Analysis
  \textbf{2} (2021), no.~4, Article 6. \doi{10.31390/josa.2.4.06}

\bibitem[FGPR23]{FriG20}
Tobias Fritz, Tom{\'a}{\v{s}} Gonda, Paolo Perrone, and Eigil~Fjeldgren
  Rischel, \emph{{Representable Markov Categories and Comparison of Statistical
  Experiments in Categorical Probability}}, Theoretical Computer Science
  \textbf{961} (2023), 113896. \doi{10.1016/j.tcs.2023.113896}

\bibitem[FM17]{FM17}
Patrick Forr{\'e} and Joris~M. Mooij, \emph{{Markov Properties for Graphical
  Models with Cycles and Latent Variables}}, Preprint, arXiv:1710.08775, 2017,
  Preprint. \doi{10.48550/arXiv.1710.08775} \eprintlink{1710.08775}

\bibitem[FM18]{FM18}
\bysame, \emph{{Constraint-based Causal Discovery for Non-linear Structural
  Causal Models with Cycles and Latent Confounders}}, Proceedings of the 34th
  Annual Conference on Uncertainty in Artificial Intelligence (UAI-2018), AUAI
  Press, 2018, pp.~269--278. \eprintlink{1807.03024}
  \urllink{http://auai.org/uai2018/proceedings/papers/117.pdf}

\bibitem[FM20]{FM20}
\bysame, \emph{{Causal Calculus in the Presence of Cycles, Latent Confounders
  and Selection Bias}}, Proceedings of the 35th Uncertainty in Artificial
  Intelligence Conference (UAI 2019), Proceedings of Machine Learning Research,
  vol. 115, PMLR, 2020, pp.~71--80. \eprintlink{1901.00433}
  \urllink{https://proceedings.mlr.press/v115/forre20a.html}

\bibitem[For21]{For21}
Patrick Forr{\'e}, \emph{{Transitional Conditional Independence}}, Preprint,
  arXiv:2104.11547, 2021, Preprint. \doi{10.48550/arXiv.2104.11547}
  \eprintlink{2104.11547}

\bibitem[FR20]{Fri19}
Tobias Fritz and Eigil~Fjeldgren Rischel, \emph{{Infinite products and zero-one
  laws in categorical probability}}, Compositionality \textbf{2} (2020), no.~3.
  \doi{10.32408/compositionality-2-3}

\bibitem[Fre15]{Fremlin}
David~H. Fremlin, \emph{{Measure Theory}}, vol. 1-6, Torres Fremlin,
  Colchester, 2000-2015.
  \urllink{https://www1.essex.ac.uk/maths/people/fremlin/mt.htm}

\bibitem[Fri20]{Fri20}
Tobias Fritz, \emph{{A synthetic approach to Markov kernels, conditional
  independence and theorems on sufficient statistics}}, Advances in Mathematics
  \textbf{370} (2020), 107239. \doi{10.1016/j.aim.2020.107239}

\bibitem[FS90]{Fre90}
Peter~J. Freyd and Andre Scedrov, \emph{{Categories, Allegories}},
  North-Holland Mathematical Library, vol.~39, North-Holland (Elsevier),
  Amsterdam, 1990, xviii + 296 pp.; ISBN 978-0-444-70368-2.
  \doi{10.1016/S0924-6509(08)X7003-5}

\bibitem[Gir82]{Gir82}
Mich{\`e}le Giry, \emph{{A Categorical Approach to Probability Theory}},
  {Categorical Aspects of Topology and Analysis}, Lecture Notes in Mathematics,
  vol. 915, Springer, Berlin, Heidelberg, 1982, Edited by B. Banaschewski,
  pp.~68--85. \doi{10.1007/BFb0092872}

\bibitem[Hey30]{Hey30}
Arend Heyting, \emph{{Die formalen Regeln der intuitionistischen Logik. I, II,
  III}}, Sitzungsberichte der Preu{\ss}ischen Akademie der Wissenschaften,
  Physikalisch-mathematische Klasse (1930), 42--56, 57--71, 158--169, Parts II
  and III are usually cited under the title ``Die formalen Regeln der
  intuitionistischen Mathematik''. \urllink{https://zbmath.org/2566292}

\bibitem[HKSY17]{Heu17}
Chris Heunen, Ohad Kammar, Sam Staton, and Hongseok Yang, \emph{A convenient
  category for higher-order probability theory}, 2017 32nd Annual ACM/IEEE
  Symposium on Logic in Computer Science (LICS), IEEE, 2017, pp.~1--12.
  \doi{10.1109/LICS.2017.8005137}

\bibitem[Jac16]{Bar16}
Bart Jacobs, \emph{{Affine Monads and Side-Effect-Freeness}}, {Coalgebraic
  Methods in Computer Science (CMCS 2016)} (Cham), Lecture Notes in Computer
  Science, vol. 9608, Springer, 2016, pp.~53--72.
  \doi{10.1007/978-3-319-40370-0_5}

\bibitem[Joh02]{Joh02}
Peter~T. Johnstone, \emph{{Sketches of an Elephant: A Topos Theory
  Compendium}}, Oxford Logic Guides, vol. 43--44, Oxford University Press,
  Oxford, 2002, 2 volumes only (Oxford Logic Guides 43 and 44), xxii + 1160
  pp.; the projected third volume has not appeared. The DOI given is the Oxford
  Academic DOI of volume 1 (ISBN 978-0-19-853425-9); volume 2 (ISBN
  978-0-19-851598-2) is 10.1093/oso/9780198515982.001.0001.
  \doi{10.1093/oso/9780198534259.001.0001}

\bibitem[Kal05]{Kal06}
Olav Kallenberg, \emph{{Probabilistic Symmetries and Invariance Principles}},
  Probability and Its Applications, Springer, New York, NY, 2005.
  \doi{10.1007/0-387-28861-9}

\bibitem[Kal17]{Kal17}
\bysame, \emph{{Random Measures, Theory and Applications}}, Probability Theory
  and Stochastic Modelling, vol.~77, Springer, Cham, 2017.
  \doi{10.1007/978-3-319-41598-7}

\bibitem[Kal21]{Kal21}
\bysame, \emph{{Foundations of Modern Probability}}, 3rd ed., Probability
  Theory and Stochastic Modelling, vol.~99, Springer, Cham, 2021.
  \doi{10.1007/978-3-030-61871-1}

\bibitem[Kan58]{Kan58}
Daniel~M. Kan, \emph{{Adjoint Functors}}, Transactions of the American
  Mathematical Society \textbf{87} (1958), no.~2, 294--329.
  \doi{10.1090/S0002-9947-1958-0131451-0}

\bibitem[Kec95]{Kec95}
Alexander~S. Kechris, \emph{{Classical Descriptive Set Theory}}, Graduate Texts
  in Mathematics, vol. 156, Springer-Verlag, New York, 1995.
  \doi{10.1007/978-1-4612-4190-4}

\bibitem[KF09]{KF09}
Daphne Koller and Nir Friedman, \emph{{Probabilistic Graphical Models:
  Principles and Techniques}}, Adaptive Computation and Machine Learning, MIT
  Press, Cambridge, MA, 2009, ISBN 978-0-262-01319-2.
  \urllink{https://mitpress.mit.edu/9780262013192/probabilistic-graphical-models/}

\bibitem[Kle65]{Kle65}
Heinrich Kleisli, \emph{{Every Standard Construction is Induced by a Pair of
  Adjoint Functors}}, Proceedings of the American Mathematical Society
  \textbf{16} (1965), no.~3, 544--546. \doi{10.1090/S0002-9939-1965-0177024-4}

\bibitem[Kle20]{Kle20}
Achim Klenke, \emph{{Probability Theory - A Comprehensive Course}}, 3rd ed.,
  Universitext, Springer, Cham, 2020. \doi{10.1007/978-3-030-56402-5}

\bibitem[Koc70]{Koc70}
Anders Kock, \emph{{Monads on symmetric monoidal closed categories}}, Archiv
  der Mathematik \textbf{21} (1970), no.~1, 1--10. \doi{10.1007/BF01220868}

\bibitem[Koc71a]{Koc71}
\bysame, \emph{{Bilinearity and cartesian closed monads}}, Mathematica
  Scandinavica \textbf{29} (1971), no.~2, 161--174.
  \doi{10.7146/math.scand.a-11042}

\bibitem[Koc71b]{Koc71c}
\bysame, \emph{{Closed categories generated by commutative monads}}, Journal of
  the Australian Mathematical Society \textbf{12} (1971), no.~4, 405--424.
  \doi{10.1017/S1446788700010272}

\bibitem[Koc72]{Koc72}
\bysame, \emph{{Strong functors and monoidal monads}}, Archiv der Mathematik
  \textbf{23} (1972), no.~1, 113--120. \doi{10.1007/BF01304852}

\bibitem[Koc12]{Koc12}
\bysame, \emph{{Commutative monads as a theory of distributions}}, Theory and
  Applications of Categories \textbf{26} (2012), no.~4, 97--131, The DOI given
  is the arXiv DOI. \doi{10.48550/arXiv.1108.5952} \eprintlink{1108.5952}

\bibitem[Kol33]{Kol33}
Andrey~Nikolaevich Kolmogorov, \emph{{Grundbegriffe der
  Wahrscheinlichkeitsrechnung}}, Ergebnisse der Mathematik und ihrer
  Grenzgebiete, vol.~2, Springer, Berlin, 1933. \doi{10.1007/978-3-642-49888-6}

\bibitem[KPS25]{KPS25}
Ohad Kammar, Seo~Jin Park, and Sam Staton, \emph{{Infinite product measures in
  quasi-Borel spaces}}, Early Ideas talk, CALCO 2025, Glasgow, 2025.

\bibitem[Lau96]{Lau96}
Steffen~L. Lauritzen, \emph{{Graphical Models}}, Oxford Statistical Science
  Series, vol.~17, Clarendon Press, Oxford, 1996.
  \doi{10.1093/oso/9780198522195.001.0001}

\bibitem[Law62]{Law62}
F.~William Lawvere, \emph{{The Category of Probabilistic Mappings}},
  Unpublished manuscript. Reproduced in the nLab and in the Lawvere Archives,
  1962. \urllink{https://ncatlab.org/nlab/files/lawvereprobability1962.pdf}

\bibitem[LDLL90]{Lau90}
S.~L. Lauritzen, A.~P. Dawid, B.~N. Larsen, and H.-G. Leimer,
  \emph{{Independence properties of directed Markov fields}}, Networks
  \textbf{20} (1990), no.~5, 491--505. \doi{10.1002/net.3230200503}

\bibitem[Lin79]{Lin79}
Harald Lindner, \emph{{Affine Parts of Monads}}, Archiv der Mathematik
  \textbf{33} (1979), no.~1, 437--443. \doi{10.1007/BF01222782}

\bibitem[McC69]{McC69}
Michael~C. McCord, \emph{{Classifying Spaces and Infinite Symmetric Products}},
  Transactions of the American Mathematical Society \textbf{146} (1969),
  273--298. \doi{10.1090/S0002-9947-1969-0251719-4}

\bibitem[ML98]{Mac98}
Saunders Mac~Lane, \emph{{Categories for the Working Mathematician}}, 2nd ed.,
  Graduate Texts in Mathematics, vol.~5, Springer, New York, 1998.
  \doi{10.1007/978-1-4757-4721-8}

\bibitem[MMS22]{MMS22}
Cristina Matache, Sean~K. Moss, and Sam Staton, \emph{{Concrete categories and
  higher-order recursion: With applications including probability,
  differentiability, and full abstraction}}, {Proceedings of the 37th Annual
  ACM/IEEE Symposium on Logic in Computer Science (LICS 2022)}, Association for
  Computing Machinery, 2022, pp.~57:1--57:14. \doi{10.1145/3531130.3533370}
  \eprintlink{2205.15917}

\bibitem[Mog91]{Mog91}
Eugenio Moggi, \emph{{Notions of Computation and Monads}}, Information and
  Computation \textbf{93} (1991), no.~1, 55--92.
  \doi{10.1016/0890-5401(91)90052-4}

\bibitem[MSR19]{Mal19}
Daniel Malinsky, Ilya Shpitser, and Thomas~S. Richardson, \emph{{A Potential
  Outcomes Calculus for Identifying Conditional Path-Specific Effects}},
  {Proceedings of the 22nd International Conference on Artificial Intelligence
  and Statistics (AISTATS 2019)}, Proceedings of Machine Learning Research,
  vol.~89, PMLR, 2019, Edited by Kamalika Chaudhuri and Masashi Sugiyama,
  pp.~3080--3088. \eprintlink{1903.03662}
  \urllink{https://proceedings.mlr.press/v89/malinsky19b.html}

\bibitem[Pac78]{Pac78}
Jan~K. Pachl, \emph{{Disintegration and Compact Measures}}, Mathematica
  Scandinavica \textbf{43} (1978), 157--168. \doi{10.7146/math.scand.a-11771}

\bibitem[Pan09]{Pra09}
Prakash Panangaden, \emph{{Labelled Markov Processes}}, Imperial College Press,
  London, 2009, Distributed by World Scientific. \doi{10.1142/p595}

\bibitem[Pea09]{Pearl09}
Judea Pearl, \emph{{Causality: Models, Reasoning, and Inference}}, 2nd ed.,
  Cambridge University Press, Cambridge, 2009. \doi{10.1017/CBO9780511803161}

\bibitem[Ram79]{Ram79}
D.~Ramachandran, \emph{{Perfect Mixtures of Perfect Measures}}, The Annals of
  Probability \textbf{7} (1979), no.~3, 444--452. \doi{10.1214/aop/1176995045}

\bibitem[Res77]{Res77}
Paul Ressel, \emph{{Some Continuity and Measurability Results on Spaces of
  Measures}}, Mathematica Scandinavica \textbf{40} (1977), no.~1, 69--78.
  \doi{10.7146/math.scand.a-11676}

\bibitem[Rey98]{Rey98}
John~C. Reynolds, \emph{{Theories of Programming Languages}}, Cambridge
  University Press, Cambridge, 1998. \doi{10.1017/CBO9780511626364}

\bibitem[Ric03]{Ric03}
Thomas~S. Richardson, \emph{{Markov Properties for Acyclic Directed Mixed
  Graphs}}, Scandinavian Journal of Statistics \textbf{30} (2003), no.~1,
  145--157. \doi{10.1111/1467-9469.00323}

\bibitem[RR81]{Rao81}
K.~P. S.~Bhaskara Rao and Bhamidi~V. Rao, \emph{{Borel Spaces}}, Dissertationes
  Mathematicae (Rozprawy Matematyczne) \textbf{190} (1981), 1--63, assigned;
  Dissertationes Mathematicae has DOIs only from vol.\ 386 (2000) on. The URL
  is the zbMATH Open record (Zbl 0472.28002).
  \urllink{https://zbmath.org/3739847}

\bibitem[RR83]{RR83}
Paul~R. Rosenbaum and Donald~B. Rubin, \emph{{The central role of the
  propensity score in observational studies for causal effects}}, Biometrika
  \textbf{70} (1983), no.~1, 41--55. \doi{10.1093/biomet/70.1.41}

\bibitem[RR13a]{RR13a}
Thomas~S. Richardson and James~M. Robins, \emph{{Single World Intervention
  Graphs: A Primer}}, {Second UAI Workshop on Causal Structure Learning,
  Bellevue, Washington}, 2013, Workshop paper.
  \urllink{https://www.stats.ox.ac.uk/~evans/uai13/Richardson.pdf}

\bibitem[RR13b]{RR13b}
\bysame, \emph{{Single World Intervention Graphs (SWIGs): A Unification of the
  Counterfactual and Graphical Approaches to Causality}}, Working Paper 128,
  Center for Statistics and the Social Sciences, University of Washington,
  2013. \urllink{https://csss.uw.edu/files/working-papers/2013/wp128.pdf}

\bibitem[Rub78]{Rub78}
Donald~B. Rubin, \emph{{Bayesian Inference for Causal Effects: The Role of
  Randomization}}, The Annals of Statistics \textbf{6} (1978), no.~1, 34--58.
  \doi{10.1214/aos/1176344064}

\bibitem[Sat18]{Sat18}
Tetsuya Sato, \emph{{The Giry Monad is not Strong for the Canonical Symmetric
  Monoidal Closed Structure on Meas}}, Journal of Pure and Applied Algebra
  \textbf{222} (2018), no.~10, 2888--2896. \doi{10.1016/j.jpaa.2017.11.004}
  \eprintlink{1612.05939}

\bibitem[Saz62]{Saz62}
Vyacheslav~Vasil'evich Sazonov, \emph{{On Perfect Measures}}, Izvestiya
  Akademii Nauk SSSR. Seriya Matematicheskaya \textbf{26} (1962), no.~3,
  391--414, In Russian. English translation: American Mathematical Society
  Translations, Series 2, vol.\ 48. \urllink{https://www.mathnet.ru/eng/im3224}

\bibitem[SC19]{SC19}
Peter Scholze and Dustin Clausen, \emph{{Lectures on Condensed Mathematics}},
  2019, Lecture notes, Universit{\"a}t Bonn, Sommersemester 2019; posted to
  arXiv in 2026. The arXiv listing gives Peter Scholze as author, with all
  results joint with Dustin Clausen. \doi{10.48550/arXiv.2605.03658}
  \eprintlink{2605.03658}

\bibitem[SC20]{SC20}
\bysame, \emph{{Lectures on Analytic Geometry}}, 2020, Lecture notes,
  Universit{\"a}t Bonn, Wintersemester 2019/20; posted to arXiv in 2026. The
  arXiv listing gives Peter Scholze as author, with all results joint with
  Dustin Clausen. \doi{10.48550/arXiv.2605.03655} \eprintlink{2605.03655}

\bibitem[See84]{See84}
Robert A.~G. Seely, \emph{{Locally cartesian closed categories and type
  theory}}, Mathematical Proceedings of the Cambridge Philosophical Society
  \textbf{95} (1984), no.~1, 33--48. \doi{10.1017/S0305004100061284}

\bibitem[Sho84]{Sho84}
Rae~M. Shortt, \emph{{Universally measurable spaces: an invariance theorem and
  diverse characterizations}}, Fundamenta Mathematicae \textbf{121} (1984),
  no.~2, 169--176. \doi{10.4064/fm-121-2-169-176}

\bibitem[{\'S}KV{\etalchar{+}}18]{SKV17}
Adam {\'S}cibior, Ohad Kammar, Matthijs V{\'a}k{\'a}r, Sam Staton, Hongseok
  Yang, Yufei Cai, Klaus Ostermann, Sean~K. Moss, Chris Heunen, and Zoubin
  Ghahramani, \emph{{Denotational validation of higher-order Bayesian
  inference}}, Proceedings of the ACM on Programming Languages \textbf{2}
  (2018), no.~POPL, 60:1--60:29, Article 60; presented at POPL 2018.
  \doi{10.1145/3158148}

\bibitem[Spa63]{Spa63}
Edwin Spanier, \emph{{Quasi-topologies}}, Duke Mathematical Journal \textbf{30}
  (1963), no.~1, 1--14. \doi{10.1215/S0012-7094-63-03001-1}

\bibitem[SSSW21]{SSSW21}
Marcin Sabok, Sam Staton, Dario Stein, and Michael Wolman, \emph{{Probabilistic
  programming semantics for name generation}}, Proceedings of the ACM on
  Programming Languages \textbf{5} (2021), no.~POPL, 1--29.
  \doi{10.1145/3434292} \eprintlink{2007.08638}

\bibitem[Ste67]{Ste67}
Norman~E. Steenrod, \emph{{A Convenient Category of Topological Spaces}},
  Michigan Mathematical Journal \textbf{14} (1967), no.~2, 133--152.
  \doi{10.1307/mmj/1028999711}

\bibitem[Str72]{Str72}
Ross Street, \emph{{The Formal Theory of Monads}}, Journal of Pure and Applied
  Algebra \textbf{2} (1972), no.~2, 149--168.
  \doi{10.1016/0022-4049(72)90019-9}

\bibitem[Str09]{Str09}
Neil Strickland, \emph{{The category of CGWH spaces}}, Unpublished lecture
  notes, University of Sheffield (the notes themselves are undated), 2009.
  \urllink{https://ncatlab.org/nlab/files/StricklandCGHWSpaces.pdf}

\bibitem[SYW{\etalchar{+}}16]{Sta16}
Sam Staton, Hongseok Yang, Frank Wood, Chris Heunen, and Ohad Kammar,
  \emph{{Semantics for probabilistic programming: higher-order functions,
  continuous distributions, and soft constraints}}, {Proceedings of the 31st
  Annual ACM/IEEE Symposium on Logic in Computer Science (LICS 2016)} (New
  York, NY, USA), Association for Computing Machinery, 2016, pp.~525--534.
  \doi{10.1145/2933575.2935313} \eprintlink{1601.04943}

\bibitem[\v{C}82]{Cen82}
Nikolai~Nikolaevich \v{C}encov, \emph{{Statistical Decision Rules and Optimal
  Inference}}, Translations of Mathematical Monographs, vol.~53, American
  Mathematical Society, Providence, RI, 1982, translated from Russian.
  \doi{10.1090/mmono/053}

\bibitem[Vig03]{Vig03}
Sebastiano Vigna, \emph{{A Guided Tour in the Topos of Graphs}}, 2003.
  \eprintlink{math/0306394}

\bibitem[VO18]{VO18}
Matthijs V\'ak\'ar and Luke Ong, \emph{{On S-Finite Measures and Kernels}},
  2018. \doi{10.48550/arXiv.1810.01837} \eprintlink{1810.01837}

\end{thebibliography}
\end{document}